\documentclass{article}
\usepackage{graphicx} % Required for inserting images

\usepackage{authblk}
\usepackage{geometry}
\usepackage[utf8]{inputenc}
\usepackage{amsmath}
\usepackage{amssymb}
\usepackage{amsthm}
\usepackage{amsfonts}
\usepackage{fancyhdr}
\usepackage{multicol}
\usepackage{graphicx}
\usepackage{longtable}
\usepackage{mleftright}
\usepackage{mathtools}
\usepackage{physics}
\usepackage{xcolor}
\usepackage[shortlabels]{enumitem}
\usepackage[square]{natbib}
\usepackage{dsfont}
\usepackage{mathrsfs}
\usepackage{bm}
\usepackage{microtype} % improved typesetting
\usepackage[citecolor=blue]{hyperref}
\usepackage[capitalise,noabbrev]{cleveref}
\usepackage{cleveref}
\usepackage{cancel}
\usepackage{pdflscape}
\usepackage{placeins}

\theoremstyle{definition}
\newtheorem{thm}{Theorem}
\newtheorem{lemm}{Lemma}
\newtheorem{prop}{Proposition}
\newtheorem{coro}{Corollary}

\newtheorem{defi}{Definition}

\theoremstyle{remark}
\newtheorem*{rem}{Remark}

\crefname{thm}{Theorem}{Theorems}
\crefname{lemm}{Lemma}{Lemmas}
\crefname{prop}{Proposition}{Propositions}
\crefname{coro}{Corollary}{Corollaries}
\crefname{ex}{Example}{Examples}
\crefname{con}{Condition}{Conditions}
\crefname{conj}{Conjecture}{Conjectures}
\crefname{assum}{Assumption}{Assumptions}
\crefname{algo}{Algorithm}{Algorithms}
\crefname{defi}{Definition}{Definitions}
\crefname{rem}{Remark}{Remarks}
\crefname{proofEnd}{proof}{proofs}
\crefname{proofEnd}{Proof}{Proofs}

\usepackage{hyperref}
\hypersetup{
    colorlinks=true,
    linkcolor=blue,
    filecolor=magenta,      
    urlcolor=cyan,
    pdftitle={Closed Form Z-Estimation under Coupled Sample-Parameter Asymptotics},
    pdfpagemode=FullScreen,
    }

\renewcommand{\left}{\mleft}
\renewcommand{\right}{\mright}

\newcommand{\R}{\mathbb{R}}

\newcommand{\N}{\mathbb{N}}

\newcommand{\E}{\mathbb{E}}
\newcommand{\Var}{\mathrm{Var}}

\newcommand{\sigmin}{\sigma_\mathrm{min}}

\newcommand{\dif}{\mathop{}\!\mathrm{d}}

\newcommand{\inP}{\stackrel{p}{\longrightarrow}}    % Convergence in P
\newcommand{\inD}{\stackrel{d}{\longrightarrow}}    % Convergence in D

\newcommand{\hattn}{\hat{\btheta}_n}
\newcommand{\tildetn}{\tilde{\btheta}_n}

\renewcommand{\log}[1]{\mathrm{log}\mleft( #1 \mright)}
\newcommand{\logtwo}[1]{\mathrm{log}_2\mleft( #1 \mright)}

\newcommand{\atantwo}[2]{\mathrm{atan2}\left( #1, #2 \right)}

\newcommand{\iid}{\overset{\mathrm{iid}}{\sim}}

\renewcommand{\op}[1]{o_{\mathbb{P}}\mleft( #1 \mright)}
\newcommand{\Op}[1]{O_{\mathbb{P}}\mleft( #1 \mright)}

\newcommand{\inv}{^{-1}}
\renewcommand{\norm}[1]{\left\lVert#1\right\rVert}

\newcommand{\ba}{{\bm{a}}}

\newcommand{\be}{{\bm{e}}}

\newcommand{\bh}{{\bm{h}}}

\newcommand{\bq}{{\bm{q}}}
\newcommand{\br}{{\bm{r}}}
\newcommand{\bs}{{\bm{s}}}
\newcommand{\bt}{{\bm{t}}}
\newcommand{\bu}{{\bm{u}}}
\newcommand{\bv}{{\bm{v}}}
\newcommand{\bw}{{\bm{w}}}
\newcommand{\bx}{{\bm{x}}}

\newcommand{\bz}{{\bm{z}}}

\newcommand{\Ab}{\mathbf{A}}
\newcommand{\Bb}{\mathbf{B}}

\newcommand{\Db}{\mathbf{D}}

\newcommand{\Gb}{\mathbf{G}}
\newcommand{\Hb}{\mathbf{H}}
\newcommand{\Ib}{\mathbf{I}}
\newcommand{\Jb}{\mathbf{J}}

\newcommand{\Mb}{\mathbf{M}}

\newcommand{\Ob}{\mathbf{O}}

\newcommand{\Rb}{\mathbf{R}}
\newcommand{\Sb}{\mathbf{S}}
\newcommand{\Tb}{\mathbf{T}}

\newcommand{\Xb}{\mathbf{X}}

\newcommand{\bF}{{\bm{F}}}
\newcommand{\bG}{{\bm{G}}}

\newcommand{\bL}{{\bm{L}}}

\newcommand{\bP}{{\bm{P}}}
\newcommand{\bQ}{{\bm{Q}}}

\newcommand{\bS}{{\bm{S}}}
\newcommand{\bT}{{\bm{T}}}
\newcommand{\bU}{{\bm{U}}}
\newcommand{\bV}{{\bm{V}}}
\newcommand{\bW}{{\bm{W}}}
\newcommand{\bX}{{\bm{X}}}
\newcommand{\bY}{{\bm{Y}}}
\newcommand{\bZ}{{\bm{Z}}}

\newcommand{\bzero}{{\bm{0}}}

\newcommand{\PP}{\mathbb{P}}

\newcommand{\bbeta}{{\bm{\beta}}}
\newcommand{\bgamma}{{\bm{\gamma}}}
\newcommand{\bepsilon}{{\bm{\varepsilon}}}

\newcommand{\btheta}{{\bm{\theta}}}

\newcommand{\bphi}{{\bm{\phi}}}

\newcommand{\bpsi}{{\bm{\psi}}}

\renewcommand{\exp}[1]{\mathrm{exp}\left(#1\right)}
\renewcommand{\sin}[1]{\mathrm{sin}\left(#1\right)}
\newcommand{\sinsq}[1]{\mathrm{sin}^2\left(#1\right)}
\renewcommand{\cos}[1]{\mathrm{cos}\left(#1\right)}
\newcommand{\cossq}[1]{\mathrm{cos}^2\left(#1\right)}

\newcommand{\dand}{\quad \text{and} \quad}

\title{An Asymptotic Framework for Closed-Form Estimators from Surrogate Estimating Equations}

\author[1]{Yigal Teplitsky%
\thanks{Present affiliation: Department of Physics, Princeton University.
Email: \href{mailto:yigal.teplitsky@princeton.edu}
{\textcolor{black}{\texttt{yigal.teplitsky@princeton.edu}}}.}}

\author[1,2]{Robert Zimmerman%
\thanks{Email: \href{mailto:robert.zimmerman@alumni.utoronto.ca}
{\textcolor{black}{\texttt{robert.zimmerman@alumni.utoronto.ca}}}.}}

\affil[1]{Department of Statistical Sciences, University of Toronto}
\affil[2]{Department of Mathematics, Imperial College London}

\date{November 5, 2026}

\begin{document}

\maketitle

\begin{abstract}
    Many statistical estimators, including maximum likelihood and method of moments estimators, are defined as roots of estimating equations. In complex models, these equations often cannot be solved explicitly; although numerical methods can provide approximate solutions, closed-form approximations can reduce computation and facilitate theoretical analysis by making an estimator's dependence on the data more transparent. In this paper, we develop a general approach to constructing such approximations and a framework for quantifying their accuracy. Our theory accommodates sequences of parameters that approach the boundary of the parameter space or diverge to infinity, and provides conditions under which a closed-form approximate estimator has the same first-order asymptotic distribution as its exact counterpart. We focus on approximations to maximum likelihood estimators, although the theory applies to general estimating equations. We apply the framework to iid models from the gamma, zeta, von Mises, and generalized Pareto families, as well as to gamma regression. Together, these examples cover scalar and vector parameters, iid and non-iid settings, and normal and nonnormal asymptotic regimes. Each example is accompanied by a simulation study, and the gamma regression example also includes a real-world data application.
\end{abstract}

\tableofcontents

\section{Introduction}

Many statistical estimators are defined implicitly as roots of data-dependent estimating equations. Given a parameter space $\Theta \subseteq \R^d$ and an estimating equation $\bq_n:\Theta \to \R^d$ of the form
\begin{equation}\label{eq:estimating_equation}
    \bq_n(\btheta \mid \bX_n) = \bzero,
\end{equation}
a root --- when it exists --- is known as a Z-estimator. When $\bq_n$ is a score function for some statistical model, for example, any maximum likelihood estimator (MLE) that is an interior point of $\Theta$ is a root of $\bq_n$. While score equations are perhaps the most prominent sources of Z-estimators, method of moments and generalized method of moments estimators under suitable regularity conditions \citep[e.g.,][]{hansenLargeSampleProperties1982} also fit into this framework. The study of such estimators dates back to \cite{huberBehaviorMaximumLikelihood1967}; modern treatments are given by \cite{vandervaartAsymptoticStatistics1998} and \cite{vandervaartWeakConvergenceEmpirical2023a} who develop asymptotic theory of these estimators.

In the simplest models, \eqref{eq:estimating_equation} can be solved explicitly, yielding a closed-form\footnote{We use the term \emph{closed form} informally to mean an explicit expression in the data obtained from finitely many algebraic operations, elementary functions, and standard special functions; we do not regard an estimator specified only as an implicit root or the output of a numerical optimization procedure as being in closed form.} function of the data $\bX_n$. More commonly, however, the estimating equation contains a nonlinear --- often transcendental --- function with no explicit inverse, couples several parameters together, or has a nearly singular Jacobian near the boundary of $\Theta$. Classical examples include the beta and Dirichlet distributions \citep{beckmanMaximumLikelihoodEstimation1978,minkaEstimatingDirichletDistribution2000}, the Student's $t$ distribution \citep{liuMLESTIMATIONDISTRIBUTION1995,hasannasabAlternativesEMAlgorithm2021}, the von Mises distribution \citep{schouEstimationConcentrationParameter1978,hornikMaximumLikelihoodEstimation2014}, and families used in extreme value analysis \citep{prescottMaximumLikelihoodEstimation1980,grimshawComputingMaximumLikelihood1993}. Although modern root-finding methods and software often make solving such equations routine, an explicit approximation avoids the need for an initial value, a stopping criterion, and repeated computations; it can also be useful in regimes where numerical root-finding is unstable, and can facilitate theoretical analysis by making the estimator's dependence on the data more transparent. Work on specific models has therefore pursued closed-form or rapidly computable alternatives when exact likelihood estimation requires iteration; examples include \cite{yeClosedFormEstimatorsGamma2017} for the gamma distribution, \cite{clausetPowerLawDistributionsEmpirical2009} for the zeta distribution, and \cite{hochangAsymptoticallyEfficientClosedform2023} for the Dirichlet distribution. 

In this paper, we present an approach to construct closed-form approximations to roots of intractable estimating equations of the form \eqref{eq:estimating_equation} by replacing the original estimating equation with a tractable surrogate that admits a closed-form solution, and we develop a general framework for determining when such approximations are accurate and for quantifying the resulting error. The framework considers sequences of data-generating parameters that approach a limiting value; importantly, this limiting value may lie on the boundary of the parameter space or be a formal value such as infinity.

The gamma model provides a simple illustration. Suppose that $X_1,X_2,\ldots,X_n$ constitute a random sample from the $\text{Gamma}(\alpha,\beta)$ distribution in the shape-rate parameterization, where the rate parameter $\beta$ is known. The MLE $\hat\alpha_n$ maximizes the log-likelihood function, and differentiating the log-likelihood shows that it is the unique solution of
\begin{equation}\label{eq:intro_gamm_orig}
    \psi(\hat\alpha_n)=\log{\beta G_n},
\end{equation}
where $\psi\coloneqq(\log{\Gamma})'$ is the digamma function and $G_n\coloneqq\left(\prod_{i=1}^n X_i\right)^{1/n}$ is the geometric mean of the sample. Because $\psi$ has no standard inverse in closed form, $\hat\alpha_n$ is usually obtained numerically, using Newton's method, for example; see \cite{greenwoodAidsFittingGamma1960,choiMaximumLikelihoodEstimation1969,minkaEstimatingGammaDistribution2002} for computational treatments. However, \eqref{eq:intro_gamm_orig} is equivalent to
\begin{equation}\label{eq:intro_gamma_EE}
    \exp{-\psi(\hat\alpha_n)}=\frac{1}{\beta G_n},
\end{equation}
and, as $\alpha\to\infty$, one can show that
\[
    \exp{-\psi(\alpha)}
    =\frac{1}{\alpha-1/2}
    +O\left(\frac{1}{\alpha^3}\right).
\]
Replacing the left side of \eqref{eq:intro_gamma_EE} by its leading approximation and solving yields the explicit estimator
\begin{equation}\label{eq:alphaapproxMLE}
    \tilde\alpha_n
    \coloneqq\frac{1}{2}
      +\beta\left(\prod_{i=1}^n X_i\right)^{1/n}.
\end{equation}
We show in \cref{sec:gamma} that, for data-generating sequences $\{\alpha_n\}_{n\geq1}$ satisfying $\alpha_n\to\infty$,
\[
    \lvert\tilde\alpha_n-\hat\alpha_n\rvert
    =o_{\mathbb P_{\alpha_n}^{\otimes n}}(1),
\]
so the two estimators are equivalent on the original parameter scale. If, in addition, $\alpha_n^3/n\to\infty$, then
\[
    \sqrt{\frac{n}{\alpha_n}}\,
    \lvert\tilde\alpha_n-\hat\alpha_n\rvert
    =o_{\mathbb P_{\alpha_n}^{\otimes n}}(1),
\]
and under the same condition, we further show that
\[
    \sqrt{\frac{n}{\alpha_n}}
    (\tilde\alpha_n-\alpha_n)
    \inD\mathcal N(0,1)
    \dand
    \sqrt{\frac{n}{\alpha_n}}
    (\hat\alpha_n-\alpha_n)
    \inD\mathcal N(0,1).
\]
Thus, the approximate estimator and the MLE have the same first-order asymptotic distribution.

The preceding calculation suggests a general strategy. First, rewrite the exact estimating equation $\bq_n(\btheta \mid \bX_n) = \bzero$ in a form that exposes the behavior of its intractable component near the limiting parameter value. Next, replace that component by a finite truncation of its asymptotic expansion, retaining sufficiently many expansion terms to control the approximation error while keeping the resulting surrogate estimating equation $\tilde\bq_n(\btheta \mid \bX_n) = \bzero$ explicitly solvable. The root $\tilde{\btheta}_n$ of the surrogate equation is then a candidate approximation to the selected exact root $\hat\btheta_n$. Whether it inherits the first-order asymptotic distribution of $\hat\btheta_n$ depends not only on convergence of the parameter to its limiting value, but on its rate of approach relative to the sample size as well.

Although deriving $\tilde{\btheta}_n$ can be straightforward, as in the gamma example above, establishing that it is close to $\hat\btheta_n$ requires a careful analysis. The essential problems are that both the equations and their selected roots are random, the sample size and the data-generating parameter move jointly, and the accuracy of the approximation depends on the scale at which the estimators will be compared. Accordingly, even if $\norm{\bq_n(\btheta \mid \bX_n) - \tilde\bq_n(\btheta \mid \bX_n)} \inP 0$ for each fixed $\btheta \in \Theta$, it does not follow that $\norm{\hat\btheta_n - \tilde\btheta_n} \inP 0$. Pointwise control need not extend uniformly to a neighborhood containing both $\hat\btheta_n$ and $\tilde\btheta_n$, and even a small perturbation can produce a large change, especially when the Jacobian of the estimating equation is nearly singular as in the general framework of \citep{riceTheoryCondition1966}; a deterministic antecedent appears in \citet[Chapter 2]{bruijnAsymptoticMethodsAnalysis1981}.

To address these difficulties, we introduce one-to-one transformations of the parameter, which we call \emph{vanishing charts}, under which the transformed data-generating parameter converges to zero along the regime under consideration. Such regimes can differ substantially from standard fixed-parameter asymptotics, exhibiting altered convergence rates, non-Gaussian limits, or even failure of MLE existence \citep{smithMaximumLikelihoodEstimation1985,selfAsymptoticPropertiesMaximum1987,geyerAsymptoticsConstrainedMEstimation1994}. After localizing the exact and surrogate roots on a common deterministic set in vanishing chart coordinates, our analysis separates three factors governing their discrepancy: the error incurred when $\bq_n$ is replaced by $\tilde\bq_n$, the sensitivity of the surrogate root to perturbations, and the magnification that occurs when errors in chart coordinates are translated back to the original parameter scale. Our main result combines uniform control of the equation error with a lower bound on the smallest singular value of the surrogate Jacobian to bound the distance between the roots in chart coordinates. Further results accommodate vector parameters with different componentwise rates, transfer these bounds to the original parameter scale, and identify when $\tilde\btheta_n$ inherits the first-order asymptotic distribution of $\hat\btheta_n$.

We study the gamma family in detail: the known-rate case is scalar, joint estimation of the shape and rate introduces different componentwise rates, and gamma regression shows that the observations need not be identically distributed (a real-world data example compares the exact and approximate gamma regression fits for measurements of black cherry trees). Our remaining examples highlight other features of the construction, with mathematical details left to the appendix. The zeta distribution is a discrete example where the MLE can fail to be finite. For the von Mises model, we preserve the explicit location estimator while approximating a ratio of Bessel functions for the concentration parameter. The generalized Pareto model has nonregular, parameter-dependent support, requiring the likelihood to be profiled prior to truncation. For all of these models, we use R simulations to assess the finite-sample accuracy of the approximate estimators and to illustrate the distinct features of each model.

The remainder of the paper develops the general theory in \cref{sec:theory}, presents the worked examples in \cref{sec:gamma,sec:additionalexamples}, and concludes with a discussion in \cref{sec:discussion}; technical details are deferred to the appendices.

\section{Theory and Practical Guidance}\label{sec:theory}

\subsection{Setup}\label{sec:defs}

For each $n\geq1$, let $\bX_n \coloneqq (X_1, \ldots, X_n)$ be the random vector of data and consider the sequence of parametric models $\{\mathbb P_{n,\btheta}:\btheta\in\Theta\subseteq \R^d\}_{n\geq1}$, where $\mathbb P_{n,\btheta}$ denotes the joint law of $\bX_n$. In the iid setting, $\bX_n=(X_1,\ldots,X_n)$ and $\mathbb P_{n,\btheta}=\mathbb P_{\btheta}^{\otimes n}$. When a joint density or mass function $p_{n,\btheta}$ exists, an interior MLE at which $\log{p_{n,\btheta}(\bX_n)}$ is differentiable is a root of the joint score function $\nabla_{\btheta}\log{p_{n,\btheta}(\bX_n)}$. For a differentiable map $\bF:\R^d\to\R^m$, let $\Jb_\bF(\bx)$ denote its Jacobian at $\bx$, and for a square matrix $\Ab$, let $\sigmin(\Ab)$ denote its smallest singular value. $\norm{\cdot}$ denotes the Euclidean norm on vectors and its induced operator norm on matrices. For $\bz\in\R^m$ and $r>0$, we write $\overline{B}_r(\bz)\coloneqq\{\bx\in\R^m:\norm{\bx-\bz}\leq r\}$ for the closed ball of radius $r$ centered at $\bz$. When it is useful to display the governing law explicitly, we write
\[
    \bZ_n\xrightarrow{\mathbb P_{n,\btheta_n}}\bz
\]
to mean that $\mathbb P_{n,\btheta_n}(\norm{\bZ_n-\bz}>\varepsilon)\to0$ for every $\varepsilon>0$; when the laws are clear, we simply write $\bZ_n\inP\bz$. We employ corresponding stochastic order notation throughout: along a sequence $\{\btheta_n\}_{n\geq1}$ and for deterministic $r_n>0$, we write $\bZ_n=o_{\mathbb P_{n,\btheta_n}}(r_n)$ when $\bZ_n/r_n\to\bzero$ in probability under the sequence of laws $\{\mathbb P_{n,\btheta_n}\}_{n\geq1}$, and $\bZ_n=O_{\mathbb P_{n,\btheta_n}}(r_n)$ when $\bZ_n/r_n$ is bounded in probability under these laws. When the laws are clear, we abbreviate these as $\op{r_n}$ and $\Op{r_n}$, respectively.

Our framework involves parameter sequences that may diverge or approach the boundary of the parameter space. We formalize this behavior through a \textit{limiting regime} and then map it to zero using a \textit{vanishing chart}:

\begin{defi}
A \emph{limiting regime} $(\mathcal R, \btheta^*)$ consists of a collection $\mathcal R$ of deterministic sequences $\{\btheta_n\}_{n\geq1}\subseteq\Theta$ and an associated \emph{limiting value} $\btheta^*$. We write $\btheta_n\to\btheta^*$ as shorthand for $\{\btheta_n\}_{n\geq1}\in\mathcal R$ and simply refer to $\mathcal R$ as the limiting regime when the limiting value is understood from context. The limiting value may be a point of $\overline\Theta$ or a formal value such as infinity.
\end{defi}

\begin{defi}
A \emph{vanishing chart} for a limiting regime $\mathcal R$ is a one-to-one map $\bphi:N\to U$, where $U\coloneqq\bphi(N)\subseteq\R^d$ and $N\subseteq\Theta$, such that, for every $\{\btheta_n\}_{n\geq1}\in\mathcal R$, we have $\btheta_n\in N$ eventually and $\bphi(\btheta_n)\to\bzero$. We write $\bphi\inv:U\to N$ for the inverse of this map. Note that neither $N$ nor $U$ need be open; when we assert that $\bphi$ or $\bphi\inv$ is differentiable, it is understood that the map admits a differentiable extension on an open neighborhood of the subset on which its derivative is used.
\end{defi}

An estimator may be naturally defined only on an event $E_n$ whose probability tends to one along the regime sequence ${\btheta_n}_{n\geq1}\in\mathcal R$ under consideration. In that case, we extend it to a measurable $\Theta$-valued random variable by assigning a fixed value in $\Theta$ on $E_n^c$ (which we refer to as an \emph{exceptional event}); this value is chosen in advance and does not depend on the unknown parameter or the regime sequence. Similarly, if an estimator falls outside the domain $N$ of a vanishing chart $\bphi\to U$, so that its image under $\bphi$ is not defined, we assign its transformed version a predetermined value in $U$. These completions ensure that the relevant random quantities are defined everywhere. Once the required high-probability containment has been established, the particular values assigned outside the corresponding events do not affect the asymptotic conclusions. We describe this containment using the following terminology:

\begin{defi}
Along a fixed regime sequence $\{\btheta_n\}_{n\geq1}\in\mathcal R$, we say that a random sequence $\{\bZ_n\}_{n\geq1}$ is \emph{localized in} deterministic sets $A_n$ if $\mathbb P_{n,\btheta_n}\left(\bZ_n\in A_n\right)\longrightarrow1$. If, in addition, $\sup_{\bz\in A_n}\norm{\bz-\bz_0}\longrightarrow0$, we say that $\bZ_n$ is \emph{localized at} $\bz_0$.
\end{defi}

\begin{defi}
\label{def:approximate_est}
Let $\bphi:N\to U$ be a vanishing chart for $\mathcal R$, and let $\{\hattn\}_{n\geq1}$ and $\{\tildetn\}_{n\geq1}$ be estimator sequences in $\Theta$. Fix $\{\btheta_n\}_{n\geq1}\in\mathcal R$. For a deterministic sequence $\rho_n\geq1$, the estimators are \emph{equivalent at scale $\rho_n$ relative to $\bphi$ along $\{\btheta_n\}_{n\geq1}$} if
\[
\mathbb P_{n,\btheta_n}\left(\hattn,\tildetn\in N\right)
\longrightarrow1
\dand
\bphi(\hattn)=o_{\mathbb P_{n,\btheta_n}}(1),
\]
and
\[
\rho_n
\norm{\bphi(\tildetn)-\bphi(\hattn)}
=o_{\mathbb P_{n,\btheta_n}}(1).
\]
These conditions imply $\bphi(\tildetn)=o_{\mathbb P_{n,\btheta_n}}(1)$. For a deterministic sequence $\rho_n'>0$, the estimators are \emph{equivalent at scale $\rho_n'$ on the original parameter scale along $\{\btheta_n\}_{n\geq1}$} if
\[
\rho_n'\norm{\tildetn-\hattn}
=o_{\mathbb P_{n,\btheta_n}}(1).
\]
When either scale sequence is identically one, we simply say that the estimators are equivalent on the corresponding scale. When the scale sequence tends to infinity, we say that the estimators are asymptotically equivalent under the corresponding normalization. Each property holds in the limiting regime $\mathcal R$ if it holds along every sequence in $\mathcal R$. If a common normalization is asserted throughout $\mathcal R$, the same deterministic normalizing sequence must work along every regime sequence. If $\hattn$ is an MLE and an equivalence property has been established at a stated chart or the original parameter scale, we call $\tildetn$ an \emph{approximate maximum likelihood estimator}, or \emph{approximate MLE}, at that scale. 
\end{defi}

Chart equivalence at scale one records simultaneous localization of the two roots near the limiting chart value, but need not quantify surrogate accuracy beyond that localization. However, if
\[
\Rb_n\bigl(\bphi(\hattn)-\bphi(\btheta_n)\bigr)\inD\bZ
\dand
\Rb_n\bigl(\bphi(\tildetn)-\bphi(\hattn)\bigr)
=o_{\mathbb P_{n,\btheta_n}}(1),
\]
then by Slutsky's theorem, the corresponding normalized error of $\tildetn$ also converges to $\bZ$; the analogous statement holds on the original parameter scale.

One may naturally ask whether the asymptotic equivalence in \cref{def:approximate_est} depends on the specific choice of vanishing chart. Our first proposition shows that both asymptotic equivalence and the corresponding normalization are preserved under a smooth mapping whose Jacobian at zero is nonsingular:

\begin{prop}\label{prop:chart_invariance}
Let $\bphi:N_1\to U_1$ and $\bpsi:N_2\to U_2$ be vanishing charts for $\mathcal R$. Suppose that $\bT\coloneqq\bpsi\circ\bphi\inv$, defined on $\bphi(N_1\cap N_2)$, admits a $C^1$ extension $\widetilde\bT$ to an open neighborhood of $\bzero$ such that $\widetilde\bT(\bzero)=\bzero$ and $\Jb_{\widetilde\bT}(\bzero)$ is nonsingular. For every $\{\btheta_n\}_{n\geq1}\in\mathcal R$, let $\{\hat{\btheta}_n^{(1)}\}_{n\geq1}$ and $\{\hat{\btheta}_n^{(2)}\}_{n\geq1}$ satisfy
\[
\mathbb P_{n,\btheta_n}\left(
\hat{\btheta}_n^{(1)},\hat{\btheta}_n^{(2)}\in N_1\cap N_2\right)\to1,
\]
and $\bphi(\hat{\btheta}_n^{(i)})=o_{\mathbb P_{n,\btheta_n}}(1)$ for $i=1,2$.
Then $\bpsi(\hat{\btheta}_n^{(i)})=o_{\mathbb P_{n,\btheta_n}}(1)$ for $i=1,2$, and, for every deterministic $\rho_n\geq1$ and every regime sequence,
\begin{equation}\label{eq:approx_equiv_chart}
\rho_n\norm{\bphi(\hat{\btheta}_n^{(1)})-
\bphi(\hat{\btheta}_n^{(2)})}=o_{\mathbb P_{n,\btheta_n}}(1)
\quad\Longleftrightarrow\quad
\rho_n\norm{\bpsi(\hat{\btheta}_n^{(1)})-
\bpsi(\hat{\btheta}_n^{(2)})}=o_{\mathbb P_{n,\btheta_n}}(1).
\end{equation}
\end{prop}
\begin{proof}
Fix an arbitrary $\{\btheta_n\}_{n\geq1}\in\mathcal R$; throughout the proof, probabilistic statements are understood under $\{\mathbb P_{n,\btheta_n}\}_{n\geq1}$. Let $\widetilde\bT:W\to\R^d$ be the assumed $C^1$ extension of $\bT=\bpsi\circ\bphi\inv$, where $W$ is an open neighborhood of $\bzero$. Since $\widetilde\bT$ is continuously differentiable at $\bzero$, there is an $r>0$ such that $\overline{B}_r(\bzero)\subseteq W$ and
\[
    \sup_{\bu\in\overline{B}_r(\bzero)}
    \norm{\Jb_{\widetilde\bT}(\bu)-\Jb_{\widetilde\bT}(\bzero)}\leq \frac{\sigmin(\Jb_{\widetilde\bT}(\bzero))}{2}.
\]
Because $\overline{B}_r(\bzero)$ is convex, the line segment joining $\bu$ and $\bv$ lies in $\overline{B}_r(\bzero)$. Applying the fundamental theorem of calculus to the path $t\longmapsto\widetilde\bT\bigl(\bu+t(\bv-\bu)\bigr)$ gives, for $\bu,\bv\in\overline{B}_r(\bzero)$,
\[
\widetilde\bT(\bv)-\widetilde\bT(\bu)
=
\int_0^1
\Jb_{\widetilde\bT}\bigl(\bu+t(\bv-\bu)\bigr)
(\bv-\bu)\dif t
\]
\citep[e.g.,][Theorem~3.2.7]{ortegaIterativeSolutionNonlinear1970}. Consequently,
\begin{equation}\label{eq:chart_bilip}
    \frac{\sigmin(\Jb_{\widetilde\bT}(\bzero))}{2}\norm{\bv-\bu}
    \leq
    \norm{\widetilde\bT(\bv)-\widetilde\bT(\bu)}
    \leq
    \left(\norm{\Jb_{\widetilde\bT}(\bzero)}+\frac{\sigmin(\Jb_{\widetilde\bT}(\bzero))}{2}\right)\norm{\bv-\bu},
\end{equation}
where the upper bound follows directly from the integral formula and the lower bound follows from subtracting $\Jb_{\widetilde\bT}(\bzero)(\bv-\bu)$ and using $\norm{\Jb_{\widetilde\bT}(\bzero)(\bv-\bu)}\geq \sigmin(\Jb_{\widetilde\bT}(\bzero))\norm{\bv-\bu}$.

For $i=1,2$, write $\bU_{i,n}\coloneqq\bphi(\hat{\btheta}_n^{(i)})$, using the fixed chart completion outside $N_1$, and define
\[
E_n\coloneqq
\{\hat{\btheta}_n^{(1)},\hat{\btheta}_n^{(2)}\in N_1\cap N_2\}
\cap
\{\bphi(\hat{\btheta}_n^{(1)}),\bphi(\hat{\btheta}_n^{(2)}) \in B_r(\bzero)\}.
\]
Our assumptions, together with $\bphi(\hat{\btheta}_n^{(1)}),\bphi(\hat{\btheta}_n^{(2)})=o_{\mathbb P_{n,\btheta_n}}(1)$, imply $\mathbb P_{n,\btheta_n}(E_n)\to1$. On the event $E_n$,
\[
\bpsi(\hat{\btheta}_n^{(1)})
=\widetilde\bT(\bphi(\hat{\btheta}_n^{(1)}))
\dand
\bpsi(\hat{\btheta}_n^{(2)})=\widetilde\bT(\bphi(\hat{\btheta}_n^{(2)})).
\]
Since $\widetilde\bT(\bzero)=\bzero$, the upper bound in \eqref{eq:chart_bilip}, applied with $\bu=\bzero$, shows that $\bpsi(\hat{\btheta}_n^{(i)})=o_{\mathbb P_{n,\btheta_n}}(1)$. Applying both bounds in \eqref{eq:chart_bilip} to $\bphi(\hat{\btheta}_n^{(1)})$ and $\bphi(\hat{\btheta}_n^{(2)})$ on $E_n$, and multiplying by any deterministic $\rho_n\geq1$, gives
\[
\begin{aligned}\frac{\sigmin(\Jb_{\widetilde\bT}(\bzero))}{2}\rho_n
 \norm{\bphi(\hat{\btheta}_n^{(1)})-\bphi(\hat{\btheta}_n^{(2)})}
&\leq
\rho_n\norm{\bpsi(\hat{\btheta}_n^{(1)})-
                 \bpsi(\hat{\btheta}_n^{(2)})}
\\
&\leq
\left(\norm{\Jb_{\widetilde\bT}(\bzero)}+
      \frac{\sigmin(\Jb_{\widetilde\bT}(\bzero))}{2}\right)
\rho_n\norm{\bphi(\hat{\btheta}_n^{(1)})-\bphi(\hat{\btheta}_n^{(2)})}.
\end{aligned}
\]
Since $\mathbb P_{n,\btheta_n}(E_n^c)\to0$, these inequalities prove \eqref{eq:approx_equiv_chart}. Since the choice of regime sequence was arbitrary, the conclusion holds throughout $\mathcal R$.
\end{proof}

\subsection{Uniform Root Stability}\label{sec:maintheorem}

We now develop results that establish the existence and localization of an exact root of the estimating equation, and bound its distance from the surrogate root. To this end, fix an arbitrary regime sequence $\{\btheta_n\}_{n\geq1}\in\mathcal R$; throughout this section, all probabilistic statements are under $\mathbb P_{n,\btheta_n}$. Let $\bphi:N\to U$ be a vanishing chart and, for random maps $\bq_n,\tilde{\bq}_n:\Theta\to\R^d$ and $\bu\in U$, write the estimating equations in chart coordinates as
\[
\bQ_n(\bu)\coloneqq\bq_n\bigl(\bphi\inv(\bu)\bigr)\dand
\tilde{\bQ}_n(\bu)\coloneqq\tilde{\bq}_n\bigl(\bphi\inv(\bu)\bigr).
\]
The corresponding roots in these coordinates are
\[
\hat{\bu}_n\coloneqq\bphi(\hattn)
\dand
\tilde{\bu}_n\coloneqq\bphi(\tildetn),
\]
using the convention for exceptional events from \cref{sec:defs}.

Our first theorem is useful when only the surrogate root has been localized. It both produces a nearby exact root and bounds every exact root in the same neighborhood:

\begin{thm}\label{thm:surrogate_root_enclosure}
Let $U\subseteq\R^d$ be open, let $\bQ:U\to\R^d$ be continuous, and suppose that $\tilde\bQ(\tilde\bu)=\bzero$ for some map $\tilde\bQ:U\to\R^d$ and some $\tilde\bu\in U$. For nonsingular $\Rb,\Jb\in\R^{d\times d}$ and $r>0$, define
\[
    D\coloneqq\left\{\bu:\norm{\Rb(\bu-\tilde\bu)}\leq r\right\}\subseteq U,
\]
and suppose that $\tilde\bQ$ is continuously differentiable on an open neighborhood of $D$. Set
\[
    \kappa\coloneqq\sup_{\bu\in D}
    \norm{\Rb\Jb\inv\Jb_{\tilde\bQ}(\bu)\Rb\inv-\Ib_d}
    \dand
    \eta\coloneqq\sup_{\bu\in D}
    \norm{\Rb\Jb\inv\left[\bQ(\bu)-\tilde\bQ(\bu)\right]}.
\]
If $\kappa<1$ and $\eta<(1-\kappa)r$, then $\bQ$ has a zero in the interior of $D$, and every zero $\hat\bu\in D$ of $\bQ$ satisfies
\begin{equation}\label{eq:surrogate_root_enclosure}
    \norm{\Rb(\hat\bu-\tilde\bu)}
    \leq\frac{\eta}{1-\kappa}.
\end{equation}
\end{thm}

\begin{proof}
For $\bx\in\overline{B}_r(\bzero)$, set $\bu(\bx)\coloneqq\tilde\bu+\Rb\inv\bx$ and define
\[
    \bF(\bx)\coloneqq\Rb\Jb\inv\tilde\bQ(\bu(\bx))
    \dand
    \bG(\bx)\coloneqq\Rb\Jb\inv\bQ(\bu(\bx)).
\]
Since $\bF(\bzero)=\bzero$, the fundamental theorem of calculus gives
\begin{equation}\label{eq:surrogate_integral_identity}
    \bF(\bx)=\int_0^1
    \Rb\Jb\inv\Jb_{\tilde\bQ}(\bu(t\bx))\Rb\inv\bx\dif t.
\end{equation}
Therefore
\[
    \bF(\bx)^\intercal \bx\geq(1-\kappa)\norm{\bx}^2
    \dand
    \norm{\bG(\bx)-\bF(\bx)}\leq\eta.
\]
On $\norm{\bx}=r$, it follows that
\[
    \bG(\bx)^\intercal \bx
    \geq(1-\kappa)r^2-\eta r>0.
\]
This condition implies that $\bG$ has a zero in the interior of $\overline{B}_r(\bzero)$. Indeed, otherwise $\bx\mapsto-r\bG(\bx)/\norm{\bG(\bx)}$ would be a continuous map from $\overline{B}_r(\bzero)$ into itself. By the Brouwer fixed point theorem, this map would have a fixed point on the boundary of $\overline{B}_r(\bzero)$, making $\bG(\bx)^\intercal \bx < 0$, a contradiction. This proves the existence of a zero of $\bQ$ in the interior of $D$.

If $\hat\bu\in D$ is any such zero and $\hat\bx=\Rb(\hat\bu-\tilde\bu)$, \eqref{eq:surrogate_integral_identity} gives
\[
    \norm{\bF(\hat\bx)}\geq(1-\kappa)\norm{\hat\bx}.
\]
Since $\bG(\hat\bx)=\bzero$ and $\norm{\bF(\hat\bx)-\bG(\hat\bx)}\leq\eta$, \eqref{eq:surrogate_root_enclosure} follows.
\end{proof}

In statistical applications, the estimating equations are random, so the conditions of \cref{thm:surrogate_root_enclosure} must be verified on events whose probabilities tend to one. The following corollary casts this verification into stochastic localization and an error bound.

\begin{coro}\label{coro:stochastic_localization_transfer}
For each $n$, suppose that the objects $\bQ_n,\tilde\bQ_n,\tilde\bu_n,\Rb_n,\Jb_n,r_n$ (which may be random) satisfy the conditions of \cref{thm:surrogate_root_enclosure} with probability tending to one. Let $D_n$, $\kappa_n$, and $\eta_n$ be the corresponding quantities in that theorem. If $\kappa_n = \op{1}$ and $\eta_n/r_n = \op{1}$, then with probability tending to one, $\bQ_n$ has a zero in $D_n$, and every zero $\bu\in D_n$ satisfies \begin{equation}\label{eq:stochastic_root_enclosure} 
    \norm{\Rb_n(\bu-\tilde\bu_n)}
    \leq
    \frac{\eta_n}{1-\kappa_n}.
\end{equation}
Let $\hat\bu_n^{\mathrm{loc}}$ denote any selected zero in $D_n$, completed on the exceptional event according to the convention of \cref{sec:defs}. For any deterministic $s_n>0$, $\eta_n=\Op{s_n}$ implies $\Rb_n(\hat\bu_n^{\mathrm{loc}}-\tilde\bu_n)=\Op{s_n}$ and $\eta_n=\op{s_n}$ implies $\Rb_n(\hat\bu_n^{\mathrm{loc}}-\tilde\bu_n)=\op{s_n}$. If, with probability tending to one, $\bQ_n$ has the unique zero $\hat\bu_n$ in the relevant domain, then $\hat\bu_n^{\mathrm{loc}}=\hat\bu_n$ with probability tending to one.
\end{coro}

\begin{proof}
Let $A_n$ be an event on which the conditions of \cref{thm:surrogate_root_enclosure} hold, with $\mathbb P_{n,\btheta_n}(A_n)\to1$, and set
\[
    B_n
    \coloneqq
    A_n\cap
    \left\{
        \kappa_n<\frac12,\ 
        \eta_n<\frac{r_n}{2}
    \right\}.
\]
The assumptions $\kappa_n=\op{1}$ and $\eta_n/r_n=\op{1}$ imply that $\mathbb P_{n,\btheta_n}(B_n)\to1$. On $B_n$, \cref{thm:surrogate_root_enclosure} gives the existence of a zero in $D_n$ and establishes \eqref{eq:stochastic_root_enclosure} for every zero in $D_n$. The stated stochastic orders follow from $1/(1-\kappa_n)=1+\op{1}$, and under the uniqueness assumption, the selected local zero $\hat\bu_n^{\mathrm{loc}}$ must equal $\hat\bu_n$ on an event whose probability tends to one.
\end{proof}

An error in an estimating equation translates into an error in its root only after we account for the sensitivity of that root to perturbations. We make that comparison in coordinates of a vanishing chart, where our main result, \cref{thm:MAIN} below, bounds the difference between the exact and surrogate roots by the uniform equation error divided by a scale set by the smallest singular value:

\begin{thm}\label{thm:MAIN}
Suppose that there are deterministic convex sets $D_n\subseteq U$ and events $E_n$ with $\mathbb P_{n,\btheta_n}(E_n)\to1$ such that, on $E_n$,
\begin{equation}\label{eq:MAIN_assumptions}
\hattn,\tildetn\in N
\dand
\bq_n(\hattn)=\bzero
\dand
\tilde{\bq}_n(\tildetn)=\bzero
\dand
\hat{\bu}_n,\tilde{\bu}_n\in D_n,
\end{equation}
and $\tilde{\bQ}_n|_{D_n}$ admits a $C^1$ extension to an open neighborhood of $D_n$. Choose one such extension and denote its Jacobian by $\Jb_{\tilde{\bQ}_n}$. Suppose also that there are possibly random matrices $\Jb_n\in\R^{d\times d}$ and deterministic $a_n>0$ such that
\[
\mathbb P_{n,\btheta_n}\left(\sigmin(\Jb_n)\geq a_n\right)\longrightarrow1
\dand
\sup_{\bu\in D_n}
\norm{\Jb_{\tilde{\bQ}_n}(\bu)-\Jb_n}
=o_{\mathbb P_{n,\btheta_n}}(a_n)
\text{ on } E_n.
\]
Define
\[
\Delta_n\coloneqq
\sup_{\bu\in D_n}
\norm{\tilde{\bQ}_n(\bu)-\bQ_n(\bu)}.
\]
Then
\begin{equation}\label{eq:MAIN_conclusion}
\mathbb P_{n,\btheta_n}\left(
\norm{\tilde{\bu}_n-\hat{\bu}_n}
\leq\frac{2\Delta_n}{a_n}
\right)
\longrightarrow1.
\end{equation}
\end{thm}

\begin{proof}
Let
\[
H_n\coloneqq E_n
\cap\{\sigmin(\Jb_n)\geq a_n\}
\cap\left\{
\sup_{\bu\in D_n}
\norm{\Jb_{\tilde{\bQ}_n}(\bu)-\Jb_n}
\leq\frac{a_n}{2}
\right\}.
\]
Our assumptions imply $\mathbb P_{n,\btheta_n}(H_n)\to1$. On $H_n$, the segment $\{\hat{\bu}_n + t(\tilde{\bu}_n - \hat{\bu}_n) : t \in [0,1]\}$ lies in $D_n$ by convexity. Moreover, \eqref{eq:MAIN_assumptions} gives $\bQ_n(\hat{\bu}_n)=\bzero$ and $\tilde{\bQ}_n(\tilde{\bu}_n)=\bzero$. Applying the fundamental theorem of calculus to the selected extension of $\tilde{\bQ}_n|_{D_n}$ along this segment therefore gives
\begin{equation}\label{eq:integral_root_identity}
\Ab_n(\tilde{\bu}_n-\hat{\bu}_n)
=
-\left(\tilde{\bQ}_n(\hat{\bu}_n)-\bQ_n(\hat{\bu}_n)\right),
\quad \text{where} \quad
\Ab_n\coloneqq
\int_0^1
\Jb_{\tilde{\bQ}_n}
\bigl(\hat{\bu}_n+t(\tilde{\bu}_n-\hat{\bu}_n)\bigr)\dif t.
\end{equation}
We also have that
\[
\norm{\Ab_n-\Jb_n}
\leq
\sup_{\bu\in D_n}
\norm{\Jb_{\tilde{\bQ}_n}(\bu)-\Jb_n}
\leq\frac{a_n}{2}.
\]
Weyl's inequality now yields $\sigmin(\Ab_n)\geq\sigmin(\Jb_n)-\norm{\Ab_n-\Jb_n}\geq a_n/2$. Thus $\Ab_n$ is nonsingular on $H_n$, with $\norm{\Ab_n\inv}\leq2/a_n$, and so
\[
\norm{\tilde{\bu}_n-\hat{\bu}_n}
\leq
\norm{\Ab_n\inv}
\norm{\tilde{\bQ}_n(\hat{\bu}_n)-\bQ_n(\hat{\bu}_n)}
\leq\frac{2\Delta_n}{a_n}.
\]
Since this bound holds on $H_n$ and $\mathbb P_{n,\btheta_n}(H_n)\to1$, \eqref{eq:MAIN_conclusion} follows.
\end{proof}

To connect the bound in \cref{thm:MAIN} with \cref{def:approximate_est}, we must localize the transformed roots near zero and compare the root error $\Delta_n/a_n$ with the desired scale. The following corollary makes this connection explicit:

\begin{coro}
\label{coro:root_rate}
Under the assumptions of \cref{thm:MAIN}, define
\[
    r_n\coloneqq\sup_{\bu\in D_n}\norm{\bu}.
\]
If $r_n\to0$, then
\begin{equation}\label{eq:root_rate}
    \hat{\bu}_n=o_{\mathbb P_{n,\btheta_n}}(1)
    \dand
    \tilde{\bu}_n=o_{\mathbb P_{n,\btheta_n}}(1)
    \dand
    \norm{\tilde{\bu}_n-\hat{\bu}_n}
    =o_{\mathbb P_{n,\btheta_n}}(1).
\end{equation}
Therefore, the roots are equivalent at scale one relative to $\bphi$ along $\{\btheta_n\}_{n\geq1}$.

If, in addition, there exists a deterministic sequence $\rho_n\geq1$ such that
\begin{equation}\label{eq:Delta_n}
    \frac{\rho_n\Delta_n}{a_n} = o_{\mathbb P_{n,\btheta_n}}(1),
\end{equation}
then
\[
    \rho_n\norm{\tilde{\bu}_n-\hat{\bu}_n} = o_{\mathbb P_{n,\btheta_n}}(1).
\]
Thus the estimators are equivalent at scale $\rho_n$ relative to $\bphi$ along $\{\btheta_n\}_{n\geq1}$; if $\rho_n\to\infty$, they are asymptotically equivalent under the normalization $\rho_n$ along that sequence. If the assumptions of \cref{thm:MAIN}, the condition $r_n\to0$, and \eqref{eq:Delta_n} hold along every sequence $\{\btheta_n\}_{n\geq1} \in \mathcal R$ with the same deterministic sequence $\rho_n$ used throughout, then the corresponding equivalence holds throughout $\mathcal R$.
\end{coro}

\begin{proof}
On $E_n$, we have $\hat\bu_n, \tilde\bu_n \in D_n$, and hence
\[
    \norm{\hat{\bu}_n}\leq r_n
    \dand
    \norm{\tilde{\bu}_n}\leq r_n
    \dand
    \norm{\tilde{\bu}_n-\hat{\bu}_n}\leq2r_n.
\]
Because $\mathbb P_{n,\btheta_n}(E_n)\to1$ and $r_n\to0$, all three quantities above have the stochastic orders asserted in \eqref{eq:root_rate}. Moreover, $E_n \subseteq \{\hat\btheta_n, \tilde\btheta_n \in N\}$, so the domain requirement in \cref{def:approximate_est} is satisfied with probability tending to one. Together with \eqref{eq:root_rate}, this proves equivalence at chart scale one.

Now, let
\[
G_n\coloneqq\left\{
\norm{\tilde{\bu}_n-\hat{\bu}_n}
\leq\frac{2\Delta_n}{a_n}
\right\}.
\]
By \cref{thm:MAIN}, $\mathbb P_{n,\btheta_n}(G_n)\to1$. For every $\varepsilon>0$,
\[
\begin{aligned}
\mathbb P_{n,\btheta_n}\left(
\rho_n\norm{\tilde{\bu}_n-\hat{\bu}_n}>\varepsilon
\right)
&\leq
\mathbb P_{n,\btheta_n}(G_n^c)
+
\mathbb P_{n,\btheta_n}\left(
\frac{2\rho_n\Delta_n}{a_n}>\varepsilon
\right)
\longrightarrow0,
\end{aligned}
\]
where the convergence follows from \eqref{eq:Delta_n}; the remaining conclusions follow from \cref{def:approximate_est}. Since the regime sequence fixed at the beginning of the subsection was arbitrary, the same argument applies throughout $\mathcal R$ whenever all the stated conditions hold along every regime sequence and the same deterministic sequence $\rho_n$ is used throughout.
\end{proof}

The following matrix-scaled version of \cref{thm:MAIN} controls the root discrepancy after a general normalization $\Rb_n$, which we will need for joint asymptotic distributions.

\begin{prop}
\label{prop:matrix_root_stability}
Suppose that there are deterministic convex sets $D_n\subseteq U$ and events $E_n$ with $\mathbb P_{n,\btheta_n}(E_n)\to1$ such that, on $E_n$,
\[
\hattn,\tildetn\in N
\dand
\bq_n(\hattn)=\bzero
\dand
\tilde{\bq}_n(\tildetn)=\bzero
\dand
\hat{\bu}_n,\tilde{\bu}_n\in D_n,
\]
  and $\tilde{\bQ}_n|_{D_n}$ admits a $C^1$ extension to an open neighborhood of $D_n$; denote the Jacobian of this extension by $\Jb_{\tilde{\bQ}_n}$. Let $\Rb_n\in\R^{d\times d}$ be deterministic and nonsingular, and let $\Jb_n\in\R^{d\times d}$ be possibly random and nonsingular with probability tending to one. Let $B_n$ be the event that $\Jb_n$ is nonsingular, and define
\[
\Gamma_n\coloneqq
\sup_{\bu\in D_n}
\norm{
\Rb_n\Jb_n\inv
\left(\tilde{\bQ}_n(\bu)-\bQ_n(\bu)\right)}.
\]
Extend $\Gamma_n$ arbitrarily on $B_n^c$. Suppose that, on $E_n\cap B_n$, 
\begin{equation}\label{eq:thm1_matrixCLTcond}
\sup_{\bu\in D_n}
\norm{
\Rb_n\Jb_n\inv
\left(\Jb_{\tilde{\bQ}_n}(\bu)-\Jb_n\right)
\Rb_n\inv}
=o_{\mathbb P_{n,\btheta_n}}(1) \dand
\Gamma_n=o_{\mathbb P_{n,\btheta_n}}(1).
\end{equation}
Then
\begin{equation}\label{eq:matrix_root_stability}
\Rb_n(\tilde{\bu}_n-\hat{\bu}_n)
=o_{\mathbb P_{n,\btheta_n}}(1).
\end{equation}
More precisely,
\begin{equation}\label{eq:matrix_root_stability_bound}
\mathbb P_{n,\btheta_n}\left(
\norm{\Rb_n(\tilde{\bu}_n-\hat{\bu}_n)}
\leq2\Gamma_n
\right)\longrightarrow1.
\end{equation}
\end{prop}

\begin{proof}
All quantities below are considered on $F_n\coloneqq E_n\cap B_n$, which satisfies $\mathbb P_{n,\btheta_n}(F_n)\to1$ by assumption. On that event, $\bQ_n(\hat{\bu}_n)=\bzero$ and $\tilde{\bQ}_n(\tilde{\bu}_n)=\bzero$. Since the segment joining $\hat{\bu}_n$ and $\tilde{\bu}_n$ lies in $D_n$, applying the fundamental theorem of calculus to the selected extension of $\tilde{\bQ}_n$ gives 
\begin{equation}\label{eq:Abn_integral}
\Ab_n(\tilde{\bu}_n-\hat{\bu}_n)
=
-\left(\tilde{\bQ}_n(\hat{\bu}_n)-\bQ_n(\hat{\bu}_n)\right)
\dand
\Ab_n\coloneqq
\int_0^1
\Jb_{\tilde{\bQ}_n}
\bigl(\hat{\bu}_n+t(\tilde{\bu}_n-\hat{\bu}_n)\bigr)\dif t.
\end{equation}
Define
\[
\Bb_n\coloneqq
\Rb_n\Jb_n\inv(\Ab_n-\Jb_n)\Rb_n\inv.
\]
Again using the fact that the segment joining $\hat{\bu}_n$ and $\tilde{\bu}_n$ lies in $D_n$,
\[
\begin{aligned}
\norm{\Bb_n}
&\leq
\int_0^1
\norm{
\Rb_n\Jb_n\inv
\left[
\Jb_{\tilde{\bQ}_n}
\bigl(\hat{\bu}_n+t(\tilde{\bu}_n-\hat{\bu}_n)\bigr)-\Jb_n
\right]
\Rb_n\inv}
\dif t\\
&\leq
\sup_{\bu\in D_n}
\norm{
\Rb_n\Jb_n\inv
\left(\Jb_{\tilde{\bQ}_n}(\bu)-\Jb_n\right)
\Rb_n\inv} \\
&=o_{\mathbb P_{n,\btheta_n}}(1),
\end{aligned}
\]
so $\Bb_n=o_{\mathbb P_{n,\btheta_n}}(1)$. Define $G_n \coloneqq F_n\cap\{\norm{\Bb_n}<1/2\}$. A direct expansion gives 
\[
\Ib_d+\Bb_n = \Rb_n\Jb_n\inv\Ab_n\Rb_n\inv,
\]
so $\Ab_n$ is nonsingular on $G_n$, and inversion yields 
\[
    (\Ib_d+\Bb_n)\inv
    =
    \Rb_n\Ab_n\inv\Jb_n\Rb_n\inv,
\]
or equivalently,
\[
    \Rb_n
    =
    (\Ib_d+\Bb_n)\inv\Rb_n\Jb_n\inv \Ab_n.
\]
On $G_n$, combining this identity with \eqref{eq:Abn_integral} gives
\[
\Rb_n(\tilde{\bu}_n-\hat{\bu}_n)
=
-(\Ib_d+\Bb_n)\inv
\Rb_n\Jb_n\inv
\left(\tilde{\bQ}_n(\hat{\bu}_n)-\bQ_n(\hat{\bu}_n)\right).
\]
With $\Gamma_n$ as defined in the statement, note that
\[
\norm{(\Ib_d+\Bb_n)\inv}
=\norm{\sum_{k=0}^\infty(-\Bb_n)^k}
\leq\frac{1}{1-\norm{\Bb_n}}
\leq2
\]
on $G_n$. Since $\hat{\bu}_n\in D_n$ on $F_n$, it follows that
\[
\norm{\Rb_n(\tilde{\bu}_n-\hat{\bu}_n)}
\leq
\norm{(\Ib_d+\Bb_n)\inv}\Gamma_n
\leq2\Gamma_n
\]
on $G_n$. Since $\PP(G_n)\to1$, this also proves \eqref{eq:matrix_root_stability_bound}. Moreover, for every $\varepsilon>0$,
\[
\mathbb P_{n,\btheta_n}\left(
\norm{\Rb_n(\tilde{\bu}_n-\hat{\bu}_n)}>\varepsilon
\right)
\leq
\mathbb P_{n,\btheta_n}(G_n^c)
+
\mathbb P_{n,\btheta_n}\left(2\Gamma_n>\varepsilon\right)
\longrightarrow0,
\]
which proves \eqref{eq:matrix_root_stability}.
\end{proof}

\begin{rem}
Let $\bu_n\coloneqq\bphi(\btheta_n)$ for all sufficiently large $n$. If $\Rb_n(\hat{\bu}_n-\bu_n)\inD\bZ$ for any weak limit $\bZ$, then by Slutsky's theorem $\Rb_n(\tilde{\bu}_n-\bu_n)\inD\bZ$, and conversely, because $\Rb_n(\tilde{\bu}_n-\hat{\bu}_n) =o_{\mathbb P_{n,\btheta_n}}(1)$. The scalar case is obtained by taking $\Rb_n\coloneqq\rho_n\Ib_d$. We will use this observation repeatedly in \cref{sec:gamma,sec:additionalexamples} to transfer asymptotic distributions from exact to approximate estimators.
\end{rem}

\cref{thm:MAIN,coro:root_rate,prop:matrix_root_stability} control the discrepancy between the roots in the coordinates of the vanishing chart. However, this control need not directly translate to the original parameter scale, since $\bphi\inv$ may magnify small coordinate errors near the limiting value. At the cost of a Lipschitz assumption on $\bphi\inv$, the following corollary accounts for this issue:

\begin{coro}
\label{coro:original_scale}
Suppose that there are deterministic convex sets $K_n\subseteq U$ and events $F_n$ with $\mathbb P_{n,\btheta_n}(F_n)\to1$ such that
\[
\hattn,\tildetn\in N
\dand
\bq_n(\hattn)=\bzero
\dand
\tilde{\bq}_n(\tildetn)=\bzero
\dand
\hat{\bu}_n,\tilde{\bu}_n\in K_n
\]
 on $F_n$, and $\tilde{\bQ}_n|_{K_n}$ admits a $C^1$ extension to an open neighborhood of $K_n$. Choose one such extension and denote its Jacobian by $\Jb_{\tilde{\bQ}_n}$. Suppose that there are possibly random matrices $\Jb_n\in\R^{d\times d}$ and deterministic $a_n,b_n>0$ such that
\[
\mathbb P_{n,\btheta_n}\left(\sigmin(\Jb_n)\geq a_n\right)\longrightarrow1\dand
\sup_{\bu\in K_n}
\norm{\Jb_{\tilde{\bQ}_n}(\bu)-\Jb_n}
=o_{\mathbb P_{n,\btheta_n}}(a_n)
\]
 on $F_n$, and $\bphi\inv$ is $b_n$-Lipschitz on $K_n$. Define
\[
\Delta_{K,n}\coloneqq
\sup_{\bu\in K_n}\norm{\tilde{\bQ}_n(\bu)-\bQ_n(\bu)}.
\]
Let $\{\rho_n'\}_{n\geq1}$ be deterministic with $\rho_n'>0$ for all $n$. If
\[
\frac{\rho_n'b_n\Delta_{K,n}}{a_n}
=o_{\mathbb P_{n,\btheta_n}}(1),
\]
then
\[
\rho_n'\norm{\tildetn-\hattn}
=o_{\mathbb P_{n,\btheta_n}}(1).
\]
\end{coro}

\begin{proof}
The assumptions of \cref{thm:MAIN} hold with $K_n$ in place of $D_n$ and $\Delta_{K,n}$ in place of $\Delta_n$. Define
\[
G_n\coloneqq\left\{
\norm{\tilde{\bu}_n-\hat{\bu}_n}
\leq\frac{2\Delta_{K,n}}{a_n}
\right\}.
\]
By \cref{thm:MAIN}, $\mathbb P_{n,\btheta_n}(G_n)\to1$. On $F_n\cap G_n$, both roots belong to $K_n$, and the fact that $\bphi\inv$ is  $b_n$-Lipschitz gives
\[
\rho_n'\norm{\tildetn-\hattn}
\leq
\rho_n'b_n\norm{\tilde{\bu}_n-\hat{\bu}_n}
\leq
\frac{2\rho_n'b_n\Delta_{K,n}}{a_n}.
\]
For every $\varepsilon>0$,
\[
\begin{aligned}
\mathbb P_{n,\btheta_n}\left(
\rho_n'\norm{\tildetn-\hattn}>\varepsilon
\right)
&\leq
\mathbb P_{n,\btheta_n}\left(F_n^c\cup G_n^c\right)
+
\mathbb P_{n,\btheta_n}\left(
\frac{2\rho_n'b_n\Delta_{K,n}}{a_n}>\varepsilon
\right)
\longrightarrow0,
\end{aligned}
\]
because $\rho_n'b_n\Delta_{K,n}/a_n=o_{\mathbb P_{n,\btheta_n}}(1)$. This also shows that the arbitrary completion outside $F_n\cap G_n$ does not affect the conclusion.
\end{proof}

\begin{rem}
The effective approximation scale is $\Delta_n/a_n$, so the admissible normalization $\rho_n$ depends jointly on the sample size, the speed at which the fixed regime path $\btheta_n$ approaches $\btheta^*$, the localization set $D_n$, and the Jacobian lower bound scale $a_n$. While enlarging $D_n$ can accommodate slower approaches to the limiting value, it typically increases $\Delta_n$ and therefore reduces the largest normalization under which equivalence can be established.
\end{rem}

\subsection{A Construction and Verification Procedure}\label{sub:construction}

Our procedure for constructing and verifying an approximate estimator can be organized into the following three broad stages:

\begin{enumerate}
    \item \emph{Derive the approximate estimator.}

    Begin by specifying the scale on which the roots are to be compared --- a scalar or matrix normalization in chart coordinates, or the original parameter scale --- since this choice can determine the appropriate chart and the number of retained terms. Identify the selected exact root $\hattn$ of $\bq_n$ and fix a regime sequence $\{\btheta_n\}_{n\geq1}\in\mathcal R$; when $\bq_n$ is a score equation, verify separately that this root is the relevant MLE. Choose a vanishing chart $\bphi:N\to U$ and write $\bu_n\coloneqq\bphi(\btheta_n)$. Within the class of charts covered by \cref{prop:chart_invariance}, select a chart that simplifies the transformed equation and the required Jacobian bounds. When useful, rewrite or profile the estimating equation without changing the selected root, and replace its intractable terms by approximations to obtain a transformed surrogate equation $\tilde{\bQ}_n(\bu)=\bzero$ whose selected root $\tilde{\bu}_n$ is explicit. The entire system must be solvable, and it must be verified that $\tilde{\bu}_n\in U$ with probability tending to one before defining $\tildetn\coloneqq\bphi\inv(\tilde{\bu}_n)$. The procedure may be iterative: if the bounds required in the next stage fail, the surrogate or localization set must be refined.

    \item \emph{Verify the approximation.}

    Begin by localizing the explicit root $\tilde{\bu}_n$. If a suitably scaled ellipsoid centered at $\tilde{\bu}_n$ remains in $U$, \cref{coro:stochastic_localization_transfer} can help to supply the localization of $\hat\bu_n$ directly from the normalized Jacobian and equation errors; otherwise, localize the exact root separately. Then choose deterministic convex sets $D_n\subseteq U$ containing both roots with probability tending to one. To apply \cref{coro:root_rate}, these sets should also satisfy
    \[
        \sup_{\bu\in D_n}\norm{\bu}\longrightarrow0.
    \]
    Next, ensure that $\tilde{\bQ}_n$ admits a $C^1$ extension to an open neighborhood of $D_n$, and find matrices $\Jb_n$ and deterministic scales $a_n>0$ such that
    \[
        \mathbb P_{n,\btheta_n}\!\left(\sigmin(\Jb_n)\geq a_n\right)\longrightarrow1 \dand
        \sup_{\bu\in D_n}
        \norm{\Jb_{\tilde{\bQ}_n}(\bu)-\Jb_n}
        =\op{a_n}.
    \]
    Define the error
    \[
        \Delta_n\coloneqq
        \sup_{\bu\in D_n}
        \norm{\tilde{\bQ}_n(\bu)-\bQ_n(\bu)}.
    \]
    Any random or diverging factor multiplying a truncated term must be included in this error. By \cref{thm:MAIN,coro:root_rate}, equivalence at a scalar chart scale $\rho_n$ follows if
    \[
        \frac{\rho_n\Delta_n}{a_n}=\op{1}.
    \]
    For a general matrix normalization, use \cref{prop:matrix_root_stability}. If comparison is required on the original parameter scale and $\bphi\inv$ is $b_n$-Lipschitz on the localization set, then by \cref{coro:original_scale}, the corresponding requirement is $b_n\Delta_n/a_n=\op{1}$.

    \item \emph{Obtain the asymptotic distribution, when required.}

    Derive the limiting distribution of whichever estimator is more tractable (usually, the explicit approximate estimator). Because the data-generating parameter $\btheta_n$ varies with $n$, generally so does the distribution of the observations. When the data are independent, the relevant central limit theorem should therefore be formulated for a triangular array, even when the observations are identically distributed within each row. After appropriate centering and normalization, verify a Lindeberg or Lyapunov condition. For example, for centered scalar summands $W_{n,i}$, a Lyapunov condition takes the form
    \[
        \frac{\sum_i\E_{n,\btheta_n}\left[|W_{n,i}|^{2+\delta}\right]}
        {\left(\sum_i\Var_{n,\btheta_n}(W_{n,i})\right)^{1+\frac{\delta}{2}}}
        \longrightarrow0
    \]
    for some $\delta>0$, where $\E_{n,\btheta_n}$ and $\Var_{n,\btheta_n}$ denote expectation and variance under $\mathbb P_{n,\btheta_n}$. In row-wise iid problems, a uniform bound on the corresponding standardized $(2+\delta)$th moments is usually sufficient. Corresponding limits for random vectors may be obtained from the multivariate Lindeberg--Feller central limit theorem \citep[Theorem 2.27]{vandervaartAsymptoticStatistics1998} or by the Cramér--Wold device. A Taylor expansion or the delta method can then transfer the limit from the underlying sample statistics to the explicit estimator.

    Once the limiting distribution of one root is established, verify that a deterministic nonsingular normalization $\Rb_n$ satisfies
    \[
        \Rb_n(\tilde{\bu}_n-\hat{\bu}_n)=\op{1}.
    \]
    Slutsky's theorem then transfers the weak limit between the two transformed roots; see the remark following \cref{prop:matrix_root_stability}. A limiting distribution on the original parameter scale may additionally require a separate expansion or the delta method for $\bphi\inv$ about $\bu_n$.
\end{enumerate}

\section{Worked Examples: the Gamma Distribution}\label{sec:gamma}

We begin our examples with the gamma distribution, which serves as a useful test case for our general framework. The distribution is simple enough that the likelihood and score functions can be written out explicitly, but the score function for the shape parameter involves the digamma function and therefore does not admit a root in closed form. The first two examples consider the regime where the gamma shape tends to infinity, while gamma regression considers the equivalent regime in terms of the precision. \cref{sub:fixed_rate} treats a known rate and gives a one-dimensional version of the construction. \cref{sub:unknown_rate} estimates shape and rate jointly, producing a two-dimensional chart and different componentwise rates. \cref{sub:glmregression} extends the construction to independent, non-identically distributed responses. \cref{sub:gamma_simulation_study} reports simulation studies for all three gamma settings, and \cref{sub:gamma_tree_data} compares gamma regression fits based on the exact MLE and our proposed closed-form estimator in a real-world data setting.

\subsection{Fixed Rate Parameter}\label{sub:fixed_rate}

To begin with, let $X_1,\ldots,X_n\iid\mathrm{Gamma}(\alpha,\beta)$ in the shape and rate parameterization, where $\beta>0$ is known. We consider the regime $\alpha\to\infty$ with $\beta$ fixed. For $x>0$, the density is
\[
    f_{\alpha,\beta}(x)
    \coloneqq
    \frac{\beta^\alpha}{\Gamma(\alpha)}
    x^{\alpha-1}e^{-\beta x}.
\]
Write $G_n\coloneqq(\prod_{i=1}^nX_i)^{1/n}$ for the sample geometric mean of the $X_i$.

\subsubsection{Derivation of the Approximate MLE}\label{subsub:fixed_rate_approx_mle}

The score equation for the shape parameter is
\[
    \psi(\alpha) = \log{\beta G_n},
\]
where $\psi\coloneqq\Gamma'/\Gamma$. Because $\psi$ is strictly increasing from $-\infty$ to $+\infty$, the equation has a unique solution $\hat\alpha_n$. Moreover, the log-likelihood $\ell_n$ is strictly concave in $\alpha$ since $\ell''_n(\alpha) = -n \psi_1(\alpha) < 0$, so $\hat\alpha_n$ is indeed the unique maximizer of $\ell_n$. Equivalently, $\hat\alpha_n$ may be defined as the root of the transformed equation
\[
    q_n(\alpha) \coloneqq \exp{-\psi(\alpha)} - \frac{1}{\beta G_n},
\]
which is convenient to use as $\alpha\to\infty$; indeed, we have
\begin{equation}\label{eq:psi_asympt}
    \psi(x)
    =
    \log{x}-\frac{1}{2x}-\frac{1}{12x^2}+O\left(\frac{1}{x^4}\right) \dand
    \exp{-\psi(x)}
    =
    \frac{1}{x-1/2}+O\left(\frac{1}{x^3}\right)
\end{equation}
as $x \to \infty$ by the standard digamma expansion \citep[Equation 6.3.18]{abramowitz1988handbook}. This suggests the vanishing chart $\phi:N\to U$, where $N\coloneqq(1/2,\infty)$ and $U\coloneqq(0,\infty)$, defined by
\[
    \phi(\alpha)\coloneqq\frac{1}{\alpha-1/2}
    \dand
    u\coloneqq\phi(\alpha)
    \dand
    \phi\inv(u)=\frac12+\frac1u,
\]
under which the exact and surrogate equations on the chart scale are
\[
    Q_n(u)
    \coloneqq
    \exp{-\psi\left(\frac12+\frac1u\right)}
    -\frac{1}{\beta G_n}
    \dand
    \tilde Q_n(u)\coloneqq u-\frac{1}{\beta G_n}.
\]
The latter function has the explicit root $\tilde u_n\coloneqq1/(\beta G_n)$, which gives the candidate approximate MLE
\[
    \tilde\alpha_n\coloneqq\phi\inv(\tilde u_n)=\frac12+\beta G_n.
\]

\begin{rem}
    In practice, of course, $\alpha$ is unknown; it is therefore useful to ask whether the observed data itself provides evidence that $\alpha$ exceeds a certain prescribed benchmark $a > 0$, chosen such that the corresponding approximate estimator is sufficiently accurate. When the rate parameter $\beta$ is known, such a check is immediate and can be formalized as a hypothesis test: consider testing $H_0: \alpha \leq a$ versus $H_A: \alpha > a$. For any $\delta \in (0,1)$, rejecting $H_0$ whenever
\[
    \beta \overline{X}_n > a + \sqrt{\frac{a(1-\delta)}{n\delta}}
\]
    defines a level-$\delta$ test, which follows from stochastic monotonicity in $\alpha$, Cantelli's inequality, and the fact that $\beta \overline{X}_n \sim \text{Gamma}(n \alpha, n)$. Similar tests for other distributions can be developed by identifying a relevant test statistic which is ancillary for all other unknown parameters in the model, such as a scale-invariant statistic when $\beta$ is unknown in the $\text{Gamma}(\alpha, \beta)$ setting; since such tests are not a focus of our work, we do not pursue them further in this paper.
\end{rem}

\subsubsection{Verification of Approximation}

Fix a deterministic sequence $\{\alpha_n\}_{n\geq1}$ satisfying $\alpha_n\to\infty$, and write
\begin{equation}\label{eq:u_n_def}
    u_n\coloneqq\phi(\alpha_n)=\frac{1}{\alpha_n-1/2}
    \dand
    \hat u_n\coloneqq\phi(\hat\alpha_n).
\end{equation}
By \cref{lemm:gamma1_localization} in \cref{app:gammadistresults},
\[
    \tilde u_n-u_n
    =
    \Op{\sqrt{\frac{u_n^3}{n}}}
    +O\left(\frac{u_n^2}{n}\right)
    +O(u_n^3)
    =
    \op{u_n}.
\]
On the event $|\tilde u_n/u_n-1|\leq1/4$, define
\[
    D_n^*\coloneqq\left\{u:|u-\tilde u_n|\leq\frac{u_n}{4}\right\}.
\]
This event has probability tending to one, and $D_n^*\subset[u_n/2,3u_n/2]$; thus, we may apply \cref{coro:stochastic_localization_transfer} with $\Rb_n\coloneqq u_n^{-1}$, $\Jb_n\coloneqq1$, and $r_n\coloneqq1/4$. Since $J_{\tilde Q_n}(u)=1$, we have $\kappa_n=0$, while
\[
    Q_n(u)-\tilde Q_n(u)
    =\exp{-\psi\left(\frac12+\frac1u\right)}-u
    =O(u^3)
\]
as $u\downarrow0$ gives
\[
    \eta_n
    =\sup_{u\in D_n^*}u_n^{-1}|Q_n(u)-\tilde Q_n(u)|
    =O(u_n^2)=o(1).
\]
The exact equation has a unique root, so \cref{coro:stochastic_localization_transfer} gives
\begin{equation}\label{eq:tilde_u_hat_u}
    |\tilde u_n-\hat u_n|=\Op{u_n^3}
    \dand
    \hat u_n-u_n=\op{u_n}.
\end{equation}
Consequently, $\tilde u_n/u_n,\hat u_n/u_n\inP1$, and the deterministic convex sets
\[
    D_n\coloneqq\left[\frac{u_n}{2},\,2u_n\right]
\]
therefore contain both roots with probability tending to one, and
$\sup_{u\in D_n}|u|=2u_n\to0$. Taking $J_n\coloneqq1$ and $a_n\coloneqq1$, we have
\[
    \Delta_n
    \coloneqq
    \sup_{u\in D_n}|Q_n(u)-\tilde Q_n(u)|
    =
    O(u_n^3).
\]
Thus the roots are equivalent at scale 1 relative to the $\phi$-chart (see \cref{def:approximate_est}). More generally, if $\rho_n\to\infty$ and $\rho_nu_n^3\to0$, then
\[
    \rho_n|\tilde u_n-\hat u_n|=\op{1},
\]
and the estimators are asymptotically equivalent under the normalization $\rho_n$ in the chart. For comparison on the original parameter scale, observe that, for $u_1,u_2\in D_n$,
\[
    |\phi\inv(u_1)-\phi\inv(u_2)|
    =
    \frac{|u_1-u_2|}{u_1u_2}
    \leq
    \frac{4}{u_n^2}|u_1-u_2|.
\]
Thus $\phi\inv$ is $b_n$-Lipschitz on $D_n$ with $b_n\coloneqq4/u_n^2$, and
\[
    \frac{b_n\Delta_n}{a_n}
    =
    O(u_n)
    \longrightarrow0.
\]
It follows from \cref{coro:original_scale} that
\[
    |\tilde\alpha_n-\hat\alpha_n|= \Op{\frac{b_n \Delta_n}{a_n}} = \Op{u_n} = \Op{\frac{1}{\alpha_n}} = \op{1},
\]
which confirms that $\tilde\alpha_n$ is an approximate MLE on the original parameter scale.

\subsubsection{Asymptotic Distribution}

Suppose now that $nu_n^3\longrightarrow0$, or equivalently, $\alpha_n^3/n\to\infty$. For each $n$, let $X_{n,1},\dots,X_{n,n} \iid \mathrm{Gamma}(\alpha_n,\beta)$, and define
\[
    Y_{n,i}\coloneqq\log{X_{n,i}}
    \dand
    \bar Y_n\coloneqq\frac1n\sum_{i=1}^nY_{n,i}.
\]
The third moment bound in \cref{lemm:gamma1_lyapunov} of \cref{app:gammadistresults} gives the Lyapunov triangular array CLT
\[
    \frac{\sqrt n\left[\bar Y_n-(\psi(\alpha_n)-\log{\beta})\right]}
    {\sqrt{\psi_1(\alpha_n)}}
    \inD
    \mathcal N(0,1).
\]
For ease of notation, define
\[
    Z_n\coloneqq\bar Y_n-\left(\psi(\alpha_n)-\log{\beta}\right)
    \dand
    c_n\coloneqq e^{-\psi(\alpha_n)}.
\]
Since
\[
    \psi_1(\alpha_n)=u_n+O(u_n^3),
\]
the preceding CLT implies
\[
    \sqrt{\frac n{u_n}}\,Z_n\inD\mathcal N(0,1).
\]
Moreover, $\tilde u_n=c_ne^{-Z_n}$ and $c_n=u_n+O(u_n^3)$. Because $Z_n=\Op{\sqrt{u_n/n}}$, an expansion of $e^{-Z_n}$ gives
\[
    \tilde u_n-u_n
    =
    -c_nZ_n+(c_n-u_n)+\Op{c_nZ_n^2}.
\]
Under the given regime,
\[
    \sqrt{\frac n{u_n^3}}(c_n-u_n)
    =
    O(\sqrt{nu_n^3})
    \longrightarrow0
    \dand
    \sqrt{\frac n{u_n^3}}\,c_nZ_n^2
    =
    \Op{\sqrt{\frac{u_n}{n}}}
    =
    \op{1}.
\]
Since $c_n/u_n\to1$, it follows that
\[
    \sqrt{\frac n{u_n^3}}
    (\tilde u_n-u_n)
    \inD
    \mathcal N(0,1).
\]
We have $|\tilde u_n-\hat u_n|=\Op{u_n^3}$ from \eqref{eq:tilde_u_hat_u}. Hence
\[
    \sqrt{\frac n{u_n^3}}
    (\tilde u_n-\hat u_n)
    =
    \Op{\sqrt{nu_n^3}}
    =
    \op{1},
\]
and Slutsky's theorem yields
\[
    \sqrt{\frac n{u_n^3}}
    (\hat u_n-u_n)
    \inD
    \mathcal N(0,1).
\]
Finally,
\cref{lemm:gamma1_inverse_chart} gives
\[
    \sqrt{nu_n}\left(\phi\inv(\tilde u_n)-\alpha_n\right)
    =
    -\sqrt{\frac n{u_n^3}}(\tilde u_n-u_n)+\op{1}
    \dand
    \sqrt{nu_n}\left(\phi\inv(\hat u_n)-\alpha_n\right)
    =
    -\sqrt{\frac n{u_n^3}}(\hat u_n-u_n)+\op{1}.
\]
Since $u_n=1/\alpha_n+O(\alpha_n^{-2})$, we conclude that
\[
    \sqrt{\frac n{\alpha_n}}
    (\tilde\alpha_n-\alpha_n)
    \inD\mathcal N(0,1)
    \dand
    \sqrt{\frac n{\alpha_n}}
    (\hat\alpha_n-\alpha_n)
    \inD\mathcal N(0,1).
\]

\subsection{Unknown Rate Parameter}\label{sub:unknown_rate}

We now suppose that $X_1,\ldots,X_n\iid\mathrm{Gamma}(\alpha,\beta)$ with both $\alpha$ and $\beta$ unknown. We again consider the regime $\alpha\to\infty$, with the true rate $\beta>0$ fixed along the regime sequence, and write
\[
    \bar X_n\coloneqq\frac1n\sum_{i=1}^nX_i
    \dand
    G_n\coloneqq\left(\prod_{i=1}^nX_i\right)^{\frac{1}{n}}.
\]

\subsubsection{Derivation of the Approximate MLE}

The score vector is
\[
    \bS_n(\alpha,\beta)
    \coloneqq
    \begin{pmatrix}
        \dfrac{\partial}{\partial\alpha}\ell_n(\alpha,\beta)\\[9pt]
        \dfrac{\partial}{\partial\beta}\ell_n(\alpha,\beta)
    \end{pmatrix}
    =
    \begin{pmatrix}
        n\left(\log{\beta}+\log{G_n}-\psi(\alpha)\right)\\[9pt]
        n\left(\dfrac{\alpha}{\beta}-\bar X_n\right)
    \end{pmatrix}.
\]
For each $\alpha>0$, the score equation for $\beta$, 
\[
    \frac{\partial}{\partial\beta}\ell_n(\alpha,\beta)
    =
    n\left(\dfrac{\alpha}{\beta}-\bar X_n\right) = 0,
\]
has the unique solution $\beta=\alpha/\bar X_n$, which maximizes the likelihood over $\beta$. Substituting this into the score equation for $\alpha$ gives the profiled score equation
\begin{equation}\label{eq:profilescore}
    \log{\hat\alpha_n}-\psi(\hat\alpha_n)
    =
    \log{\frac{\bar X_n}{G_n}},
\end{equation}
with $\hat\beta_n=\hat\alpha_n/\bar X_n$. For $n\geq2$, the AM--GM inequality gives $\log{\bar X_n/G_n}>0$ almost surely (we use the completion convention of \cref{sec:defs} on the null event of a constant sample). Moreover, the function
\[
    g(\alpha)\coloneqq\log{\alpha}-\psi(\alpha)
\]
is strictly decreasing from $\infty$ to zero, so the profiled score equation \eqref{eq:profilescore} changes sign exactly once. Its unique root $\hat\alpha_n$, together with $\hat\beta_n=\hat\alpha_n/\bar X_n$, is therefore the unique joint MLE.

An equivalent system defining the exact root is
\[
    \bq_n(\alpha,\beta)
    \coloneqq
    \begin{pmatrix}
        \log{\alpha}-\psi(\alpha)-\log{\dfrac{\bar X_n}{G_n}}\\[9pt]
        \dfrac{\beta}{\alpha}-\dfrac{1}{\bar X_n}
    \end{pmatrix}
    =
    \bzero.
\]
By \eqref{eq:psi_asympt}, we have
\[
    \log{\alpha}-\psi(\alpha)
    =
    \frac{1}{2\alpha}
    +\frac{1}{12\alpha^2}
    +O\left(\frac{1}{\alpha^4}\right).
\]
We retain the first two terms of the expansion and use the surrogate
\[
    \tilde\bq_n(\alpha,\beta)
    \coloneqq
    \begin{pmatrix}
        \dfrac{1}{(2\alpha)}+\dfrac{1}{(12\alpha^2)}-\log{\dfrac{\bar X_n}{G_n}}\\[9pt]
        \dfrac{\beta}{\alpha}-\dfrac{1}{\bar X_n}
    \end{pmatrix},
\]
which suggests the chart
\begin{equation}\label{eq:gamma2_phi}
    \bphi(\alpha,\beta)
    \coloneqq
    \begin{pmatrix}
        \dfrac{1}{(2\alpha)}+\dfrac{1}{(12\alpha^2)}\\[9pt]
        \dfrac{\beta}{\alpha}
    \end{pmatrix}
    \dand
    \bu\coloneqq\bphi(\alpha,\beta)
    =\begin{pmatrix}u_1\\u_2\end{pmatrix}.
\end{equation}
This is a bijection from $(0,\infty)^2$ to $U\coloneqq(0,\infty)^2$, with inverse
\begin{equation}\label{eq:gamma2_phiinv}
    \bphi\inv(\bu)
    =
    \begin{pmatrix}
        \dfrac{1}{-3+\sqrt{9+12u_1}}\\[9pt]
        \dfrac{u_2}{-3+\sqrt{9+12u_1}}
    \end{pmatrix}.
\end{equation}
For fixed $\beta>0$, both chart coordinates tend to zero as $\alpha\to\infty$, so $\bphi$ is a vanishing chart for the regime under consideration. Define
\[
    h(u_1)
    \coloneqq
    \log{\frac{1}{-3+\sqrt{9+12u_1}}}
    -
    \psi\left(\frac{1}{-3+\sqrt{9+12u_1}}\right).
\]
The exact and surrogate equations in chart coordinates are, respectively, 
\[
    \bQ_n(\bu)
    \coloneqq
    \begin{pmatrix}
        h(u_1)-\log{\dfrac{\bar X_n}{G_n}}\\[9pt]
        u_2-\dfrac{1}{\bar X_n}
    \end{pmatrix}
    \dand
    \tilde\bQ_n(\bu)
    \coloneqq
    \begin{pmatrix}
        u_1-\log{\dfrac{\bar X_n}{G_n}}\\[9pt]
        u_2-\dfrac{1}{\bar X_n}
    \end{pmatrix},
\]
and the surrogate has the explicit root
\[
    \tilde\bu_n
    \coloneqq
    \begin{pmatrix}
        \log{\dfrac{\bar X_n}{G_n}}\\[9pt]
        \dfrac{1}{\bar X_n}
    \end{pmatrix}.
\]
Because $\bar X_n>G_n>0$ almost surely for $n\geq2$, we have $\tilde\bu_n \in U$ almost surely. The resulting candidate approximate MLE is
\[
    \tilde\alpha_n
    \coloneqq
    \frac{1}{-3+\sqrt{9+12\log{\bar X_n/G_n}}}
    \dand
    \tilde\beta_n
    \coloneqq
    \frac{\tilde\alpha_n}{\bar X_n}.
\]

\begin{rem}
    \begin{enumerate}
        \item It is essential to retain the term $1/(12\alpha^2)$ in \eqref{eq:gamma2_phi} for approximation on the original parameter scale. If only $1/(2\alpha)$ were retained, the neglected term would have order $u_1^2$, while the inverse chart has Lipschitz constant $u_1^{-2}$; the resulting amplified error would be only $O(1)$. Retaining the second-order term reduces the equation error to $O(u_1^4)$ and the amplified error to $O(u_1^2)$.
        \item The formula $1/(-3+\sqrt{9+12 u})$ suffers from cancellation as $u \to 0$, which is the regime of interest for the vanishing chart. In our simulations, we use the more numerically stable equivalent form $(3+\sqrt{9+12 u})/(12u)$.
    \end{enumerate}

\end{rem}

\subsubsection{Verification of Approximation}

For notational convenience, write $\btheta\coloneqq(\alpha,\beta)^\intercal$, $\hat\btheta_n\coloneqq(\hat\alpha_n,\hat\beta_n)^\intercal$, and $\tilde\btheta_n\coloneqq(\tilde\alpha_n,\tilde\beta_n)^\intercal$. Fix $\beta>0$ and a deterministic sequence $\{\alpha_n\}_{n\geq1}$ satisfying $\alpha_n\to\infty$. For each $n$, let $X_{n,1},\ldots,X_{n,n} \iid\mathrm{Gamma}(\alpha_n,\beta)$, with $\bar X_n$ and $G_n$ defined row-wise as before, and write
\[
    \bu_n
    \coloneqq
    \bphi(\alpha_n,\beta)
    =
    \begin{pmatrix}u_{1,n}\\u_{2,n}\end{pmatrix}
    =
    \begin{pmatrix}
        \dfrac{1}{(2\alpha_n)}+\dfrac{1}{(12\alpha_n^2)}\\[9pt]
        \dfrac{\beta}{\alpha_n}
    \end{pmatrix}.
\]
In particular, $u_{2,n}\sim2\beta u_{1,n}$. The calculations in \cref{lemm:gamma2_approx_root} of
\cref{app:gammadistresults} give the coordinatewise bounds
\[
\begin{aligned}
    \tilde u_{1,n}-u_{1,n}
    &=
    \Op{\frac{u_{1,n}}{\sqrt n}}
    +O\left(\frac{u_{1,n}}n\right)
    +O(u_{1,n}^4),\\
    \tilde u_{2,n}-u_{2,n}
    &=
    \Op{\frac{u_{1,n}^{3/2}}{\sqrt n}}
    +O\left(\frac{u_{1,n}^2}{n}\right),
\end{aligned}
\]
and so $\tilde\bu_n=\bu_n+\op{u_{1,n}}$. Choose $r_0>0$ so that $\overline{B}_{r_0}((1,2\beta)^\intercal)\subseteq(0,\infty)^2$. On the event that $u_{1,n}^{-1}\tilde\bu_n\in\overline{B}_{r_0/2}((1,2\beta)^\intercal)$, define
\[
    D_n^*\coloneqq\left\{\bu:
    \norm{\frac{1}{u_{1,n}}(\bu-\tilde\bu_n)}\leq\frac{r_0}{2}\right\}.
\]
This event has probability tending to one and $D_n^*\subset U$. Apply \cref{coro:stochastic_localization_transfer} with $\Rb_n\coloneqq u_{1,n}^{-1}\Ib_2$, $\Jb_n\coloneqq\Ib_2$, and $r_n\coloneqq r_0/2$. Since $\Jb_{\tilde\bQ_n}=\Ib_2$, we have $\kappa_n=0$. The two equations differ only in their first coordinates, and $h(u_1)-u_1=O(u_1^4)$ as $u_1\downarrow0$; hence
\[
    \eta_n
    =\sup_{\bu\in D_n^*}
    \norm{\Rb_n[\bQ_n(\bu)-\tilde\bQ_n(\bu)]}
    =O(u_{1,n}^3)=o(1).
\]
Because the joint MLE is unique, the equation $\bQ_n(\bu)=\bzero$ has the unique root $\hat\bu_n\coloneqq\bphi(\hat\alpha_n,\hat\beta_n)$. Hence, \cref{coro:stochastic_localization_transfer} gives
\begin{equation}\label{eq:gamma2_chart_root_difference}
    \norm{\tilde\bu_n-\hat\bu_n}=\Op{u_{1,n}^4}\dand
    \hat\bu_n=\bu_n+\op{u_{1,n}}.
\end{equation}
Since
$u_{2,n}\asymp u_{1,n}$, these relations imply, for $j=1,2$,
\[
    \frac{\hat u_{j,n}}{u_{j,n}}\inP1
    \dand
    \frac{\tilde u_{j,n}}{u_{j,n}}\inP1.
\]
It follows that the deterministic convex sets
\[
    D_n
    \coloneqq
    \left[\frac{u_{1,n}}2,\,2u_{1,n}\right]
    \times
    \left[\frac{u_{2,n}}2,\,2u_{2,n}\right]
\]
contain both transformed roots with probability tending to one. Moreover,
\[
    \sup_{\bu\in D_n}\norm{\bu}
    =
    O(u_{1,n})
    \longrightarrow0.
\]
Because $\Jb_{\tilde\bQ_n}=\Ib_2$, we may take $\Jb_n\coloneqq\Ib_2$ and $a_n\coloneqq1$. Since the two equations differ only in their first coordinates and $h(u_1)-u_1=O(u_1^4)$ as $u_1\downarrow0$, we have
\[
    \Delta_n
    \coloneqq
    \sup_{\bu\in D_n}
    \norm{\bQ_n(\bu)-\tilde\bQ_n(\bu)}
    =
    O(u_{1,n}^4).
\]
Now, \eqref{eq:gamma2_chart_root_difference} already establishes approximation relative to the vanishing chart. On the original parameter scale, \cref{lemm:gamma2_inverse_chart} gives
\[
    \sup_{\bu\in D_n}
    \norm{\Jb_{\bphi\inv}(\bu)}
    =
    O\left(\frac{1}{u_{1,n}^2}\right),
\]
so that $\bphi\inv$ is $b_n$-Lipschitz on $D_n$ for some
$b_n=O(u_{1,n}^{-2})$, and
\[
    \frac{b_n\Delta_n}{a_n}
    =
    O(u_{1,n}^2)
    \longrightarrow0.
\]
It follows from \cref{coro:original_scale} that
\[
    \norm{\tilde\btheta_n-\hat\btheta_n}
    =
    \op{1}.
\]
Hence $(\tilde\alpha_n,\tilde\beta_n)$ is an approximate MLE both relative to the chart and on the original parameter scale.

\subsubsection{Joint Asymptotic Distribution}

We assume henceforth that $\sqrt n\,u_{1,n}^3\longrightarrow0$, or equivalently, that $\alpha_n^3/\sqrt n\to\infty$. This condition makes the $O(u_{1,n}^4)$ deterministic bias negligible on the natural scale of the first coordinate. Define
\[
    \Rb_n
    \coloneqq
    \begin{pmatrix}
        \dfrac{\sqrt n}{u_{1,n}}&0\\[9pt]
        0&\dfrac{\sqrt n}{u_{1,n}^{3/2}}
    \end{pmatrix}.
\]
By \cref{lemm:gamma2_approx_root}, we have
\[
    \frac{\sqrt n}{u_{1,n}}
    (\tilde u_{1,n}-u_{1,n})
    \inD\mathcal N(0,2)
    \dand
    \frac{\sqrt n}{u_{1,n}^{3/2}}
    (\tilde u_{2,n}-u_{2,n})
    \inD\mathcal N(0,8\beta^2).
\]
To identify their joint limit, let
\[
    W_n\coloneqq\sum_{i=1}^nX_{n,i}
    \dand
    P_{n,i}\coloneqq\frac{X_{n,i}}{W_n}.
\]
By the standard gamma--Dirichlet decomposition \citep[e.g.,][Chapter~49]{kotzContinuousMultivariateDistributions2000}, $W_n$ is independent of $\bP_n\coloneqq(P_{n,1},\ldots,P_{n,n})$. Moreover,
\[
    \tilde u_{1,n}
    =
    -\log{n}-\frac1n\sum_{i=1}^n\log{P_{n,i}}
    \dand
    \tilde u_{2,n}
    =
    \frac n{W_n}.
\]
Thus $\tilde u_{1,n}$ and $\tilde u_{2,n}$ are independent for every $n$. Their marginal limits therefore combine to give
\begin{equation}\label{eq:gamma2_tildeujoint}
    \Rb_n(\tilde\bu_n-\bu_n)
    \inD
    \mathcal N\left(
        \bzero,
        \begin{pmatrix}
            2&0\\
            0&8\beta^2
        \end{pmatrix}
    \right).
\end{equation}

We transfer this limit to the exact root. Taking $\Jb_n=\Ib_2$, the normalized Jacobian error in \cref{prop:matrix_root_stability} is zero. Because only the first equation is approximated,
\[
    \tilde\bQ_n(\bu)-\bQ_n(\bu)
    =
    \begin{pmatrix}
        u_1-h(u_1)\\
        0
    \end{pmatrix}.
\]
On the sets $D_n$ from the preceding subsection,
\[
    \sup_{\bu\in D_n}
    \norm{
        \Rb_n\left(\tilde\bQ_n(\bu)-\bQ_n(\bu)\right)
    }
    =
    O(\sqrt n\,u_{1,n}^3)
    =
    o(1).
\]
Thus, \cref{prop:matrix_root_stability} gives
\[
    \Rb_n(\tilde\bu_n-\hat\bu_n)=\op{1}.
\]
By Slutsky's theorem, we therefore have
\[
    \Rb_n(\hat\bu_n-\bu_n)
    \inD
    \mathcal N\left(
        \bzero,
        \begin{pmatrix}
            2&0\\
            0&8\beta^2
        \end{pmatrix}
    \right).
\]
It remains to transform these limits back to the original parameter scale. Let $\bF\coloneqq\bphi\inv$ and define
\[
    \Tb_n
    \coloneqq
    \begin{pmatrix}
        \sqrt n\,u_{1,n}&0\\
        0&\sqrt n
    \end{pmatrix}.
\]
By \cref{lemm:gamma2_inverse_chart}\ref{lemm:gamma2_inverse_chart_scalings},
\[
    \Tb_n\Jb_{\bF}(\bu_n)\Rb_n\inv
    \longrightarrow
    \Ab
    \coloneqq
    \begin{pmatrix}
        -\dfrac{1}{2}&0\\[9pt]
        -\beta&0
    \end{pmatrix}.
\]
Moreover, \cref{lemm:gamma2_inverse_chart}\ref{lemm:gamma2_inverse_chart_remainder} gives
\[
    \Tb_n
    \left[\bF(\tilde\bu_n)-\bF(\bu_n)
      -\Jb_{\bF}(\bu_n)(\tilde\bu_n-\bu_n)\right]
    =
    \op{1}
    \dand
    \Tb_n
    \left[\bF(\hat\bu_n)-\bF(\bu_n)
      -\Jb_{\bF}(\bu_n)(\hat\bu_n-\bu_n)\right]
    =
    \op{1}.
\]
Writing
\[
    \boldsymbol{\Sigma}_u
    \coloneqq
    \begin{pmatrix}
        2&0\\
        0&8\beta^2
    \end{pmatrix},
\]
the limiting covariance under the normalization $\Tb_n$ is
\[
    \Ab\boldsymbol{\Sigma}_u\Ab^\intercal
    =
    \begin{pmatrix}
        \dfrac{1}{2}&\beta\\[9pt]
        \beta&2\beta^2
    \end{pmatrix}.
\]
Since $2u_{1,n}\sim1/\alpha_n$, rescaling the first coordinate gives the joint limits
\begin{equation}\label{eq:gamma2_jointlimit}
    \begin{pmatrix}
        \sqrt n\dfrac{(\tilde\alpha_n-\alpha_n)}{\alpha_n}\\[9pt]
        \sqrt n(\tilde\beta_n-\beta)
    \end{pmatrix}
    \inD
    \mathcal N\left(
        \bzero,
        \begin{pmatrix}
            2&2\beta\\
            2\beta&2\beta^2
        \end{pmatrix}
    \right)
    \dand
    \begin{pmatrix}
        \sqrt n\dfrac{(\hat\alpha_n-\alpha_n)}{\alpha_n}\\[9pt]
        \sqrt n(\hat\beta_n-\beta)
    \end{pmatrix}
    \inD
    \mathcal N\left(
        \bzero,
        \begin{pmatrix}
            2&2\beta\\
            2\beta&2\beta^2
        \end{pmatrix}
    \right).
\end{equation}

\begin{rem}
The common covariance matrix in the asymptotic distributions in \eqref{eq:gamma2_jointlimit} has rank one because
\[
    \sqrt n(\tilde\beta_n-\beta)
    =
    \beta\frac{\sqrt n}{\alpha_n}
    (\tilde\alpha_n-\alpha_n)
    +\op{1}.
\]
This singular limit does not indicate nonidentifiability. It merely arises because under the normalization $\Tb_n$, the variation in the second chart coordinate is of smaller order after transforming back to the original parameter scale.
\end{rem}

\subsection{Gamma Regression}\label{sub:glmregression}

Gamma regression extends the preceding example to non-identically distributed observations. For convenience, we turn to the mean--dispersion parameterization of the gamma distribution, whose density for $y>0$ is
\begin{equation}\label{eq:glm_density}
    f_{\mu,\varphi}(y)
    \coloneqq
    \frac{y^{1/\varphi-1}
    \exp{-\frac{y}{(\mu\varphi)}}}
    {\Gamma\left(\dfrac{1}{\varphi}\right)(\mu\varphi)^{1/\varphi}}.
\end{equation}
Let $k\coloneqq p+1$ be fixed, and let $\Xb_n\in\R^{n\times k}$ be a deterministic design matrix with rows $\bx_{n,i}^\intercal$, whose first column is $\vec{\bm 1}_n$. We assume that
\begin{equation}\label{eq:glm_design}
    \sup_{n\geq1}\max_{1\leq i\leq n}\norm{\bx_{n,i}}<\infty
    \dand
    \frac1n\Xb_n^\intercal\Xb_n\longrightarrow\Sb,
\end{equation}
where $\Sb$ is positive definite. For a fixed vector of coefficients $\bbeta_0\in\R^k$, suppose that the response is $\bY_n \coloneqq (Y_{n,1},\ldots,Y_{n,n})^\intercal$, where $Y_{n,1},\ldots,Y_{n,n}$ are independent and gamma distributed with means
\[
    \mu_{n,i}\coloneqq\exp{\bx_{n,i}^\intercal\bbeta_0}
\]
and common dispersion $\varphi_n\downarrow0$, so that
\[
    Z_{n,i}\coloneqq Y_{n,i}\exp{-\bx_{n,i}^\intercal\bbeta_0}
    \iid \mathrm{Gamma}(1,\varphi_n),
\]
where the notation on the right refers to the parameterization in \eqref{eq:glm_density}.

The limiting regime $\varphi_n\downarrow0$ can be interpreted in the context of generalized linear models (GLMs), where in the notation of \citet{mccullaghGeneralizedLinearModels2019}, $\nu\coloneqq\varphi^{-1}$ is the precision parameter and $\varphi$ is the squared coefficient of variation. Under the log link, the GLM weight matrix is the identity matrix, and when $\varphi$ is small, the first-order covariance of the MLE of $\bbeta$ agrees with that obtained by applying least squares to the log responses \citep[Sections 8.2 and 8.3.4]{mccullaghGeneralizedLinearModels2019}. Below, we show that the least squares estimator based on the log responses, paired with an explicit estimator of $\varphi$, forms an approximate MLE obtained from a tractable surrogate for the gamma likelihood equations. Our results then show that, under suitable asymptotic conditions, the least squares estimator is asymptotically equivalent to the exact MLE of $\bbeta$ and inherits its first-order asymptotic distribution.

\subsubsection{Derivation of the Approximate MLE}\label{subsub:glm_approx_mle}

For a candidate $\bbeta\in\R^k$, define
\[
    Z_{n,i}(\bbeta)
    \coloneqq
    Y_{n,i}\exp{-\bx_{n,i}^\intercal\bbeta}\dand
    \bZ_n(\bbeta)
    \coloneqq
    \bigl(
        Z_{n,1}(\bbeta),\ldots,Z_{n,n}(\bbeta)
    \bigr)^\intercal
    \dand
    \bZ_n\coloneqq\bZ_n(\bbeta_0),
\]
and let
\[
    h(z)\coloneqq z-1-\log{z}
    \dand
    g(\varphi)\coloneqq-\log{\varphi}-\psi\left(\dfrac{1}{\varphi}\right).
\]
Differentiating the log likelihood and multiplying the two score equations by nonzero constant factors gives the equivalent system
\begin{equation}\label{eq:glm_score}
    \bq_n(\varphi,\bbeta)
    \coloneqq
    \begin{pmatrix}
        g(\varphi)-\dfrac1n\displaystyle\sum_{i=1}^n
        h\left(Z_{n,i}(\bbeta)\right)\\[9pt]
        \dfrac{\varphi}{2n}\Xb_n^\intercal
        \left(\bZ_n(\bbeta)-\vec{\bm 1}_n\right)
    \end{pmatrix}
    =\bzero.
\end{equation}
Define the log-response vector and the log residuals by
\[
    \bL_n
    \coloneqq
    \bigl(\log{Y_{n,1}},\ldots,\log{Y_{n,n}}\bigr)^\intercal
    \dand
    r_{n,i}(\bbeta)
    \coloneqq
    \log{Y_{n,i}}-\bx_{n,i}^\intercal\bbeta.
\]
For each fixed $i$,
\[
    r_{n,i}(\bbeta_0)=\log{Z_{n,i}}\inP0
\]
under the limiting regime. The expansions
\[
    g(\varphi)=\frac{\varphi}{2}+O(\varphi^2)
    \dand
    h(e^r)=\frac{r^2}{2}+O(|r|^3)
    \dand
    e^r-1=r+O(r^2)
\]
therefore suggest the surrogate equations
\[
    \Xb_n^\intercal
    \left(\bL_n-\Xb_n\bbeta\right)
    =
    \bzero\dand
    \varphi
    =
    \frac1n
    \norm{\bL_n-\Xb_n\bbeta}^2.
\]
For sufficiently large $n$, our assumptions imply that $\Xb_n^\intercal\Xb_n$ is nonsingular. Let
\[
    \Hb_n
    \coloneqq
    \Xb_n(\Xb_n^\intercal\Xb_n)^{-1}\Xb_n^\intercal\dand
    \Mb_n
    \coloneqq
    \Ib_n-\Hb_n.
\]
Solving the surrogate equations gives the candidate estimators
\begin{equation}\label{eq:glm_estimators}
    \tilde\bbeta_n
    \coloneqq
    (\Xb_n^\intercal\Xb_n)^{-1}
    \Xb_n^\intercal\bL_n\dand
    \tilde\varphi_n
    \coloneqq
    \frac1n\norm{\Mb_n\bL_n}^2.
\end{equation}
Thus $\tilde\bbeta_n$ is the least squares estimator based on the log responses, while $\tilde\varphi_n$ is the residual sum of squares divided by $n$. Note that $n > k$ eventually, and since $\bL_n$ has a continuous distribution, $\tilde\varphi_n>0$ almost surely.

\begin{rem}
The estimator $\tilde\varphi_n$ can also be viewed as an explicit version of the deviance approximation to the gamma dispersion MLE given by \cite{mccullaghGeneralizedLinearModels2019}. Define
\[
    D_n(\bbeta)
    \coloneqq
    2\sum_{i=1}^n h\left(Z_{n,i}(\bbeta)\right)
    =
    -2\sum_{i=1}^n
    \left[
        \log{
            \frac{Y_{n,i}}{\mu_{n,i}(\bbeta)}}
        -\frac{Y_{n,i} - \mu_{n,i}(\bbeta)}{\mu_{n,i}(\bbeta)}
    \right].
\]
If $\hat\nu_n=1/\hat\varphi_n$, then the first coordinate of
\eqref{eq:glm_score}, evaluated at the MLE, becomes
\[
    2n\left(\log{\hat\nu_n}-\psi(\hat\nu_n)\right)
    =
    D_n(\hat\bbeta_n),
\]
which is the likelihood equation given by \citet[Equation (8.1)]{mccullaghGeneralizedLinearModels2019}. Their large precision approximation has leading term $\hat\varphi_n \approx (1/n)D_n(\hat\bbeta_n)$. In our construction, $Z_{n,i}(\bbeta)=\exp{r_{n,i}(\bbeta)}$ and $2h(e^r)=r^2+O(|r|^3)$. Consequently, the leading quadratic approximation to $(1/n)D_n(\bbeta)$ is the mean squared log residual. At $\bbeta=\tilde\bbeta_n$, this is exactly
\[
    \frac1n
    \sum_{i=1}^n r_{n,i}(\tilde\bbeta_n)^2
    =
    \tilde\varphi_n.
\]
Thus we obtain $\tilde\varphi_n$ by replacing the gamma GLM deviance with its quadratic approximation, while the corresponding linear approximation to the $\bbeta$ score function yields the least squares estimator $\tilde\bbeta_n$.
\end{rem}

We use the vanishing chart
\[
    \bphi(\varphi,\bbeta)
    \coloneqq
    \begin{pmatrix}
        \dfrac{\varphi}{2}\\[9pt]
        \dfrac{\varphi}{2}\bbeta
    \end{pmatrix}
    \dand
    \bu\coloneqq\bphi(\varphi,\bbeta)
    =
    \begin{pmatrix}
        u_1\\
        \bu_2
    \end{pmatrix},
\]
which is a bijection from $N\coloneqq(0,\infty)\times\R^k$ to $U\coloneqq(0,\infty)\times\R^k$ with inverse
\begin{equation}\label{eq:glm_phiinv}
    \bphi\inv(\bu)
    =
    \begin{pmatrix}
        2u_1\\[9pt]
        \dfrac{\bu_2}{u_1}
    \end{pmatrix}.
\end{equation}
For fixed $\bbeta_0$, this chart satisfies $\bphi(\varphi_n,\bbeta_0)\longrightarrow\bzero$. For $\bu=(u_1,\bu_2^\intercal)^\intercal$, set
\[
    \br_n(\bu)
    \coloneqq
    \bL_n-\Xb_n\frac{\bu_2}{u_1},
\]
and let $r_{n,i}(\bu)$ denote its $i$th component. The exact equation \eqref{eq:glm_score} in chart coordinates is
\begin{equation}\label{eq:glm_exact_Q}
    \bQ_n(\bu)
    \coloneqq
    \begin{pmatrix}
        g(2u_1)-\dfrac1n\displaystyle\sum_{i=1}^n
        h\left(\exp{r_{n,i}(\bu)}\right)\\[9pt]
        \dfrac{u_1}{n}\Xb_n^\intercal
        \left[\bZ_n\left(\dfrac{\bu_2}{u_1}\right)-\vec{\bm 1}_n\right]
    \end{pmatrix}.
\end{equation}
The following surrogate for $\bQ_n$ has the same root as this system:
\begin{equation}\label{eq:glm_surrogate_Q}
    \tilde\bQ_n(\bu)
    \coloneqq
    \begin{pmatrix}
        u_1-\tilde u_{1,n}\\[9pt]
        \dfrac1n\Xb_n^\intercal
        \left(u_1\bL_n-\Xb_n\bu_2\right)
    \end{pmatrix}
    \dand
    \tilde u_{1,n}\coloneqq\frac{\tilde\varphi_n}{2},
\end{equation}
whose unique root is
\[
    \tilde\bu_n
    \coloneqq
    \begin{pmatrix}
        \tilde u_{1,n}\\
        \tilde u_{1,n}\tilde\bbeta_n
    \end{pmatrix},
\]
which belongs to $U$ almost surely for sufficiently large $n$.

\subsubsection{Verification of Approximation}\label{subsub:glm_verification}

For notational convenience, write
\[
    \btheta_n
    \coloneqq
    (\varphi_n,\bbeta_0^\intercal)^\intercal
    \dand
    \hat\btheta_n
    \coloneqq
    (\hat\varphi_n,\hat\bbeta_n^\intercal)^\intercal
    \dand
    \tilde\btheta_n
    \coloneqq
    (\tilde\varphi_n,\tilde\bbeta_n^\intercal)^\intercal,
\]
and let
\[
    \bu_n\coloneqq\bphi(\btheta_n)
    \dand
    \hat\bu_n\coloneqq\bphi(\hat\btheta_n)
    \dand
    \tilde\bu_n\coloneqq\bphi(\tilde\btheta_n).
\]
Write $\mathbb P_n\coloneqq\mathbb P_{n,\btheta_n}$ for the true joint law. In particular,
\[
    \bu_n
    \coloneqq
    \begin{pmatrix}
        u_{1,n}\\
        \bu_{2,n}
    \end{pmatrix}
    =
    \begin{pmatrix}
        \dfrac{\varphi_n}{2}\\[9pt]
        \dfrac{\varphi_n}{2}\bbeta_0
    \end{pmatrix}.
\]

By \cref{lemm:glm_estimator_rates}, the likelihood has a unique maximizer with probability tending to one; the required localization and equation error bounds are established in \cref{lemm:glm_localization_remainder}. In particular, there exist deterministic convex sets $D_n$ such that
\begin{equation}\label{eq:glm_localization}
    \mathbb P_n
    \left(\hat\bu_n,\tilde\bu_n\in D_n\right)
    \longrightarrow1
    \dand
    \sup_{\bu\in D_n}\norm{\bu}
    =
    O(\varphi_n)
    \longrightarrow0,
\end{equation}
and
\begin{equation}\label{eq:glm_uniform_remainder}
    \Delta_{D,n}
    \coloneqq
    \sup_{\bu\in D_n}
    \norm{\tilde\bQ_n(\bu)-\bQ_n(\bu)}
    =
    \op{\varphi_n}.
\end{equation}
The surrogate has constant Jacobian
\[
    \Jb_n
    \coloneqq
    \Jb_{\tilde\bQ_n}
    =
    \begin{pmatrix}
        1&\bzero^\intercal\\[9pt]
        \dfrac{1}{n}\Xb_n^\intercal\bL_n&-\dfrac{1}{n}\Xb_n^\intercal\Xb_n
    \end{pmatrix},
\]
whose inverse is computed via block matrix inversion as
\[
    \Jb_n\inv
    =
    \begin{pmatrix}
        1&\bzero^\intercal\\[0.3em]
        \tilde\bbeta_n&-n(\Xb_n^\intercal\Xb_n)^{-1}
    \end{pmatrix}.
\]
Since $\tilde\bbeta_n\inP\bbeta_0$ and $n^{-1}\Xb_n^\intercal\Xb_n\to\Sb$, we have that $\norm{\Jb_n\inv}=\Op{1}$. Consequently, there is a constant $a>0$ such that
\[
    \mathbb P_n
    \left(\sigmin(\Jb_n)\geq a\right)
    \longrightarrow1.
\]
Moreover, because the surrogate Jacobian is identically equal to $\Jb_n$, the Jacobian approximation condition in \cref{thm:MAIN} holds with zero error. Taking $\rho_n=1/\varphi_n$ in \cref{coro:root_rate}, 
\eqref{eq:glm_localization} and \eqref{eq:glm_uniform_remainder} give
\begin{equation}\label{eq:glm_chart_approximation}
    \frac1{\varphi_n}
    \norm{\tilde\bu_n-\hat\bu_n}
    =
    \op{1}.
\end{equation}
Thus $\tilde\bu_n$ and $\hat\bu_n$ are asymptotically equivalent relative to the vanishing chart at scale $1/\varphi_n$.

It remains to translate this comparison back to the original parameter scale. On $D_n$, the inverse chart has Jacobian
\[
    \Jb_{\bphi\inv}(\bu)
    =
    \begin{pmatrix}
        2&\bzero^\intercal\\[9pt]
        -\dfrac{\bu_2}{u_1^2}&\dfrac{\Ib_k}{u_1}
    \end{pmatrix},
\]
and the construction of $D_n$ ensures that $u_1\asymp\varphi_n$ and that $\bu_2/u_1$ remains bounded uniformly on $D_n$. Hence
\[
    \sup_{\bu\in D_n}
    \norm{\Jb_{\bphi\inv}(\bu)}
    =
    O\left(\frac{1}{\varphi_n}\right),
\]
so $\bphi\inv$ is $b_n$-Lipschitz on $D_n$ for some
$b_n\coloneqq C/\varphi_n$. Since
\[
    \Delta_{D,n}
    =
    \op{\varphi_n}
    =
    \op{\frac{a}{b_n}},
\]
\cref{coro:original_scale} yields
\[
    \norm{\tilde\btheta_n-\hat\btheta_n}
    =
    \op{1}.
\]
Together with \eqref{eq:glm_chart_approximation}, we obtain the explicit componentwise differences
\[
    |\tilde\varphi_n-\hat\varphi_n|
    =
    \op{\varphi_n}
    \dand
    \norm{\tilde\bbeta_n-\hat\bbeta_n}
    =
    \op{1}.
\]
The estimators in \eqref{eq:glm_estimators} are therefore approximate MLEs both relative to the chart and on the original parameter scale.

\subsubsection{Joint Asymptotic Distribution}

The approximate estimator of $\bbeta$ in \eqref{eq:glm_estimators} can have a first-order bias in its first component, which corresponds to the intercept. To see this, define $\be_0 \coloneqq(1,0,\ldots,0)^\intercal\in\R^k$, 
\[
    m_n \coloneqq \E\left[\log{Z_{n,i}}\right]
      = \psi\left(\dfrac{1}{\varphi_n}\right)+\log{\varphi_n}
    \dand
    v_n\coloneqq \Var(\log{Z_{n,i}})
      =\psi_1\left(\dfrac{1}{\varphi_n}\right).
\]
and set
\[
    c_n\coloneqq\left(1-\frac{k}{n}\right)v_n.
\]
Since the first column of $\Xb_n$ is $\vec{\bm 1}_n$ and $\Mb_n$ has rank $n-k$, we have that
\[
    \E\left[\tilde\bbeta_n\right]
    =
    \bbeta_0+m_n\be_0
    \dand
    \E\left[\tilde\varphi_n\right]
    =
    c_n,
\]
so $m_n\be_0$ and $c_n-\varphi_n$ are deterministic biases of $\tilde\bbeta_n$ and $\tilde\varphi_n$, respectively.

The joint limits proved in \cref{lemm:glm_joint_limits} are
\begin{equation}\label{eq:glm_joint_approx_centered}
    \begin{pmatrix}
        \dfrac{\sqrt n}{\varphi_n}
        (\tilde\varphi_n-c_n)\\[9pt]
        \sqrt{\dfrac n{\varphi_n}}
        (\tilde\bbeta_n-\bbeta_0-m_n\be_0)
    \end{pmatrix}
    \inD
    \mathcal N_{k+1}\left(
        \bzero,
        \begin{pmatrix}
            2&\bzero^\intercal\\
            \bzero&\Sb\inv
        \end{pmatrix}
    \right)
\end{equation}
and
\begin{equation}\label{eq:glm_joint_exact}
    \begin{pmatrix}
        \dfrac{\sqrt n}{\varphi_n}
        (\hat\varphi_n-\varphi_n)\\[9pt]
        \sqrt{\dfrac n{\varphi_n}}
        (\hat\bbeta_n-\bbeta_0)
    \end{pmatrix}
    \inD
    \mathcal N_{k+1}\left(
        \bzero,
        \begin{pmatrix}
            2&\bzero^\intercal\\
            \bzero&\Sb\inv
        \end{pmatrix}
    \right).
\end{equation}
Thus the dispersion error has order $\varphi_n/\sqrt n$, whereas the coefficient error has order $\sqrt{\varphi_n/n}$. The two blocks in each limiting normal vector are independent. For the exact MLE, the zero cross-covariance term follows from the fact that the coefficient and dispersion score functions are orthogonal:
\[
    \operatorname{Cov}
    \left(Z_{n,i}-1,h(Z_{n,i})\right)
    =
    0.
\]
Note that the deterministic centering terms $c_n$ and $\bbeta_0 + m_n\be_0$ in \eqref{eq:glm_joint_approx_centered} cannot always be replaced by the true parameters. Indeed,
\[
    m_n
    =
    -\frac{\varphi_n}{2}
    +O(\varphi_n^2)
    \dand
    c_n-\varphi_n
    =
    \frac{\varphi_n^2}{2}
    +O\left(
        \varphi_n^3+\frac{\varphi_n}{n}
    \right).
\]
If $n\varphi_n\to\lambda\in[0,\infty)$, then
\eqref{eq:glm_joint_approx_centered} becomes
\begin{equation}\label{eq:glm_joint_approx_shifted}
    \begin{pmatrix}
        \dfrac{\sqrt n}{\varphi_n}
        (\tilde\varphi_n-\varphi_n)\\[9pt]
        \sqrt{\dfrac n{\varphi_n}}
        (\tilde\bbeta_n-\bbeta_0)
    \end{pmatrix}
    \inD
    \mathcal N_{k+1}\left(
        \begin{pmatrix}
            0\\[9pt]
            -\sqrt\lambda\,\dfrac{\be_0}{2}
        \end{pmatrix},
        \begin{pmatrix}
            2&\bzero^\intercal\\
            \bzero&\Sb\inv
        \end{pmatrix}
    \right).
\end{equation}
When $n\varphi_n\to0$, the explicit estimators and the exact MLE therefore have the same first-order distribution when centered at the true parameter. If $n\varphi_n\to\lambda>0$, their limiting covariance matrices still agree, but the approximate intercept has the nonzero mean shift in \eqref{eq:glm_joint_approx_shifted}. If $n\varphi_n\to\infty$, however, the intercept bias dominates the normalization. This restriction arises from deterministic bias rather than random approximation error. 

Removing the leading bias is useful because it relaxes the sufficient condition for the explicit and exact estimators to have the same first-order distribution from $n\varphi_n\to0$ to $n\varphi_n^4\to0$. Note that the correction need not reduce the realized distance from the exact MLE in every sample; rather, it removes the leading deterministic centering terms. As we show below, the correction is asymptotically negligible on both the unscaled original parameter scale and the unscaled chart scale, so it preserves the equivalence already established in \cref{subsub:glm_verification}. Because
\[
    \E\left[\tilde\varphi_n\right]
    =
    \left(1-\frac{k}{n}\right)v_n,
\]
first define
\[
    \check v_n
    \coloneqq
    \frac{\tilde\varphi_n}{1-k/n},
\]
which corrects $\tilde\varphi_n$ for the degrees of freedom. Expansions of the polygamma functions give
\[
    \varphi_n
    =
    v_n-\frac12v_n^2+O(v_n^3)
    \dand
    -m_n
    =
    \frac12v_n-\frac16v_n^2+O(v_n^3),
\]
which lead to the corrections
\begin{equation}\label{eq:glm_bias_corrected_estimators}
    \tilde\varphi_n^{\mathrm{bc}}
    \coloneqq
    \check v_n-\frac12\check v_n^2
    \dand
    \tilde\bbeta_n^{\mathrm{bc}}
    \coloneqq
    \tilde\bbeta_n
    +
    \left(
        \frac12\check v_n-\frac16\check v_n^2
    \right)\be_0.
\end{equation}
Only the intercept is corrected because the common mean $m_n$ of the log errors lies in the intercept column of the design matrix. Since $\check v_n\inP0$, the bias-corrected dispersion estimator is positive with probability tending to one; on the complementary event, we use the completion convention of \cref{sec:defs}. For the joint limit of these bias-corrected estimators, if $n\varphi_n^4\longrightarrow0$ (or equivalently $\sqrt n\,\varphi_n^2\to0$), then
\begin{equation}\label{eq:glm_joint_bias_corrected}
    \begin{pmatrix}
        \dfrac{\sqrt n}{\varphi_n}
        (\tilde\varphi_n^{\mathrm{bc}}-\varphi_n)\\[9pt]
        \sqrt{\dfrac n{\varphi_n}}
        (\tilde\bbeta_n^{\mathrm{bc}}-\bbeta_0)
    \end{pmatrix}
    \inD
    \mathcal N_{k+1}\left(
        \bzero,
        \begin{pmatrix}
            2&\bzero^\intercal\\
            \bzero&\Sb\inv
        \end{pmatrix}
    \right).
\end{equation}
In fact, the corrected explicit estimator and the exact MLE are jointly asymptotically equivalent:
\begin{equation}\label{eq:glm_bias_corrected_equivalence}
    \begin{pmatrix}
        \dfrac{\sqrt n}{\varphi_n}
        (\tilde\varphi_n^{\mathrm{bc}}-\hat\varphi_n)\\[9pt]
        \sqrt{\dfrac n{\varphi_n}}
        (\tilde\bbeta_n^{\mathrm{bc}}-\hat\bbeta_n)
    \end{pmatrix}
    =\op{1}.
\end{equation}
Since $\check v_n=\Op{\varphi_n}$, the corrections are $o_{\mathbb P}(1)$ on the original parameter scale and $o_{\mathbb P}(\varphi_n)$ in the chart coordinates. The bias-corrected estimators $\tilde\varphi_n^{\mathrm{bc}}$ and $\tilde\bbeta_n^{\mathrm{bc}}$ are therefore also approximate MLEs on both scales. Provided that $n\varphi_n^4\longrightarrow0$, their joint limit agrees with the exact MLE limit in \eqref{eq:glm_joint_exact}.

\subsection{Simulation Study}\label{sub:gamma_simulation_study}

We compared the three gamma models under the same simulation design. Let $a_n\coloneqq n^\gamma$, where
\[
    \gamma\in\{0.1,0.4,0.7\}\dand
    n\in\{250,500,1000,2500,5000,10000,25000\}.
\]
For the two iid models of \cref{sub:fixed_rate,sub:unknown_rate}, we set $\alpha_n\coloneqq a_n$ and $\beta\coloneqq2$. For the gamma regression model of \cref{sub:glmregression}, we set $\varphi_n\coloneqq a_n^{-1}$ and used a design matrix with rows $\bx_{n,i}\coloneqq(1,x_{n,i})^\intercal$, where the $x_{n,i}$ were equally spaced and centered on $[-1,1]$, and took $\bbeta_0\coloneqq(0.5,-0.75)^\intercal$. We replicated each design point $2000$ times. Our setup placed points on both sides of the sufficient rate conditions derived above: $\gamma>1/3$ for the known-rate model, $\gamma>1/6$ for the unknown-rate model, and $\gamma>1/4$ for the bias-corrected gamma regression estimator.

In each replication, we computed the exact MLE using Newton's method and the corresponding approximate estimators from \cref{sub:fixed_rate,sub:unknown_rate,sub:glmregression}. We divided the absolute difference between an approximate estimator and the exact MLE by the asymptotic standard error of the corresponding MLE coordinate, which compared the approximation error with sampling variation. We also divided the RMSE of the approximate estimator by that of the exact MLE to measure any finite-sample loss in accuracy.

\begin{figure}[!t]
    \centering
    \includegraphics[width=0.94\textwidth]{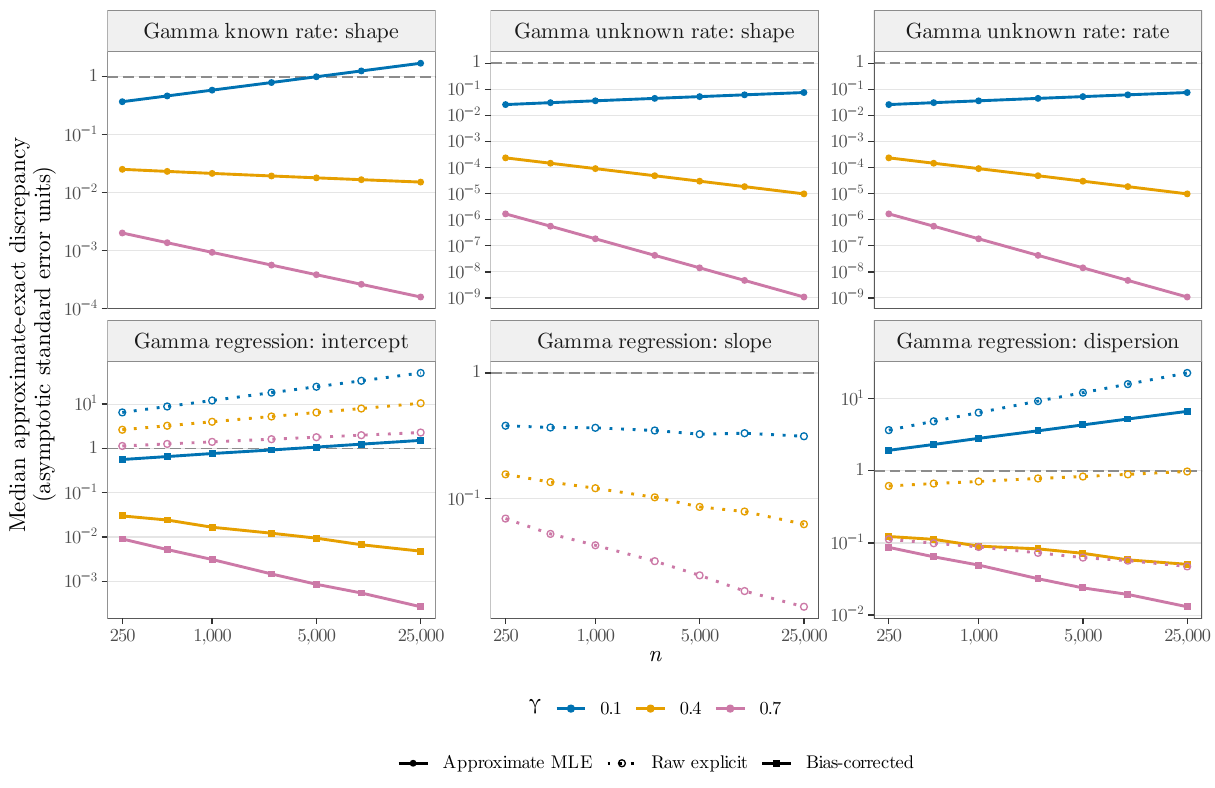}
    \caption{Median absolute difference between each explicit estimator and the exact MLE, divided by the asymptotic standard error of the corresponding MLE coordinate. The horizontal dashed line marks one asymptotic standard error. The vertical scales differ among panels. Each point is based on $2000$ replications.}
    \label{fig:gamma-unified-discrepancy}
\end{figure}

\cref{fig:gamma-unified-discrepancy} shows a clear transition across the rate conditions. In the known-rate model, the normalized discrepancy increases with $n$ for $\gamma=0.1$, but decreases for $\gamma=0.4$ and $0.7$. At $n=25000$, its median is $1.71$ for $\gamma=0.1$, compared with $0.0153$ for $\gamma=0.4$ and $1.59\times10^{-4}$ for $\gamma=0.7$. In the unknown-rate model, the discrepancy likewise increases only for $\gamma=0.1$, which is below the threshold of $1/6$; even there, its median at $n=25000$ is $0.0753$ for both coordinates. In the regression panels, the raw intercept bias remains visible throughout the grid, while the raw dispersion discrepancy is substantial in the slower regimes but is already small for $\gamma=0.7$. For $\gamma>1/4$, the bias-corrected discrepancy decreases with $n$ in every coordinate.

\FloatBarrier

\begin{figure}[!t]
    \centering
    \includegraphics[width=0.94\textwidth]{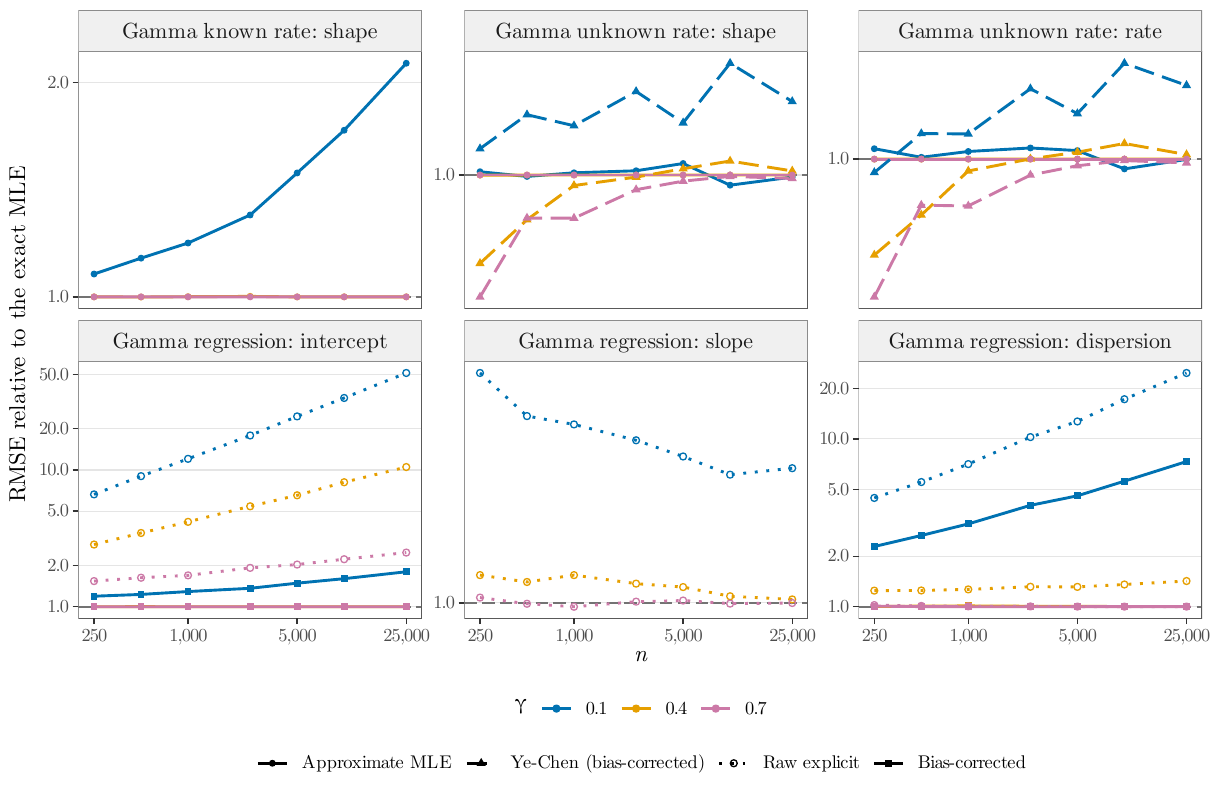}
    \caption{RMSEs of each explicit estimator relative to that of the exact MLE. The curves for the estimators of \cite{yeClosedFormEstimatorsGamma2017} appear only for the unknown-rate model, and the raw and bias-corrected curves appear only for gamma regression; note that the vertical scales differ across panels. Each point is based on $2000$ replications.}
    \label{fig:gamma-unified-rmse}
\end{figure}

For the unknown-rate model, we also included the bias-corrected mixed log moment estimators of \citet{yeClosedFormEstimatorsGamma2017} as a closed-form benchmark, which were derived from two likelihood equations for the generalized gamma family as a noniterative alternative to the gamma MLE. In our shape-rate parameterization, these estimators are given by 
\[ 
    \check\alpha_n \coloneqq \left(\frac{n-1}{n+2}\right)\frac{\bar X_n}{\check\tau_n}
    \dand
    \check\beta_n \coloneqq \left(\frac{n-1}{n+2}\right)\frac{1}{\check\tau_n}
    \dand
    \check\tau_n\coloneqq\frac1n\sum_{i=1}^n(X_i - \bar{X}_n)\log{X_i}.
\]
The RMSE comparisons in \cref{fig:gamma-unified-rmse} lead to the same conclusions. For known rate and $n=25000$, the RMSE ratio is $2.13$ when $\gamma=0.10$ and within $0.01\%$ of 1 when $\gamma=0.40$ or $0.70$. With unknown rate, our approximate estimator is nearly indistinguishable from the MLE throughout the grid. The bias-corrected estimators of \cite{yeClosedFormEstimatorsGamma2017} are similarly competitive: their RMSE ratios range from $0.978$ to $1.020$ over all design points and both coordinates. In gamma regression, the bias correction brings the intercept and dispersion RMSE ratios close to one for $\gamma=0.40$ and $0.70$. The raw intercept retains substantial bias over the entire grid, while the slope requires no correction and approaches the MLE as $\varphi_n$ decreases.

\FloatBarrier

\begin{figure}[!t]
    \centering
    \includegraphics[width=0.90\textwidth]{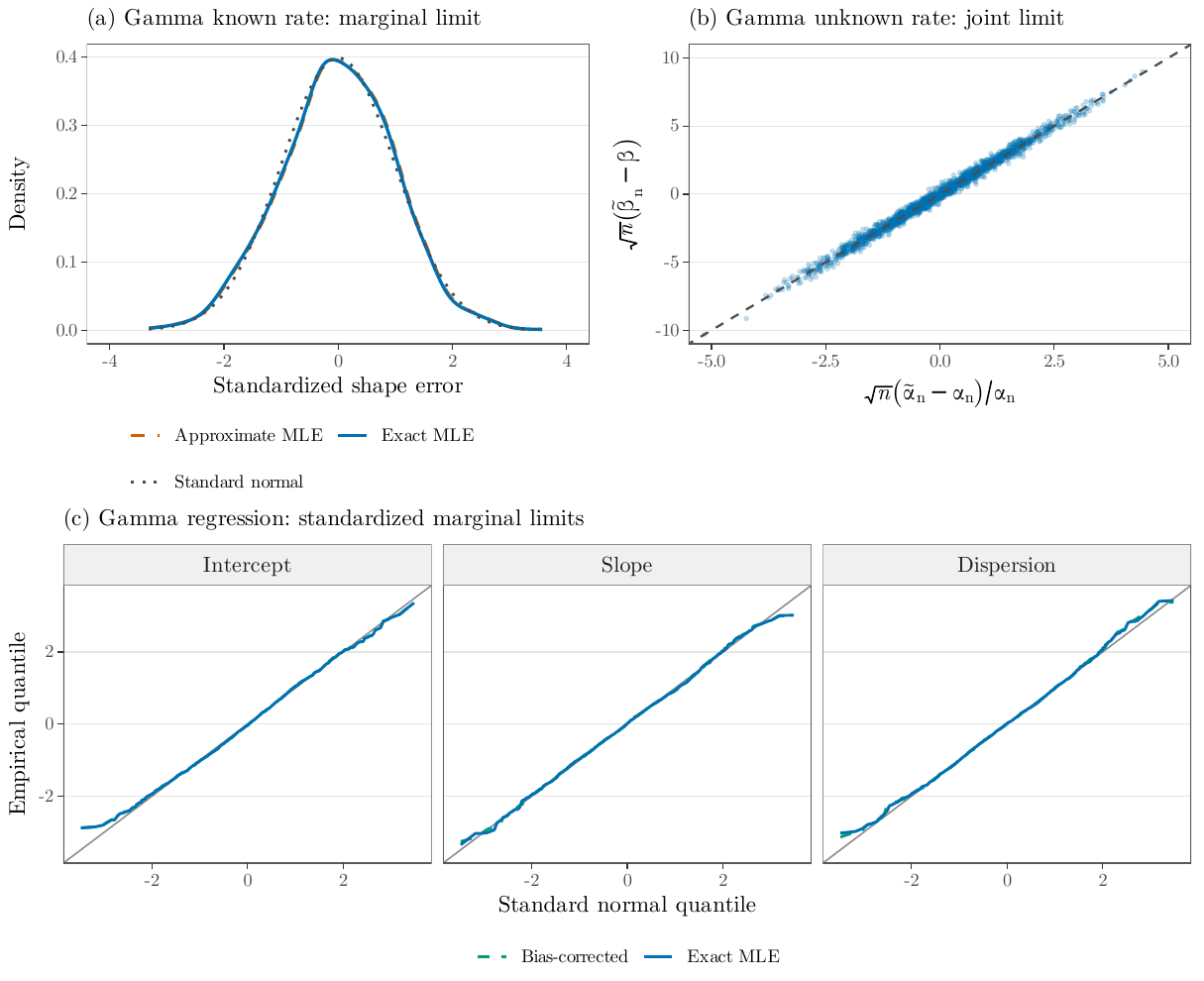}
    \caption{Distributional diagnostics at $n=25000$ and $\gamma=0.4$. Panel (a) compares the standardized shape errors in the known-rate model with the standard normal density. Panel (b) shows the joint approximate-estimator errors in the unknown-rate model; the dashed line is $z_\beta=2z_\alpha$. Panel (c) gives normal quantile plots for unit-variance standardizations of the exact and bias-corrected regression estimators.}
    \label{fig:gamma-unified-diagnostics}
\end{figure}

Finally, \cref{fig:gamma-unified-diagnostics} examines the limiting distributions at $(n,\gamma)=(25000,0.4)$, which satisfies all three sufficient rate conditions. In the known-rate model, the approximate and exact standardized shape errors are visually indistinguishable and both closely follow the standard normal density. In the unknown-rate model, the empirical variances of $\sqrt n(\tilde\alpha_n-\alpha_n)/\alpha_n$ and $\sqrt n(\tilde\beta_n-\beta)$ are $2.02$ and $8.12$, respectively, compared with the limiting values $2$ and $8$, and their correlation is $0.996$, compared with the limiting value $1$. The regression quantile plots use versions of the normalizations in \eqref{eq:glm_joint_exact} and \eqref{eq:glm_joint_bias_corrected} with unit variance. The empirical variances of the corrected intercept, slope, and dispersion errors are $0.994$, $0.977$, and $1.016$, respectively. Thus the finite-sample results reproduce both the marginal limits and the singular joint behavior derived above.

\subsection{Application: Black Cherry Tree Data}\label{sub:gamma_tree_data}

To illustrate the practical use of our framework, we conclude this section with a brief application to the \texttt{trees} data distributed with R, which record the diameter, height, and timber volume of 31 felled black cherry trees \citep{meyerForestMensuration1953}. These data have been used to illustrate gamma regression with a log link by \citet[Chapter~4]{woodGeneralizedAdditiveModels2017}. For $i=1,\dots,31$, we fit the parametric model
\[
    \mathrm{Volume}_i\mid d_i,h_i \sim \text{Gamma}(\exp{\beta_0+\beta_1d_i+\beta_2h_i}, \varphi),
\]
in the mean-dispersion parameterization given by \eqref{eq:glm_density}, where $d_i$ and $h_i$ denote the diameter and height of the $i$th tree, respectively.

The raw explicit estimator provides a direct way to assess the plausibility of the small-$\varphi$ regime without numerically maximizing the likelihood. It gives $\tilde\varphi_n=0.008456$ and $\sqrt{\tilde\varphi_n}=0.0920$, so the fitted conditional coefficient of variation is $9.20\%$. At the corresponding coefficient estimate $\tilde\bbeta_n$, the log residuals $\tilde r_i \coloneqq \log{Y_i}-\bx_i^\intercal\tilde\bbeta_n$ are also small, with $\max_i|\tilde r_i|=0.1773$. Thus, the quadratic approximation $h(e^r)\approx r^2/2$ is highly accurate for this fit: the exact average deviance term is
\[
    \frac{2}{n}\sum_i h\left(e^{\tilde r_i}\right)
    =0.008486,
\]
whereas its quadratic approximation is
\[
    \frac1n\sum_i\tilde r_i^2
    =\tilde\varphi_n
    =0.008456,
\]
a relative difference of only $0.36\%$.

The usual gamma GLM, fitted numerically using the \texttt{glm} function in R with the log link, converged in six iterations. \cref{tab:gamma-tree-data} compares it with both the raw estimators in \eqref{eq:glm_estimators} and the bias-corrected estimators in \eqref{eq:glm_bias_corrected_estimators}. For the regression coefficients, the usual GLM uses $\hat\varphi_n^{\mathrm{P}}(\Xb_n^\intercal\Xb_n)^{-1}$, where $\hat\varphi_n^{\mathrm{P}}$ is the Pearson statistic divided by the residual degrees of freedom; the raw and corrected fits replace $\hat\varphi_n^{\mathrm{P}}$ by $\tilde\varphi_n$ and $\tilde\varphi_n^{\mathrm{bc}}$, respectively. The dispersion intervals use the first-order variance $2\varphi^2/n$, with the corresponding estimate substituted for $\varphi$.

\begin{table}[!htbp]
\centering
\begingroup
\renewcommand{\arraystretch}{1.2}
\newcommand{\estci}[2]{%
    \shortstack[c]{\rule{0pt}{2.5ex}#1\\[2pt]#2}%
}

\begin{tabular}{l|ccc}
\hline
\rule{0pt}{2.8ex}
parameter
& \texttt{glm}
& raw explicit
& bias-corrected explicit
\\[1pt]
\hline

intercept
& \estci{$0.0923$}{$(-0.3308,\,0.5154)$}
& \estci{$0.1026$}{$(-0.2985,\,0.5037)$}
& \estci{$0.1073$}{$(-0.3138,\,0.5283)$}
\\[3pt]

diameter
& \estci{$0.1453$}{$(0.1323,\,0.1582)$}
& \estci{$0.1453$}{$(0.1330,\,0.1576)$}
& \estci{$0.1453$}{$(0.1324,\,0.1582)$}
\\[3pt]

height
& \estci{$0.0166$}{$(0.0102,\,0.0230)$}
& \estci{$0.0164$}{$(0.0103,\,0.0224)$}
& \estci{$0.0164$}{$(0.0100,\,0.0227)$}
\\[3pt]

dispersion
& \estci{$0.009410$}{$(0.004726,\,0.014095)$}
& \estci{$0.008456$}{$(0.004246,\,0.012666)$}
& \estci{$0.009318$}{$(0.004679,\,0.013957)$}
\\[1pt]
\hline
\end{tabular}
\endgroup
\caption{Gamma regression estimates for the black cherry tree data, with 95\% Wald intervals in parentheses. The \texttt{glm} dispersion is the Pearson estimate divided by the residual degrees of freedom; the full likelihood estimate of dispersion is $0.008455$.}
\label{tab:gamma-tree-data}
\end{table}

Both explicit coefficient vectors differ from the usual GLM estimate by less than $0.07$ of the GLM standard error, and all three sets of intervals lead to the same conclusions. The raw coefficient standard errors are $5.2\%$ smaller than the GLM values, whereas the corrected standard errors differ by only $0.49\%$. The raw dispersion estimate $0.008456$ is nearly identical to the full likelihood estimate $0.008455$, while the corrected estimate $0.009318$ is closer to the Pearson estimate $0.009410$ reported by \texttt{glm}. Relative to the GLM fitted means, the mean and maximum absolute relative differences decrease from $0.42\%$ and $0.63\%$ for the raw fit to $0.10\%$ and $0.29\%$ after correction.

Thus the correction improves overall agreement with the conventional GLM fit, although it does not necessarily reduce the realized distance from the exact MLE in every coordinate or every sample. Its purpose is to remove deterministic bias relative to the true parameter. The correction is negligible for the unscaled comparison with the MLE, so the corrected and raw estimators retain the same unscaled equivalence properties; however, it removes terms that may be non-negligible under the first-order normalization, thereby weakening the sufficient condition for sharing the MLE's first-order distribution from $n\varphi_n\to0$ to $n\varphi_n^4\to0$.

\FloatBarrier
\section{Additional Examples}\label{sec:additionalexamples}

In this section, we illustrate the generality of our approach by applying it to several additional parametric families, focusing here on the main constructions and resulting approximate MLEs. Complete proofs and calculations, including verification of the assumptions of \cref{thm:MAIN}, the derivation of corresponding asymptotic distributions, and additional simulation results, are deferred to the appendices.

\subsection{Zeta Distribution}\label{sub:zeta}

For $s>1$ and $k\in\{1,2,\ldots\}$, the $\mathrm{Zeta}(s)$ distribution has mass function
\[
    \mathbb P(X=k)=\frac{k^{-s}}{\zeta(s)}.
\]
As $s \to \infty$, the distribution concentrates on the value $k=1$, with 
\[
    \PP_s(X > 1) = 1 - \frac{1}{\zeta(s)} = \frac{1}{2^s}(1 + o(1)).
\]
For $X_1,\ldots,X_n\iid\mathrm{Zeta}(s)$, let $G_n\coloneqq(\prod_{i=1}^nX_i)^{1/n}$. The score function and its derivative are, respectively,
\[
    q_n(s)\coloneqq-n\left(\frac{\zeta'(s)}{\zeta(s)}+\log{G_n}\right) \dand q_n'(s)=-n\Var_s(\log{X})<0.
\]
As $s\downarrow1$, $q_n(s)\to\infty$, while as $s\to\infty$, $q_n(s)\to-n\log{G_n}$; therefore, a unique finite MLE exists exactly when $G_n>1$. If $G_n=1$ (or equivalently $X_1=\cdots=X_n=1$), the likelihood approaches its supremum as $s\to\infty$ but has no finite maximizer; on this event, we use the completion convention of \cref{sec:defs}.

The function $\zeta'/\zeta$ is the logarithmic derivative of the Riemann zeta function, which admits a well-known expansion as a Dirichlet series:
\begin{equation}\label{eq:zeta_dirichlet}
    -\frac{\zeta'(s)}{\zeta(s)}
    =\sum_{m=2}^\infty\frac{\Lambda(m)}{m^s},
\end{equation}
valid for all $s > 1$ \citep[Theorem 294]{hardyIntroductionTheoryNumbers2008}, where $\Lambda$ is the von Mangoldt function.\footnote{This function, of fundamental importance in number theory, maps a positive integer $k$ to $\log{p}$ if $k$ is an integer power of some prime $p$, and $0$ otherwise.} As $s\to\infty$, the series is dominated by its first term, giving
\begin{equation}\label{eq:zeta_dirichlet_asymptotic}
    -\frac{\zeta'(s)}{\zeta(s)}
    =\frac{\log{2}}{2^s}+O\left(\frac{1}{3^s}\right).
\end{equation}
Thus, in contrast to the gamma examples of \cref{sec:gamma}, this approximation arises from a Dirichlet series rather than an expansion in integer powers of the parameter. Define the vanishing chart by $\phi(s)\coloneqq2^{-s}$, set $u\coloneqq\phi(s)$ and $U\coloneqq(0,1/2)$, and note that $\phi^{-1}(u)=-\logtwo{u}$. The exact and surrogate equations become
\[
    Q_n(u)\coloneqq n\left(m(u)-\log{G_n}\right)\dand
    \widetilde Q_n(u)\coloneqq n\left(\log{2}u-\log{G_n}\right),
    \quad \text{where} \quad
    m(u)\coloneqq\sum_{m=2}^\infty\Lambda(m)u^{\logtwo{m}}.
\]
The surrogate root and corresponding approximate MLE are, respectively,
\[
    \widetilde u_n\coloneqq\frac{\log{G_n}}{\log{2}}\dand
    \widetilde s_n
    \coloneqq-\logtwo{\frac{\log{G_n}}{\log{2}}}
    =-\logtwo{\logtwo{G_n}}.
\]

Fix $s_n\to\infty$, set $u_n\coloneqq2^{-s_n}$ and $p\coloneqq\logtwo{3}$, and write $\mathbb P_n\coloneqq\mathbb P_{s_n}^{\otimes n}$. If $nu_n\to\infty$, then we have that
\[ \hat u_n = u_n + o_{\mathbb P_n}(u_n)\dand
\tilde u_n = u_n + o_{\mathbb P_n}(u_n).
\]
Moreover,
\[  \mathbb P_n(G_n = 1) = \frac{1}{\zeta(s_n)^n} = \exp{-nu_n(1+o(1))} \longrightarrow 0, \]
so both $\hat u_n$ and $\tilde u_n$ belong to $U$ with probability tending to one.

The verification in \cref{app:sub:zeta} gives
\[ 
    |\widetilde u_n-\widehat u_n|
    =O_{\mathbb P_n}(u_n^p)
    \dand
    |\widetilde s_n-\widehat s_n|
    =O_{\mathbb P_n}\left(u_n^{p-1}\right)
    =O_{\mathbb P_n}\left(2^{-s_n(p-1)}\right)
    =o_{\mathbb P_n}(1).
\]
The exact MLE also satisfies
\[
    \sqrt{\frac{n}{2^{s_n}}}(\widehat s_n-s_n)
    \inD\mathcal N\left(0,\frac1{\log{2}^2}\right).
\]
If, in addition, $\sqrt n\,u_n^{p-1/2}\to0$, the verification in \cref{sec:zeta} shows that the discrepancy between $\hat{s}_n$ and $\tilde{s}_n$ is negligible under the relevant normalization. Slutsky's theorem then gives
\[
    \sqrt{\frac{n}{2^{s_n}}}(\widetilde s_n-s_n)
    \inD\mathcal N\left(0,\frac1{\log{2}^2}\right).
\]
Other closed-form estimators of $s$ are available in the literature. For example, \citet[Equation 3.7]{clausetPowerLawDistributionsEmpirical2009} obtain one by approximating a lower-truncated discrete power law by its continuous analogue and applying a continuity correction. The accuracy of this approximation depends on a sufficiently large lower truncation $x_{\min}$, whereas the model we consider here has $x_{\min}=1$; our construction instead retains the exact likelihood and provides an asymptotic comparison with the exact MLE.

The full verification and simulation study are provided in \cref{sec:zeta}, along with supporting moment calculations and a verification of the relevant Lyapunov condition.

\subsection{von Mises Distribution}\label{sub:vonmises}

For $\kappa>0$, the $\mathrm{von\,Mises}(\mu,\kappa)$ density on $[0,2\pi)$ is
\[
    f_{\mu,\kappa}(x)
    \coloneqq\frac{\exp{\kappa\cos{x-\mu}}}{2\pi I_0(\kappa)},
\]
where $\mu$ is understood modulo $2\pi$ and $I_m$ denotes the modified Bessel function of the first kind of order $m$. The location parameter $\mu$ determines the central direction of the distribution, while the concentration parameter $\kappa$ controls the spread around $\mu$. While the MLE of $\mu$ is explicit, the MLE of $\kappa$ is not available in closed form and requires numerical inversion of a ratio of Bessel functions.

For $X_1, \ldots, X_n \iid \mathrm{von\,Mises}(\mu,\kappa)$, define
\[
    \overline C_n\coloneqq\frac1n\sum_{i=1}^n\cos{X_i}
    \dand
    \overline S_n\coloneqq\frac1n\sum_{i=1}^n\sin{X_i}
    \dand
    R_n\coloneqq\sqrt{\overline C_n^2+\overline S_n^2},
\]
and let $\widehat\mu_n\coloneqq\operatorname{atan2}(\overline S_n,\overline C_n)$.

Writing $A(\kappa)\coloneqq I_1(\kappa)/I_0(\kappa)$, the MLE $(\widehat\mu_n,\widehat\kappa_n)$ is characterized by
\[
    \widehat\mu_n=\operatorname{atan2}(\overline S_n,\overline C_n)\dand
    A(\widehat\kappa_n)=R_n.
\]
The function $A$ increases strictly from $0$ to $1$ on $(0,\infty)$. For $n\geq2$, the events $R_n=0$ and $R_n=1$ have probability zero, so the MLE pair is unique (modulo $2\pi$) almost surely; we use the completion convention of \cref{sec:defs} on exceptional events.

To work with an ordinary coordinate for $\mu$, fix an open interval $M\subset\R$ of length less than $2\pi$ that contains a chosen representative of the true value $\mu$. Whenever the sample mean direction belongs to the corresponding arc, we take $\widehat\mu_n$ to be its unique representative in $M$. The argument given in \cref{sec:vonmises} shows that this representation is available with probability tending to one as $\kappa\to\infty$.

To preserve the explicit location estimator while placing both parameters in a vanishing chart, we combine the location and concentration equations as follows:
\begin{equation}\label{eq:vonmises-system}
    \bq_n(\mu,\kappa)
    \coloneqq
    \begin{pmatrix}
        \mu\left(A(\kappa)-1\right)-\widehat\mu_n(R_n-1)\\
        A(\kappa)-R_n
    \end{pmatrix}
    =\bzero.
\end{equation}
This modification preserves the selected root. Indeed, if $q_{2,n}=0$, then
\[
    q_{1,n}=(R_n-1)(\mu-\widehat\mu_n),
\]
which forces $\mu=\widehat\mu_n$ in this representation since $R_n<1$ almost surely. The factor $A(\kappa)-1$ tends to zero as $\kappa\to\infty$, which allows both coordinates of the chart below to vanish.

The Bessel ratio has the expansion
\[
    A(\kappa)
    =1-\frac1{2\kappa}-\frac1{8\kappa^2}
    -\frac1{8\kappa^3}+O\left(\frac{1}{\kappa^4}\right)
\]
\citep[Formulas 9.7.1 and 9.7.3]{abramowitz1988handbook}. Retaining the first two non-constant terms gives $g(\kappa)\coloneqq-1/(2\kappa)-1/(8\kappa^2)$; replacing $A(\kappa)-1$ by $g(\kappa)$ in \eqref{eq:vonmises-system} leaves the location root unchanged and gives the candidate estimators
\[
    \widetilde\mu_n\coloneqq\widehat\mu_n
    \dand
    \widetilde\kappa_n
    \coloneqq\frac{1+\sqrt{3-2R_n}}{4(1-R_n)}.
\]
Thus only the concentration equation is approximated. The leading term of
$\widetilde\kappa_n$ coincides with the classical approximation
for the concentration parameter of a $p$-dimensional von Mises--Fisher
distribution,
\[
    \widehat\kappa
    \approx
    \frac{p-1}{2(1-\overline R)}
\]
\citep[Equation~(10.3.7)]{mardiaDirectionalStatistics1999}, which is also used
by \citet[Equation~(4.1)]{banerjeeClusteringUnitHypersphere2005}. In the present setting, set $p\coloneqq2$ and $\overline R\coloneqq R_n$; writing $\delta_n\coloneqq1-R_n$, we have
\[
    \widetilde\kappa_n
    =
    \frac{1}{2\delta_n}
    +\frac14
    -\frac{\delta_n}{8}
    +O(\delta_n^2).
\]
Thus the classical approximation is the leading term of
$\widetilde\kappa_n$, while the explicit estimator retains an additional
constant-order correction. Indeed, inversion of the expansion of $A$ above
gives
\[
    \widehat\kappa_n
    =
    \frac{1}{2\delta_n}
    +\frac14
    +\frac{3\delta_n}{8}
    +O(\delta_n^2),
\]
so the two estimators agree through their constant terms. The following analysis
makes this formal comparison rigorous under the coupled asymptotic regime. 

Set $N\coloneqq M\times(0,\infty)$. The vanishing chart
\[
    \bphi(\mu,\kappa)
    \coloneqq\left(\mu g(\kappa),g(\kappa)\right)^{\mathsf T}
    \dand
    \bu\coloneqq\bphi(\mu,\kappa)=(u_1,u_2)^{\mathsf T}
    \dand
    U\coloneqq\{(u_1,u_2):u_2<0,\ \frac{u_1}{u_2}\in M\},
\]
has inverse
\[
    \bphi^{-1}(u_1,u_2)
    =
    \left(\frac{u_1}{u_2},
    -\frac{1+\sqrt{1-2u_2}}{4u_2}\right)^{\mathsf T}.
\]
For fixed $\mu$, both coordinates vanish as $\kappa\to\infty$, while the location parameter is encoded in the ratio $u_1/u_2$.

Fix $\kappa_n\to\infty$ and write $\mathbb P_n\coloneqq\mathbb P_{\mu,\kappa_n}^{\otimes n}$. The verification in \cref{app:vonMisesdistresults} first localizes the explicit transformed root and then uses \cref{coro:stochastic_localization_transfer} to locate the exact root. With probability tending to one, both roots consequently belong to a common deterministic set $K_n\subset U$ satisfying $\sup_{\bu\in K_n}\norm{\bu}\to0$. On these sets,
\[ 
\Jb_{\widetilde{\bQ}_n}=\Ib_2 \dand
\sup_{\bu\in K_n}\norm{\bQ_n(\bu)-\widetilde{\bQ}_n(\bu)}
    =O\left(\frac{1}{\kappa_n^3}\right),
\]
and $\bphi^{-1}$ is $O(\kappa_n^2)$-Lipschitz. The enclosure gives
\[
    \norm{\widetilde{\bu}_n-\widehat{\bu}_n}
    =O_{\mathbb P_n}\left(\frac{1}{\kappa_n^3}\right),
\]
while \cref{coro:original_scale} gives
\[
    \widetilde\mu_n=\widehat\mu_n \dand
    \widetilde\kappa_n-\widehat\kappa_n
    =O_{\mathbb P_n}\left(\frac{1}{\kappa_n}\right)
    =o_{\mathbb P_n}(1).
\]
The Fisher information for a one observation is
\[
    \mathcal I_1(\mu,\kappa)
    =
    \begin{pmatrix}
        \kappa A(\kappa)&0\\
        0&A'(\kappa)
    \end{pmatrix},
\]
where $A'(\kappa)\sim1/(2\kappa^2)$; thus, the two coordinates have different natural rates. The exact MLE satisfies
\[
    \begin{pmatrix}
        \sqrt{n\kappa_n}(\widehat\mu_n-\mu)\\[9pt]
        \sqrt n\,\dfrac{(\widehat\kappa_n-\kappa_n)}{\kappa_n}
    \end{pmatrix}
    \inD
    \mathcal N\left(
        \bzero,\begin{pmatrix}1&0\\0&2\end{pmatrix}
    \right).
\]
If in addition $\sqrt n/\kappa_n^2\to0$ (equivalently, $\kappa_n/n^{1/4}\to\infty$), then
\[ \frac{\sqrt n}{\kappa_n}
(\widetilde\kappa_n-\widehat\kappa_n)
=o_{\mathbb P_n}(1).
\]
Since the exact and approximate location MLEs coincide, Slutsky's theorem gives the same joint limit for $(\widetilde\mu_n,\widetilde\kappa_n)$.

For comparison, \citet{abeyasekeraEstimationParametersMises1982} use a normal approximation to a highly concentrated von Mises distribution and, for $n>3$, propose the closed-form estimator
\[
    \widehat\kappa_{\mathrm{lin}}
    \coloneqq
\frac{1}{
\frac{1}{n-3}
\sum_{i=1}^n(X_i-\overline X_n)^2
}.
\]
Because this estimator treats the angles as ordinary real observations, an unwrapping convention is needed when the sample straddles the circular discontinuity. Their construction approximates the distribution by a normal law, whereas ours retains the exact von Mises likelihood and approximates only the Bessel ratio equation. The detailed verification and a simulation study that includes a comparison with the estimator of \citet{abeyasekeraEstimationParametersMises1982} are given in \cref{sec:vonmises}.

\subsection{Generalized Pareto Distribution}\label{sub:GPD}

Finally, we consider the generalized Pareto distribution, which presents two additional complications: its support depends on one parameter, and the score function for another parameter has a singularity at the boundary of the parameter space.

Fix $0<\bar\xi<1/6$ and consider the restricted parameter space
\[
    \Theta_{\bar\xi}
    \coloneqq
    \{(\mu,\sigma,\xi):\mu\in\R,\ \sigma>0,\ 0\leq\xi\leq\bar\xi\}.
\]
For $\xi>0$ and $x\geq\mu$, the generalized Pareto distribution has density
\[
    f_{\mu,\sigma,\xi}(x)
    \coloneqq
    \frac1{\sigma}
    \left(1+\xi\frac{x-\mu}{\sigma}\right)^{-1-\frac{1}{\xi}}.
\]
At $\xi=0$ and for $x\geq\mu$, its continuous limit is the shifted exponential density
\[
    f_{\mu,\sigma,0}(x)
    \coloneqq
    \frac1{\sigma}\exp{-\frac{x-\mu}{\sigma}}.
\]
We observe $X_{n,1},\ldots,X_{n,n}\iid\mathrm{GPD}(\mu,\sigma,\xi_n)$, where $\mu\in\R$ and $\sigma>0$ are fixed and $\xi_n\downarrow0$, and write $\mathbb P_n\coloneqq\mathbb P_{\mu,\sigma,\xi_n}^{\otimes n}$. The log-likelihood function with all three parameters unrestricted can be unbounded, so we impose the upper bound $\xi \leq \bar\xi$, guaranteeing the existence of a global maximizer by \cref{prop:gpd_global_mle}. The restriction $\bar\xi < 1/6$ is a convenient choice rather than a sharp bound; under the conditions given below, it is never approached in the limit. This restriction is consistent with the need for care in computing GPD likelihood estimators as investigated by \cite{grimshawComputingMaximumLikelihood1993}.

Writing $X_{(1)}$ for the sample minimum of the $X_i$, the log-likelihood is increasing in the location parameter $\mu$ throughout the feasible region, and so every global maximizer has $\widehat\mu_n = X_{(1)}$. Set 
\[
    Y_i\coloneqq X_{n,i}-X_{(1)}
    \dand
    M_{k,n}\coloneqq\frac1n\sum_{i=1}^nY_i^k.
\]
We profile $\mu$ from the likelihood at $X_{(1)}$ and approximate only the scale and shape coordinates. The profile scores for the remaining parameters are given by
\[
    S_\sigma(\sigma,\xi)
    \coloneqq \frac{1}{\sigma}\left(\sum_{i=1}^n
       \frac{(1+\xi)t_i(\sigma)}{1+\xi t_i(\sigma)} - n\right)
    \dand
    S_\xi(\sigma,\xi) \coloneqq \frac{1}{\xi}\left(\sum_{i=1}^n
    \frac{\log{1+\xi t_i(\sigma)}}{\xi} - \sum_{i=1}^n
       \frac{(1+\xi)t_i(\sigma)}{1+\xi t_i(\sigma)} \right)
\]
where $t_i(s)\coloneqq Y_i/s$. Writing
\[
 B_n(\sigma,\xi)\coloneqq\frac1n\sum_{i=1}^n
    \frac{\log{1+\xi t_i(\sigma)}}{\xi}\dand
 A_n(\sigma,\xi)\coloneqq\frac1n\sum_{i=1}^n
       \frac{(1+\xi)t_i(\sigma)}{1+\xi t_i(\sigma)}, 
\]
the equivalent system
\[
    \bq_n(\sigma,\xi)
    \coloneqq
    \begin{pmatrix}
        B_n(\sigma,\xi)-1\\[9pt]
        \dfrac{(B_n(\sigma,\xi)-A_n(\sigma,\xi))}{\xi}
    \end{pmatrix}
    =
    \frac1n\begin{pmatrix}
        \sigma&\xi\\
        0&1
    \end{pmatrix}
    \begin{pmatrix}S_\sigma(\sigma,\xi)\\S_\xi(\sigma,\xi)\end{pmatrix}
\]
has the same interior roots because the matrix is nonsingular for $\sigma>0$. This rewriting exposes the removable singularities at $\xi=0$ and makes the truncation transparent. Expanding at $\xi=0$ and retaining the constant and linear terms gives
\[
    \widetilde\bq_n(\sigma,\xi)
    \coloneqq
    \begin{pmatrix}
        \dfrac{M_{1,n}}{\sigma}-1-\dfrac{\xi M_{2,n}}{(2\sigma^2)}\\[9pt]
        -\dfrac{M_{1,n}}{\sigma}+\dfrac{M_{2,n}}{(2\sigma^2)}
        +\xi\left(\dfrac{M_{2,n}}{\sigma^2}-\dfrac{2M_{3,n}}{(3\sigma^3)}\right)
    \end{pmatrix}.
\]
Setting
\[
 C_n\coloneqq6M_{1,n}M_{2,n}+8M_{3,n} \dand
 E_n\coloneqq8M_{1,n}M_{3,n}-3M_{2,n}^2,
\]
we select the root $(\widetilde\sigma_n, \widetilde\xi_n)$, where
\[
 \widetilde\sigma_n
 \coloneqq
 \frac{C_n-\sqrt{C_n^2-48M_{2,n}E_n}}{24M_{2,n}}\dand
 \widetilde\xi_n
 \coloneqq
 \frac{2\widetilde\sigma_n(M_{1,n}-\widetilde\sigma_n)}{M_{2,n}}.
\]
At the population moments of the shifted exponential distribution, the two scale roots are $\sigma$ and $3\sigma/2$, so we select the smaller branch. We also set $\widetilde\mu_n\coloneqq X_{(1)}$. The selection and local uniqueness of this root are established in \cref{lemm:gpd_surrogate}. If the selected raw branch is not real or does not belong to $(0,\infty)\times[0,\bar\xi]$, we define $(\widetilde\sigma_n,\widetilde\xi_n)\coloneqq(M_{1,n},0)$, the shifted-exponential boundary estimator. The probability of this replacement tends to zero under the localization regime below, so it leaves the asymptotic arguments for the interior branch unchanged while producing a constrained estimator in every sample.

The first surrogate equation also gives a useful interpretation of the scale estimator. On the event where the selected root is well defined, we have
\[
    \widetilde\sigma_n
    =
    \frac{
        M_{1,n}
        +\sqrt{M_{1,n}^2-2M_{2,n}\widetilde\xi_n}
    }{2}.
\]
Under the regime considered below, $\widetilde\xi_n=o_{\mathbb P_n}(1)$, so expansion about $\widetilde\xi_n=0$ gives
\[
    \widetilde\sigma_n
    =
    M_{1,n}
    -\frac{M_{2,n}}{2M_{1,n}}\widetilde\xi_n
    +O_{\mathbb P_n}\left(\widetilde\xi_n^2\right).
\]
At $\xi=0$, the GPD reduces to the shifted exponential distribution, for which the profile MLE of the scale parameter after setting $\widehat\mu_n=X_{(1)}$ is precisely
\[
    M_{1,n}=\frac1n\sum_{i=1}^n\left(X_{n,i}-X_{(1)}\right).
\]
Thus $\widetilde\sigma_n$ may be viewed as adding a first-order shape correction to the shifted-exponential scale estimator.

The population target of the selected branch satisfies
\[
    \sigma_*(\xi)
    =\sigma\left(1+8\xi^2-68\xi^3+O(\xi^4)\right)
    \dand
    \xi_*(\xi)
    =\xi-10\xi^2+100\xi^3+O(\xi^4),
\]
so the population discrepancy begins at order $\xi^2$ in both coordinates. After profiling out the common location estimate, we use the vanishing chart
\[
    \bphi(\sigma,\xi)
    \coloneqq
    \begin{pmatrix}
        \dfrac{\xi}{\sigma}\\[9pt]
        \dfrac{\xi}{\sigma^2}
    \end{pmatrix}
    \dand
    \bF(u_2,u_3)
    \coloneqq\bphi\inv(u_2,u_3)
    =
    \begin{pmatrix}
        \dfrac{u_2}{u_3}\\[9pt]
        \dfrac{u_2^2}{u_3}
    \end{pmatrix},
\]
on $0<\xi<\bar\xi$, whose image is
\[
 U\coloneqq\left\{(u_2,u_3):u_2,u_3>0,\ \frac{u_2^2}{u_3}<\bar\xi\right\}.
\]
Suppose first that we are in the regime where $\sqrt n\,\xi_n\to\infty$, which keeps the true shape parameter away from zero at the $\sqrt n$ scale. The verification in \cref{prop:gpd_global_mle,prop:gpd_root_stability} shows that the restricted global MLE exists uniquely with probability tending to one and
\[
    \norm{
        \begin{pmatrix}
            \widetilde\sigma_n\\
            \widetilde\xi_n
        \end{pmatrix}
        -
        \begin{pmatrix}
            \widehat\sigma_n\\
            \widehat\xi_n
        \end{pmatrix}
    }
    =
    O_{\mathbb P_n}(\xi_n^2)
    =
    o_{\mathbb P_n}(1).
\]
Because the two estimators have the same location coordinate, $\widetilde\btheta_n$ is an approximate MLE on the original parameter scale.

If in addition $\sqrt n\,\xi_n^2\to0$, then the approximation error is negligible on the $\sqrt n$ scale:
\[
    \sqrt n
    \begin{pmatrix}
        \widetilde\sigma_n-\widehat\sigma_n\\
        \widetilde\xi_n-\widehat\xi_n
    \end{pmatrix}
    =
    o_{\mathbb P_n}(1).
\]
Together, the conditions $\sqrt n\,\xi_n\to\infty$ and $\sqrt n\,\xi_n^2\to0$ are equivalent to $n^{-1/2}\ll\xi_n\ll n^{-1/4}$. By \cref{thm:gpd_limits}, both estimators then satisfy
\[
    \sqrt n
    \begin{pmatrix}
        \widehat\sigma_n-\sigma\\
        \widehat\xi_n-\xi_n
    \end{pmatrix}
    \inD
    \mathcal N\left(
        \bzero,
        \begin{pmatrix}
            2\sigma^2&-\sigma\\
            -\sigma&1
        \end{pmatrix}
    \right)
    \dand
    \sqrt n
    \begin{pmatrix}
        \widetilde\sigma_n-\sigma\\
        \widetilde\xi_n-\xi_n
    \end{pmatrix}
    \inD
    \mathcal N\left(
        \bzero,
        \begin{pmatrix}
            2\sigma^2&-\sigma\\
            -\sigma&1
        \end{pmatrix}
    \right).
\]
The standard location estimator $\widehat\mu_n=\widetilde\mu_n =X_{(1)}$ has the separate non-Gaussian limit
\[
    \frac n\sigma(\widehat\mu_n-\mu)
    =
    \frac n\sigma(\widetilde\mu_n-\mu)
    \inD\mathrm{Exp}(1).
\]
For comparison, \citet{hoskingParameterQuantileEstimation1987a} give an explicit estimator based on probability weighted moments for the two parameter GPD with a known threshold. In the present three-parameter setting, we adapt their formulas to the excesses over $X_{(1)}$. The detailed proofs and results for finite samples are given in \cref{sec:Pareto}.

\section{Discussion}\label{sec:discussion}

It is often desirable to obtain closed-form estimators of statistical parameters, and it can therefore be worthwhile to replace an implicitly defined root of an estimating equation by an explicit root of a suitable surrogate equation. Such a replacement immediately prompts the question of when such a root remains close to the selected root of the original equation. Our framework answers this question by distinguishing three sources of error: the approximation of the surrogate equation itself, the sensitivity of the surrogate root to perturbations, and the amplification of errors through the inverse chart. While we focus on MLEs, our results apply to general Z-estimators.

Our examples show that no single approximation rule works in all cases. The gamma models of \cref{sec:gamma} illustrate how increasing the complexity of the underlying statistical model leads to an increase in the complexity of the approximate estimator. The examples in \cref{sec:additionalexamples} require different constructions: the zeta example approximates a Dirichlet series and requires us to consider a non-null event on which a finite MLE does not exist, the von Mises example preserves the exact MLE in one coordinate while approximating the other, and the GPD example requires that we carefully rewrite the profiled score equation prior to truncation, and has a separate non-Gaussian limit in one component. Our simulations and real-world data application support our theoretical developments and illustrate the importance of the rate condition stated in each example.

Our framework has some limitations. The construction of an approximate estimator is inherently model dependent: choosing a useful chart, reexpressing the estimating equation in an appropriate form, and selecting the truncation order can each require substantial work. Our framework requires a prespecified limiting regime, so in applications we require some basis for assuming that the data-generating parameter is sufficiently close to the corresponding limiting value. The plausibility of the required rate is model- and application-dependent; our simulations and data example illustrate both successful and visibly inadequate regimes. Although \cref{thm:surrogate_root_enclosure} can establish the existence of a nearby exact root, global existence, uniqueness, and identification of that root as the relevant MLE can remain model-specific, as in the zeta and GPD examples. Finally, while our handling of exceptional events is harmless as $n \to \infty$, such events can matter at moderate sample sizes.

There are a number of avenues for future work. Instead of focusing solely on fixed limiting regimes, allowing for uniform control over broad families of parameter sequences would be useful when the rate of approach to the limiting value is unknown. Because our results are first order, a natural question is whether the exact and surrogate roots agree at higher orders. Edgeworth and saddlepoint approximations for roots of estimating equations may provide suitable frameworks for investigating this question \citep{barndorff-nielsenEdgeworthSaddlePointApproximations1979,fieldSmallSampleAsymptotic1982,danielsSaddlepointApproximationsEstimating1983}. Showing that the exact and approximate estimators have the same distributional expansions beyond first order would require sharper control of the distance between them.

The bias correction in the gamma regression example of \cref{sub:glmregression} suggests another direction. By removing the leading deterministic centering terms, the correction relaxes the rate condition required for first-order equivalence while preserving both the estimator's closed form and the properties of an approximate MLE. It would be useful to develop general methods for constructing such corrections and to determine how far the rate restrictions can be relaxed without sacrificing equivalence to the exact estimator.

One may also create a limiting regime by embedding the model under consideration within a larger family, similar to exponential tilting. The GPD example serves as a prototype, since at $\xi=0$ it reduces to a shifted exponential family. More generally, one might embed $\PP_{\btheta}$ within $\PP_{\btheta,\bgamma}$ so that $\PP_{\btheta,\bgamma}$ approaches $\PP_{\btheta}$ as $\bgamma\to\bgamma^*$ (making this statement rigorous is also left for future work). The added parameter could then provide the limiting direction even when the model of interest has no natural limiting regime in $\btheta$, thereby allowing us to seek an approximate estimator in the enlarged family and compare its $\btheta$ component with the estimator in the embedded limiting model. Determining what makes an embedding useful, when this comparison is valid, and how much efficiency is lost or gained through embedding are open questions.

R code is available at
\href{https://github.com/rob-zimmerman/closed-form-z-estimation}{%
\textcolor{black}{\texttt{https://github.com/rob-zimmerman/closed-form-z-estimation}}}.

\section*{Acknowledgments}
The authors thank Nancy Reid for helpful discussions.

The authors acknowledge the use of ChatGPT to audit a number of mathematical arguments and derivations in the preparation of this paper and to assist in the setup and presentation of our simulation studies, and the use of Mathematica to verify several algebraic manipulations in the appendices.

\bibliography{References}
\appendix

\numberwithin{equation}{section}
\numberwithin{thm}{section}
\numberwithin{lemm}{section}
\numberwithin{prop}{section}
\numberwithin{coro}{section}
\numberwithin{assum}{section}
\numberwithin{algo}{section}
\numberwithin{defi}{section}
\numberwithin{conj}{section}
\numberwithin{theorem}{section}
\numberwithin{lemma}{section}
\numberwithin{proposition}{section}
\numberwithin{corollary}{section}
\numberwithin{definition}{section}
\numberwithin{remark}{section}

\newpage
\section{Technical Results for the Gamma Models}\label{app:gammadistresults}

\subsection{Gamma Moment Identities}

\begin{lemm}\label{lemm:gammaexps}
    Let $X \sim \mathrm{Gamma}(\alpha,\beta)$. If $s>-\alpha$, then
    \begin{equation}\label{eq:gamma_moments}
        \E\left[X^s\right]
        =
        \beta^{-s}\frac{\Gamma(\alpha+s)}{\Gamma(\alpha)}.
    \end{equation}
    Moreover, for $k\in\N$,
    \[
        \E\left[\log{X}^k\right]
        =
        \sum_{r=0}^k
        \binom{k}{r}
        (-\log{\beta})^r
        \frac{\Gamma^{(k-r)}(\alpha)}{\Gamma(\alpha)}.
    \]
    Equivalently,
    \begin{equation}\label{eq:logXk_bell}
        \E\left[\log{X}^k\right]
        =
        \sum_{r=0}^k
        \binom{k}{r}
        (-\log{\beta})^r
        Y_{k-r}\!\left(
            \psi(\alpha),\psi_1(\alpha),\ldots,
            \psi_{k-r-1}(\alpha)
        \right),
    \end{equation}
    where $Y_m$ denotes the $m$th complete exponential Bell polynomial, with the convention $Y_0\coloneqq1$.
\end{lemm}

\begin{proof}
    The Gamma integral gives \eqref{eq:gamma_moments} directly:
    \[
        \E\left[X^s\right]
        =
        \frac{\beta^\alpha}{\Gamma(\alpha)}
        \int_0^\infty x^{\alpha+s-1}e^{-\beta x}\dif x
        =
        \beta^{-s}\frac{\Gamma(\alpha+s)}{\Gamma(\alpha)}.
    \]
For the logarithmic moments, let $y\coloneqq\beta x$. Then
    \begin{align*}
        \E\left[\log{X}^k\right]
        &=
        \frac{1}{\Gamma(\alpha)}
        \int_0^\infty
        (\log{y}-\log{\beta})^k
        y^{\alpha-1}e^{-y}\dif y\\
        &=
        \frac{1}{\Gamma(\alpha)}
        \sum_{r=0}^k
        \binom{k}{r}
        (-\log{\beta})^r
        \int_0^\infty
        (\log{y})^{k-r}
        y^{\alpha-1}e^{-y}\dif y\\
        &=
        \sum_{r=0}^k
        \binom{k}{r}
        (-\log{\beta})^r
        \frac{\Gamma^{(k-r)}(\alpha)}{\Gamma(\alpha)}.
    \end{align*}
    Finally, applying the Fa\`a di Bruno formula \citep[Theorem 3.4 C]{comtetAdvancedCombinatoricsArt1974} to $\Gamma(\alpha)=\exp{\log{\Gamma(\alpha)}}$ gives, for $m\geq1$,
    \[
        \Gamma^{(m)}(\alpha)
        =
        \Gamma(\alpha)
        Y_m\!\left(
            \psi(\alpha),\psi_1(\alpha),\ldots,\psi_{m-1}(\alpha)
        \right),
    \]
    which yields \eqref{eq:logXk_bell}.
\end{proof}

\subsection{Known Rate}

The following results establish the localization, inverse chart expansion, and triangular array central limit theorem used in \cref{sub:fixed_rate}.

\begin{lemm}\label{lemm:gammaexps1}
    Let $X_1, X_2, \ldots, X_n \iid \mathrm{Gamma}(\alpha, \beta)$, and let $G_n \coloneqq \left(\prod_{i=1}^n X_i\right)^{1/n}$ be the geometric mean. Then, for every $k$ satisfying $k/n<\alpha$,
    \[
        \E\left[(\beta G_n)^{-k}\right]
        =
        \left(
            \frac{\Gamma(\alpha-k/n)}{\Gamma(\alpha)}
        \right)^{n}.
    \]
    In particular, whenever $2/n<\alpha$,
    \[
        \E\left[\frac{1}{\beta G_n}\right]
        =
        \left(
            \frac{\Gamma(\alpha-1/n)}{\Gamma(\alpha)}
        \right)^{n},
    \]
    and
    \[
        \Var\left(\frac{1}{\beta G_n}\right)
        =
        \left(
            \frac{\Gamma(\alpha-2/n)}{\Gamma(\alpha)}
        \right)^{n}
        -
        \left(
            \frac{\Gamma(\alpha-1/n)}{\Gamma(\alpha)}
        \right)^{2n}.
    \]
\end{lemm}
\begin{proof}
    Since $G_n^{-k}=\prod_{i=1}^nX_i^{-k/n}$, independence and \cref{lemm:gammaexps} give, whenever $k/n<\alpha$,
    \[
        \begin{aligned}
        \E\left[(\beta G_n)^{-k}\right]
        &=\beta^{-k}\prod_{i=1}^n\E\left[X_i^{-\frac{k}{n}}\right]
        =\beta^{-k}
          \left(
          \beta^{\frac{k}{n}}
          \frac{\Gamma(\alpha-k/n)}{\Gamma(\alpha)}
          \right)^{n}
          =\left(
          \frac{\Gamma(\alpha-k/n)}{\Gamma(\alpha)}
          \right)^{n},
        \end{aligned}
    \]
    from which both the first moment and variance identities follow by taking $k=1$ and $k=2$, respectively.
\end{proof}

\begin{lemm}
\label{lemm:gamma1_localization}
Let $X_{n,1},\ldots,X_{n,n}\iid\mathrm{Gamma}(\alpha_n,\beta)$, where $\beta>0$ is fixed and $\alpha_n\to\infty$, and set
\[
    u_n\coloneqq\frac{1}{\alpha_n-1/2}\dand
    \tilde u_n\coloneqq\frac{1}{\beta G_n}.
\]
Then
\begin{equation}\label{eq:gamma1_localization_rate}
    \E_{\alpha_n}\left[\tilde u_n\right]
    =
    u_n+O\left(\frac{u_n^2}{n}\right)+O(u_n^3)\dand
    \Var_{\alpha_n}(\tilde u_n)
    =O\left(\frac{u_n^3}{n}\right),
\end{equation}
and hence
\begin{equation}\label{eq:gamma1_localization}
    \tilde u_n-u_n
    =
    \Op{\sqrt{\frac{u_n^3}{n}}}
    +O\left(\frac{u_n^2}{n}\right)
    +O(u_n^3)
    =\op{u_n}.
\end{equation}
\end{lemm}

\begin{proof}
For fixed $k\in\{1,2\}$, define
\[
    L_{k,n}
    \coloneqq
    n\left(\log{\Gamma\left(\alpha_n-\dfrac{k}{n}\right)}-\log{\Gamma(\alpha_n)}\right).
\]
Applying Taylor's theorem to the function
$x\mapsto\log{\Gamma(\alpha_n-x)}$ about $x = 0$ and evaluating at $x=k/n$, gives
\begin{equation}\label{eq:gamma1_log_ratio}
    L_{k,n}
    =
    -k\psi(\alpha_n)
    +\frac{k^2}{2n}\psi_1(\alpha_n)
    +O\left(\frac{1}{n^2\alpha_n^2}\right).
\end{equation}
Before multiplication by $n$, the remainder is bounded by
\[
    C\left(\frac{k}{n}\right)^3
    \sup_{x\in[\alpha_n-\frac{k}{n},\alpha_n]}|\psi_2(x)|.
\]
Now,
\[
    \psi_2(x)
    =
    -2\sum_{j=0}^\infty \frac{1}{(x+j)^{3}} \quad \text{ and so } \quad
    |\psi_2(x)|
    =
    2\sum_{j=0}^\infty \frac{1}{(x+j)^{3}}
    \leq \frac{2}{x^3} + 2 \int_0^\infty \frac{\dif t}{(x+t)^{3}}
    =
    \frac{2}{x^3} + \frac{1}{x^2}
    =
    O\left(\frac{1}{x^2}\right),
\]
where the first equality is Equation 6.4.10 of \citet{abramowitz1988handbook}. Thus $\psi_2(x)=O(x^{-2})$ uniformly on the interval $[\alpha_n-k/n,\alpha_n]$, and the remainder in \eqref{eq:gamma1_log_ratio} stays valid when $n$ and $\alpha_n$ vary together. Next, standard digamma and trigamma expansions \citep[Equations 6.3.18 and 6.4.12]{abramowitz1988handbook}, together with $u_n = ( \alpha_n - 1/2)^{-1}$, give both
\[
    e^{-\psi(\alpha_n)}=u_n+O(u_n^3)\dand
    \psi_1(\alpha_n)=u_n+O(u_n^3).
\]
Now, by \cref{lemm:gammaexps1}, $\E_{\alpha_n}\left[\tilde u_n^k\right]=e^{L_{k,n}}$. We first take $k=1$ in \eqref{eq:gamma1_log_ratio} and exponentiate, which yields
\begin{equation}\label{eq:gamma1_expectation}
    \E_{\alpha_n}\left[\tilde u_n\right]
    =
    e^{-\psi(\alpha_n)}
    \left(1+O\left(\frac{u_n}{n}\right)\right)
    =
    u_n+O\left(\frac{u_n^2}{n}\right)+O(u_n^3).
\end{equation}
Moreover,
\[
    L_{2,n}-2L_{1,n}
    =
    \frac{\psi_1(\alpha_n)}{n}
    +O\left(\frac{1}{n^2\alpha_n^2}\right)
    =O\left(\frac{u_n}{n}\right),
\]
and $e^{2L_{1,n}}=O(u_n^2)$. Thus,
\[
    \Var_{\alpha_n}(\tilde u_n)
    =
    e^{2L_{1,n}}
    \left(e^{L_{2,n}-2L_{1,n}}-1\right)
    =O\left(\frac{u_n^3}{n}\right).
\]
Chebyshev's inequality then gives, for sufficiently large $n$ and every $M>0$,
\[
    \mathbb P_{\alpha_n}^{\otimes n}
    \left(
        |\tilde u_n-\E_{\alpha_n}\left[\tilde u_n\right]|
        >M\sqrt{\frac{u_n^3}{n}}
    \right)
    \leq\frac{C}{M^2},
\]
which, in combination with \eqref{eq:gamma1_expectation}, yields the first equality in \eqref{eq:gamma1_localization}. Since
\[ \sqrt{\frac{u_n^3}{n}}=o(u_n)
\dand
\frac{u_n^2}{n}=o(u_n)
\dand
u_n^3=o(u_n)  \]
the second equality follows.
\end{proof}

\begin{lemm}\label{lemm:gamma1_lyapunov}
    Let $X_\alpha\sim\mathrm{Gamma}(\alpha,\beta)$ with $\beta>0$ fixed, and set
    \[
        Y_\alpha\coloneqq\log{X_\alpha},
        \quad
        \mu_\alpha\coloneqq\psi(\alpha)-\log{\beta}\dand
        \sigma_\alpha^2\coloneqq\psi_1(\alpha).
    \]
    There exist $\alpha_0>0$ and $C<\infty$ such that
    \begin{equation}\label{eq:gamma1_lyapunov}
        \sup_{\alpha\geq\alpha_0}
        \E\left[
            \left|
            \frac{Y_\alpha-\mu_\alpha}{\sigma_\alpha}
            \right|^3
        \right]
        \leq C.
    \end{equation}
    Consequently, if $X_{n,1},\ldots,X_{n,n}\iid\mathrm{Gamma}(\alpha_n,\beta)$ with $\alpha_n\to\infty$, then the Lyapunov condition above holds for $Y_{n,i}\coloneqq\log{X_{n,i}}$.
\end{lemm}

\begin{proof}
    Without loss of generality, we may assume that $\beta=1$ since changing $\beta$ only translates $\log{X_\alpha}$. Define
    \[
        Z_\alpha
        \coloneqq
        \frac{Y_\alpha-\psi(\alpha)}
        {\sqrt{\psi_1(\alpha)}}.
    \]
    For fixed $t$ and
    \[
        h_\alpha(t)\coloneqq\frac{t}{\sqrt{\psi_1(\alpha)}},
    \]
    \cref{lemm:gammaexps} gives
    \[
        \E\left[e^{tZ_\alpha}\right]
        =
        \frac{\Gamma(\alpha+h_\alpha(t))}
        {\Gamma(\alpha)}
        \exp{-h_\alpha(t)\psi(\alpha)},
    \]
    so
    \[
        \log{\E\left[e^{tZ_\alpha}\right]}
        =
        \log{\Gamma(\alpha+h_\alpha(t))}
        -
        \log{\Gamma(\alpha)}
        -
        h_\alpha(t)\psi(\alpha).
    \]
    Since $\psi_1(\alpha)\sim\alpha^{-1}$ as $\alpha\to\infty$, we have $h_\alpha(t)=O(\sqrt{\alpha})$, uniformly for $t$ in any fixed bounded interval. Taylor's theorem therefore gives
    \begin{equation}\label{eq:log_mgf_Z_alpha}
        \log{\E\left[e^{tZ_\alpha}\right]}
        =
        \frac{h_\alpha(t)^2}{2}\psi_1(\alpha)
        +
        \frac{h_\alpha(t)^3}{6}\psi_2(\xi_{\alpha,t})
        =
        \frac{t^2}{2}
        +
        O(\alpha^{-\frac{1}{2}}),
    \end{equation}
    where $\xi_{\alpha,t}$ lies between $\alpha$ and $\alpha+h_\alpha(t)$ and we have used $\psi_2(x)=O(x^{-2})$. For sufficiently large $\alpha$, we have $\xi_{\alpha,t}\geq\alpha/2$. Exponentiating both sides of \eqref{eq:log_mgf_Z_alpha}, we find that there exists $\eta>0$ such that
    \[
        \sup_{\alpha\geq\alpha_0}
        \left(
            \E\left[e^{\eta Z_\alpha}\right]
            +
            \E\left[e^{-\eta Z_\alpha}\right]
        \right)
        <\infty.
    \]
    Since
    \[
        e^{\eta|x|}
        \leq
        e^{\eta x}+e^{-\eta x},
    \]
    we thus have
    \[
        \sup_{\alpha\geq\alpha_0}
        \E\left[e^{\eta|Z_\alpha|}\right]
        <\infty.
    \]
    Using the elementary bound
    \begin{equation}\label{eq:x3_bound}
        x^3\leq\frac{6}{\eta^3}e^{\eta x},
    \end{equation}
    valid for all $x \geq 0$, we have
    \[
        \sup_{\alpha\geq\alpha_0}\E\left[|Z_\alpha|^3\right] \leq \frac{6}{\eta^3} \sup_{\alpha\geq\alpha_0} \E\left[e^{\eta|Z_\alpha|}\right] <\infty,
    \]
    which proves \eqref{eq:gamma1_lyapunov}. Finally, writing
    \[
        Z_{n,i}
        \coloneqq
        \frac{Y_{n,i}-\mu_{\alpha_n}}
        {\sqrt{\psi_1(\alpha_n)}},
    \]
    we have that
    \[
        \frac{
            n\E\left[|Y_{n,1}-\mu_{\alpha_n}|^3\right]
        }{
            \left(n\psi_1(\alpha_n)\right)^{3/2}
        }
        =
        \frac{\E\left[|Z_{n,1}|^3\right]}{\sqrt n}
        \longrightarrow0,
    \]
    so the Lyapunov condition holds.
\end{proof}

\begin{lemm}
\label{lemm:gamma1_inverse_chart}
Let $u_n\downarrow0$ and suppose that a positive random sequence $V_n$ satisfies
\[
    V_n-u_n=\Op{\sqrt{\frac{u_n^3}{n}}}
    \dand
    \mathbb P\left(V_n\in[\frac{u_n}{2},2u_n]\right)\longrightarrow1.
\]
Then for $\phi\inv(u)=1/2+1/u$, we have
\[
    \sqrt{nu_n}\left(\phi\inv(V_n)-\phi\inv(u_n)\right)
    =
    -\sqrt{\frac{n}{u_n^3}}(V_n-u_n)+\op{1}.
\]
\end{lemm}

\begin{proof}
The exact identity
\[
    \frac1{V_n}-\frac1{u_n}
    =
    -\frac{V_n-u_n}{u_n^2}
    +\frac{(V_n-u_n)^2}{u_n^2V_n}
\]
shows that, on $\{V_n\in[u_n/2,2u_n]\}$, the remainder is bounded by $2(V_n-u_n)^2/u_n^3=\Op{1/n}$. Multiplying by $\sqrt{nu_n}$ therefore gives
\[
    \Op{\sqrt{\frac{u_n}{n}}}=\op{1}.
\]
Since $\mathbb P\left(V_n\in[u_n/2,2u_n]\right)\longrightarrow1$, this bound establishes $\op{1}$ remainder as claimed.
\end{proof}

\subsection{Unknown Rate}

\begin{lemm}\label{lemm:gammaexps2}
    Let $X_1,\ldots,X_n\iid\mathrm{Gamma}(\alpha,\beta)$, and suppose $n\alpha>2$. Then the following hold:
    \begin{enumerate}[label=(\roman*)]
        \item\label{lemm:EXinv}
        \[
            \E\left[\frac{1}{\bar X_n}\right]
            =
            \frac{n\beta}{n\alpha-1}.
        \]

        \item\label{lemm:Varinv}
        \[
            \Var\left(\frac{1}{\bar X_n}\right)
            =
            \frac{n^2\beta^2}
            {(n\alpha-1)^2(n\alpha-2)}.
        \]

        \item\label{lemm:ElogXG}
        \[
            \E\left[\log{\frac{\bar X_n}{G_n}}\right]
            =
            \psi(n\alpha)-\psi(\alpha)-\log{n}.
        \]

        \item\label{lemm:VarlogXG}
        \[
            \Var\left(\log{\frac{\bar X_n}{G_n}}\right)
            =
            \frac{\psi_1(\alpha)}{n}
            -
            \psi_1(n\alpha).
        \]
    \end{enumerate}
\end{lemm}

\begin{proof}
    Let
    \[
        S_n\coloneqq\sum_{i=1}^n X_i.
    \]
    Since $S_n\sim\mathrm{Gamma}(n\alpha,\beta)$,
    \[
        \bar X_n=\frac{S_n}{n}
        \sim
        \mathrm{Gamma}(n\alpha,n\beta).
    \]
    The negative gamma moments in \cref{lemm:gammaexps} give
    \[
        \E\left[\frac{1}{\bar X_n}\right]
        =\frac{n\beta}{n\alpha-1}\dand
        \E\left[\frac{1}{\bar X_n^2}\right]
        =\frac{n^2\beta^2}{(n\alpha-1)(n\alpha-2)}.
    \]
    Subtracting the square of the first expectation from the second proves parts \ref{lemm:EXinv} and \ref{lemm:Varinv}.

    For the remaining two claims, \eqref{eq:logXk_bell} in \cref{lemm:gammaexps} with $k=1,2$ gives 
    \[
        \log{\frac{\bar X_n}{G_n}}
        =
        \log{S_n}-\log{ n}-\frac1n\sum_{i=1}^n\log{X_i}.
    \]
    It follows that
    \begin{align*}
        \E\left[\log{\frac{\bar X_n}{G_n}}\right]
        &=
        \left(\psi(n\alpha)-\log{\beta}\right)
        -\log{n}
        -\left(\psi(\alpha)-\log{\beta}\right)
        =
        \psi(n\alpha)-\psi(\alpha)-\log{n},
    \end{align*}
which is \ref{lemm:ElogXG}. For \ref{lemm:VarlogXG}, write $P_i\coloneqq X_i/S_n$. By the standard gamma--Dirichlet decomposition \citep[Chapter~49]{kotzContinuousMultivariateDistributions2000}, $S_n$ is independent of $(P_1,\ldots,P_n)$, and
    \[
        \log{X_i}=\log{S_n}+\log{P_i}.
    \]
    It follows that
    \[
        \mathrm{Cov}(\log{S_n},\log{X_i})
        =
        \Var(\log{S_n})
        =
        \psi_1(n\alpha).
    \]
    Since the $X_i$ are independent,
    \[
        \Var\left(\frac1n\sum_{i=1}^n\log{X_i}\right)
        =
        \frac{\psi_1(\alpha)}{n}.
    \]
    Therefore
    \begin{align*}
        \Var\left(\log{\frac{\bar X_n}{G_n}}\right)
        &=
        \psi_1(n\alpha)
        +
        \frac{\psi_1(\alpha)}{n}
        -
        2\psi_1(n\alpha)
        =
        \frac{\psi_1(\alpha)}{n}
        -
        \psi_1(n\alpha).
    \end{align*}
\end{proof}

\begin{lemm}\label{lemm:gamma2CLT}
    Let
    \[
        X_{n,1},\ldots,X_{n,n}
        \iid
        \mathrm{Gamma}(\alpha_n,\beta)
        \dand
        \alpha_n\to\infty,
    \]
    where $\beta>0$ is fixed, and define the score vector evaluated at $(\alpha_n,\beta)$ for one observation by
    \[
        \bS_{n,i}
        \coloneqq
        \begin{pmatrix}
            \log{\beta}+\log{X_{n,i}}-\psi(\alpha_n)\\[9pt]
            \dfrac{\alpha_n}{\beta}-X_{n,i}
        \end{pmatrix}.
    \]
    Then its covariance matrix is given by
    \begin{equation}\label{eq:score_covariance}
        \boldsymbol{\mathcal I}_n
        \coloneqq
        \Var(\bS_{n,i})
        =
        \begin{pmatrix}
            \psi_1(\alpha_n) & -\dfrac{1}{\beta}\\[9pt]
            -\dfrac{1}{\beta} & \dfrac{\alpha_n}{\beta^2}
        \end{pmatrix},
    \end{equation}
    and there exist $\alpha_0>0$ and $C<\infty$ such that
    \begin{equation}\label{eq:uniform_third_moment_bound}
        \sup_{\alpha\geq\alpha_0}
        \sup_{\norm{\bt}=1}
        \E_{\alpha,\beta}
        \left[
            \left|
                \bt^\intercal
                \boldsymbol{\mathcal I}_{\alpha}^{-\frac{1}{2}}
                \bS_\alpha
            \right|^3
        \right]
        \leq C,
    \end{equation}
    where $\bS_\alpha$ and $\boldsymbol{\mathcal I}_\alpha$ are obtained by replacing $\alpha_n$ with $\alpha$ in the preceding definitions. Consequently,
    \[
        \frac{1}{\sqrt n}
        \boldsymbol{\mathcal I}_n^{-\frac{1}{2}}
        \sum_{i=1}^n\bS_{n,i}
        \inD
        \mathcal N(\bzero,\Ib_2).
    \]
\end{lemm}

\begin{proof}
    Let
    \[
        W_\alpha\coloneqq\beta X_\alpha
        \sim\mathrm{Gamma}(\alpha,1)
    \]
    and define
    \[
        v_\alpha
        \coloneqq
        \psi_1(\alpha)-\frac1\alpha.
    \]
    Since
    \[
        \psi_1(\alpha)
        =
        \sum_{k=0}^\infty\frac{1}{(\alpha+k)^2}
        >
        \int_0^\infty\frac{\dif x}{(\alpha+x)^2}
        =
        \frac1\alpha,
    \]
    we have $v_\alpha>0$. Moreover,
    \[
        v_\alpha
        =
        \frac{1}{2\alpha^2}
        +
        \frac{1}{6\alpha^3}
        +
        O\left(\frac{1}{\alpha^5}\right).
    \]

    Set
    \[
        T_\alpha
        \coloneqq
        \frac{
            \log{W_\alpha}-\psi(\alpha)
            -(W_\alpha-\alpha)/\alpha
        }{
            \sqrt{v_\alpha}
        }\dand
        Z_\alpha
        \coloneqq
        \frac{W_\alpha-\alpha}{\sqrt\alpha}.
    \]
    To verify the covariance identity \eqref{eq:score_covariance}, let
    $\widetilde W_\alpha\sim\mathrm{Gamma}(\alpha+1,1)$. Multiplying the
    $\mathrm{Gamma}(\alpha,1)$ density by $w$ changes the shape parameter
    from $\alpha$ to $\alpha+1$ and contributes the normalizing factor
    $\alpha$, so
    \[
        \E\left[W_\alpha\log{W_\alpha}\right]
        =
        \alpha\E\left[\log{\widetilde W_\alpha}\right]
        =
        \alpha\psi(\alpha+1).
    \]
    Since $\E\left[W_\alpha\right]=\alpha$ and
    $\E\left[\log{W_\alpha}\right]=\psi(\alpha)$, the recurrence
    $\psi(\alpha+1)=\psi(\alpha)+1/\alpha$ gives
    \[
        \operatorname{Cov}(\log{W_\alpha},W_\alpha)=1.
    \]
    Together with
    \[
        \Var(\log{W_\alpha})=\psi_1(\alpha)
        \dand
        \Var(W_\alpha)=\alpha,
    \]
    which follow from \cref{lemm:gammaexps}, this proves
    \eqref{eq:score_covariance}. The same identities also give
    \[
        \Var(T_\alpha)=1
        \dand
        \Var(Z_\alpha)=1
        \dand
        \operatorname{Cov}(T_\alpha,Z_\alpha)=0.
    \]
    Equivalently, if
    \[
        \Db_\alpha
        \coloneqq
        \begin{pmatrix}
            v_\alpha^{-\frac{1}{2}}
            &
            \dfrac{\beta}{\alpha\sqrt{v_\alpha}}
            \\[9pt]
            0
            &
            -\dfrac{\beta}{\sqrt\alpha}
        \end{pmatrix},
    \]
    then
    \[
        \Db_\alpha\bS_\alpha
        =
        \begin{pmatrix}
            T_\alpha\\
            Z_\alpha
        \end{pmatrix}
        \dand
        \Db_\alpha
        \boldsymbol{\mathcal I}_\alpha
        \Db_\alpha^\intercal
        =
        \Ib_2.
    \]
    Hence
    \[
        \Ob_\alpha
        \coloneqq
        \Db_\alpha
        \boldsymbol{\mathcal I}_\alpha^{\frac{1}{2}}
    \]
    is orthogonal, and therefore
    \[
        \boldsymbol{\mathcal I}_\alpha^{-\frac{1}{2}}\bS_\alpha
        =
        \Ob_\alpha^\intercal
        \begin{pmatrix}
            T_\alpha\\
            Z_\alpha
        \end{pmatrix}.
    \]
    Next, we establish uniform bounds for the third moments of $T_\alpha$ and $Z_\alpha$. For the latter, the fourth centered moment of a $\mathrm{Gamma}(\alpha,1)$ variable \citep{willinkRelationshipsCentralMoments2003} and H\"older's inequality give
    \[
        \sup_{\alpha\geq\alpha_0}\E\left[|Z_\alpha|^3\right]
        \leq
        \sup_{\alpha\geq\alpha_0}
        \left(\E\left[Z_\alpha^4\right]\right)^{\frac{3}{4}}
        =
        \left(3+\frac6{\alpha_0}\right)^{\frac{3}{4}}
        <\infty.
    \]
    For $T_\alpha$, write
    \[
        \xi_\alpha
        \coloneqq
        \log{W_\alpha}-\psi(\alpha)
        -\frac{W_\alpha-\alpha}{\alpha}
        \dand
        T_\alpha=\frac{\xi_\alpha}{\sqrt{v_\alpha}}.
    \]
    For $r>-\alpha$, direct integration gives
    \[
        \E\left[e^{r\xi_\alpha}\right]
        =
        e^{r\left(1-\psi(\alpha)\right)}
        \frac{\Gamma(\alpha+r)}{\Gamma(\alpha)}
        \left(1+\frac r\alpha\right)^{-(\alpha+r)}.
    \]
    For $s$ in a neighborhood of zero, set
    \[
        r_\alpha(s)
        \coloneqq
        \frac{s}{\sqrt{v_\alpha}}
        \dand
        c_\alpha(s)
        \coloneqq
        \frac{r_\alpha(s)}{\alpha}
        =
        \frac{s}{\alpha\sqrt{v_\alpha}}.
    \]
    Since
    \[
        \alpha\sqrt{v_\alpha}\longrightarrow\frac1{\sqrt2},
    \]
    we may choose $\eta>0$ and $\alpha_0$ sufficiently large that
    \[
        |c_\alpha(s)|\leq\frac12
        \quad
        \text{ for all }\alpha\geq\alpha_0
        \text{ and }|s|\leq\eta.
    \]
    Applying Stirling's formula to $\alpha(1+c)$ and $\alpha$, and noting that $\alpha(1+c)\geq\alpha/2$ when $|c|\leq1/2$, gives, uniformly over this range,
    \begin{align*}
        \log{\Gamma(\alpha(1+c))}
        -\log{\Gamma(\alpha)}
        &=
        \alpha c\log{\alpha}
        +
        \alpha(1+c)\log{1+c}
        -
        \alpha c
        -\frac12\log{1+c}
        +
        O\left(\frac{1}{\alpha}\right).
    \end{align*}
    Substituting $c=c_\alpha(s)$ into the preceding display yields, uniformly for $|s|\leq\eta$,
    \[
        \log{\E\left[e^{sT_\alpha}\right]}
        =
        \alpha c_\alpha(s)
        \left(\log{\alpha}-\psi(\alpha)\right)
        -
        \frac12\log{1+c_\alpha(s)}
        +
        O\left(\frac{1}{\alpha}\right).
    \]
    Since
    \[
        \alpha\left(\log{\alpha}-\psi(\alpha)\right)
        =
        \frac12+O\left(\frac{1}{\alpha}\right),
    \]
    it follows that
    \[
        \sup_{\alpha\geq\alpha_0}
        \left(
            \E\left[e^{\eta T_\alpha}\right]
            +
            \E\left[e^{-\eta T_\alpha}\right]
        \right)
        <\infty.
    \]
    Since
    \[
        e^{\eta|x|}
        \leq
        e^{\eta x}+e^{-\eta x},
    \]
    the preceding bound implies
    \[
        \sup_{\alpha\geq\alpha_0}
        \E\left[e^{\eta|T_\alpha|}\right]
        <\infty.
    \]
    The inequality \eqref{eq:x3_bound} then gives
    \[
        \sup_{\alpha\geq\alpha_0}\E\left[|T_\alpha|^3\right] \leq \frac{6}{\eta^3} \sup_{\alpha\geq\alpha_0} \E\left[e^{\eta|T_\alpha|}\right] <\infty,
    \]
    proving the uniform third moment bound for $T_\alpha$. 

    Let $\norm{\bt}=1$. Since $\Ob_\alpha$ is orthogonal, there is a unit vector $\bq_\alpha\coloneqq\Ob_\alpha\bt$ such that
    \[
        \bt^\intercal
        \boldsymbol{\mathcal I}_\alpha^{-\frac{1}{2}}\bS_\alpha
        =
        q_{\alpha,1}T_\alpha+q_{\alpha,2}Z_\alpha.
    \]
    Therefore
    \[
        \left|
            \bt^\intercal
            \boldsymbol{\mathcal I}_\alpha^{-\frac{1}{2}}\bS_\alpha
        \right|^3
        \leq
        4\left(|T_\alpha|^3+|Z_\alpha|^3\right),
    \]
    which proves \eqref{eq:uniform_third_moment_bound}. Moreover, for each fixed nonzero $\br\in\R^2$, \eqref{eq:uniform_third_moment_bound} gives
    \[
        \frac{
            n\E\left[
                \left|
                    \br^\intercal
                    \boldsymbol{\mathcal I}_n^{-1/2}
                    \bS_{n,1}
                \right|^3
            \right]
        }{
            \left(n\norm{\br}^2\right)^{3/2}
        }
        \leq
        \frac{C}{\sqrt n}
        \longrightarrow0.
    \]
    Therefore,
    \[
        \frac1{\sqrt n}
        \br^\intercal
        \boldsymbol{\mathcal I}_n^{-\frac{1}{2}}
        \sum_{i=1}^n\bS_{n,i}
        \inD
        \mathcal N(0,\norm{\br}^2).
    \]
    The Cram\'er--Wold device completes the proof.
\end{proof}

\begin{lemm}
\label{lemm:gamma2_approx_root}
Let $X_{n,1},\ldots,X_{n,n}\iid\mathrm{Gamma}(\alpha_n,\beta)$, where $\beta>0$ is fixed and $\alpha_n\to\infty$, and write
\[
    \bar X_n\coloneqq\frac1n\sum_{i=1}^nX_{n,i}
    \dand
    G_n\coloneqq\left(\prod_{i=1}^nX_{n,i}\right)^{\frac{1}{n}}.
\]
Define
\[
    u_{1,n}
    \coloneqq
    \frac{1}{2\alpha_n}+\frac{1}{12\alpha_n^2},
    \quad
    u_{2,n}\coloneqq\frac{\beta}{\alpha_n},
    \quad
    \tilde u_{1,n}\coloneqq\log{\frac{\bar X_n}{G_n}}\dand
    \tilde u_{2,n}\coloneqq\frac1{\bar X_n}.
\]
Then the following hold:
\begin{enumerate}[label=(\roman*)]
    \item \label{lemm:gamma2_approx_root_first_rate} The first coordinate satisfies
    \[
        \tilde u_{1,n}-u_{1,n}
        =
        \Op{\frac{u_{1,n}}{\sqrt n}}
        +O\left(\frac{u_{1,n}}n\right)
        +O(u_{1,n}^4).
    \]

    \item \label{lemm:gamma2_approx_root_second_rate} The second coordinate satisfies
    \[
        \tilde u_{2,n}-u_{2,n}
        =
        \Op{\frac{u_{1,n}^{3/2}}{\sqrt n}}
        +O\left(\frac{u_{1,n}^2}{n}\right).
    \]

    \item \label{lemm:gamma2_approx_root_first_limit} If $\sqrt n\,u_{1,n}^3\to0$, then
    \[
        \frac{\sqrt n}{u_{1,n}}(\tilde u_{1,n}-u_{1,n})
        \inD\mathcal N(0,2).
    \]

    \item \label{lemm:gamma2_approx_root_second_limit} The second coordinate also satisfies
    \[
        \frac{\sqrt n}{u_{1,n}^{3/2}}(\tilde u_{2,n}-u_{2,n})
        \inD\mathcal N(0,8\beta^2).
    \]
\end{enumerate}
\end{lemm}

\begin{proof}
First note that $u_{1,n}\sim(2\alpha_n)^{-1}$; in particular, $u_{1,n}\asymp\alpha_n^{-1}$, and $n\alpha_n>2$ for all sufficiently large $n$, so \cref{lemm:gammaexps2} applies.

For \ref{lemm:gamma2_approx_root_first_rate}, parts \ref{lemm:ElogXG} and \ref{lemm:VarlogXG} of \cref{lemm:gammaexps2} give
\[
    \E\left[\tilde u_{1,n}\right]
    =
    \psi(n\alpha_n)-\psi(\alpha_n)-\log{n} \dand
    \Var(\tilde u_{1,n})
    =
    \frac{\psi_1(\alpha_n)}{n}-\psi_1(n\alpha_n).
\]
Substituting in the standard digamma expansion in \eqref{eq:psi_asympt} yields
\[
    \psi(n\alpha_n)-\psi(\alpha_n)-\log{n}
    =
    u_{1,n}
    -\frac{1}{2n\alpha_n}
    -\frac{1}{12n^2\alpha_n^2}
    +O\left(\frac{1}{\alpha_n^4}\right),
\]
where the remainder is uniform because $(n\alpha_n)^{-4}\leq\alpha_n^{-4}$. Hence
\begin{equation}\label{eq:E_tilde_u_1_minus_u_1}
    \E\left[\tilde u_{1,n}\right]-u_{1,n}
    =
    O\left(\frac{u_{1,n}}{n}\right)
    +O(u_{1,n}^4).
\end{equation}
The trigamma expansion \citep[Equation 6.4.12]{abramowitz1988handbook}
\[
    \psi_1(x)
    =
    \frac1x+\frac{1}{2x^2}+\frac{1}{6x^3}+O\left(\frac{1}{x^5}\right)
\]
gives
\[
\begin{aligned}
    \Var(\tilde u_{1,n})
    &=
    \frac{1-1/n}{2n\alpha_n^2}
    +
    \frac{1-1/n^2}{6n\alpha_n^3}
    +O\left(\frac{1}{n\alpha_n^5}\right)=
    O\left(\frac{u_{1,n}^2}{n}\right).
\end{aligned}
\]
Chebyshev's inequality therefore implies that
\[
    \tilde u_{1,n}-\E\left[\tilde u_{1,n}\right]
    =
    \Op{\frac{u_{1,n}}{\sqrt n}}.
\]
Combining this with \eqref{eq:E_tilde_u_1_minus_u_1} proves \ref{lemm:gamma2_approx_root_first_rate}.

For \ref{lemm:gamma2_approx_root_second_rate}, parts \ref{lemm:EXinv} and \ref{lemm:Varinv} of \cref{lemm:gammaexps2} give
\begin{equation}\label{eq:E_tilde_u_2_minus_u_2}
    \E\left[\tilde u_{2,n}\right]-u_{2,n}
    =
    \frac{n\beta}{n\alpha_n-1}-\frac{\beta}{\alpha_n}
    =
    \frac{\beta}{\alpha_n(n\alpha_n-1)}
    =
    O\left(\frac{u_{1,n}^2}{n}\right)
\end{equation}
and
\[
\begin{aligned}
    \Var(\tilde u_{2,n})
    &=
    \frac{n^2\beta^2}
    {(n\alpha_n-1)^2(n\alpha_n-2)}
    &=
    O\left(\frac{1}{n\alpha_n^3}\right)
    =
    O\left(\frac{u_{1,n}^3}{n}\right).
\end{aligned}
\]
Another application of Chebyshev's inequality yields
\[
    \tilde u_{2,n}-\E\left[\tilde u_{2,n}\right]
    =
    \Op{\frac{u_{1,n}^{3/2}}{\sqrt n}},
\]
and combining this with \eqref{eq:E_tilde_u_2_minus_u_2} proves \ref{lemm:gamma2_approx_root_second_rate}.

For \ref{lemm:gamma2_approx_root_first_limit}, let
\[
    \xi_{n,i}
    \coloneqq
    \log{\beta X_{n,i}}-\psi(\alpha_n)
    +1-\frac{\beta X_{n,i}}{\alpha_n} \dand
    v_n\coloneqq\Var(\xi_{n,i})
    =
    \psi_1(\alpha_n)-\frac1{\alpha_n}.
\]
where the variance identity follows from
\[
\begin{aligned}
\Var\left(\log{\beta X_{n,i}}\right)&=\psi_1(\alpha_n) \dand
\Var(\beta X_{n,i}/\alpha_n)=\operatorname{Cov}\left(\log{\beta X_{n,i}},
                            \beta X_{n,i}/\alpha_n\right)
=\frac1{\alpha_n}.
\end{aligned}
\]
Note that $\E\left[\xi_{n,i}\right]=0$. The uniform third moment bound \eqref{eq:uniform_third_moment_bound} in \cref{lemm:gamma2CLT} gives the Lyapunov condition for $\xi_{n,i}/\sqrt{v_n}$, and hence
\[
    \frac{1}{\sqrt{nv_n}}\sum_{i=1}^n\xi_{n,i}
    \inD\mathcal N(0,1).
\]
Now, set
\[
    Z_n
    =
    \frac{\beta\bar X_n}{\alpha_n}-1.
\]
Since $\E\left[Z_n\right]=0$ and $\Var(Z_n)=1/(n\alpha_n)$, Chebyshev's inequality gives $Z_n=\Op{(n\alpha_n)^{-1/2}}=\op{1}$. Moreover, the identities
\[
    \log{1+Z_n}
    -\frac1n\sum_{i=1}^n
      \left(\log{\beta X_{n,i}}-\psi(\alpha_n)\right)
    =
    \tilde u_{1,n}-\left(\log{\alpha_n}-\psi(\alpha_n)\right)
\]
and
\[
    Z_n
    -\frac1n\sum_{i=1}^n
      \left(\log{\beta X_{n,i}}-\psi(\alpha_n)\right)
    =
    -\frac1n\sum_{i=1}^n\xi_{n,i}
\]
yield that
\[
    \tilde u_{1,n}-\left(\log{\alpha_n}-\psi(\alpha_n)\right)
    =
    -\frac1n\sum_{i=1}^n\xi_{n,i}
    +\log{1+Z_n}-Z_n.
\]
On the event $|Z_n|\leq1/2$, Taylor's theorem gives $|\log{1+Z_n}-Z_n|\leq C Z_n^2$. This event has probability tending to one, and hence
\[
    \log{1+Z_n}-Z_n
    =
    \Op{\frac{1}{n\alpha_n}}.
\]
The expansions of $\psi$ and $\psi_1$ used previously also give
\[
    \frac{v_n}{u_{1,n}^2}\longrightarrow2
    \dand
    \log{\alpha_n}-\psi(\alpha_n)-u_{1,n}
    =
    O(u_{1,n}^4).
\]
Therefore,
\[
\begin{aligned}
    \frac{\sqrt n}{u_{1,n}}(\tilde u_{1,n}-u_{1,n})
    &=-\frac{\sqrt{v_n}}{u_{1,n}}
      \frac{1}{\sqrt{nv_n}}\sum_{i=1}^n\xi_{n,i}
      +\frac{\sqrt n}{u_{1,n}}
       \left(\log{\alpha_n}-\psi(\alpha_n)-u_{1,n}\right)
      +\frac{\sqrt n}{u_{1,n}}
       \left(\log{1+Z_n}-Z_n\right).
\end{aligned}
\]
Since $\sqrt n\,u_{1,n}^3\to0$ by assumption, the second term is $o(1)$. Since $\alpha_n u_{1,n}\to1/2$, the third term is
\[
    \Op{\frac{1}{\sqrt n\,\alpha_n u_{1,n}}}
    =
    \op{1}.
\]
Slutsky's theorem then proves \ref{lemm:gamma2_approx_root_first_limit}.

For \ref{lemm:gamma2_approx_root_second_limit}, let $\mu_n\coloneqq\alpha_n/\beta$ and $\Delta_n\coloneqq\bar X_n-\mu_n$. The corresponding central limit theorem for $(\beta X_{n,i}-\alpha_n)/\sqrt{\alpha_n}$, also established in \cref{lemm:gamma2CLT}, gives
\[
    \frac{\sqrt n\,\beta}{\sqrt{\alpha_n}}
    \Delta_n
    \inD\mathcal N(0,1).
\]
The central limit theorem implies
\[
    \Delta_n
    =\Op{\frac{\sqrt{\alpha_n}}{\beta\sqrt n}}.
\]
Since $\mu_n=\alpha_n/\beta$, it follows that
\[
    \frac{\Delta_n}{\mu_n}
    =\Op{\frac{1}{\sqrt{n\alpha_n}}}
    =\op{1},
\]
and hence $\PP(|\Delta_n|\leq\mu_n/2)\to1$. The identity
\[
    \frac1{\bar X_n}-\frac1{\mu_n}
    =
    -\frac{\Delta_n}{\mu_n^2}
    +\frac{\Delta_n^2}{\mu_n^2\bar X_n}
\]
shows that, on $\{|\Delta_n|\leq\mu_n/2\}$, we have $\bar X_n=\mu_n+\Delta_n\geq\mu_n/2$ and therefore
\[
    \left|\frac{\Delta_n^2}{\mu_n^2\bar X_n}\right|
    \leq
    \frac{2\Delta_n^2}{\mu_n^3}.
\]
Therefore
\[
    \frac{\Delta_n^2}{\mu_n^2\bar X_n}
    =
    \Op{\frac{1}{n\alpha_n^2}}.
\]
Since $\alpha_n^{3/2}u_{1,n}^{3/2}\to2^{-3/2}$, the scaled remainder satisfies
\[
    \frac{\sqrt n}{u_{1,n}^{3/2}}
    \Op{\frac{1}{n\alpha_n^2}}
    =
    \Op{\frac{1}{\sqrt{n\alpha_n}}}
    =
    \op{1}.
\]
It follows that
\[
    \frac{\sqrt n}{u_{1,n}^{3/2}}
    (\tilde u_{2,n}-u_{2,n})
    =
    -\frac{\beta}{\alpha_n^{3/2}u_{1,n}^{3/2}}
    \frac{\sqrt n\,\beta}{\sqrt{\alpha_n}}\Delta_n
    +\op{1}.
\]
The deterministic term converges to $-2^{3/2}\beta$, so Slutsky's theorem gives the limiting variance $8\beta^2$, thereby proving \ref{lemm:gamma2_approx_root_second_limit}.
\end{proof}

\begin{lemm}
\label{lemm:gamma2_inverse_chart}
Let $\beta>0$, $\alpha_n\to\infty$,
\[
    \bu_n
    =
    \begin{pmatrix}
        \dfrac{1}{(2\alpha_n)}+\dfrac{1}{(12\alpha_n^2)}\\[9pt]
        \dfrac{\beta}{\alpha_n}
    \end{pmatrix}
    =
    \begin{pmatrix}
        u_{1,n}\\
        u_{2,n}
    \end{pmatrix},
\]
and let
\[
    K_n
    \coloneqq
    \left[\frac{u_{1,n}}{2},2u_{1,n}\right]
    \times
    \left[\frac{u_{2,n}}{2},2u_{2,n}\right].
\]
For $\bF\coloneqq\bphi\inv$ as defined in \eqref{eq:gamma2_phiinv}, the following hold:
\begin{enumerate}[label=(\roman*)]
    \item \label{lemm:gamma2_inverse_chart_jacobian} The Jacobian satisfies
    \begin{equation}\label{eq:gamma2_inverse_chart_jacobian}
        \sup_{\bu\in K_n}\norm{\Jb_{\bF}(\bu)}
        =
        O\left(\frac{1}{u_{1,n}^2}\right),
    \end{equation}
    and therefore $\bF$ is $b_n$-Lipschitz on $K_n$ for some $b_n=O(u_{1,n}^{-2})$.

    \item \label{lemm:gamma2_inverse_chart_scalings} Define
    \[
        \Rb_n
        \coloneqq
        \operatorname{diag}
        \left(\frac{\sqrt n}{u_{1,n}},
        \frac{\sqrt n}{u_{1,n}^{3/2}}\right)\dand
        \Tb_n
        \coloneqq
        \operatorname{diag}(\sqrt n\,u_{1,n},\sqrt n).
    \]
    Then
    \[
        \Tb_n\Jb_{\bF}(\bu_n)\Rb_n\inv
        \longrightarrow
        \begin{pmatrix}
            -\dfrac{1}{2}&0\\[9pt]
            -\beta&0
        \end{pmatrix}.
    \]

    \item \label{lemm:gamma2_inverse_chart_remainder} If a random sequence $\bV_n$ taking values in $(0, \infty)^2$ satisfies $\PP(\bV_n\in K_n)\to1$ and $\Rb_n(\bV_n-\bu_n)=\Op{1}$, then
    \[
        \Tb_n
        \left[\bF(\bV_n)-\bF(\bu_n)
        -\Jb_{\bF}(\bu_n)(\bV_n-\bu_n)\right]
        =
        \op{1}.
    \]
\end{enumerate}
\end{lemm}

\begin{proof}
Write
\[
    d(u_1)\coloneqq-3+\sqrt{9+12u_1} \dand
    f(u_1)\coloneqq\frac1{d(u_1)},
\]
so $\bF(\bu)=\left(f(u_1),u_2f(u_1)\right)^\intercal$. Rationalizing $d$ gives
\[
    d(u_1)
    =
    \frac{12u_1}{3+\sqrt{9+12u_1}}.
\]
Thus we have $d(u_1)\asymp u_{1,n}$ and $\sqrt{9+12u_1}\asymp1$ uniformly on $K_n$. Moreover,
\[
    f'(u_1)
    =
    -\frac{6}{\sqrt{9+12u_1}\,d(u_1)^2}
\]
and
\[
    f''(u_1)
    =
    \frac{36}{(9+12u_1)^{3/2}d(u_1)^2}
    +
    \frac{72}{(9+12u_1)d(u_1)^3}.
\]
It follows that, uniformly on $K_n$,
\[
    f(u_1)=O\left(\frac{1}{u_{1,n}}\right)
    \dand
    f'(u_1)=O\left(\frac{1}{u_{1,n}^2}\right)
    \dand
    f''(u_1)=O\left(\frac{1}{u_{1,n}^3}\right).
\]
Also, $u_2=O(u_{1,n})$ because $u_{2,n}\asymp u_{1,n}$. Hence
\begin{equation}\label{eq:gamma2_inverse_chart_jacobian_expression}
    \Jb_{\bF}(\bu)
    =
    \begin{pmatrix}
        f'(u_1)&0\\
        u_2f'(u_1)&f(u_1)
    \end{pmatrix}
    =
    \begin{pmatrix}
        O\left(\dfrac{1}{u_{1,n}^2}\right)&0\\[9pt]
        O\left(\dfrac{1}{u_{1,n}}\right)&O\left(\dfrac{1}{u_{1,n}}\right)
    \end{pmatrix}
\end{equation}
uniformly on $K_n$, which proves \eqref{eq:gamma2_inverse_chart_jacobian}. Since $K_n$ is convex, the fundamental theorem of calculus gives, for $\bv,\bw\in K_n$,
\[
    \bF(\bv)-\bF(\bw)
    =
    \int_0^1
    \Jb_{\bF}\left[\bw+t(\bv-\bw)\right](\bv-\bw)\dif t,
\]
and therefore the Lipschitz constant for $\bF$ on $K_n$ can be taken of order $u_{1,n}^{-2}$, which proves \ref{lemm:gamma2_inverse_chart_jacobian}.

At $\bu_n$, we have $d(u_{1,n})=1/\alpha_n$ and $u_{2,n}=\beta d(u_{1,n})$. Since $\sqrt{9+12u_{1,n}}=3+1/\alpha_n$, direct substitution into the Jacobian expression in \eqref{eq:gamma2_inverse_chart_jacobian_expression} gives
\[
    \Tb_n\Jb_{\bF}(\bu_n)\Rb_n\inv
    =
    \begin{pmatrix}
        -\dfrac{6\alpha_n^2u_{1,n}^2}{3+1/\alpha_n}&0\\[9pt]
        -\dfrac{6\beta\alpha_n u_{1,n}}{3+1/\alpha_n}
        &\alpha_n u_{1,n}^{3/2}
    \end{pmatrix}
\]
and
\[
    \alpha_n u_{1,n}
    =
    \frac12+\frac{1}{12\alpha_n}\dand
    \alpha_n u_{1,n}^{3/2}
    =
    (\alpha_n u_{1,n})u_{1,n}^{1/2}
    \longrightarrow0,
\]
which together prove \ref{lemm:gamma2_inverse_chart_scalings}.

For \ref{lemm:gamma2_inverse_chart_remainder}, write $\delta_{j,n}\coloneqq V_{j,n}-u_{j,n}$ and $A_n\coloneqq\{\bV_n\in K_n\}$. The condition $\Rb_n(\bV_n-\bu_n)=\Op{1}$ gives
\[
    \delta_{1,n}
    =
    \Op{\frac{u_{1,n}}{\sqrt n}}\dand
    \delta_{2,n}
    =
    \Op{\frac{u_{1,n}^{3/2}}{\sqrt n}}.
\]
On $A_n$, the segment joining $\bu_n$ and $\bV_n$ lies in $K_n$. The integral form of Taylor's theorem therefore shows that the first component remainder equals
\[
    \delta_{1,n}^2
    \int_0^1(1-t)
    f''(u_{1,n}+t\delta_{1,n})\dif t
    =
    \Op{\frac{\delta_{1,n}^2}{u_{1,n}^3}}
    =
    \Op{\frac{1}{nu_{1,n}}}.
\]
After multiplying by $\sqrt n\,u_{1,n}$, this quantity is $\Op{n^{-1/2}}=\op{1}$. As for the second component, $F_2(\bu)=u_2f(u_1)$, so uniformly on $K_n$,
\[
    \partial_{11}F_2(\bu)
    =u_2f''(u_1)=O\left(\frac{1}{u_{1,n}^2}\right)
    \dand
    \partial_{12}F_2(\bu)
    =f'(u_1)=O\left(\frac{1}{u_{1,n}^2}\right)
    \dand
    \partial_{22}F_2(\bu)=0.
\]
Applying the same integral expansion, its remainder is
\[
    \Op{\frac{\delta_{1,n}^2}{u_{1,n}^2}}
    +
    \Op{\frac{|\delta_{1,n}\delta_{2,n}|}{u_{1,n}^2}}
    =
    \Op{\frac1n}
    +
    \Op{\frac{\sqrt{u_{1,n}}}{n}},
\]
which, after multiplying by $\sqrt n$, is
\[
    \Op{n^{-\frac{1}{2}}}
    +
    \Op{\sqrt{\frac{u_{1,n}}{n}}}
    =
    \op{1}.
\]
Finally, $\PP(A_n^c)\to0$, so the two scaled componentwise bounds both hold in probability without conditioning on $A_n$, which proves \ref{lemm:gamma2_inverse_chart_remainder}.
\end{proof}

\subsection{Gamma Regression}

\begin{lemm}\label{lemm:glm_estimator_rates}
Assume the gamma regression model and design conditions of \cref{sub:glmregression}, where $\varphi_n\downarrow0$ and $\bbeta_0\in\R^k$ is fixed. Then the following hold:
\begin{enumerate}[label=(\roman*)]
    \item The explicit estimators in \eqref{eq:glm_estimators} satisfy
    \begin{equation}\label{eq:glm_surrogate_rates}
        \tilde\bbeta_n-\bbeta_0
        =
        O(\varphi_n)
        +
        \Op{\sqrt{\frac{\varphi_n}{n}}},
    \end{equation}
    and
    \begin{equation}\label{eq:glm_surrogate_phi_rate}
        \tilde\varphi_n-\varphi_n
        =
        O\left(\varphi_n^2+\frac{\varphi_n}{n}\right)
        +
        \Op{\frac{\varphi_n}{\sqrt n}}.
    \end{equation}

    \item The exact MLE exists uniquely with probability tending to one and satisfies
    \begin{equation}\label{eq:glm_exact_rate1}
        \hat\bbeta_n-\bbeta_0
        =
        \Op{\sqrt{\frac{\varphi_n}{n}}}
    \end{equation}
    and
    \begin{equation}\label{eq:glm_exact_rate2}
        \hat\varphi_n-\varphi_n
        =
        \Op{\frac{\varphi_n}{\sqrt n}}.
    \end{equation}
\end{enumerate}
\end{lemm}

\begin{proof}
Write
\[
    W_{n,i}\coloneqq\log{Z_{n,i}}
    \dand
    \bW_n\coloneqq\bigl(W_{n,1},\ldots,W_{n,n}\bigr)^\intercal,
\]
and let
\[
    m_n\coloneqq\E\left[W_{n,i}\right]
    =
    \psi\left(\dfrac{1}{\varphi_n}\right)+\log{\varphi_n}.
\]
The standard digamma and trigamma expansions \citep[Formulas 6.3.18 and 6.4.12]{abramowitz1988handbook} give
\begin{equation}\label{eq:glm_loggamma_moments}
    m_n=-\frac{\varphi_n}{2}+O(\varphi_n^2)\dand
    \Var(W_{n,i})
    =
    \psi_1\left(\dfrac{1}{\varphi_n}\right)
    =
    \varphi_n+O(\varphi_n^2).
\end{equation}
By \cref{lemm:gammaexps}, the cumulants of $W_{n,i}$ of orders $r\geq2$ are $\psi_{r-1}(1/\varphi_n)$. The general identity $\mu_4=\kappa_4+3\kappa_2^2$ relating central moments and cumulants \citep{willinkRelationshipsCentralMoments2003} therefore gives
\[
    \E\left[(W_{n,i}-m_n)^4\right]
    =
    \psi_3\left(\dfrac{1}{\varphi_n}\right)
    +3\psi_1\left(\dfrac{1}{\varphi_n}\right)^2
    =
    O(\varphi_n^2).
\]
Hence, H\"older's inequality gives
\[
    \E\left[|W_{n,i}-m_n|^3\right]
    \leq
    \left(
        \E\left[(W_{n,i}-m_n)^4\right]
    \right)^{\frac{3}{4}}
    =
    O(\varphi_n^{\frac{3}{2}}).
\]
Since $m_n=O(\varphi_n)$ and
$|x+y|^3\leq 4(|x|^3+|y|^3)$, it follows that
\[
    \E\left[|W_{n,i}|^3\right]
    =
    O(\varphi_n^{\frac{3}{2}}).
\]

It remains to control the corresponding moment weighted by $Z_{n,i}$, which we do by exponentially tilting the gamma distribution. Set $a_n\coloneqq1/\varphi_n$, so that $Z_{n,i}\sim\mathrm{Gamma}(a_n,a_n)$ in the shape-rate parameterization. Because $\E\left[Z_{n,i}\right]=1$, multiplication of its density by $z$ defines a probability density. Moreover,
\[
    z\frac{a_n^{a_n}}{\Gamma(a_n)}
    z^{a_n-1}\exp{-a_nz}
    =
    \frac{a_n^{a_n+1}}{\Gamma(a_n+1)}
    z^{a_n}\exp{-a_nz},
\]
which is the density of a $\mathrm{Gamma}(a_n+1,a_n)$ random variable. Thus, if $Z_{n,i}^{\star}\sim\mathrm{Gamma}(a_n+1,a_n)$, then
\[
    \E\left[Z_{n,i}|W_{n,i}|^3\right]
    =
    \E\left[|\log{Z_{n,i}^{\star}}|^3\right].
\]
The mean of $\log{Z_{n,i}^{\star}}$ is
\[
    \psi(a_n+1)-\log{a_n}
    =
    m_n+\varphi_n
    =
    O(\varphi_n),
\]
and its fourth centered moment is
\[
    \psi_3(a_n+1)+3\psi_1(a_n+1)^2
    =
    O(\varphi_n^2).
\]
Applying the same Lyapunov bound to the centered logarithm and then using its $O(\varphi_n)$ mean yields
\[
    \E\left[|\log{Z_{n,i}^{\star}}|^3\right]
    =
    O(\varphi_n^{\frac{3}{2}}).
\]
Combining the weighted and unweighted bounds gives
\begin{equation}\label{eq:glm_weighted_log_moment}
    \E\left[(1+Z_{n,i})|W_{n,i}|^3\right]
    =
    O(\varphi_n^{\frac{3}{2}}).
\end{equation}

Next, we establish the rates of the explicit estimators $\tilde\bbeta_n$ and $\tilde\varphi_n$. Since our regression setup includes an intercept,
\begin{align*}
    \tilde\bbeta_n-\bbeta_0
    &=
    (\Xb_n^\intercal\Xb_n)^{-1}
    \Xb_n^\intercal\bW_n
    =
    m_n\be_0
    +
    (\Xb_n^\intercal\Xb_n)^{-1}\Xb_n^\intercal
    (\bW_n-m_n\vec{\bm 1}_n).
\end{align*}
Since the components of $\bW_n-m_n\vec{\bm 1}_n$ are independent with variance
$\psi_1(1/\varphi_n)$ and $n^{-1}\Xb_n^\intercal\Xb_n$, the second term has covariance
\begin{align*}
    &(\Xb_n^\intercal\Xb_n)^{-1}\Xb_n^\intercal
    \left(\psi_1\left(\dfrac{1}{\varphi_n}\right)\Ib_n\right)
    \Xb_n(\Xb_n^\intercal\Xb_n)^{-1}
    =
    \psi_1(1/\varphi_n)(\Xb_n^\intercal\Xb_n)^{-1}.
\end{align*}
The design assumptions and \eqref{eq:glm_loggamma_moments} therefore give
\[
    (\Xb_n^\intercal\Xb_n)^{-1}\Xb_n^\intercal
    (\bW_n-m_n\vec{\bm 1}_n)
    =
    \Op{\sqrt{\frac{\varphi_n}{n}}}.
\]
Together with $m_n=O(\varphi_n)$, this proves \eqref{eq:glm_surrogate_rates}. Moreover,
\[
    \tilde\varphi_n
    =
    \frac1n
    (\bW_n-m_n\vec{\bm 1}_n)^\intercal
    \Mb_n
    (\bW_n-m_n\vec{\bm 1}_n).
\]
Because $\Mb_n$ is an idempotent matrix of rank $n-k$, we have $\Tr(\Mb_n)=n-k$, and \citet[Theorem~1.5]{seberLinearRegressionAnalysis2012} gives
\[
    \E\left[\tilde\varphi_n\right]
    =
    \psi_1\left(\dfrac{1}{\varphi_n}\right)\left(1-\frac{k}{n}\right)
    =
    \varphi_n+O\left(\varphi_n^2+\frac{\varphi_n}{n}\right).
\]
Now, if $\bepsilon$ has iid centered components with variance $\sigma^2$ and fourth cumulant $\kappa_4 = \E\left[\varepsilon_i^4\right]-3\sigma^4$, then for any symmetric matrix $\Ab$, \citet[Theorem~1.6]{seberLinearRegressionAnalysis2012} gives
\begin{equation}\label{eq:quad_form_variance}
    \Var(\bepsilon^\intercal\Ab\bepsilon)
    =
    2\sigma^4\Tr(\Ab^2)
    +
    \kappa_4\sum_{i=1}^n A_{ii}^2,
\end{equation}
where $\sigma^2\coloneqq\psi_1(1/\varphi_n)$, $\kappa_4\coloneqq\psi_3(1/\varphi_n)$, and $\Ab\coloneqq n^{-1}\Mb_n$. Therefore,
\begin{align*}
    \Var(\tilde\varphi_n)
    &=
    \frac{2\psi_1(1/\varphi_n)^2}{n^2}
    \Tr(\Mb_n^2)
    +
    \frac{\psi_3(1/\varphi_n)}{n^2}
    \sum_{i=1}^n(\Mb_n)_{ii}^2.
\end{align*}
Since $\Mb_n^2=\Mb_n$, we have $\Tr(\Mb_n^2)=n-k$. Moreover, $0\leq(\Mb_n)_{ii}\leq1$, and hence
\[
    \sum_{i=1}^n(\Mb_n)_{ii}^2\leq n.
\]
Together with the polygamma expansions, these bounds give
\[
    \Var(\tilde\varphi_n)
    =
    O\left(\frac{\varphi_n^2}{n}\right),
\]
which proves \eqref{eq:glm_surrogate_phi_rate}.

Next, recall from \cref{sub:glmregression} that $h(z)=z-1-\log{z}$. For $\boldsymbol{\delta}\in\R^k$, define
\[
    A_n(\boldsymbol{\delta})
    \coloneqq
    \frac1n\sum_{i=1}^n
    h\left(Y_{n,i}\exp{-\bx_{n,i}^\intercal(\bbeta_0+\boldsymbol{\delta})}\right)
    =
    \frac1n\sum_{i=1}^n
    h\left(Z_{n,i}\exp{-\bx_{n,i}^\intercal\boldsymbol{\delta}}\right).
\]
For each fixed $\varphi>0$, apart from terms that do not depend on $\boldsymbol{\delta}$, $(n/\varphi)A_n(\boldsymbol{\delta})$ is the negative log likelihood evaluated at $\bbeta=\bbeta_0+\boldsymbol{\delta}$. Consequently, $\nabla A_n(\bzero)$ is the negative coefficient score at $\bbeta_0$, scaled by $\varphi/n$. More generally, the gradient and Hessian are, respectively,
\begin{align*}
    \nabla A_n(\boldsymbol{\delta})
    &=
    -\frac1n\Xb_n^\intercal
    \left(\bZ_n(\bbeta_0+\boldsymbol{\delta})-\vec{\bm 1}_n\right) \dand
    \nabla^2 A_n(\boldsymbol{\delta})
    =
    \frac1n\Xb_n^\intercal
    \operatorname{diag}\left(\bZ_n(\bbeta_0+\boldsymbol{\delta})\right)
    \Xb_n.
\end{align*}
Since $n^{-1}\Xb_n^\intercal\Xb_n\to\Sb$ and $\Sb$ is positive definite, $\Xb_n$ has full column rank for all sufficiently large $n$. For every nonzero $\bv\in\R^k$,
\[
    \bv^\intercal\nabla^2 A_n(\boldsymbol{\delta})\bv
    =
    \frac1n\sum_{i=1}^n
    Z_{n,i}\exp{-\bx_{n,i}^\intercal\boldsymbol{\delta}}
    (\bx_{n,i}^\intercal\bv)^2
    >0
\]
whenever $\Xb_n$ has full column rank; thus, $A_n$ is strictly convex. It is also \emph{coercive}, meaning that $A_n(\boldsymbol{\delta})\longrightarrow\infty$ as $\norm{\boldsymbol{\delta}}\longrightarrow\infty$. Indeed, if it were not coercive, then there would exist a sequence $\{\boldsymbol{\delta}_j\}$ with $\norm{\boldsymbol{\delta}_j}\to\infty$ along which $A_n$ remained
bounded. After passing to a subsequence, suppose that
\[
    \frac{\boldsymbol{\delta}_j}
         {\norm{\boldsymbol{\delta}_j}}
    \longrightarrow\bv
\]
for some unit vector $\bv$. Full column rank implies that $\bx_{n,i}^\intercal\bv\neq0$ for at least one $i$, and hence $|\bx_{n,i}^\intercal\boldsymbol{\delta}_j|\to\infty$. Since $h(t)\to\infty$ as either $t\downarrow0$ or $t\to\infty$, the corresponding summand in $A_n(\boldsymbol{\delta}_j)$ diverges. Because $h$ is nonnegative, this forces $A_n(\boldsymbol{\delta}_j)\to\infty$, a contradiction. Therefore, $A_n$ is coercive. A continuous coercive function on $\R^k$ attains its minimum by Weierstrass' theorem \citep[Proposition~3.2.1]{bertsekasConvexOptimizationTheory2009}, and strict convexity makes the minimizer $\widehat{\boldsymbol{\delta}}_n$ unique.

At zero,
\begin{equation}\label{eq:gradA0}
    \nabla A_n(\bzero)
    =
    -\frac1n\Xb_n^\intercal
    (\bZ_n-\vec{\bm 1}_n)
    =
    \Op{\sqrt{\frac{\varphi_n}{n}}}.
\end{equation}
Moreover,
\[
    \nabla^2 A_n(\bzero)-\frac1n\Xb_n^\intercal\Xb_n
    =
    \frac1n\sum_{i=1}^n
    (Z_{n,i}-1)\bx_{n,i}\bx_{n,i}^\intercal.
\]
By independence, $\E\left[Z_{n,i}-1\right]=0$, $\Var(Z_{n,i})=\varphi_n$, and
\begin{align*}
    \E\left[
        \norm{
            \nabla^2 A_n(\bzero)-\frac1n\Xb_n^\intercal\Xb_n
        }_{\mathrm F}^2
    \right]
    &=
    \frac{\varphi_n}{n^2}
    \sum_{i=1}^n
    \norm{\bx_{n,i}\bx_{n,i}^\intercal}_{\mathrm F}^2 
    =
    \frac{\varphi_n}{n^2}
    \sum_{i=1}^n\norm{\bx_{n,i}}^4
    \leq
    \frac{C^4\varphi_n}{n},
\end{align*}
where the constant $C$ may be chosen so that $\max_{1\leq i\leq n}\norm{\bx_{n,i}}\leq C$ for every $n$ by \eqref{eq:glm_design}. It follows that
\[
    \nabla^2 A_n(\bzero)
    =
    \frac1n\Xb_n^\intercal\Xb_n
    +
    \Op{\sqrt{\frac{\varphi_n}{n}}}
    =
    \frac1n\Xb_n^\intercal\Xb_n+\op{1}.
\]

Now suppose that
$\norm{\boldsymbol{\delta}}\leq\varepsilon$. Then
\[
    |\bx_{n,i}^\intercal\boldsymbol{\delta}|
    \leq
    \norm{\bx_{n,i}}\norm{\boldsymbol{\delta}}
    \leq C\varepsilon,
\]
and hence
\[
    \exp{-\bx_{n,i}^\intercal\boldsymbol{\delta}}
    \geq
    \exp{-C\varepsilon}.
\]
Consequently, for every $\bv\in\R^k$,
\begin{align*}
    \bv^\intercal
    \nabla^2 A_n(\boldsymbol{\delta})\bv
    &=
    \frac1n\sum_{i=1}^n
    Z_{n,i}
    \exp{-\bx_{n,i}^\intercal\boldsymbol{\delta}}
    (\bx_{n,i}^\intercal\bv)^2 \\
    &\geq
    \exp{-C\varepsilon}
    \frac1n\sum_{i=1}^n
    Z_{n,i}(\bx_{n,i}^\intercal\bv)^2 \\
    &=
    \exp{-C\varepsilon}
    \bv^\intercal\nabla^2 A_n(\bzero)\bv.
\end{align*}
Therefore,
\[
    \nabla^2 A_n(\boldsymbol{\delta})
    \succeq
    \exp{-C\varepsilon}\nabla^2 A_n(\bzero)
\]
uniformly over $\norm{\boldsymbol{\delta}}\leq\varepsilon$. Since $n^{-1}\Xb_n^\intercal\Xb_n\to\Sb$ and $\Sb$ is positive definite, there is a constant $c>0$ such that, with probability tending to one,
\[
    \inf_{\norm{\boldsymbol{\delta}}\leq\varepsilon}
    \lambda_{\min}
    \left(\nabla^2 A_n(\boldsymbol{\delta})\right)
    \geq c.
\]
On this event, for $\norm{\boldsymbol{\delta}}=\varepsilon$, Taylor's theorem gives
\begin{align*}
    A_n(\boldsymbol{\delta})-A_n(\bzero)
    &=
    \nabla A_n(\bzero)^\intercal\boldsymbol{\delta}
    +
    \int_0^1(1-t)
    \boldsymbol{\delta}^\intercal
    \nabla^2 A_n(t\boldsymbol{\delta})
    \boldsymbol{\delta}\dif t 
    \geq
    -\varepsilon\norm{\nabla A_n(\bzero)}
    +\frac{c\varepsilon^2}{2}.
\end{align*}
By \eqref{eq:gradA0}, the right side is positive with probability tending to one. Convexity therefore implies that $\norm{\widehat{\boldsymbol{\delta}}_n}<\varepsilon$ with probability tending to one. On this event, the first-order condition gives
\[
    \bzero
    =
    \nabla A_n(\bzero)
    +
    \left(
        \int_0^1
        \nabla^2 A_n(t\widehat{\boldsymbol{\delta}}_n)\dif t
    \right)
    \widehat{\boldsymbol{\delta}}_n.
\]
The matrix in braces has smallest eigenvalue at least $c$, so
\[
    \norm{\widehat{\boldsymbol{\delta}}_n}
    \leq
    \frac{1}{c}\norm{\nabla A_n(\bzero)}
    =
    \Op{\sqrt{\frac{\varphi_n}{n}}}.
\]
Consequently,
\begin{equation}\label{eq:glm_beta_exact_localization}
    \hat\bbeta_n-\bbeta_0
    =
    \widehat{\boldsymbol{\delta}}_n
    =
    \Op{\sqrt{\frac{\varphi_n}{n}}},
\end{equation}
which is \eqref{eq:glm_exact_rate1}.

Finally, since $h(t)\geq0$ with equality only when $t=1$, we have that $A_n(\widehat{\boldsymbol{\delta}}_n)=0$ only if, for every $i=1,\ldots,n$,
\[
    Z_{n,i}
    \exp{-\bx_{n,i}^\intercal
    \widehat{\boldsymbol{\delta}}_n}
    =1.
\]
Equivalently, $\bL_n$ must lie in the column space of $\Xb_n$, and for $n>k$, this event has probability zero because $\bL_n$ has a continuous distribution and $\operatorname{col}(\Xb_n)$ is a proper linear subspace of $\R^n$. Therefore, $A_n(\widehat{\boldsymbol{\delta}}_n)>0$ almost surely for all sufficiently large $n$. Now, the dispersion score is
\begin{equation}\label{eq:glm_profile_dispersion}
    g(\hat\varphi_n)
    =
    A_n(\widehat{\boldsymbol{\delta}}_n) \quad \text{ where } \quad
    g(\varphi)=-\log{\varphi}-\psi\left(\dfrac{1}{\varphi}\right).
\end{equation}
The function $g$ is continuous and strictly increasing on $(0,\infty)$; indeed, for $x>0$,
\[
    \psi_1(x)
    =
    \sum_{j=0}^{\infty}\frac{1}{(x+j)^2}
    >
    \int_0^\infty\frac{\dif t}{(x+t)^2}
    =
    \frac1x,
\]
so
\[
    g'(\varphi)
    =
    \frac{\psi_1(1/\varphi)-\varphi}{\varphi^2}
    >0.
\]
Moreover, $g(\varphi)\to0$ as $\varphi\downarrow0$ by the expansion $g(\varphi)=\varphi/2+O(\varphi^2)$. The digamma recurrence
\[
    \psi\left(\dfrac{1}{\varphi}\right)
    =
    \psi\left(1+\dfrac{1}{\varphi}\right)-\varphi
\]
\citep[Equation~6.3.5]{abramowitz1988handbook}, together with $\psi(1)=-\gamma$ \citep[Equation~6.3.2]{abramowitz1988handbook}, gives
\[
    g(\varphi)
    =
    \varphi-\log{\varphi}-\psi\left(1+\dfrac{1}{\varphi}\right)
    =
    \varphi-\log{\varphi}+\gamma+o(1)
\]
as $\varphi\to\infty$. Thus $g(\varphi)\to\infty$, so \eqref{eq:glm_profile_dispersion} has a unique solution. Equivalently, writing $a\coloneqq1/\varphi$, the derivative of the profile log likelihood with respect to $a$ is $n\left(g(1/a)-A_n(\widehat{\boldsymbol{\delta}}_n)\right)$, which changes from positive to negative at this solution; it is therefore the unique maximizer of the profile log likelihood. At the true coefficient,
\[
    \E\left[A_n(\bzero)\right]=g(\varphi_n).
\]
Using $\operatorname{Cov}(Z_{n,i},\log{Z_{n,i}})=\varphi_n$, we also obtain
\[
    \Var\left(h(Z_{n,i})\right)
    =
    \psi_1\left(\dfrac{1}{\varphi_n}\right)-\varphi_n
    =
    \frac{\varphi_n^2}{2}+O(\varphi_n^3).
\]
It follows that
\[
    A_n(\bzero)
    =
    g(\varphi_n)
    +
    \Op{\frac{\varphi_n}{\sqrt n}}.
\]
A Taylor expansion of $A_n$ between zero and $\widehat{\boldsymbol{\delta}}_n$, together with \eqref{eq:glm_beta_exact_localization}, then gives
\begin{align*}
    A_n(\widehat{\boldsymbol{\delta}}_n)-A_n(\bzero)
    &=
    \nabla A_n(\bzero)^\intercal
    \widehat{\boldsymbol{\delta}}_n 
    +
    \int_0^1(1-t)
    \widehat{\boldsymbol{\delta}}_n^\intercal
    \nabla^2 A_n(t\widehat{\boldsymbol{\delta}}_n)
    \widehat{\boldsymbol{\delta}}_n\dif t 
    =
    \Op{\frac{\varphi_n}{n}},
\end{align*}
where we used the fact that the Hessian is uniformly $\Op{1}$ along the line segment joining zero and $\widehat{\boldsymbol{\delta}}_n$, which follows from the bounded design condition and \eqref{eq:glm_beta_exact_localization}. Consequently,
\[
    g(\hat\varphi_n)-g(\varphi_n)
    =
    \Op{\frac{\varphi_n}{\sqrt n}}.
\]
The preceding bounds also imply $g(\hat\varphi_n)=\Op{\varphi_n}$, and hence $\hat\varphi_n\inP0$. The trigamma expansion gives $g'(\varphi)=1/2+O(\varphi)$ as $\varphi\downarrow0$. The mean value theorem applied to \eqref{eq:glm_profile_dispersion} therefore yields
\[
    \hat\varphi_n-\varphi_n
    =
    \Op{\frac{\varphi_n}{\sqrt n}},
\]
which is \eqref{eq:glm_exact_rate2}.
\end{proof}

\begin{lemm}\label{lemm:glm_localization_remainder}
Under the assumptions of \cref{lemm:glm_estimator_rates} and \cref{sub:glmregression}, let
\[
    \btheta_n=(\varphi_n,\bbeta_0^\intercal)^\intercal
    \dand
    \bu_n=\bphi(\btheta_n).
\]
Write $\mathbb P_n\coloneqq\mathbb P_{n,\btheta_n}$ for the true joint law.
Then there are deterministic convex sets $D_n$ containing $\bu_n$ such that
\[
    \mathbb P_n
    \left(\hat\bu_n,\tilde\bu_n\in D_n\right)
    \longrightarrow1\dand
    \sup_{\bu\in D_n}\norm{\bu}=O(\varphi_n).
\]
Moreover, for $\bQ_n$ and $\tilde\bQ_n$ as defined in \eqref{eq:glm_exact_Q} and \eqref{eq:glm_surrogate_Q},
\begin{equation}\label{eq:glm_appendix_uniform_remainder}
    \sup_{\bu\in D_n}
    \norm{\tilde\bQ_n(\bu)-\bQ_n(\bu)}
    =
    \op{\varphi_n}.
\end{equation}
\end{lemm}

\begin{proof}

Put
\[
    \eta_{\varphi,n}\coloneqq\varphi_n+n^{-\frac{1}{2}}
    \dand
    \eta_{\beta,n}\coloneqq\varphi_n+\sqrt{\frac{\varphi_n}{n}}.
\]
Because
\[
    \frac{\eta_{\beta,n}^2}{\varphi_n}
    =
    O\left(\varphi_n+\sqrt{\frac{\varphi_n}{n}}+\frac1n\right)
    \longrightarrow0,
\]
we may choose a deterministic $t_n\to\infty$ sufficiently slowly that
\[
    r_{\varphi,n}\coloneqq t_n\eta_{\varphi,n}\longrightarrow0
    \dand
    r_{\beta,n}\coloneqq t_n\eta_{\beta,n}\longrightarrow0
    \dand
    r_{\beta,n}^2=o(\varphi_n).
\]
With $u_{1,n}\coloneqq\varphi_n/2$, define
\begin{equation}\label{eq:glm_Dn}
    D_n
    \coloneqq
    \left\{
        \bu=(u_1,\bu_2^\intercal)^\intercal:
        \left|\frac{u_1}{u_{1,n}}-1\right|
        \leq r_{\varphi,n},
        \quad
        \norm{\bu_2-u_1\bbeta_0}
        \leq r_{\beta,n}u_1
    \right\}.
\end{equation}
For all sufficiently large $n$, $r_{\varphi,n}<1$, so $D_n\subset U=(0,\infty)\times\R^k$ and $\bu_n\in D_n$. To see that $D_n$ is convex, let $\bv,\bw\in D_n$, write $\bv=(v_1,\bv_2^\intercal)^\intercal$ and $\bw=(w_1,\bw_2^\intercal)^\intercal$, and take $\lambda\in[0,1]$. The first coordinate of $\lambda\bv+(1-\lambda)\bw$ satisfies the first constraint in \eqref{eq:glm_Dn}, and the triangle inequality gives
\begin{align*}
    &\norm{
        \lambda\bv_2+(1-\lambda)\bw_2
        -\left(\lambda v_1+(1-\lambda)w_1\right)\bbeta_0
    } 
    \leq
    \lambda r_{\beta,n}v_1
    +(1-\lambda)r_{\beta,n}w_1
    =
    r_{\beta,n}\left(\lambda v_1+(1-\lambda)w_1\right).
\end{align*}
Thus $\lambda\bv+(1-\lambda)\bw\in D_n$.

For either estimator, denote its original coordinates by $(\varphi_n^\star,\bbeta_n^\star)$ and write
\[
    \bu_n^\star
    \coloneqq
    \bphi(\varphi_n^\star,\bbeta_n^\star)
    =
    \begin{pmatrix}
        u_{1,n}^\star\\
        \bu_{2,n}^\star
    \end{pmatrix}.
\]
These chart coordinates satisfy
\[
    \left|
        \frac{u_{1,n}^\star}{u_{1,n}}-1
    \right|
    =
    \left|
        \frac{\varphi_n^\star}{\varphi_n}-1
    \right|
    \dand
    \frac{
        \norm{\bu_{2,n}^\star-u_{1,n}^\star\bbeta_0}
    }{u_{1,n}^\star}
    =
    \norm{\bbeta_n^\star-\bbeta_0}.
\]
The rates in \cref{lemm:glm_estimator_rates} show that the two quantities on the right are $\Op{\eta_{\varphi,n}}$ and $\Op{\eta_{\beta,n}}$, respectively. Since $\eta_{\varphi,n}=o(r_{\varphi,n})$ and $\eta_{\beta,n}=o(r_{\beta,n})$, the chart coordinates of both estimators belong to $D_n$ with probability tending to one. Finally, the definition of $D_n$ gives
\[
    \sup_{\bu\in D_n}\norm{\bu}=O(\varphi_n).
\]

Finally, for $\bu\in D_n$, set
\[
    \bbeta(\bu)\coloneqq\frac{\bu_2}{u_1}
    \dand
    \boldsymbol{\delta}(\bu)\coloneqq\bbeta(\bu)-\bbeta_0,
\]
so that $\sup_{\bu\in D_n}\norm{\boldsymbol{\delta}(\bu)}\leq r_{\beta,n}$. Recall that
\[
    \br_n(\bu)
    =
    \bL_n-\Xb_n\frac{\bu_2}{u_1}
\]
is the log-residual vector. The least-squares decomposition gives
\begin{align*}
    \frac1{2n}\norm{\br_n(\bu)}^2
    &=
    \tilde u_{1,n}
    +
    \frac12
    \left(\bbeta(\bu)-\tilde\bbeta_n\right)^\intercal
    \frac1n\Xb_n^\intercal\Xb_n
    \left(\bbeta(\bu)-\tilde\bbeta_n\right).
\end{align*}
Since $\tilde\bbeta_n-\bbeta_0=\Op{\eta_{\beta,n}}=\op{r_{\beta,n}}$ and $n^{-1}\Xb_n^\intercal\Xb_n\to\Sb$, the second term is $\Op{r_{\beta,n}^2}$ uniformly on $D_n$.

Recalling that $h(z)=z-1-\log{z}$, we have
\[
    h(e^x)-\frac{x^2}{2}
    =
    e^x-1-x-\frac{x^2}{2}.
\]
Taylor's theorem gives, for some fixed $C>1$ and every $x\in\R$, the global bound
\[
    \left|h(e^x)-\frac{x^2}{2}\right|
    \leq C(1+e^x)|x|^3.
\]
Writing $W_{n,i}\coloneqq\log{Z_{n,i}}$ and increasing $C$ if necessary, boundedness of the design matrix and $r_{\beta,n}\to0$ give, for all sufficiently large $n$ and every $\bu\in D_n$,
\[
    |r_{n,i}(\bu)-W_{n,i}|
    =
    |\bx_{n,i}^\intercal\boldsymbol{\delta}(\bu)|
    \leq Cr_{\beta,n}
    \dand
    \exp{r_{n,i}(\bu)}
    =
    Z_{n,i}\exp{-\bx_{n,i}^\intercal\boldsymbol{\delta}(\bu)}
    \leq CZ_{n,i}.
\]
Since $(a+b)^3\leq4(a^3+b^3)$ for $a,b\geq0$, these bounds imply
\begin{align*}
    &\sup_{\bu\in D_n}
    \frac1n\sum_{i=1}^n
    (1+\exp{r_{n,i}(\bu)})|r_{n,i}(\bu)|^3 \leq
    4C^4\left(
        \frac1n\sum_{i=1}^n
        (1+Z_{n,i})|W_{n,i}|^3
        +
        r_{\beta,n}^3
        \frac1n\sum_{i=1}^n(1+Z_{n,i})
    \right).
\end{align*}
By \eqref{eq:glm_weighted_log_moment} and Markov's inequality, the first average on the right is $\Op{\varphi_n^{3/2}}$, while the second is $\Op{1}$. Therefore,
\begin{equation}\label{eq:glm_h_uniform}
    \sup_{\bu\in D_n}
    \frac1n\sum_{i=1}^n
    \left|
        h\left(\exp{r_{n,i}(\bu)}\right)
        -\frac12r_{n,i}(\bu)^2
    \right|
    =
    \Op{\varphi_n^{\frac{3}{2}}+r_{\beta,n}^3}.
\end{equation}
The first coordinate difference can be written as
\begin{align*}
    \tilde Q_{n,1}(\bu)-Q_{n,1}(\bu)
    &=
    \left(u_1-g(2u_1)\right) 
    +
    \left[
        \frac1n\sum_{i=1}^n
        h\left(\exp{r_{n,i}(\bu)}\right)
        -\frac1{2n}\norm{\br_n(\bu)}^2
    \right] 
    +
    \left[ \vphantom{\frac1n\sum_{i=1}^n}
        \frac1{2n}\norm{\br_n(\bu)}^2
        -\tilde u_{1,n}
    \right].
\end{align*}
Because $u_1=O(\varphi_n)$ uniformly on $D_n$, the expansion $g(2u_1)=u_1+O(u_1^2)$ shows that these three terms are respectively $O(\varphi_n^2)$, $\Op{\varphi_n^{3/2}+r_{\beta,n}^3}$, and $\Op{r_{\beta,n}^2}$, uniformly on $D_n$. Since $r_{\beta,n}\to0$,
\[
    \sup_{\bu\in D_n}
    |\tilde Q_{n,1}(\bu)-Q_{n,1}(\bu)|
    =
    \Op{\varphi_n^{\frac{3}{2}}+r_{\beta,n}^2}.
\]
For the remaining coordinates, the identity $e^x-1-x=h(e^x)$ gives
\[
    \tilde\bQ_{n,2}(\bu)-\bQ_{n,2}(\bu)
    =
    -\frac{u_1}{n}\Xb_n^\intercal
    \left[
        \bZ_n\left(\bbeta(\bu)\right)
        -\vec{\bm 1}_n
        -\br_n(\bu)
    \right],
\]
where the $i$th component of the bracketed vector is $h(\exp{r_{n,i}(\bu)})$. The least-squares decomposition and \eqref{eq:glm_h_uniform} give
\[
    \sup_{\bu\in D_n}
    \frac1n\sum_{i=1}^n
    h\left(\exp{r_{n,i}(\bu)}\right)
    =
    \Op{\varphi_n+r_{\beta,n}^2}.
\]
Since $u_1=O(\varphi_n)$ uniformly on $D_n$ and the design rows are bounded,
\[
    \sup_{\bu\in D_n}
    \norm{\tilde\bQ_{n,2}(\bu)-\bQ_{n,2}(\bu)}
    =
    \Op{\varphi_n(\varphi_n+r_{\beta,n}^2)}.
\]
Therefore
\[
    \sup_{\bu\in D_n}
    \norm{\tilde\bQ_n(\bu)-\bQ_n(\bu)}
    =
    \Op{\varphi_n^{\frac{3}{2}}+r_{\beta,n}^2}
    =
    \op{\varphi_n},
\]
which is \eqref{eq:glm_appendix_uniform_remainder}.
\end{proof}

\begin{lemm}\label{lemm:glm_joint_limits}
Under the assumptions of \cref{sub:glmregression}, the approximate and exact estimators satisfy \eqref{eq:glm_joint_approx_centered} and \eqref{eq:glm_joint_exact}, respectively. Moreover, if $n\varphi_n^4\to0$, the bias-corrected estimators in \eqref{eq:glm_bias_corrected_estimators} satisfy \eqref{eq:glm_joint_bias_corrected} and are asymptotically equivalent to the exact MLE as in \eqref{eq:glm_bias_corrected_equivalence}.
\end{lemm}

\begin{proof}
Continue to write
\[
    W_{n,i}=\log{Z_{n,i}}
    \dand
    m_n=\E\left[W_{n,i}\right]
    \dand
    \epsilon_{n,i}\coloneqq W_{n,i}-m_n
    \dand
    v_n\coloneqq\Var(W_{n,i}).
\]
We first collect the moments needed to compute the limiting covariance matrices and verify the Lyapunov conditions. By \cref{lemm:gammaexps}, the cumulant generating function of $W_{n,i}$ is
\[
    K_n(t)
    \coloneqq
    \log{\E\left[e^{tW_{n,i}}\right]}
    =
    \log{\Gamma\left(\frac{1}{\varphi_n}+t\right)}
    -\log{\Gamma\left(\frac{1}{\varphi_n}\right)}
    +t\log{\varphi_n}.
\]
Thus, its cumulants of orders $r\geq2$ are $\psi_{r-1}(1/\varphi_n)$. Standard polygamma expansions give the following four expansions:
\begin{equation}\label{eq:glm_joint_log_moments}
\begin{aligned}
    m_n
    &=\psi\left(\dfrac{1}{\varphi_n}\right)+\log{\varphi_n}
      =-\frac{\varphi_n}{2}
      -\frac{\varphi_n^2}{12}
      +O(\varphi_n^4),\\
    v_n
    &=\psi_1\left(\dfrac{1}{\varphi_n}\right)
      =\varphi_n
      +\frac{\varphi_n^2}{2}
      +\frac{\varphi_n^3}{6}
      +O(\varphi_n^5),\\
    \E\left[\epsilon_{n,i}^3\right]
    &=\psi_2\left(\dfrac{1}{\varphi_n}\right)
      =-\varphi_n^2+O(\varphi_n^3),\\
    \E\left[\epsilon_{n,i}^4\right]
    &=\psi_3\left(\dfrac{1}{\varphi_n}\right)+3v_n^2
      =3\varphi_n^2+O(\varphi_n^3).
\end{aligned}
\end{equation}
These expansions will determine the covariance matrices below. For the Lyapunov conditions, we will use the three bounds
\begin{equation}\label{eq:glm_epsilon_sixth_bound}
    \E\left[|\epsilon_{n,i}|^6\right]
    =O(\varphi_n^3),
\end{equation}
\begin{equation}\label{eq:glm_Z_third_bound}
    \E\left[|Z_{n,i}-1|^3\right]
    =O(\varphi_n^{\frac{3}{2}}),
\end{equation}
and
\begin{equation}\label{eq:glm_h_centered_third_bound}
    \E\left[
        |h(Z_{n,i})-g(\varphi_n)|^3
    \right]
    =O(\varphi_n^3).
\end{equation}
The identity \citep{willinkRelationshipsCentralMoments2003}
\begin{equation}\label{eq:glm_sixth_moment_cumulant_identity}
    \mu_6
    =
    \kappa_6+15\kappa_4\kappa_2+10\kappa_3^2+15\kappa_2^3
\end{equation}
expresses the sixth centered moment of $W_{n,i}$ as
\begin{equation}\label{eq:glm_epsilon_sixth_moment}
    \E\left[\epsilon_{n,i}^6\right]
    =
    \psi_5\left(\dfrac{1}{\varphi_n}\right)
    +15\psi_3\left(\dfrac{1}{\varphi_n}\right)v_n
    +10\psi_2\left(\dfrac{1}{\varphi_n}\right)^2
    +15v_n^3.
\end{equation}
The polygamma expansions \eqref{eq:glm_joint_log_moments} therefore give $\E\left[|\epsilon_{n,i}|^6\right]=O(\varphi_n^3)$, which proves \eqref{eq:glm_epsilon_sixth_bound}. Next, the fourth centered moment of $Z_{n,i}$ is
\[
    \E\left[(Z_{n,i}-1)^4\right]
    =
    3\varphi_n^2+6\varphi_n^3.
\]
H\"older's inequality gives
\[
    \E\left[|Z_{n,i}-1|^3\right]
    \leq
    \left(
        \E\left[(Z_{n,i}-1)^4\right]
    \right)^{\frac{3}{4}}
    =
    O(\varphi_n^{\frac{3}{2}}),
\]
which proves \eqref{eq:glm_Z_third_bound}.

It remains to control the third absolute moment involving $h$. Since
\[
    |h(e^w)|^3
    \leq
    C(1+e^{3w})|w|^6,
\]
we have
\begin{equation}\label{eq:glm_h_third_moment_bound}
    \E\left[|h(Z_{n,i})|^3\right]
    \leq
    C\left(
        \E\left[|W_{n,i}|^6\right]
        +
        \E\left[Z_{n,i}^3|W_{n,i}|^6\right]
    \right).
\end{equation}
The bound in \eqref{eq:glm_epsilon_sixth_bound}, together with $m_n=O(\varphi_n)$, gives
$\E\left[|W_{n,i}|^6\right]=O(\varphi_n^3)$. To treat the weighted term, set $a_n\coloneqq1/\varphi_n$ and let
$Z_{n,i}^{\star}\sim\mathrm{Gamma}(a_n+3,a_n)$. Multiplying the $\mathrm{Gamma}(a_n,a_n)$ density by $z^3$ and dividing by $\E[Z_{n,i}^3]$ yields the $\mathrm{Gamma}(a_n+3,a_n)$ density. Hence,
\[
    \E\left[Z_{n,i}^3|W_{n,i}|^6\right]
    =
    \E\left[Z_{n,i}^3\right]
    \E\left[|\log{Z_{n,i}^{\star}}|^6\right].
\]
Here
\[
    \E\left[Z_{n,i}^3\right]
    =
    (1+\varphi_n)(1+2\varphi_n)
    =O(1),
\]
while the mean of $\log{Z_{n,i}^{\star}}$ is
$\psi(a_n+3)-\log{a_n}=O(\varphi_n)$. Applying the identity \eqref{eq:glm_sixth_moment_cumulant_identity}, with the cumulants of orders $r\geq2$ now given by $\psi_{r-1}(a_n+3)$, yields
\[
    \E\left[
        \left|
            \log{Z_{n,i}^{\star}}
            -\E\left[\log{Z_{n,i}^{\star}}\right]
        \right|^6
    \right]
    =
    O(\varphi_n^3).
\]
Together with the $O(\varphi_n)$ mean, this gives $\E\left[|\log{Z_{n,i}^{\star}}|^6\right]=O(\varphi_n^3)$. It follows from \eqref{eq:glm_h_third_moment_bound} that
$\E\left[|h(Z_{n,i})|^3\right]=O(\varphi_n^3)$. Finally,
$g(\varphi_n)=\E\left[h(Z_{n,i})\right]=O(\varphi_n)$, so
\[
    \E\left[
        |h(Z_{n,i})-g(\varphi_n)|^3
    \right]
    =
    O(\varphi_n^3),
\]
which proves \eqref{eq:glm_h_centered_third_bound}. To connect these bounds with the arrays below, note that \eqref{eq:glm_epsilon_sixth_bound} also implies
\begin{equation}\label{eq:glm_epsilon_lyapunov_bounds}
\begin{aligned}
    \E\left[|\epsilon_{n,i}|^3\right]
    =O(\varphi_n^{\frac{3}{2}}) \dand
    \E\left[|\epsilon_{n,i}^2-v_n|^3\right]
    =O(\varphi_n^3).
\end{aligned}
\end{equation}
The bounds in \eqref{eq:glm_epsilon_lyapunov_bounds} apply to the approximate estimator array, while \eqref{eq:glm_Z_third_bound} and \eqref{eq:glm_h_centered_third_bound} apply to the exact score array.

We first prove the approximate estimator limit. Since the first column of $\Xb_n$ is $\vec{\bm 1}_n$, we have
\begin{equation}\label{eq:glm_joint_beta_decomp}
    \tilde\bbeta_n-\bbeta_0-m_n\be_0
    =
    (\Xb_n^\intercal\Xb_n)^{-1}\Xb_n^\intercal\bepsilon_n
    \dand
    \bepsilon_n=(\epsilon_{n,1},\ldots,\epsilon_{n,n})^\intercal
\end{equation}
and
\[
    \tilde\varphi_n
    =\frac1n\bepsilon_n^\intercal\Mb_n\bepsilon_n.
\]
Since $\E[\bepsilon_n]=\bzero$, $\Var(\bepsilon_n)=v_n\Ib_n$, and $\Tr(\Mb_n)=n-k$, \citet[Theorem~1.5]{seberLinearRegressionAnalysis2012} gives
\[
    \E\left[\tilde\varphi_n\right]
    =c_n=\left(1-\frac{k}{n}\right)v_n.
\]
The design assumptions imply
\[
    \max_{1\leq i\leq n}(\Hb_n)_{ii}=O\left(\frac{1}{n}\right)\dand
    \Tr(\Hb_n)=\Tr(\Hb_n^2)=k.
\]
The variance formula for quadratic forms \eqref{eq:quad_form_variance} and \eqref{eq:glm_joint_log_moments} give
\[
\begin{aligned}
    \Var(\bepsilon_n^\intercal\Hb_n\bepsilon_n)
    &=
    2v_n^2\Tr(\Hb_n^2)
    +\psi_3\left(\dfrac{1}{\varphi_n}\right)
    \sum_{i=1}^n(\Hb_n)_{ii}^2
    =
    O(\varphi_n^2)
    +O\left(\frac{\varphi_n^3}{n}\right),
\end{aligned}
\]
where
$\sum_i(\Hb_n)_{ii}^2
\leq k\max_i(\Hb_n)_{ii}=O(n^{-1})$.
Consequently,
\[
    \Var\left(
        \frac{
            \bepsilon_n^\intercal\Hb_n\bepsilon_n-kv_n
        }{\sqrt n\,\varphi_n}
    \right)
    =
    O\left(\frac{1}{n}\right)+O\left(\frac{\varphi_n}{n^2}\right)
    =o(1),
\]
and Chebyshev's inequality yields
\begin{equation}\label{eq:glm_projection_negligible}
    \frac{1}{\sqrt n\,\varphi_n}
    \left(\bepsilon_n^\intercal\Hb_n\bepsilon_n-kv_n\right)
    =\op{1}.
\end{equation}

Now, consider the independent triangular array
\[
    \frac1{\sqrt n}
    \sum_{i=1}^n
    \begin{pmatrix}
        \dfrac{(\epsilon_{n,i}^2-v_n)}{\varphi_n}\\[9pt]
        \dfrac{\bx_{n,i}\epsilon_{n,i}}{\sqrt{\varphi_n}}
    \end{pmatrix}.
\]
Write
$\bar\bx_n\coloneqq n^{-1}\sum_{i=1}^n\bx_{n,i}$. Independence gives the covariance matrix
\[
    \begin{pmatrix}
        \dfrac{\psi_3(1/\varphi_n)+2v_n^2}{\varphi_n^2}
        &
        \dfrac{\psi_2(1/\varphi_n)}{\varphi_n^{3/2}}
        \bar\bx_n^\intercal\\[9pt]
        \dfrac{\psi_2(1/\varphi_n)}{\varphi_n^{3/2}}
        \bar\bx_n
        &
        \dfrac{v_n}{n\varphi_n}\Xb_n^\intercal\Xb_n
    \end{pmatrix}.
\]
By \eqref{eq:glm_joint_log_moments}, the upper-left entry converges to $2$ and
\[
    \frac{\psi_2(1/\varphi_n)}{\varphi_n^{3/2}}
    =O(\sqrt{\varphi_n}).
\]
The bounded design condition gives $\norm{\bar\bx_n}=O(1)$, while
$\dfrac{v_n}{n\varphi_n}\Xb_n^\intercal\Xb_n\to\Sb$. Hence the covariance matrix converges to
\[
    \begin{pmatrix}
        2&\bzero^\intercal\\
        \bzero&\Sb
    \end{pmatrix}.
\]
In addition, the bounded design condition and \eqref{eq:glm_epsilon_lyapunov_bounds} give
\[
    \frac1{n^{\frac{3}{2}}}
    \sum_{i=1}^n
    \E\left[
        \norm{
            \begin{pmatrix}
                \dfrac{(\epsilon_{n,i}^2-v_n)}{\varphi_n}\\[9pt]
                \dfrac{\bx_{n,i}\epsilon_{n,i}}{\sqrt{\varphi_n}}
            \end{pmatrix}
        }^3
    \right]
    =O(n^{-1/2}).
\]
Thus the Lyapunov condition holds. The multivariate Lindeberg--Feller theorem \citep[Theorem~2.27]{vandervaartAsymptoticStatistics1998}, \eqref{eq:glm_joint_beta_decomp}, \eqref{eq:glm_projection_negligible}, and $n^{-1}\Xb_n^\intercal\Xb_n\to\Sb$ now prove \eqref{eq:glm_joint_approx_centered}.

We next treat the exact MLE. Consider the normalized sum of independent vectors
\[
    \frac1{\sqrt n}
    \sum_{i=1}^n
    \begin{pmatrix}
        2\dfrac{\left(h(Z_{n,i})-g(\varphi_n)\right)}{\varphi_n}\\[9pt]
        \bx_{n,i}\dfrac{(Z_{n,i}-1)}{\sqrt{\varphi_n}}
    \end{pmatrix}
\]
whose summands have mean zero. The moment identities
\[
    \Var\left(h(Z_{n,i})\right)
    =\psi_1\left(\dfrac{1}{\varphi_n}\right)-\varphi_n
    =\frac{\varphi_n^2}{2}+O(\varphi_n^3)
\]
and
\begin{equation}\label{eq:glm_score_orthogonality}
    \operatorname{Cov}\left(Z_{n,i}-1,h(Z_{n,i})\right)
    =\Var(Z_{n,i})
     -\operatorname{Cov}(Z_{n,i},\log{Z_{n,i}})
    =0
\end{equation}
show that its covariance matrix is
\[
    \begin{pmatrix}
        \dfrac{4\left(\psi_1(1/\varphi_n)-\varphi_n\right)}
        {\varphi_n^2}
        &
        \bzero^\intercal\\[9pt]
        \bzero
        &
        \dfrac{1}{n}\Xb_n^\intercal\Xb_n
    \end{pmatrix},
\]
which converges to
\[
    \begin{pmatrix}
        2&\bzero^\intercal\\
        \bzero&\Sb
    \end{pmatrix}.
\]
The bounded design condition, \eqref{eq:glm_Z_third_bound}, and \eqref{eq:glm_h_centered_third_bound} give
\[
    \frac1{n^{\frac{3}{2}}}
    \sum_{i=1}^n
    \E\left[
        \norm{
            \begin{pmatrix}
                2\dfrac{\left(h(Z_{n,i})-g(\varphi_n)\right)}{\varphi_n}\\[9pt]
                \bx_{n,i}\dfrac{(Z_{n,i}-1)}{\sqrt{\varphi_n}}
            \end{pmatrix}
        }^3
    \right]
    =
    O(n^{-1/2})
    \longrightarrow0.
\]
Thus the Lyapunov condition holds, and the normalized sum converges to the corresponding normal law.

It remains to relate this limit to the exact roots. The coefficient score function and \eqref{eq:glm_beta_exact_localization} give
\begin{equation}\label{eq:glm_exact_beta_linearization}
    \sqrt{\frac n{\varphi_n}}
    (\hat\bbeta_n-\bbeta_0)
    =
    \sqrt{\dfrac{n}{\varphi_n}}
    (\Xb_n^\intercal\Xb_n)^{-1}
    \Xb_n^\intercal(\bZ_n-\vec{\bm 1}_n)
    +\op{1}.
\end{equation}
To see this, define the integrated Hessian
\[
    \overline{\Sb}_n
    \coloneqq
    \int_0^1
    \frac1n\Xb_n^\intercal
    \operatorname{diag}
    \left(
        \bZ_n\left(\bbeta_0+t
        \widehat{\boldsymbol{\delta}}_n\right)
    \right)
    \Xb_n\dif t.
\]
The coefficient score equation gives
\[
    \widehat{\boldsymbol{\delta}}_n
    =
    \overline{\Sb}_n\inv
    \frac1n\Xb_n^\intercal
    (\bZ_n-\vec{\bm 1}_n).
\]
The calculation for $\nabla^2A_n(\bzero)$ in the proof of \cref{lemm:glm_estimator_rates} gives
\[
    \frac1n\Xb_n^\intercal
    \operatorname{diag}(\bZ_n)\Xb_n
    -\frac1n\Xb_n^\intercal\Xb_n
    =\op{1}.
\]
Moreover, the bounded design condition and
$\widehat{\boldsymbol{\delta}}_n=\op{1}$ imply, uniformly in $t\in[0,1]$, that
\[
    \max_{1\leq i\leq n}
    \left|
        \exp{-t\bx_{n,i}^\intercal
        \widehat{\boldsymbol{\delta}}_n}-1
    \right|
    =\op{1}.
\]
Since $n^{-1}\sum_iZ_{n,i}=\Op{1}$, it follows that
$\overline{\Sb}_n-\frac1n\Xb_n^\intercal\Xb_n=\op{1}$ and hence
$\overline{\Sb}_n\inv-n(\Xb_n^\intercal\Xb_n)^{-1}=\op{1}$.
Together with
$n^{-1/2}\varphi_n^{-1/2}
\Xb_n^\intercal(\bZ_n-\vec{\bm 1}_n)=\Op{1}$,
this proves \eqref{eq:glm_exact_beta_linearization}.

For the dispersion coordinate, the proof of \cref{lemm:glm_estimator_rates} showed that
\[
    A_n(\widehat{\boldsymbol{\delta}}_n)-A_n(\bzero)
    =\Op{\frac{\varphi_n}{n}}.
\]
Therefore,
\[
    \frac{\sqrt n}{\varphi_n}
    \left(
        A_n(\widehat{\boldsymbol{\delta}}_n)-A_n(\bzero)
    \right)
    =
    \Op{n^{-\frac{1}{2}}}
    =\op{1},
\]
so estimating the coefficient does not contribute to the first-order dispersion limit.
Since $g'(\varphi_n)\to1/2$, a Taylor expansion of \eqref{eq:glm_profile_dispersion} gives
\begin{equation}\label{eq:glm_exact_phi_linearization}
    \frac{\sqrt n}{\varphi_n}
    (\hat\varphi_n-\varphi_n)
    =
    \frac{2}{\sqrt n\,\varphi_n}
    \sum_{i=1}^n
    \left(h(Z_{n,i})-g(\varphi_n)\right)
    +\op{1}.
\end{equation}
Combining the joint score CLT with \eqref{eq:glm_exact_beta_linearization} and \eqref{eq:glm_exact_phi_linearization} proves \eqref{eq:glm_joint_exact}.

Finally, consider the bias correction. From \eqref{eq:glm_joint_approx_centered},
\[
    \check v_n-v_n=\Op{\frac{\varphi_n}{\sqrt n}}.
\]
The expansions in \eqref{eq:glm_joint_log_moments} imply
\begin{equation}\label{eq:glm_bias_correction_expansions}
    v_n-\frac12v_n^2
    =\varphi_n+O(\varphi_n^3)\dand
    \frac12v_n-\frac16v_n^2
    =g(\varphi_n)+O(\varphi_n^3).
\end{equation}
Set
\[
    q_\varphi(v)\coloneqq v-\frac12v^2
    \dand
    q_\beta(v)\coloneqq\frac12v-\frac16v^2.
\]
Since $q_\varphi'(v_n)\to1$,
\[
    q_\varphi(\check v_n)-q_\varphi(v_n)
    =
    \left(1+\op{1}\right)(\check v_n-v_n).
\]
The deterministic remainder in the first expansion of \eqref{eq:glm_bias_correction_expansions} is negligible on the dispersion scale exactly when
\[
    \frac{\sqrt n}{\varphi_n}O(\varphi_n^3)
    =O(\sqrt n\,\varphi_n^2)
    =o(1).
\]
This follows from $n\varphi_n^4\to0$. For the coefficient bias correction, $m_n=-g(\varphi_n)$ and the second expansion in \eqref{eq:glm_bias_correction_expansions} give
\[
    m_n+q_\beta(v_n)=O(\varphi_n^3),
\]
while
\[
    q_\beta(\check v_n)-q_\beta(v_n)
    =
    \Op{\frac{\varphi_n}{\sqrt n}}.
\]
The latter term is $\op{\sqrt{\varphi_n/n}}$. The deterministic remainder is also negligible on the coefficient scale because
\[
    \sqrt{\frac n{\varphi_n}}O(\varphi_n^3)
    =
    O(\sqrt n\,\varphi_n^{\frac{5}{2}})
    =
    O(\sqrt n\,\varphi_n^2)\sqrt{\varphi_n}
    =o(1).
\]
Thus the coefficient bias correction removes the intercept bias without altering its random first-order term. Combining these conclusions with \eqref{eq:glm_joint_approx_centered} proves \eqref{eq:glm_joint_bias_corrected}.

It remains to compare the corrected and exact estimators on the same probability space. Set
\[
    r_{n,i}\coloneqq\epsilon_{n,i}-(Z_{n,i}-1).
\]
The identities $\E[\epsilon_{n,i}]=\E[Z_{n,i}-1]=0$ and
$\E[Z_{n,i}\epsilon_{n,i}]=\varphi_n$ give
\[
    \E[r_{n,i}]=0
    \dand
    \E[r_{n,i}^2]=v_n-\varphi_n=O(\varphi_n^2).
\]
Writing $\Sb_n\coloneqq n^{-1}\Xb_n^\intercal\Xb_n$, the bounded design assumption therefore yields
\[
\begin{aligned}
    \E\left[
        \left\|
        \frac{1}{\sqrt{n\varphi_n}}\Xb_n^\intercal\boldsymbol r_n
        \right\|^2
    \right]
    &=
    \frac{v_n-\varphi_n}{n\varphi_n}
    \Tr(\Xb_n^\intercal\Xb_n)
    =O(\varphi_n),
\end{aligned}
\]
where $\boldsymbol r_n\coloneqq(r_{n,1},\ldots,r_{n,n})^\intercal$. Hence
$\Sb_n^{-1}(n\varphi_n)^{-1/2}\Xb_n^\intercal\boldsymbol r_n=\op{1}$. The preceding bias-correction calculation and \eqref{eq:glm_joint_beta_decomp} give
\[
    \sqrt{\frac n{\varphi_n}}
    (\tilde\bbeta_n^{\mathrm{bc}}-\bbeta_0)
    =
    \sqrt{\frac n{\varphi_n}}
    (\Xb_n^\intercal\Xb_n)^{-1}
    \Xb_n^\intercal\bepsilon_n
    +\op{1}.
\]
Subtracting \eqref{eq:glm_exact_beta_linearization} now implies
\begin{equation}\label{eq:glm_corrected_beta_equivalence}
    \sqrt{\frac n{\varphi_n}}
    (\tilde\bbeta_n^{\mathrm{bc}}-\hat\bbeta_n)
    =\op{1}.
\end{equation}

For the dispersion coordinate, \eqref{eq:glm_projection_negligible} and $k$ fixed give
\[
    \frac{\sqrt n}{\varphi_n}(\check v_n-v_n)
    =
    \frac{1}{\sqrt n\,\varphi_n}
    \sum_{i=1}^n(\epsilon_{n,i}^2-v_n)
    +\op{1}.
\]
The Taylor expansion of $q_\varphi(\check v_n)$ used above consequently gives
\begin{equation}\label{eq:glm_corrected_phi_linearization}
    \frac{\sqrt n}{\varphi_n}
    (\tilde\varphi_n^{\mathrm{bc}}-\varphi_n)
    =
    \frac{1}{\sqrt n\,\varphi_n}
    \sum_{i=1}^n(\epsilon_{n,i}^2-v_n)
    +\op{1}.
\end{equation}
Define
\[
    d_{n,i}\coloneqq
    (\epsilon_{n,i}^2-v_n)
    -2\{h(Z_{n,i})-g(\varphi_n)\}.
\]
Since $m_n=-g(\varphi_n)$, we have
$h(Z_{n,i})-g(\varphi_n)=Z_{n,i}-1-\epsilon_{n,i}$, so $\E[d_{n,i}]=0$. Exponential tilting of the gamma density, together with the polygamma recurrence formulas, gives
\[
    \E[Z_{n,i}\epsilon_{n,i}]=\varphi_n
    \dand
    \E[Z_{n,i}\epsilon_{n,i}^2]=v_n.
\]
Using these identities and the centered log-gamma moments collected above,
\begin{equation}\label{eq:glm_coupled_dispersion_moment}
\begin{aligned}
    \E[d_{n,i}^2]
    &=
    \psi_3(1/\varphi_n)+2v_n^2
    +4(v_n-\varphi_n)+4\psi_2(1/\varphi_n)
    =O(\varphi_n^3).
\end{aligned}
\end{equation}
It follows by Chebyshev's inequality that
\[
    \frac{1}{\sqrt n\,\varphi_n}\sum_{i=1}^n d_{n,i}
    =\op{1}.
\]
Subtracting \eqref{eq:glm_exact_phi_linearization} from
\eqref{eq:glm_corrected_phi_linearization} proves
\[
    \frac{\sqrt n}{\varphi_n}
    (\tilde\varphi_n^{\mathrm{bc}}-\hat\varphi_n)
    =\op{1}.
\]
Together with \eqref{eq:glm_corrected_beta_equivalence}, this is
\eqref{eq:glm_bias_corrected_equivalence}.
\end{proof}

\FloatBarrier

\section{Zeta Distribution}\label{sec:zeta}

\subsection{Preliminary Results}\label{app:zetadistresults}

The first lemma gives the moment expansions for $\log{X}$ used throughout \cref{sec:zeta}. The second uses these expansions to verify the triangular array central limit theorem needed for the limiting distributions below.

\begin{lemm}\label{lemm:zetaexps}
    If $X\sim\mathrm{Zeta}(s)$, then, for $r=1,2,3$,
    \[
        \E\left[(\log{X})^r\right]
        =\frac{(-1)^r\zeta^{(r)}(s)}{\zeta(s)}.
    \]
    In particular,
    \[
        \E\left[\log{X}\right]=-\frac{\zeta'(s)}{\zeta(s)}
        \dand
        \Var(\log{X})
        =\frac{\zeta''(s)}{\zeta(s)}
        -\left(\frac{\zeta'(s)}{\zeta(s)}\right)^2.
    \]
    Moreover, as $s\to\infty$, the raw moments for $r=1,2,3$ and the variance satisfy, respectively,
    \begin{equation}\label{eq:zeta_log_moment_expansions}
        \begin{aligned}
            \E\left[(\log{X})^r\right]
            &= \frac{(\log{2})^r}{2^s}+O\left(\frac{1}{3^s}\right),\\
            \Var(\log{X})
            &= \frac{(\log{2})^2}{2^s}+O\left(\frac{1}{3^s}\right).
        \end{aligned}
    \end{equation}
\end{lemm}
\begin{proof}
    The zeta function is represented on $s>1$ by the absolutely convergent Dirichlet series
    \[
        \zeta(s)=\sum_{k=1}^\infty k^{-s}.
    \]
    A Dirichlet series converges uniformly on compact subsets of its half-plane of convergence and may be differentiated termwise there \citep[Theorems~11.11--11.12]{apostolIntroductionAnalyticNumber1976}; therefore the differentiated series through order three converge absolutely and uniformly on every such interval. Since
    \[
        \frac{\mathrm{d}^r}{\mathrm{d}s^r}k^{-s}
        =(-\log{k})^rk^{-s},
    \]
    termwise differentiation gives, for $r=1,2,3$,
    \begin{equation}\label{eq:zeta_derivative_series}
        \zeta^{(r)}(s)
        =(-1)^r\sum_{k=1}^\infty\frac{(\log{k})^r}{k^s}.
    \end{equation}
    On the other hand, the mass function of $X$ gives
    \begin{equation}\label{eq:zeta_log_moments}
        \E\left[(\log{X})^r\right]
        =\sum_{k=1}^\infty(\log{k})^r\frac{k^{-s}}{\zeta(s)}
        =\frac1{\zeta(s)}\sum_{k=1}^\infty
          \frac{(\log{k})^r}{k^s}.
    \end{equation}
    Combining \eqref{eq:zeta_derivative_series} and \eqref{eq:zeta_log_moments} proves the first identity. Because the $k=1$ terms vanish, the sums may equivalently begin at $k=2$.

    The stated mean and variance identities follow immediately. For every fixed $r\in\N$, isolating the term $k=2$ and bounding the remaining tail gives, as $s\to\infty$,
    \[
        \sum_{k=2}^\infty\frac{(\log{k})^r}{k^s}
        =\frac{(\log{2})^r}{2^s}+O\left(\frac{1}{3^s}\right),
    \]
    where the implied constant may depend on $r$.
    Since $\zeta(s)=1+2^{-s}+O(3^{-s})$, division by $\zeta(s)$ does not change the displayed expansion. Taking $r=1,2$ and observing that $\left(\E\left[\log{X}\right]\right)^2=O(4^{-s})$ proves \eqref{eq:zeta_log_moment_expansions}.
\end{proof}

\begin{lemm}\label{lemm:zeta_Lyapunov}
    Let $X_s\sim\mathrm{Zeta}(s)$, $Y_s\coloneqq\log{X_s}$, $\mu_s\coloneqq\E\left[Y_s\right]$, and $\sigma_s^2\coloneqq\Var(Y_s)$. Then there exists a constant $C<\infty$ such that, for all sufficiently large $s$,
    \[
        \E\left[\left|\frac{Y_s-\mu_s}{\sigma_s}\right|^3\right]
        \leq C2^{\frac{s}{2}}.
    \]
    Consequently, if $X_{n,1},\ldots,X_{n,n}\iid\mathrm{Zeta}(s_n)$, $s_n\to\infty$, and $n2^{-s_n}\to\infty$, then, with $Y_{n,i}\coloneqq\log{X_{n,i}}$ and $\overline Y_n\coloneqq n^{-1}\sum_{i=1}^nY_{n,i}$, the Lyapunov condition
    \[
        \frac{\sum_{i=1}^n\E\left[|Y_{n,i}-\mu_{s_n}|^3\right]}
        {\left(\sum_{i=1}^n\Var(Y_{n,i})\right)^{3/2}}
        \longrightarrow0
    \]
    holds, and so
    \[
        \frac{\sqrt n(\overline Y_n-\mu_{s_n})}{\sigma_{s_n}}
        \inD\mathcal N(0,1).
    \]
\end{lemm}

\begin{proof}
    By \cref{lemm:zetaexps}, $\E\left[Y_s^3\right]=O(2^{-s})$, $\mu_s=O(2^{-s})$, and $\sigma_s^2=\log{2}^2 2^{-s}+O(3^{-s})$. Since $Y_s\geq0$,
    \[
        |Y_s-\mu_s|^3\leq4(Y_s^3+\mu_s^3),
    \]
    and therefore $\E\left[|Y_s-\mu_s|^3\right]=O(2^{-s})$. The leading term in $\sigma_s^2$ is positive, so
    \[
        \E\left[\left|\frac{Y_s-\mu_s}{\sigma_s}\right|^3\right]
        =O(2^{\frac{s}{2}}).
    \]
    Within each row, the observations are identically distributed. The Lyapunov ratio is therefore
    \[
        \frac{n\E\left[|Y_{n,1}-\mu_{s_n}|^3\right]}{(n\sigma_{s_n}^2)^{3/2}}
        =\frac1{\sqrt n}
        \E\left[\left|\frac{Y_{s_n}-\mu_{s_n}}{\sigma_{s_n}}\right|^3\right]
        =O\left(\sqrt{\frac{2^{s_n}}n}\right)\longrightarrow0,
    \]
    as required. The Lyapunov condition implies the Lindeberg condition, so the Lindeberg--Feller theorem \citep[Theorem~2.27]{vandervaartAsymptoticStatistics1998} gives the stated central limit theorem.
\end{proof}

\subsection{Root Enclosure and Approximation on the Original Scale}\label{app:sub:zeta}

Fix $s_n\to\infty$, set $u_n\coloneqq2^{-s_n}$ and $p\coloneqq\logtwo{3}$, and let $X_{n,1},\ldots,X_{n,n}\iid\mathrm{Zeta}(s_n)$. Assume throughout this subsection that $nu_n\to\infty$. All stochastic orders are understood to be under $\mathbb P_n\coloneqq\mathbb P_{s_n}^{\otimes n}$. Define
\[
    m(u)\coloneqq\sum_{k=2}^\infty\Lambda(k)u^{\logtwo{k}}
    =\E_{-\logtwo{u}}\left[\log{X}\right]
    \dand
    v(u)\coloneqq\Var_{-\logtwo{u}}(\log{X}),
\]
where $\Lambda$ is the von Mangoldt function. The moment expansion in \cref{lemm:zetaexps} gives the leading terms of $m(u)$ and $v(u)$. To obtain the corresponding expansion of $m'(u)$, set $s(u)\coloneqq-\logtwo{u}$ and differentiate the series for $m$ termwise, as justified by the same compact uniform convergence argument used in the proof of \cref{lemm:zetaexps}. The chain rule gives
\[
    m'(u)
    =\log{2}
    +\frac{1}{u\log{2}}
      \sum_{k=3}^\infty\frac{\Lambda(k)\log{k}}{k^{s(u)}}.
\]
Since $s(u)=-\logtwo{u}$ and $p=\logtwo{3}$, we have $3^{-s(u)}=u^p$. Because $0\leq\Lambda(k)\leq\log{k}$, the tail estimate used in the proof of \cref{lemm:zetaexps} gives
\[
    0
    \leq
    \sum_{k=3}^\infty\frac{\Lambda(k)\log{k}}{k^{s(u)}}
    \leq
    \sum_{k=3}^\infty\frac{(\log{k})^2}{k^{s(u)}}
    =O\left(3^{-s(u)}\right)
    =O(u^p).
\]
The prefactor $1/(u\log{2})$ therefore makes the remainder in $m'(u)$ of order $O(u^{p-1})$. Thus, as $u\searrow0$,
\begin{equation}\label{eq:Zeta-mun}
    m(u)=(\log{2})u+O(u^p)
    \dand
    m'(u)=\log{2}+O(u^{p-1})
    \dand
    v(u)=(\log{2})^2u+O(u^p).
\end{equation}
Because $\log{G_n}=n^{-1}\sum_i\log{X_{n,i}}$, we have
\[
    \Var_{s_n}(\log{G_n})
    =\frac{v(u_n)}{n}
    =O\left(\frac{u_n}{n}\right).
\]
Chebyshev's inequality therefore gives
\begin{equation}\label{eq:Zeta-logGn-approx}
    \log{G_n}=m(u_n)+O_{\mathbb P_n}\!\left(\sqrt{\frac{u_n}{n}}\right).
\end{equation}
Since $p>1$ and $nu_n\to\infty$, the deterministic remainder $u_n^p$ and the stochastic remainder $\sqrt{u_n/n}$ are both $o(u_n)$. It follows that
\begin{equation}\label{eq:Zeta-logGn}
    \log{G_n}=\log{2}u_n+o_{\mathbb P_n}(u_n) \dand
    \widetilde u_n=u_n+o_{\mathbb P_n}(u_n).
\end{equation}
The coefficients in the series for $m$ are nonnegative, and its $k=2$ term is $\log{2}u$, so $m$ is strictly increasing. The exceptional event $G_n=1$ has probability
\[
    \mathbb P_n(G_n=1)
    =\zeta(s_n)^{-n}
    =\exp{-nu_n(1+o(1))}\longrightarrow0.
\]
On the event $|\widetilde u_n/u_n-1|\leq1/4$, define
\[
    D_n^*\coloneqq\left\{u:|u-\widetilde u_n|\leq\frac{u_n}{4}\right\}.
\]
This event has probability tending to one and $D_n^*\subset[u_n/2,3u_n/2]\subset U$ for all sufficiently large $n$. Apply \cref{coro:stochastic_localization_transfer} with $\Rb_n\coloneqq u_n^{-1}$, $\Jb_n\coloneqq n\log{2}$, and $r_n\coloneqq1/4$. The surrogate is linear, so $\kappa_n=0$, while \eqref{eq:Zeta-mun} gives
\[
\begin{aligned}
    \eta_n
    &\coloneqq\sup_{u\in D_n^*}
    \left|u_n^{-1}(n\log{2})^{-1}[Q_n(u)-\widetilde Q_n(u)]\right|
    =O(u_n^{p-1})=o(1).
\end{aligned}
\]
On $G_n>1$, the exact equation has a unique root because $m$ is strictly increasing. \cref{coro:stochastic_localization_transfer} therefore gives
\begin{equation}\label{eq:zeta_root_difference}
    |\widetilde u_n-\widehat u_n|
    =O_{\mathbb P_n}(u_n^p)
    \dand
    \frac{\widehat u_n}{u_n}\inP1.
\end{equation}
Consequently, after applying the completion convention, the deterministic interval $K_n\coloneqq[u_n/2,2u_n]$ contains both roots with probability tending to one. On this interval we may take $a_n\coloneqq n\log{2}$.

Returning to the original scale, the same calculation on $K_n$ gives $\Delta_{K,n}=O(nu_n^p)$. Moreover, for $u\in K_n$,
\[
    |(\phi^{-1})'(u)|
    =\frac{1}{u\log{2}}
    \leq\frac{2}{u_n\log{2}}.
\]
Thus $\phi^{-1}(u)=-\logtwo{u}$ is $b_n\coloneqq2/(u_n\log{2})$-Lipschitz on $K_n$, and
\[
    \frac{b_n\Delta_{K,n}}{a_n}=O(u_n^{p-1})=o(1).
\]
Therefore, the quantitative bound in the proof of \cref{coro:original_scale} yields
\begin{equation}\label{eq:zeta_original_scale_difference}
    |\widetilde s_n-\widehat s_n|
    =O_{\mathbb P_n}\left(u_n^{p-1}\right)
    =O_{\mathbb P_n}\left(2^{-s_n(p-1)}\right)
    =o_{\mathbb P_n}(1).
\end{equation}

\subsection{Limiting Distributions}

Suppose that $nu_n\to\infty$. By \cref{lemm:zeta_Lyapunov}, the triangular array $Y_{n,i}=\log{X_{n,i}}$ satisfies
\[
    \frac{\sqrt n\left(\log{G_n}-m(u_n)\right)}{\sqrt{v(u_n)}}
    \inD\mathcal N(0,1).
\]
On the event $G_n>1$, whose probability tends to one by \cref{app:sub:zeta}, the exact root satisfies $m(\widehat u_n)=\log{G_n}$. The mean value theorem therefore gives, for some $u_n^*$ between $u_n$ and $\widehat u_n$,
\[
    \sqrt{\frac n{u_n}}(\widehat u_n-u_n)
    =
    \frac{\sqrt n\left(\log{G_n}-m(u_n)\right)}{\sqrt{v(u_n)}}
    \frac{\sqrt{v(u_n)/u_n}}{m'(u_n^*)}.
\]
The localization already proved gives $u_n^*/u_n\inP1$, and therefore $m'(u_n^*)\inP\log{2}$; at the same time, \eqref{eq:Zeta-mun} gives $v(u_n)/u_n\to\log{2}^2$. Slutsky's theorem now yields
\[
    \sqrt{\frac n{u_n}}(\widehat u_n-u_n)\inD\mathcal N(0,1).
\]

To make the discrepancy between the roots negligible under the normalization $\sqrt{n/u_n}$, assume in addition that
\begin{equation}\label{eq:zeta_additional_condition}
    \sqrt n\,u_n^{p-\frac{1}{2}}\longrightarrow0.
\end{equation}
Then \eqref{eq:zeta_root_difference} gives directly
\begin{equation}\label{eq:Zeta-dif-op1}
    \sqrt{\frac n{u_n}}(\widetilde u_n-\widehat u_n)
    =O_{\mathbb P_n}\left(\sqrt n\,u_n^{p-\frac12}\right)
    =o_{\mathbb P_n}(1).
\end{equation}
Combining this result with the limit for $\widehat u_n$ and applying Slutsky's theorem gives
\[
    \sqrt{\frac n{u_n}}(\widetilde u_n-u_n)
    \inD\mathcal N(0,1).
\]

It remains to translate these limits back to the original parameter scale. Let $g(u)\coloneqq-\logtwo{u}$, so
\[
    g'(u)=-\frac{1}{u\log{2}}
    \dand
    g''(u)=\frac{1}{u^2\log{2}}.
\]
For $w_n \in \{\widehat u_n,\widetilde u_n\}$, under the additional condition \eqref{eq:zeta_additional_condition}, the corresponding limit above gives
\[
    w_n-u_n
    =O_{\mathbb P_n}\left(\sqrt{\frac{u_n}{n}}\right).
\]
Both $u_n$ and $w_n$ belong to $K_n=[u_n/2,2u_n]$ with probability tending to one. On this event, Taylor's theorem gives, for some $\xi_n$ between $u_n$ and $w_n$,
\begin{equation}\label{eq:Zeta-Taylor-exp}
    g(w_n)-g(u_n)
    =-\frac{w_n-u_n}{u_n\log{2}}
      +\frac{(w_n-u_n)^2}{2\xi_n^2\log{2}}.
\end{equation}
Since $\xi_n\in K_n$, the remainder in \eqref{eq:Zeta-Taylor-exp} is $O_{\mathbb P_n}(1/(nu_n))$. Therefore,
\[
    \sqrt{nu_n}\left(g(w_n)-g(u_n)\right)
    =-\frac{1}{\log{2}}
      \sqrt{\frac n{u_n}}(w_n-u_n)
      +O_{\mathbb P_n}\left(\frac{1}{\sqrt{nu_n}}\right).
\]
Because $nu_n\to\infty$, the final term is $o_{\mathbb P_n}(1)$. Taking $w_n=\widehat u_n$ gives
\[
    \sqrt{nu_n}(\widehat s_n-s_n)
    \inD\mathcal N\left(0,\frac1{\log{2}^2}\right),
\]
while taking $w_n=\widetilde u_n$ under the additional condition \eqref{eq:zeta_additional_condition} gives
\[
    \sqrt{nu_n}(\widetilde s_n-s_n)
    =\sqrt{\frac n{2^{s_n}}}(\widetilde s_n-s_n)
    \inD\mathcal N\left(0,\frac1{\log{2}^2}\right).
\]

\subsection{Simulation Study}

We used the same sample size grid and number of replications as in \cref{sub:gamma_simulation_study}. We set
\[
    u_n\coloneqq n^{-\gamma}
    \dand
    s_n\coloneqq-\logtwo{u_n}=\gamma\logtwo{n},
\]
where
\[
    \gamma\in\{0.4,0.5,0.7\}
    \dand
    n\in\{250,500,1000,2500,5000,10000,25000\}.
\]
For each design point, we generated $2000$ samples from the $\mathrm{Zeta}(s_n)$ distribution. Substituting $u_n=n^{-\gamma}$ and $p=\logtwo{3}$ into the sufficient condition $\sqrt n\,u_n^{p-1/2}\to0$ gives
\[
    \sqrt n\,u_n^{p-1/2}
    =n^{1/2-\gamma(p-1/2)}\longrightarrow0
    \quad\Longleftrightarrow\quad
    \gamma>\frac{1}{2\logtwo{3}-1}\approx0.461.
\]
Thus, the grid straddled the sufficient threshold for the approximate estimator to share the first-order limiting distribution of the exact MLE. All three choices also satisfied $nu_n\to\infty$, as required for the probability of a finite MLE to tend to one.

For each sample with $G_n>1$, we computed the exact MLE numerically and the approximate estimator $\widetilde s_n$ given above. As in the gamma simulations, we divided the absolute difference between the approximate and exact estimators on the original parameter scale by the asymptotic standard error $1/\left(\log{2}\sqrt{nu_n}\right)$ of the exact MLE. We also compared the RMSE of the approximate estimator with that of the exact MLE.

\begin{figure}[!t]
    \centering
    \includegraphics[width=0.80\textwidth]{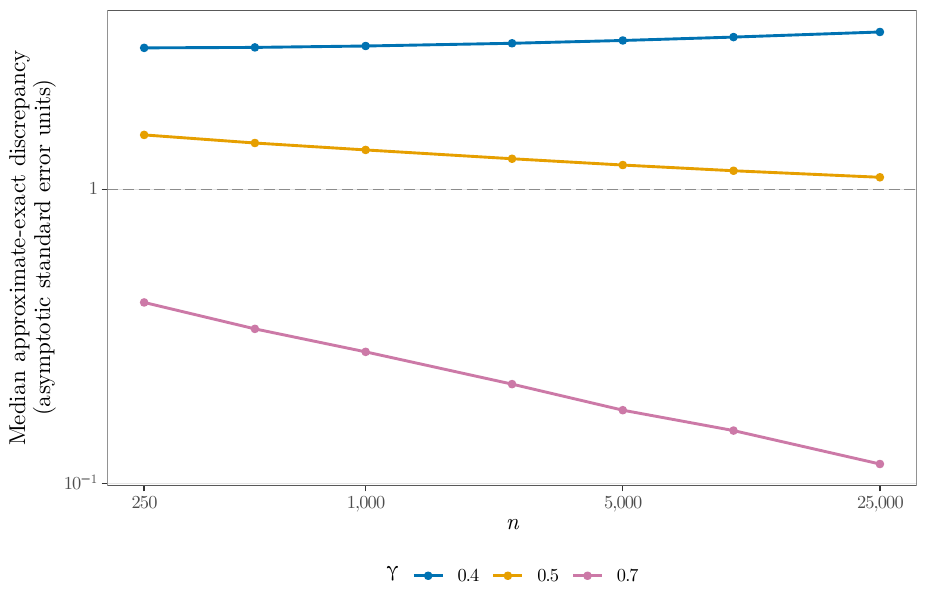}
    \caption{Median absolute difference between the approximate and exact MLEs on the original parameter scale, divided by the asymptotic standard error of the exact MLE. The horizontal dashed line marks one asymptotic standard error. Each point is based on $2000$ replications.}
    \label{fig:zeta-discrepancy}
\end{figure}

\cref{fig:zeta-discrepancy} shows the transition predicted by the rate condition. When $\gamma=0.4$, which is below the threshold, the approximate-exact difference remains between three and four asymptotic standard errors throughout the sample size grid. It decreases steadily for the other two regimes, although the decrease is slow for $\gamma=0.5$, which is only slightly above the threshold. At $n=25000$, the median discrepancy on the $s$ scale is $1.10$ standard errors for $\gamma=0.5$ and $0.116$ for $\gamma=0.7$.

\begin{figure}[!ht]
    \centering
    \includegraphics[width=0.80\textwidth]{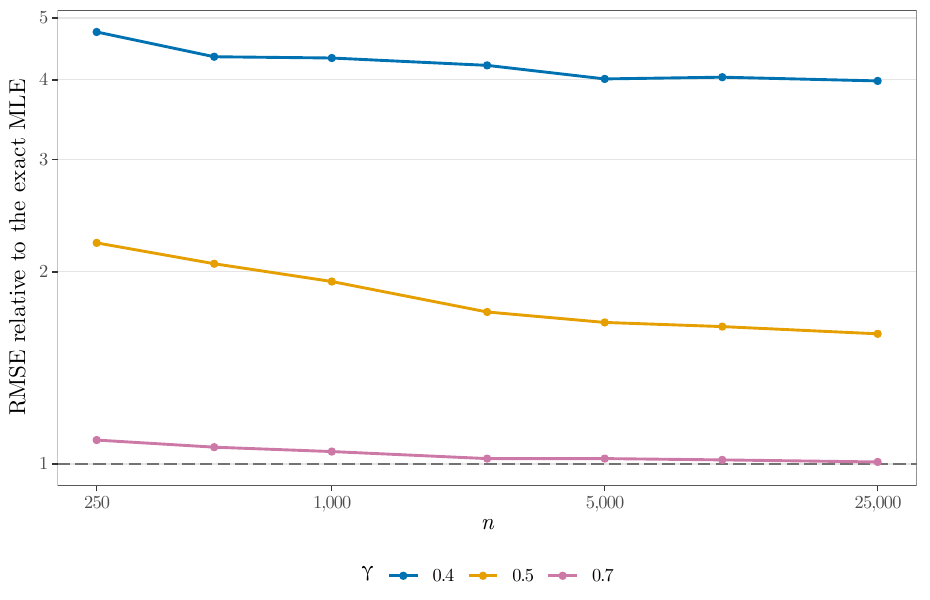}
    \caption{RMSE of the approximate estimator relative to that of the exact MLE on the original parameter scale. The dashed horizontal line marks equal RMSE. Each point is based on $2000$ replications with a finite MLE.}
    \label{fig:zeta-rmse}
\end{figure}

\cref{fig:zeta-rmse} leads to the same conclusion. At $n=25000$, the RMSE ratio for the approximate MLE is $3.98$ when $\gamma=0.4$, $1.60$ when $\gamma=0.5$, and $1.01$ when $\gamma=0.7$. Thus the ratio approaches one rapidly for $\gamma=0.7$, decreases more slowly for $\gamma=0.5$, and remains far above one below the threshold despite declining modestly over the displayed grid.

\FloatBarrier

\begin{figure}[!t]
    \centering
    \includegraphics[width=0.82\textwidth]{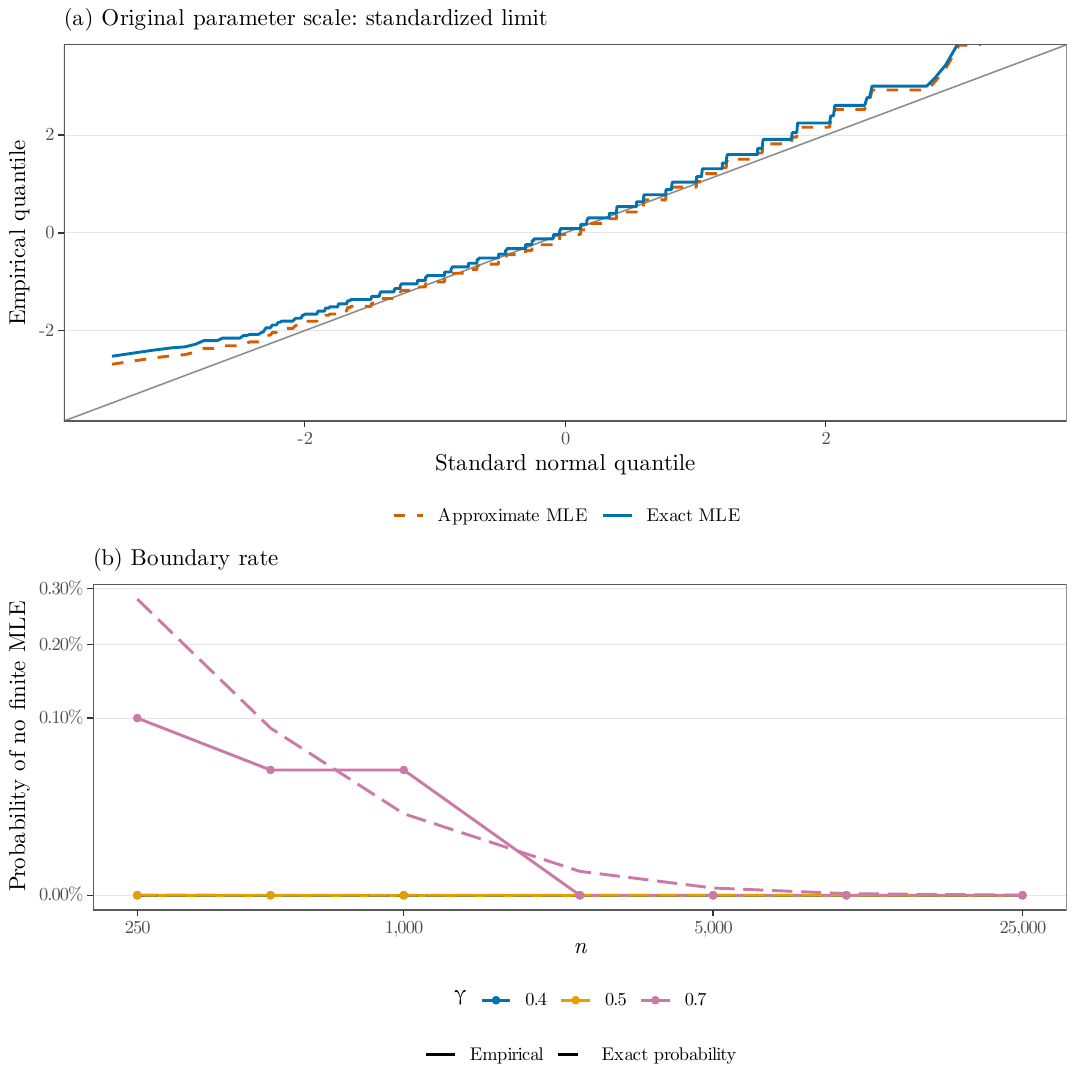}
    \caption{Distributional and boundary diagnostics. Panel (a) compares the exact and approximate standardized errors on the original parameter scale with standard normal quantiles at $n=25000$ and $\gamma=0.7$. Panel (b) compares the empirical frequency of samples with no finite MLE with its exact probability $\zeta(s_n)^{-n}$.}
    \label{fig:zeta-diagnostics}
\end{figure}

The distributional diagnostics in \cref{fig:zeta-diagnostics} use $(n,\gamma)=(25000,0.7)$. On the original parameter scale, the empirical mean and variance of the unit variance normalization for the approximate estimator are $0.008$ and $1.073$, respectively, and the exact and approximate quantiles are close. Some discreteness remains in the tails because the expected number of observations greater than one grows only at rate $nu_n=n^{0.3}$.

When $G_n=1$, the finite MLE does not exist. Such samples occurred only for $\gamma=0.7$ and $n\leq1000$ in our experiment, and the empirical frequencies are consistent with the exact probability $\zeta(s_n)^{-n}$ shown in the lower panel of \cref{fig:zeta-diagnostics}. We excluded these samples only from summaries requiring a finite MLE. Among the remaining samples, the approximate estimator always belongs to the parameter space. The numerical score inversion is also stable: its largest absolute residual over all $42000$ samples is less than $3\times10^{-12}$.

\FloatBarrier

\section{von Mises Distribution}\label{sec:vonmises}

\subsection{Preliminary Results}\label{app:vonMisesdistresults}

\subsubsection{Bessel Expansions and Moment Identities}

\begin{lemm}\label{lemm:vonMisesBessel}
As $\kappa\to\infty$,
\[
    A(\kappa)
    \coloneqq\frac{I_1(\kappa)}{I_0(\kappa)}
        =1-\frac1{2\kappa}-\frac1{8\kappa^2}
            -\frac1{8\kappa^3}+O\left(\frac{1}{\kappa^4}\right),
\]
and
\[
    \frac{I_2(\kappa)}{I_0(\kappa)}
        =1-\frac2\kappa+\frac1{\kappa^2}
            +\frac1{4\kappa^3}+O\left(\frac{1}{\kappa^4}\right).
\]
Moreover,
\[
    A'(\kappa)
    =\frac1{2\kappa^2}+\frac1{4\kappa^3}+O\left(\frac{1}{\kappa^4}\right)
    \dand
    A''(\kappa)=O\left(\frac{1}{\kappa^3}\right).
\]
\end{lemm}

\begin{proof}
For fixed $\nu$, the standard expansion is
\[
I_\nu(\kappa)
=\frac{e^\kappa}{\sqrt{2\pi\kappa}}
\left[
1-\frac{m_\nu-1}{8\kappa}
+\frac{(m_\nu-1)(m_\nu-9)}{2!(8\kappa)^2}
-\frac{(m_\nu-1)(m_\nu-9)(m_\nu-25)}{3!(8\kappa)^3}
+O\left(\frac{1}{\kappa^4}\right)
\right],
\]
where $m_\nu\coloneqq4\nu^2$ \citep[Formula 9.7.1]{abramowitz1988handbook}. For $\nu=0,1,2$, the bracketed factors are, respectively,
\[
\begin{aligned}
1&+\frac1{8\kappa}+\frac9{128\kappa^2}
    +\frac{75}{1024\kappa^3}+O\left(\frac{1}{\kappa^4}\right),\\
1&-\frac3{8\kappa}-\frac{15}{128\kappa^2}
    -\frac{105}{1024\kappa^3}+O\left(\frac{1}{\kappa^4}\right),\\
1&-\frac{15}{8\kappa}+\frac{105}{128\kappa^2}
    +\frac{315}{1024\kappa^3}+O\left(\frac{1}{\kappa^4}\right).
\end{aligned}
\]
Dividing the second and third series by the first gives the expansions of $A$ and $I_2/I_0$; as a check, the latter can also be recovered from $I_0(\kappa)-I_2(\kappa)=2I_1(\kappa)/\kappa$.

Differentiating the Bessel ratio gives
\[
    A'(\kappa)
    =\frac{1+I_2(\kappa)/I_0(\kappa)}{2}-A(\kappa)^2
\]
and substituting the two ratio expansions gives the stated formula for $A'$. Using $I_2(\kappa)/I_0(\kappa)=1-2A(\kappa)/\kappa$, the exact derivative identity can also be written as
\[
    A'(\kappa)
    =
    1-\frac{A(\kappa)}{\kappa}-A(\kappa)^2.
\]
Differentiating this identity gives
\[
    A''(\kappa)
    =
    -\frac{A'(\kappa)}{\kappa}
    +\frac{A(\kappa)}{\kappa^2}
    -2A(\kappa)A'(\kappa).
\]
Substitution of the established expansions of $A$ and $A'$ into this exact identity is now legitimate. Indeed,
\[
    \frac{A'(\kappa)}{\kappa}=O\left(\frac{1}{\kappa^3}\right)\dand
    \frac{A(\kappa)}{\kappa^2}-2A(\kappa)A'(\kappa)
    =O\left(\frac{1}{\kappa^3}\right),
\]
which proves $A''(\kappa)=O\left(\frac{1}{\kappa^3}\right)$.
\end{proof}

\begin{lemm}\label{lemm:vonMisesMoments}
If $W=X-\mu\sim\mathrm{von\,Mises}(0,\kappa)$ and $B(\kappa)=I_2(\kappa)/I_0(\kappa)$, then the first moments satisfy
\[
\begin{gathered}
    \E\left[\sin{W}\right]=0
    \dand \E\left[\cos{W}\right]=A(\kappa)
    \dand
    \E\left[\sin{W}\cos{W}\right]=0,
\end{gathered}
\]
while the variances satisfy
\[
\begin{gathered}
    \Var(\sin{W})=\frac{1-B(\kappa)}2=\frac{A(\kappa)}{\kappa} \dand
    \Var(\cos{W})=\frac{1+B(\kappa)}2-A(\kappa)^2=A'(\kappa).
\end{gathered}
\]
\end{lemm}

\begin{proof}
For every integer $k\geq1$, the integral representation of $I_k$ gives
\[
\begin{aligned}
    \E\left[\cos{kW}\right]
    &=\frac{1}{2\pi I_0(\kappa)}
      \int_{-\pi}^{\pi}
    \exp{\kappa\cos{w}}\cos{kw}\dif w
      =\frac{I_k(\kappa)}{I_0(\kappa)} \dand
    \E\left[\sin{kW}\right]&=0,
\end{aligned}
\]
where the second identity follows by symmetry. Taking $k=1$ gives the first two expectations in the lemma, and the same symmetry gives $\E\left[\sin{W}\cos{W}\right]=0$. The double angle identities
\[
    \sinsq{w}=\frac{1-\cos{2w}}2
    \dand
    \cossq{w}=\frac{1+\cos{2w}}2
\]
then yield
\[
    \Var(\sin{W})=\frac{1-B(\kappa)}2
    \dand
    \Var(\cos{W})=\frac{1+B(\kappa)}2-A(\kappa)^2.
\]
The recurrence $1-B(\kappa)=2A(\kappa)/\kappa$ gives the alternative form of the first variance; differentiating $I_1/I_0$ identifies the second with $A'(\kappa)$.
\end{proof}

\subsubsection{MLE Characterization and Root Preservation}

For $X_1,\ldots,X_n\iid\mathrm{von\,Mises}(\mu,\kappa)$, the log likelihood, up to an additive constant, is
\[
    \ell_n(\mu,\kappa)
    \coloneqq-n\log{I_0(\kappa)}+\kappa\sum_{i=1}^n\cos{X_i-\mu},
\]
with score
\[
    \partial_\mu\ell_n
    =\kappa\sum_{i=1}^n\sin{X_i-\mu}
    \dand
    \partial_\kappa\ell_n
    =\sum_{i=1}^n\cos{X_i-\mu}-nA(\kappa)
    \dand A(\kappa)\coloneqq\frac{I_1(\kappa)}{I_0(\kappa)}.
\]
Let $\overline C_n$, $\overline S_n$, $R_n$, and $\widehat\mu_n$ be as in \cref{sub:vonmises}. The identity
\[
    \sum_{i=1}^n\cos{X_i-\mu}
    =nR_n\cos{\mu-\widehat\mu_n}
\]
shows that, whenever $R_n>0$, the likelihood is uniquely maximized in the angular coordinate at $\widehat\mu_n$ modulo $2\pi$; the resulting profile score for $\kappa$ is $n\left(R_n-A(\kappa)\right)$.

Since $A(\kappa)=(\log{I_0(\kappa)})'$, differentiation under the integral sign gives
\[
    A'(\kappa)=\Var_{\mu,\kappa}\left(\cos{X-\mu}\right)>0.
\]
Moreover, $A(\kappa)\downarrow0$ as $\kappa\downarrow0$ and $A(\kappa)\uparrow1$ as $\kappa\to\infty$, so $A(\kappa)=R_n$ has a unique positive solution whenever $0<R_n<1$. For $n\geq2$, the equality $R_n=1$ requires all sample directions to coincide, whereas $R_n=0$ requires their vector sum to vanish; both are null events under a continuous von Mises law. Thus, for $n\geq2$, we have that $0<R_n<1$ almost surely and the MLE pair is unique almost surely. If $R_n=1$, there is no finite concentration maximizer. If $R_n=0$, the profile likelihood approaches its supremum as $\kappa\downarrow0$, while the direction is unidentified. As elsewhere, we complete the estimators according to \cref{sec:defs} on these exceptional events.

The modified system in \cref{sub:vonmises} preserves this MLE root exactly. To see this, write
\[
    \bq_n(\mu,\kappa)
    \coloneqq
    \begin{pmatrix}
        \mu\left(A(\kappa)-1\right)-\widehat\mu_n(R_n-1)\\
        A(\kappa)-R_n
    \end{pmatrix}.
\]
If $q_{2,n}=0$, then
\[
    q_{1,n}=(R_n-1)(\mu-\widehat\mu_n).
\]
Since $R_n<1$ almost surely, the common root is precisely $(\widehat\mu_n,\widehat\kappa_n)$.

\subsection{Uniform Root Stability}

The chart is built from the approximation
\[
    g(\kappa)\coloneqq-\frac1{2\kappa}-\frac1{8\kappa^2}.
\]
The function $g$ is strictly increasing from $-\infty$ to $0$ on $(0,\infty)$. Thus, for an open interval $M$ of length less than $2\pi$, the map
\[
    \bphi(\mu,\kappa)\coloneqq\left(\mu g(\kappa),g(\kappa)\right)^{\mathsf T}
\]
is a bijection from $N\coloneqq M\times(0,\infty)$ onto
\[
    U\coloneqq\{(u_1,u_2):u_2<0,\ \frac{u_1}{u_2}\in M\},
\]
with inverse
\[
    \bphi^{-1}(u_1,u_2)
    =\left(\frac{u_1}{u_2},k(u_2)\right)^{\mathsf T}
    \dand
    k(u_2)\coloneqq-\frac{1+\sqrt{1-2u_2}}{4u_2}.
\]
The remaining discrepancy between the exact Bessel ratio and the chart approximation is captured by $h(u_2)\coloneqq A\left(k(u_2)\right)-1$. As $u_2\nearrow0$, \cref{lemm:vonMisesBessel} gives
\[
    h(u_2)=u_2+O(|u_2|^3) \dand
    h'(u_2)=1+O(u_2^2).
\]
Because $k=g^{-1}$ and both $k$ and $A$ are strictly increasing, $h$ is a strictly increasing bijection from $(-\infty,0)$ onto $(-1,0)$.
The transformed exact and surrogate systems are given by
\[
    \bQ_n(\bu)
    \coloneqq
    \begin{pmatrix}
        (\dfrac{u_1}{u_2})h(u_2)-\widehat\mu_n(R_n-1)\\[9pt]
        h(u_2)-(R_n-1)
    \end{pmatrix}\dand
    \widetilde{\bQ}_n(\bu)
    \coloneqq
    \begin{pmatrix}
        u_1-\widehat\mu_n(R_n-1)\\
        u_2-(R_n-1)
    \end{pmatrix}.
\]
In particular, $\Jb_{\widetilde{\bQ}_n}=\Ib_2$, and the surrogate root is
\[
    \widetilde{\bu}_n
    \coloneqq\left(\widehat\mu_n(R_n-1),R_n-1\right)^{\mathsf T}.
\]

Fix $\kappa_n\to\infty$, keep $\mu$ fixed, and write
\[
    u_{2,n}\coloneqq g(\kappa_n)
    \dand
    \bu_n\coloneqq(\mu u_{2,n},u_{2,n})^{\mathsf T}
    \dand
    |u_{2,n}|\asymp\frac{1}{\kappa_n}.
\]
All stochastic orders in this subsection are under $\mathbb P_n\coloneqq\mathbb P_{\mu,\kappa_n}^{\otimes n}$. With $W_{n,i}\coloneqq X_{n,i}-\mu$ represented in $[-\pi,\pi)$, set
\[
    C_n\coloneqq\frac1n\sum_{i=1}^n\cos{W_{n,i}}
    \dand
    S_n\coloneqq\frac1n\sum_{i=1}^n\sin{W_{n,i}}.
\]
The moment identities in \cref{lemm:vonMisesMoments}, together with the expansions in \cref{lemm:vonMisesBessel} and Chebyshev's inequality, give
\[
    C_n=A(\kappa_n)
        +\Op{\frac1{\kappa_n\sqrt n}}\dand
    S_n=\Op{\frac1{\sqrt{n\kappa_n}}}.
\]
Rotating all observations by $-\mu$ leaves the resultant length unchanged, so
\[
    R_n^2
    =\left(\frac1n\sum_{i=1}^n\cos{X_{n,i}}\right)^2
     +\left(\frac1n\sum_{i=1}^n\sin{X_{n,i}}\right)^2
    =C_n^2+S_n^2.
\]
On the event $C_n>1/2$, whose probability tends to one, the representative of the sample mean direction in $M$ therefore satisfies
\[
    R_n-C_n=\frac{S_n^2}{R_n+C_n}
            =\Op{\frac1{n\kappa_n}}\dand
    \widehat\mu_n-\mu
            =\atantwo{S_n}{C_n}
            =\Op{\frac1{\sqrt{n\kappa_n}}}.
\]
Consequently,
\begin{equation}\label{eq:vonMises_Rn_bound}
    R_n-A(\kappa_n)
    =\Op{\frac1{\kappa_n\sqrt n}
              +\frac1{n\kappa_n}}.
\end{equation}
Combining this bound with $A(\kappa_n)-1=u_{2,n}+O(|u_{2,n}|^3)$ gives
\[
    \widetilde u_{2,n}=R_n-1=u_{2,n}+\op{|u_{2,n}|}.
\]
Moreover,
\[
\begin{aligned}
    \widetilde u_{1,n}-\mu u_{2,n}
    &=(\widehat\mu_n-\mu)(R_n-1)
      +\mu\left((R_n-1)-u_{2,n}\right)
      =\op{|u_{2,n}|},
\end{aligned}
\]
so $\widetilde u_{1,n}=\mu u_{2,n}+\op{|u_{2,n}|}$. We now show that localizing the surrogate root is enough to locate the exact root. Choose $r_0>0$ so that $\overline{B}_{r_0}((-\mu,-1)^\mathsf T)\subseteq\{(v_1,v_2):v_2<0,\ v_1/v_2\in M\}$. On the event that $|u_{2,n}|^{-1}\widetilde{\bu}_n\in\overline{B}_{r_0/2}((-\mu,-1)^\mathsf T)$, define
\[
    D_n^*\coloneqq\left\{\bu:
    \norm{|u_{2,n}|^{-1}(\bu-\widetilde{\bu}_n)}\leq\frac{r_0}{2}\right\}.
\]
This event has probability tending to one and $D_n^*\subset U$. Apply \cref{coro:stochastic_localization_transfer} with $\Rb_n\coloneqq|u_{2,n}|^{-1}\Ib_2$, $\Jb_n\coloneqq\Ib_2$, and $r_n\coloneqq r_0/2$. The surrogate Jacobian is the identity, so in the notation of \cref{coro:stochastic_localization_transfer} we have $\kappa_n=0$. Moreover,
\[
    \bQ_n(\bu)-\widetilde{\bQ}_n(\bu)
    =\left(h(u_2)-u_2\right)
      \begin{pmatrix}
\frac{u_1}{u_2}\\[9pt]
1
\end{pmatrix}.
\]
On $D_n^*$, $|u_2|\asymp|u_{2,n}|$ and $u_1/u_2$ is bounded. Since $h(u_2)-u_2=O(|u_2|^3)$,
\[
    \eta_n
    \coloneqq\sup_{\bu\in D_n^*}
    \norm{\Rb_n[\bQ_n(\bu)-\widetilde{\bQ}_n(\bu)]}
    =O(|u_{2,n}|^2)=o(1).
\]
The exact system has a unique root in $U$, so \cref{coro:stochastic_localization_transfer} gives
\begin{equation}\label{eq:vonMises_chart_difference}
    \norm{\widetilde{\bu}_n-\widehat{\bu}_n}
    =\Op{|u_{2,n}|^3}
    =\Op{\kappa_n^{-3}}.
\end{equation}
In particular, $\norm{\widehat{\bu}_n-\bu_n}=\op{|u_{2,n}|}$ follows from localizing the surrogate root.

Choose $\delta>0$ so that $[\mu-\delta,\mu+\delta]\subset M$, and define the deterministic convex sets
\[
    K_n\coloneqq
    \left\{\bu\in U:
        \frac32u_{2,n}\leq u_2\leq\frac12u_{2,n},\
        \left|\frac{u_1}{u_2}-\mu\right|\leq\delta
    \right\}.
\]
Since $u_2<0$, the ratio restriction is equivalent to the linear inequalities
\[
    (\mu+\delta)u_2\leq u_1\leq(\mu-\delta)u_2,
\]
so $K_n$ is convex. Let $E_n$ be the event on which both estimators lie in $N$, satisfy their respective estimating equations, and their transformed roots lie in $K_n$. The localization results above and the null exceptional events give $\mathbb P_n(E_n)\to1$. Moreover, $\sup_{\bu\in K_n}\norm{\bu}\to0$. On $K_n$,
\[
    \bQ_n(\bu)-\widetilde{\bQ}_n(\bu)
    =\left(h(u_2)-u_2\right)
      \begin{pmatrix}\dfrac{u_1}{u_2}\\[9pt]1\end{pmatrix},
\]
and therefore
\[
    \Delta_{K,n}
    \coloneqq\sup_{\bu\in K_n}
      \norm{\bQ_n(\bu)-\widetilde{\bQ}_n(\bu)}
    =O(|u_{2,n}|^3)=O\left(\frac{1}{\kappa_n^3}\right).
\]

For the inverse chart,
\[
    \Jb_{\bphi^{-1}}(\bu)
    =
    \begin{pmatrix}
        \dfrac{1}{u_2}&-\dfrac{u_1}{u_2^2}\\[9pt]
        0&\dfrac{1}{g'\left(k(u_2)\right)}
    \end{pmatrix}.
\]
On $K_n$, we have $|u_2|\asymp\kappa_n^{-1}$, $|u_1/u_2|=O(1)$, and $g'\left(k(u_2)\right)\asymp\kappa_n^{-2}$. Consequently, $\bphi^{-1}$ is $b_n$-Lipschitz there for some $b_n$ satisfying $b_n=O(\kappa_n^2)$, and $b_n\Delta_{K,n}=O\left(\frac{1}{\kappa_n}\right)$. Applying \cref{coro:original_scale} gives
\[
    \widetilde\mu_n=\widehat\mu_n \dand
    \widetilde\kappa_n-\widehat\kappa_n
    =\Op{\frac{1}{\kappa_n}}=\op{1}.
\]
This approximation on the original parameter scale requires no additional relation between $n$ and $\kappa_n$.

\subsection{Limiting Distributions}

\begin{lemm}\label{lemm:vonMisesCLT}
Let $X_{n,1},\ldots,X_{n,n}\iid \mathrm{von\,Mises}(\mu,\kappa_n)$, where $\kappa_n\to\infty$, and represent $W_{n,i}\coloneqq X_{n,i}-\mu$ in $[-\pi,\pi)$. Then
\[
    \frac1{\sqrt n}\sum_{i=1}^n
    \begin{pmatrix}
        \sqrt{\dfrac{\kappa_n}{A(\kappa_n)}}\,\sin{W_{n,i}}\\[9pt]
        \dfrac{\left(\cos{W_{n,i}}-A(\kappa_n)\right)}{\sqrt{A'(\kappa_n)}}
    \end{pmatrix}
    \inD\mathcal N(\bzero,\Ib_2).
\]
\end{lemm}

\begin{proof}
We first establish a uniform moment bound. Let $W_\kappa\sim\mathrm{von\,Mises}(0,\kappa)$ and define $V_\kappa\coloneqq\sqrt\kappa W_\kappa$. The elementary inequality
\[
    1-\cos{w}\geq\frac{2w^2}{\pi^2}
\]
holds for $|w|\leq\pi$, and the standard Bessel asymptotic gives $e^{-\kappa}I_0(\kappa)\asymp\kappa^{-1/2}$. For $|v|<\pi\sqrt\kappa$, the density of $V_\kappa$ is
\[
    f_{V_\kappa}(v)
    =
    \frac{
        \exp{-\kappa\left(1-\cos{v/\sqrt\kappa}\right)}
    }{
        2\pi\sqrt\kappa\,e^{-\kappa}I_0(\kappa)
    },
\]
and it is zero otherwise. It follows that $f_{V_\kappa}(v)\leq C\exp{-2v^2/\pi^2}$ for all sufficiently large $\kappa$. Hence, for every fixed $r>0$,
\[
    \sup_{\kappa\geq\kappa_0}\E\left[|V_\kappa|^r\right]<\infty.
\]
Define
\[
    \bZ_\kappa
    =
    \begin{pmatrix}
        \sqrt{\dfrac{\kappa}{A(\kappa)}}\,\sin{W_\kappa}\\[9pt]
        \dfrac{\left(\cos{W_\kappa}-A(\kappa)\right)}{\sqrt{A'(\kappa)}}
    \end{pmatrix}.
\]
By \cref{lemm:vonMisesBessel,lemm:vonMisesMoments}, $\E\left[\bZ_\kappa\right]=\bzero$ and $\Var(\bZ_\kappa)=\Ib_2$. Moreover, $A(\kappa)\to1$, so $A(\kappa)$ is bounded away from zero for all sufficiently large $\kappa$. Consequently,
\[
    \sqrt{\frac{\kappa}{A(\kappa)}}\,|\sin{W_\kappa}|
    \leq
    C|V_\kappa|.
\]
For the second component,
\[
    \kappa|\cos{W_\kappa}-A(\kappa)|
    \leq\frac12V_\kappa^2+\kappa\left(1-A(\kappa)\right).
\]
Since $\kappa\left(1-A(\kappa)\right)$ is bounded and $\kappa\sqrt{A'(\kappa)}\to1/\sqrt2$, the preceding moment bound implies
\begin{equation}\label{eq:vonMiseslyapu}
    \sup_{\kappa\geq\kappa_0}\E\left[\norm{\bZ_\kappa}^3\right]<\infty.
\end{equation}
Let $\bZ_{\kappa_n,i}$ be the vector constructed from $W_{n,i}$. For every fixed nonzero $\bt\in\mathbb R^2$, the Lyapunov ratio is
\[
    \frac{
        \sum_{i=1}^n
        \E\left[|\bt^{\mathsf T}\bZ_{\kappa_n,i}|^3\right]
    }{
            \left(
                \sum_{i=1}^n
                \Var(\bt^{\mathsf T}\bZ_{\kappa_n,i})
            \right)^{3/2}
    }
    =
    \frac{
        \E\left[|\bt^{\mathsf T}\bZ_{\kappa_n}|^3\right]
    }{
        \sqrt n\,\norm{\bt}^3
    }
    =O(n^{-\frac{1}{2}}),
\]
where the last equality follows from \eqref{eq:vonMiseslyapu}. The Lyapunov central limit theorem and the Cram\'er--Wold device complete the proof.
\end{proof}

For one observation, the score is
\[
    \bs_{\mu,\kappa}(X)
    \coloneqq
    \left(\kappa\sin{X-\mu},\ \cos{X-\mu}-A(\kappa)\right)^{\mathsf T}.
\]
Writing $W\coloneqq X-\mu$, the density of $W$ is even, whereas both $\sin{W}$ and
$\sin{W}\left(\cos{W}-A(\kappa)\right)$ are odd. Hence
\[
    \E\left[\sin{W}\right]=0
    \dand
    \E\left[\sin{W}\left(\cos{W}-A(\kappa)\right)\right]=0,
\]
so the off-diagonal covariance vanishes. The variance identities in \cref{lemm:vonMisesMoments} then give
\[
    \mathcal I_1(\mu,\kappa)
    \coloneqq\Var\left(\bs_{\mu,\kappa}(X)\right)
    =
    \begin{pmatrix}
        \kappa A(\kappa)&0\\
        0&A'(\kappa)
    \end{pmatrix}.
\]
Here the information function for the concentration coordinate can also be written through the exact identity
\[
    A'(\kappa)
    =\frac{1+I_2(\kappa)/I_0(\kappa)}{2}-A(\kappa)^2
    =\Var\left(\cos{X-\mu}\right).
\]
The expansions in \cref{lemm:vonMisesBessel} imply
\[
    \kappa A(\kappa)\sim\kappa \dand
    A'(\kappa)\sim\frac1{2\kappa^2}.
\]
Thus the concentration entry of the inverse Fisher information for $n$ observations is
\[
    \frac1{nA'(\kappa_n)}
    \sim\frac{2\kappa_n^2}{n}.
\]

We derive the joint MLE limit directly from the empirical sine and cosine, retaining the definitions of $C_n$ and $S_n$ from the preceding subsection. The triangular array Lyapunov result in \cref{lemm:vonMisesCLT} gives
\[
    \begin{pmatrix}
        \sqrt{\dfrac{n\kappa_n}{A(\kappa_n)}}\,S_n\\[9pt]
        \sqrt n\,\dfrac{\left(C_n-A(\kappa_n)\right)}{\sqrt{A'(\kappa_n)}}
    \end{pmatrix}
    \inD\mathcal N(\bzero,\Ib_2).
\]
The identity
\[
    \widehat\mu_n-\mu
    =\arctan{\frac{S_n}{C_n}},
\]
together with $R_n-C_n=S_n^2/(R_n+C_n)$ and Taylor's theorem, gives
\[
    \widehat\mu_n-\mu
    =\frac{S_n}{A(\kappa_n)}
      +\op{(n\kappa_n)^{-\frac{1}{2}}}
    \dand
    R_n-A(\kappa_n)
    =C_n-A(\kappa_n)
    +\op{\left(\frac{1}{\kappa_n\sqrt n}\right)}.
\]
Since $A(\kappa_n)\to1$ and $\kappa_n^2A'(\kappa_n)\to1/2$,
\[
    \begin{pmatrix}
        \sqrt{n\kappa_n}(\widehat\mu_n-\mu)\\
        \sqrt n\,\kappa_n\left(R_n-A(\kappa_n)\right)
    \end{pmatrix}
    \inD
    \mathcal N\left(
        \bzero,\begin{pmatrix}1&0\\[9pt]0&\dfrac{1}{2}\end{pmatrix}
    \right).
\]

The bound \eqref{eq:vonMises_Rn_bound} and $1-A(\kappa_n)\sim1/(2\kappa_n)$ give
\[
    \frac{1-R_n}{1-A(\kappa_n)}
    =
    1-\frac{R_n-A(\kappa_n)}{1-A(\kappa_n)}
    =
    1+\op{1}.
\]
Since $1-A(\widehat\kappa_n)=1-R_n$, the preceding display is equivalently
\[
    \frac{1-A(\widehat\kappa_n)}{1-A(\kappa_n)}\inP1.
\]
Because $\kappa_n\to\infty$, we have $1-A(\kappa_n)\to0$, and the last display therefore implies $1-A(\widehat\kappa_n)\inP0$. Since $A$ is strictly increasing and approaches one only as its argument tends to infinity, it follows that $\widehat\kappa_n\to\infty$ in probability. Applying $1-A(x)\sim(2x)^{-1}$ at both $\kappa_n$ and $\widehat\kappa_n$ now gives
\[
    \frac{1-A(\widehat\kappa_n)}{1-A(\kappa_n)}
    =
    \frac{\kappa_n}{\widehat\kappa_n}\left(1+\op{1}\right),
\]
and hence
\[
    \frac{\widehat\kappa_n}{\kappa_n}\inP1.
\]

Let $\overline\kappa_n$ lie between $\kappa_n$ and $\widehat\kappa_n$. The mean value theorem gives
\[
    R_n-A(\kappa_n)
    =
    A'(\overline\kappa_n)(\widehat\kappa_n-\kappa_n).
\]
Because $\widehat\kappa_n/\kappa_n\inP1$, we also have $\overline\kappa_n/\kappa_n\inP1$, and therefore
\[
    \kappa_n^2A'(\overline\kappa_n)\inP\frac12.
\]
Combining this with the previous bound on $R_n-A(\kappa_n)$ yields
\[
    \widehat\kappa_n-\kappa_n
    =
    \Op{\frac{\kappa_n}{\sqrt n}}.
\]

Taylor's theorem then gives
\[
\begin{aligned}
    R_n-A(\kappa_n)
    &=
    A'(\kappa_n)(\widehat\kappa_n-\kappa_n)
    +\frac12A''(\kappa_n^\dagger)
    (\widehat\kappa_n-\kappa_n)^2,
\end{aligned}
\]
where $\kappa_n^\dagger$ lies between $\kappa_n$ and $\widehat\kappa_n$. The convergence $\widehat\kappa_n/\kappa_n\inP1$ also gives $\kappa_n^\dagger/\kappa_n\inP1$; using $A''(\kappa)=O(\kappa^{-3})$ and the preliminary rate, the quadratic term is $\Op{1/(n\kappa_n)}$. Since $A'(\kappa_n)\sim1/(2\kappa_n^2)$, rearranging the expansion gives
\[
    \widehat\kappa_n-\kappa_n
    =\frac{R_n-A(\kappa_n)}{A'(\kappa_n)}
      +\Op{\frac{\kappa_n}{n}}.
\]
It follows that
\[
    \begin{pmatrix}
        \sqrt{n\kappa_n}(\widehat\mu_n-\mu)\\[9pt]
        \sqrt n\,\dfrac{(\widehat\kappa_n-\kappa_n)}{\kappa_n}
    \end{pmatrix}
    \inD
    \mathcal N\left(
        \bzero,\begin{pmatrix}1&0\\0&2\end{pmatrix}
    \right).
\]

The location coordinates of the exact and approximate MLEs agree exactly; for the concentration coordinates, root stability gives
\[
    \frac{\sqrt n}{\kappa_n}
    (\widetilde\kappa_n-\widehat\kappa_n)
    =\Op{\frac{\sqrt n}{\kappa_n^2}}.
\]
If $\sqrt n/\kappa_n^2\to0$ (equivalently, $\kappa_n/n^{1/4}\to\infty$), this discrepancy is $\op{1}$, and Slutsky's theorem gives the approximate MLE the same joint limit displayed above.

\subsection{Simulation Study}

We used the same sample size grid and number of replications as in the gamma and zeta studies. We set $\mu\coloneqq0$ and $\kappa_n\coloneqq n^\gamma$, where
\[
    \gamma\in\{0.15,0.25,0.40\}
    \dand
    n\in\{250,500,1000,2500,5000,10000,25000\},
\]
and generated $2000$ samples at each design point. These regimes lay below, at, and above the sufficient condition $\gamma>1/4$ for the approximate MLE to have the same limiting distribution as the exact MLE.

In each replication, we obtained the exact MLE by numerically solving $A(\widehat\kappa_n)=R_n$ and computed the approximate MLE from the explicit formula in \cref{sub:vonmises}. We also included the estimator of \citet{abeyasekeraEstimationParametersMises1982} as an explicit benchmark because, like our estimator, it is designed for the high-concentration, large-$\kappa$ regime. To apply their estimator without ambiguity at the circular discontinuity, we represented the observations in the interval of length $2\pi$ centered at $\widehat\mu_n$. Writing $X_{n,i}^{\circ}$ for this representation and $\overline X_n^{\circ}$ for its arithmetic mean, we used the benchmark
\[
    \widehat\kappa_n^{\mathrm{lin}}
    \coloneqq
    \frac{1}{
        (n-3)^{-1}
        \sum_{i=1}^n
        (X_{n,i}^{\circ}-\overline X_n^{\circ})^2
    }.
\]

\begin{figure}[!htbp]
    \centering
    \includegraphics[width=0.88\textwidth]{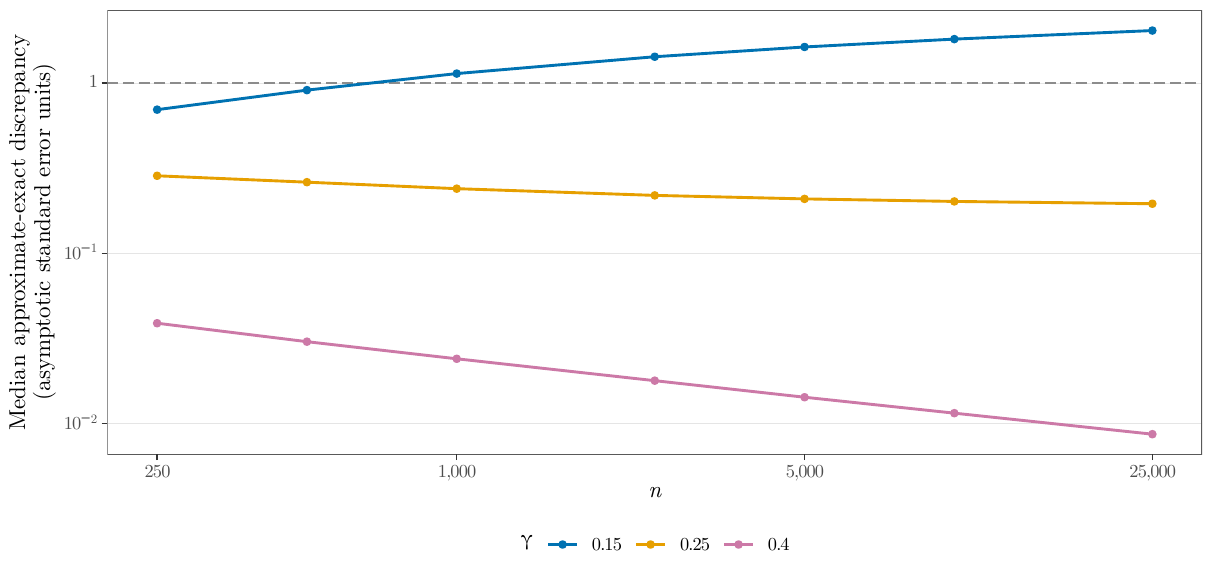}
    \caption{Median absolute difference between the approximate and exact concentration estimators, divided by the asymptotic standard error of the exact MLE. The horizontal dashed line marks one asymptotic standard error. Each point is based on $2000$ replications.}
    \label{fig:vonmises-discrepancy}
\end{figure}

For the comparison in \cref{fig:vonmises-discrepancy}, we divided $|\widetilde\kappa_n-\widehat\kappa_n|$ by the asymptotic standard error $\sqrt{2}\kappa_n/\sqrt n$. The median difference increases with $n$ when $\gamma=0.15$, remains of constant order when $\gamma=0.25$, and decreases rapidly when $\gamma=0.40$. At $n=25000$, the respective medians are $2.04$, $0.195$, and $0.00863$ standard errors. This pattern agrees with the factor $\sqrt n/\kappa_n^2=n^{1/2-2\gamma}$ in the theoretical bound.

\begin{figure}[!htbp]
    \centering
    \includegraphics[width=0.88\textwidth]{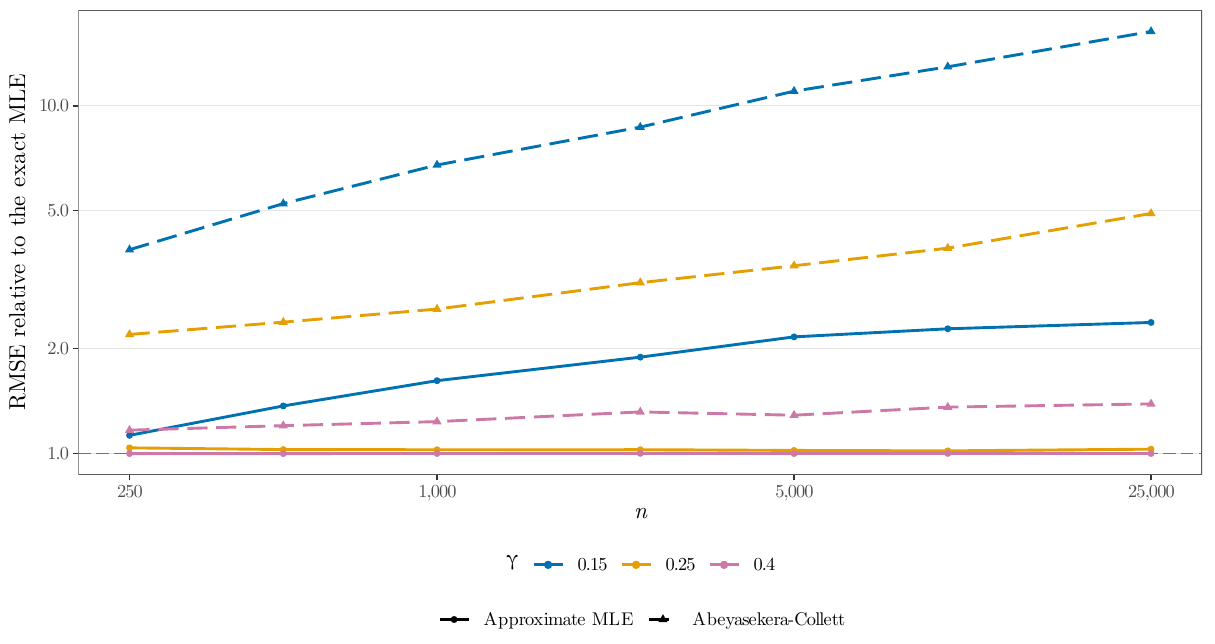}
    \caption{RMSE of each explicit concentration estimator relative to that of the exact MLE. The dashed horizontal line marks equal RMSE. Each point is based on $2000$ replications.}
    \label{fig:vonmises-rmse}
\end{figure}

The RMSE comparison in \cref{fig:vonmises-rmse} gives a similar picture. At $n=25000$, the RMSE of the approximate MLE relative to that of the exact MLE is $2.38$ for $\gamma=0.15$, $1.03$ for $\gamma=0.25$, and $1.000$ for $\gamma=0.40$. The Abeyasekera--Collett estimator has larger RMSE throughout the grid. Its RMSE ratios at the same sample size are $16.4$, $4.90$, and $1.39$, respectively, reflecting a negative bias that remains visible even in the largest concentration regime considered here.

\begin{figure}[!htbp]
    \centering
    \includegraphics[width=0.88\textwidth]{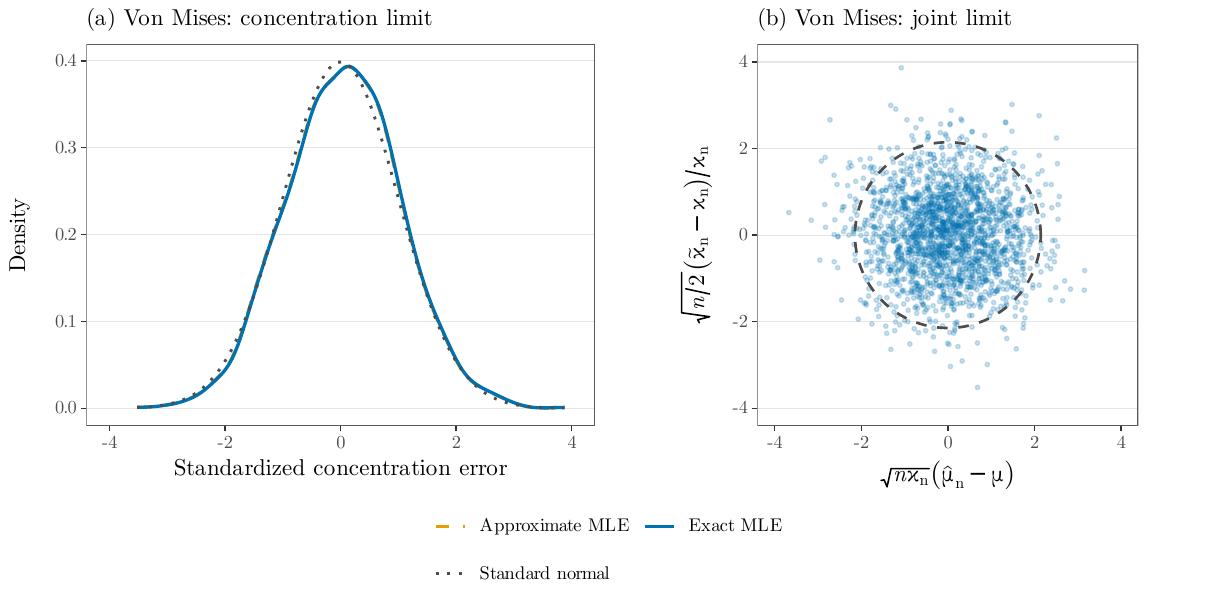}
    \caption{Distributional diagnostics at $n=25000$ and $\gamma=0.40$. Panel (a) compares the standardized concentration errors of the exact and approximate MLEs with the standard normal density. Panel (b) shows the joint standardized errors of the location and approximate concentration estimators; the dashed circle encloses $90\%$ probability under the limiting bivariate standard normal distribution.}
    \label{fig:vonmises-diagnostics}
\end{figure}

Finally, \cref{fig:vonmises-diagnostics} examines the joint limit at $(n,\gamma)=(25000,0.40)$. After division by its limiting standard deviation, the exact concentration error has empirical mean $0.051$ and variance $0.959$, while the approximate error has empirical mean $0.043$ and variance $0.959$. The standardized location error has empirical mean $-0.001$ and variance $1.023$, and its empirical correlation with the approximate concentration error is $-0.047$. The dashed circle in panel (b), which contains $90\%$ probability under the limiting bivariate standard normal distribution, contains $90.1\%$ of the simulated points.

\FloatBarrier

\section{Generalized Pareto Distribution}\label{sec:Pareto}

Throughout this section, fix $0<\bar\xi<1/6$, $\mu\in\R$, and $\sigma>0$. Let $\xi_n\downarrow0$ and suppose that
\[
 X_{n,1},\ldots,X_{n,n}\iid\mathrm{GPD}(\mu,\sigma,\xi_n).
\]
Throughout the GPD appendix, stochastic orders are understood under $\mathbb P_n\coloneqq\mathbb P_{\mu,\sigma,\xi_n}^{\otimes n}$.

\subsection{Preliminary Results}\label{app:GPDdistresults}

\subsubsection{Moments and the Sample Minimum}

\begin{lemm}\label{lemm:gpd_moments}
For any integer $k\geq1$ and any $n$ for which $k\xi_n<1$,
\begin{equation}\label{eq:gpd_moments}
m_k(\xi_n)\coloneqq\E_{\mathbb P_n}\left[X_{n,1}-\mu\right]^k
 =
 \frac{k!\sigma^k}{\prod_{j=1}^k(1-j\xi_n)}.
\end{equation}
Moreover,
\begin{equation}\label{eq:gpd_min_limit}
\frac n\sigma(X_{(1)}-\mu)\inD\mathrm{Exp}(1)
 \dand X_{(1)}-\mu=O_{\mathbb P_n}\left(\frac{1}{n}\right).
\end{equation}
For every fixed $k\in\N$, define
\[
    M_{k,n}\coloneqq\frac1n\sum_{i=1}^n(X_{n,i}-X_{(1)})^k.
\]
Then
\begin{equation}\label{eq:gpd_sample_moments}
M_{k,n}=m_k(\xi_n)+O_{\mathbb P_n}(n^{-\frac{1}{2}}).
\end{equation}
\end{lemm}

\begin{proof}
Let $Z\coloneqq X_{n,1}-\mu$. For $k\xi_n<1$, the tail integral formula gives
\[
\begin{aligned}
\E_{\mathbb P_n}\left[Z^k\right]
&=k\int_0^\infty
  y^{k-1}\left(1+\frac{\xi_n y}{\sigma}\right)^{-\frac{1}{\xi_n}}
  \dif y
  =\frac{k\sigma^k}{\xi_n^k}
  B\left(k,\frac1{\xi_n}-k\right)
=\frac{k!\sigma^k}{\prod_{j=1}^k(1-j\xi_n)},
\end{aligned}
\]
where the formula at $\xi_n=0$ is obtained by continuity. This proves \eqref{eq:gpd_moments}. In particular, the restriction $\bar\xi<1/6$ uniformly bounds moments through order six over the restricted model. More generally, because $\xi_n\to0$, for every fixed $q\in\N$ the moments of order $q$ are uniformly bounded along the regime sequence for all sufficiently large $n$. For the sample minimum, the survival function is
\[
 \mathbb P_n\left(\frac n\sigma(X_{(1)}-\mu)>x\right)
 =
 \left(1+\frac{\xi_n x}{n}\right)^{-\frac{n}{\xi_n}}
 \longrightarrow e^{-x},
\]
with the exponential interpretation when $\xi_n=0$, proving \eqref{eq:gpd_min_limit}.

Put $\delta_n\coloneqq X_{(1)}-\mu$ and $Z_{n,i}\coloneqq X_{n,i}-\mu$, and fix $k\in\N$. The binomial formula and \eqref{eq:gpd_min_limit} give
\[
 \left|M_{k,n}-\frac1n\sum_{i=1}^nZ_{n,i}^k\right|
 \leq
 C_k\delta_n\left(\frac1n\sum_{i=1}^n
       (1+Z_{n,i}^{k-1})+\delta_n^{k-1}\right)
 =O_{\mathbb P_n}\left(\frac{1}{n}\right).
\]
The eventual $(2k)$th-moment bound makes the variance of the triangular array $n^{-1}\sum_{i=1}^nZ_{n,i}^k$ of order $O(n^{-1})$, while its expectation is $m_k(\xi_n)$. Chebyshev's inequality therefore gives
\[
    \frac1n\sum_{i=1}^nZ_{n,i}^k
    =m_k(\xi_n)+O_{\mathbb P_n}(n^{-\frac{1}{2}}).
\]
Combining this estimate with the $O_{\mathbb P_n}(n^{-1})$ effect of replacing $\mu$ by $X_{(1)}$ proves \eqref{eq:gpd_sample_moments}.
\end{proof}

\subsubsection{Restricted Global MLE}

\begin{prop}\label{prop:gpd_global_mle}
For every $n\geq2$, the likelihood over
\[
 \Theta_{\bar\xi}
 =\{(\mu',s,\eta):\mu'\in\R,\ s>0,\ 0\leq\eta\leq\bar\xi\}
\]
has a global maximizer almost surely. If $\sqrt n\,\xi_n\to\infty$, this maximizer is unique with probability tending to one and satisfies
\[
 \widehat\mu_n=X_{(1)}
 \dand
 \widehat\sigma_n-\sigma=O_{\mathbb P_n}(n^{-\frac{1}{2}})
 \dand
 \widehat\xi_n-\xi_n=O_{\mathbb P_n}(n^{-\frac{1}{2}}),
\]
and $(\widehat\sigma_n,\widehat\xi_n)\in(0,\infty)\times(0,\bar\xi)$. It is the unique root of the profile score in a fixed neighborhood of $(\sigma,0)$, where the neighborhood is understood relative to $(0,\infty)\times[0,\bar\xi]$.
\end{prop}

\begin{proof}
For $\eta>0$, differentiation of the log likelihood over the feasible region $\mu'\leq X_{(1)}$ gives
\[
 \frac{\partial\ell_n}{\partial\mu'}
 =
 \frac{1+\eta}{s}\sum_{i=1}^n
 \frac{1}{1+\eta(X_{n,i}-\mu')/s}>0.
\]
At $\eta=0$, the derivative is $n/s>0$ as well. Thus every global maximizer must have $\mu'=m_n\coloneqq X_{(1)}$. Writing $Y_i\coloneqq X_{n,i}-m_n$, we obtain the profile criterion
\begin{equation}\label{eq:gpd_profile_lik}
\ell_n^p(s,\eta)
 \coloneqq
 -n\log{s}-\left(1+\frac{1}{\eta}\right)\sum_{i=1}^n
     \log{1+\frac{\eta Y_i}{s}},
\end{equation}
with the continuous extension
\[
 \ell_n^p(s,0)=-n\log{s}-\frac1s\sum_{i=1}^nY_i.
\]

The bound $\ell_n^p(s,\eta)\leq-n\log{s}$ holds uniformly in $\eta\in[0,\bar\xi]$ and tends to $-\infty$ as $s\to\infty$. For $r\geq0$, the function $\eta\mapsto\log{1+\eta r}/\eta$, continuously extended at zero, is decreasing. Since exactly one $Y_i$ is zero almost surely,
\begin{equation}\label{eq:gpd_small_scale_bound}
\ell_n^p(s,\eta)
 \leq
 -n\log{s}-\frac1{\bar\xi}\sum_{i=1}^n
    \log{1+\frac{\bar\xi Y_i}{s}}.
\end{equation}
Fix $n$ and a sample for which exactly one $Y_i$ is zero. For each of the remaining observations, as $s\downarrow0$,
\[
    \log{1+\frac{\bar\xi Y_i}{s}}
    =-\log{s}+\log{\bar\xi Y_i}+o(1).
\]
Thus, pathwise, the right side of \eqref{eq:gpd_small_scale_bound} is
\[
    \left(\frac{n-1}{\bar\xi}-n\right)\log{s}+O(1),
\]
which tends to $-\infty$ because $\bar\xi<1/6$ and $n\geq2$. Continuity therefore gives a finite global maximizer.

We next identify the maximizer. Let $Z_{n,i}\coloneqq X_{n,i}-\mu$. The sample minimum result in \cref{lemm:gpd_moments} gives $m_n-\mu=O_{\mathbb P_n}(n^{-1})$. For $z\geq0$, define
\[
 \lambda(z;s,\eta)
 \coloneqq-\log{s} - \left(1+\frac{1}{\eta}\right)\log{1+\frac{\eta z}{s}},
\]
with $\lambda(z;s,0)\coloneqq-\log{s}-z/s$. Let $Z_0\sim\mathrm{GPD}(0,\sigma,0)$ and define
\[
 L(s,\eta)\coloneqq\E\left[\lambda(Z_0;s,\eta)\right].
\]
Fix $\mathcal K_{\varepsilon,M}\coloneqq[\varepsilon,M]\times[0,\bar\xi]$ and write $\lambda_{\btheta}(z)\coloneqq\lambda(z;s,\eta)$ for $\btheta\coloneqq(s,\eta)$. A Taylor expansion in $\eta$ at zero gives
\[
 \sup_{\btheta\in\mathcal K_{\varepsilon,M}}
 \left(|\lambda_{\btheta}(z)|+
 \norm{\nabla_{\btheta}\lambda_{\btheta}(z)}\right)
 \leq G(z)\coloneqq C_{\varepsilon,M}(1+z^2).
\]
By \cref{lemm:gpd_moments}, $\sup_{n\geq n_0}\E_{\mathbb P_n}\left[G(Z_{n,1})^2\right]<\infty$ for some $n_0$. To couple these variables, let $W\sim\mathrm{Exp}(1)$ and define
\[
    Z_n^\star
    \coloneqq\sigma\frac{\exp{\xi_nW}-1}{\xi_n},
\]
with $Z_n^\star\coloneqq\sigma W$ when $\xi_n=0$. For $z\geq0$ and $\xi_n>0$,
\[
\begin{aligned}
    \mathbb P\left(Z_n^\star>z\right)
    =\mathbb P\left(
        W>\frac1{\xi_n}\log{1+\frac{\xi_nz}{\sigma}}
      \right)
      =\left(1+\frac{\xi_nz}{\sigma}\right)^{-1/\xi_n},
\end{aligned}
\]
which is the survival function of $Z_{n,1}=X_{n,1}-\mu$. Hence $Z_n^\star\stackrel{d}{=}Z_{n,1}$. Under this coupling, the following convergence holds almost surely:
\[
    Z_n^\star\longrightarrow\sigma W.
\]
Uniform integrability therefore gives, for each $\btheta\in\mathcal K_{\varepsilon,M}$,
\[
 \E_{\mathbb P_n}\left[\lambda_{\btheta}(Z_{n,1})\right]
 \longrightarrow L(\btheta).
\]

To obtain uniform convergence, fix $\delta>0$ and choose a finite $\delta$-grid $\btheta_1,\ldots,\btheta_N$ of $\mathcal K_{\varepsilon,M}$. Chebyshev's inequality gives, for every $a>0$,
\[
\mathbb P_n\left(
 \max_{1\leq j\leq N}
 \left|
 \frac1n\sum_{i=1}^n
 \left(\lambda_{\btheta_j}(Z_{n,i})
 -\E_{\mathbb P_n}\left[\lambda_{\btheta_j}(Z_{n,1})\right]\right)
 \right|>a\right)
 \leq \frac{NC}{na^2}
 \]
 and
\[
\frac1n\sum_{i=1}^nG(Z_{n,i})
 =\E_{\mathbb P_n}\left[G(Z_{n,1})\right]
   +O_{\mathbb P_n}(n^{-\frac{1}{2}})
 =O_{\mathbb P_n}(1).
\]
For every $\btheta$, choose $\btheta_j$ with $\norm{\btheta-\btheta_j}\leq\delta$. The displayed gradient bound then yields
\[
\begin{aligned}
&\sup_{\btheta\in\mathcal K_{\varepsilon,M}}
\left|
\frac1n\sum_{i=1}^n\lambda_{\btheta}(Z_{n,i})-L(\btheta)
\right|
\leq
\max_{1\leq j\leq N}
\left|
\frac1n\sum_{i=1}^n\lambda_{\btheta_j}(Z_{n,i})
-L(\btheta_j)
\right|
+\delta\left(\frac1n\sum_{i=1}^nG(Z_{n,i})+\E\left[G(Z_0)\right]\right).
\end{aligned}
\]
Letting $n\to\infty$ and then $\delta\downarrow0$ proves
\begin{equation}\label{eq:gpd_uniform_lln}
\sup_{(s,\eta)\in[\varepsilon,M]\times[0,\bar\xi]}
 \left|
 \frac{1}{n}\ell_{n,\mu}(s,\eta)-L(s,\eta)
 \right|
 \xrightarrow{\mathbb P_n}0,
\end{equation}
where $\ell_{n,\mu}$ denotes the likelihood with location fixed at $\mu$. Since
\[
 \sup_{\mathcal K_{\varepsilon,M}}
 |\partial_z\lambda(z;s,\eta)|
 \leq \frac{1+\bar\xi}{\varepsilon},
\]
replacing $Z_{n,i}$ by $Y_i=Z_{n,i}-(m_n-\mu)$ changes the normalized criterion by at most $C(m_n-\mu)=O_{\mathbb P_n}(n^{-1})$. Thus \eqref{eq:gpd_uniform_lln} also holds with $\ell_{n,\mu}$ replaced by $\ell_n^p$.

Let $f_{s,\eta}$ denote the GPD density with location $\mu$, scale $s$, and shape $\eta$. Jensen's inequality gives
\[
\begin{aligned}
L(s,\eta)-L(\sigma,0)
&=\E\left[\log{\frac{f_{s,\eta}(Z_0+\mu)}
                         {f_{\sigma,0}(Z_0+\mu)}}\right]\leq
\log{\E\left[\frac{f_{s,\eta}(Z_0+\mu)}
                         {f_{\sigma,0}(Z_0+\mu)}\right]}
=0.
\end{aligned}
\]
Equality requires the two densities to coincide. When $\eta>0$,
\[
 \frac{\partial^2}{\partial z^2}\log{f_{s,\eta}(z+\mu)}
 =\frac{\eta(1+\eta)}{(s+\eta z)^2}>0,
\]
whereas $\log{f_{\sigma,0}(z+\mu)}$ is affine in $z$. When $\eta=0$, equality requires $s=\sigma$. Hence $L$ has the unique maximizer $(\sigma,0)$.

It remains to control the tails of the scale space. Choose $c_0>0$ so that $q\coloneqq\mathbb P(Z_0>2c_0)>\bar\xi$, and fix $q'\in(\bar\xi,q)$. If
\[
 p_n\coloneqq\mathbb P_n(Z_{n,1}>2c_0),
\]
then $p_n\to q$ and
\[
 \frac1n\sum_{i=1}^n\mathds{1}\left(Z_{n,i}>2c_0\right)
 =p_n+O_{\mathbb P_n}(n^{-\frac{1}{2}})
 \xrightarrow{\mathbb P_n}q,
\]
because the variance of the average is at most $1/(4n)$. Since $m_n-\mu=o_{\mathbb P_n}(1)$, with probability tending to one at least $q'n$ of the excesses $Y_i$ exceed $c_0$. On this event, \eqref{eq:gpd_small_scale_bound} gives
\[
 \sup_{\eta\in[0,\bar\xi]}\frac1n\ell_n^p(s,\eta)
 \leq
 -\log{s}-\frac{q'}{\bar\xi}\log{1+\frac{\bar\xi c_0}{s}},
\]
which tends to $-\infty$ as $s\downarrow0$ because $q'/\bar\xi>1$. The bound $\ell_n^p/n\leq-\log{s}$ handles $s\to\infty$. Moreover, \eqref{eq:gpd_uniform_lln} at $(\sigma,0)$ gives
\[
 \frac1n\ell_n^p(\sigma,0)
 \xrightarrow{\mathbb P_n}L(\sigma,0).
\]
The two tail bounds therefore allow fixed $0<\varepsilon<M<\infty$ to be chosen so that, with probability tending to one, no global maximizer has $s\notin[\varepsilon,M]$. On this rectangle, \eqref{eq:gpd_uniform_lln} and the unique maximization of $L$ at $(\sigma,0)$ give, for every global maximizer,
\begin{equation}\label{eq:gpd_global_consistency}
\widehat\sigma_n\xrightarrow{\mathbb P_n}\sigma
\dand
 \widehat\xi_n\xrightarrow{\mathbb P_n}0.
\end{equation}

For $Z_0\sim\mathrm{GPD}(0,\sigma,0)$, differentiating $\lambda(z;s,\eta)$ with respect to $(s,\eta)$ and evaluating at $(\sigma,0)$ gives
\[
 \nabla_{s,\eta}\lambda(Z_0;\sigma,0)
 =
 \begin{pmatrix}
 \dfrac{(Z_0-\sigma)}{\sigma^2}\\[9pt]
 \dfrac{Z_0^2}{(2\sigma^2)}-\dfrac{Z_0}{\sigma}
 \end{pmatrix}.
\]
Differentiating once more and using $\E\left[Z_0^k\right]=k!\sigma^k$ gives
\begin{equation}\label{eq:gpd_I0}
\mathcal I_0\coloneqq-\E\left[\nabla_{s,\eta}^2\lambda(Z_0;\sigma,0)\right]
 =
 \begin{pmatrix}\dfrac{1}{\sigma^2}&\dfrac{1}{\sigma}\\[9pt]
                 \dfrac{1}{\sigma}&2\end{pmatrix},
\end{equation}
whose determinant is $\frac{1}{\sigma^2}>0$. On every fixed compact rectangle $\mathcal K\subset(0,\infty)\times[0,\bar\xi]$,
\[
 \sup_{(s,\eta)\in\mathcal K}
 \left(\norm{\nabla_{s,\eta}^2\lambda(z;s,\eta)}
       +\norm{\nabla_{s,\eta}^3\lambda(z;s,\eta)}\right)
 \leq C(1+z^4).
\]
Applying the argument leading to \eqref{eq:gpd_uniform_lln} componentwise to the Hessian therefore gives, on a sufficiently small fixed neighborhood $\mathcal V$ of $(\sigma,0)$,
\begin{equation}\label{eq:gpd_uniform_hessian}
 \sup_{(s,\eta)\in\mathcal V}
 \norm{
 \frac1n\nabla^2\ell_{n,\mu}(s,\eta)
 -\E\left[\nabla_{s,\eta}^2\lambda(Z_0;s,\eta)\right]
 }
 \xrightarrow{\mathbb P_n}0.
\end{equation}

To compare the fixed-location and profile likelihoods, set
\[
 \bs(z;s,\eta)\coloneqq\nabla_{s,\eta}\lambda(z;s,\eta)\dand
 \Hb(z;s,\eta)\coloneqq\nabla_{s,\eta}^2\lambda(z;s,\eta).
\]
On $\mathcal V$,
\begin{equation}\label{eq:gpd_derivative_envelope}
 \sup_{(s,\eta)\in\mathcal V}
 \left(\norm{\partial_z\bs(z;s,\eta)}
       +\norm{\partial_z\Hb(z;s,\eta)}\right)
 \leq C(1+z^3).
\end{equation}
Put $\delta_n\coloneqq m_n-\mu$, so that $Y_i\coloneqq Z_{n,i}-\delta_n$ and $0\leq Y_i\leq Z_{n,i}$. The mean value theorem and \eqref{eq:gpd_derivative_envelope} give, for every $(s,\eta)\in\mathcal V$,
\[
\norm{\bs(Y_i;s,\eta)-\bs(Z_{n,i};s,\eta)}
+\norm{\Hb(Y_i;s,\eta)-\Hb(Z_{n,i};s,\eta)}
\leq C\delta_n(1+Z_{n,i}^3),
\]
because each intermediate point lies between $Y_i$ and $Z_{n,i}$. By \eqref{eq:gpd_min_limit}, $\delta_n=O_{\mathbb P_n}(n^{-1})$, while \eqref{eq:gpd_moments} and Markov's inequality give
\[
    \frac1n\sum_{i=1}^n(1+Z_{n,i}^3)=O_{\mathbb P_n}(1).
\]
Consequently,
\begin{equation}\label{eq:gpd_profile_derivative_comparison}
\begin{aligned}
&\sup_{\mathcal V}
 \norm{\nabla\ell_n^p-\nabla\ell_{n,\mu}}
 =O_{\mathbb P_n}(1) \dand
 \frac1n\sup_{\mathcal V}
 \norm{\nabla^2\ell_n^p-\nabla^2\ell_{n,\mu}}
 =O_{\mathbb P_n}\left(\frac{1}{n}\right).
\end{aligned}
\end{equation}
By \eqref{eq:gpd_I0} and continuity, there is a fixed convex neighborhood $\mathcal V_0\subseteq\mathcal V$, understood relative to $(0,\infty)\times[0,\bar\xi]$, on which $\ell_n^p$ is strictly concave with probability tending to one. At $(\sigma,\xi_n)$, the fixed-location score has mean zero and
\[
 \E_{\mathbb P_n}\left[
 \norm{\nabla\ell_{n,\mu}(\sigma,\xi_n)}^2
 \right]
 =n\E_{\mathbb P_n}\left[\norm{\bs(Z_{n,1};\sigma,\xi_n)}^2\right]
 =O(n).
\]
Consequently, the first comparison above gives
\[
 \nabla\ell_n^p(\sigma,\xi_n)=O_{\mathbb P_n}(\sqrt n).
\]
For some $\lambda_0>0$, with probability tending to one the smallest eigenvalue of $-\nabla^2\ell_n^p$ is at least $\lambda_0n$ throughout $\mathcal V_0$. If $\norm{\bh}=M/\sqrt n$ and the segment from $(\sigma,\xi_n)$ to $(\sigma,\xi_n)+\bh$ lies in $\mathcal V_0$, Taylor's theorem gives
\[
 \ell_n^p\left((\sigma,\xi_n)+\bh\right)-\ell_n^p(\sigma,\xi_n)
 \leq
 \norm{\nabla\ell_n^p(\sigma,\xi_n)}\norm{\bh}
 -\frac{\lambda_0n}{2}\norm{\bh}^2
 =M O_{\mathbb P_n}(1)-\frac{\lambda_0M^2}{2}.
\]
Choose a deterministic sequence $M_n\to\infty$ sufficiently slowly that
\[
    \frac{M_n}{\sqrt n\,\xi_n}\longrightarrow0.
\]
Then $M_n/\sqrt n\to0$, and $\overline{B}_{M_n/\sqrt n}((\sigma,\xi_n)^\intercal)\subseteq\mathcal V_0\cap((0,\infty)\times(0,\bar\xi))$ for all sufficiently large $n$.
Taking $M=M_n$ in the preceding bound and using $\norm{\nabla\ell_n^p(\sigma,\xi_n)}/\sqrt n=O_{\mathbb P_n}(1)$ shows that, with probability tending to one, the likelihood on the boundary of $\overline{B}_{M_n/\sqrt n}((\sigma,\xi_n)^\intercal)$ is smaller than at its center. The interior of this ball therefore contains a local maximizer, and strict concavity makes it the unique score root in $\mathcal V_0$. By \eqref{eq:gpd_global_consistency}, every global maximizer lies in $\mathcal V_0$ with probability tending to one, so a global maximizer exists and coincides with this local maximizer with probability tending to one.

To obtain the rate, fix $\varepsilon>0$ and choose a constant $M$ large enough that the preceding boundary inequality holds with probability at least $1-\varepsilon/2$ for all sufficiently large $n$. Together with strict concavity on $\mathcal V_0$ and \eqref{eq:gpd_global_consistency}, this implies that the unique global maximizer lies within distance $M/\sqrt n$ of $(\sigma,\xi_n)$ with probability at least $1-\varepsilon$ for all sufficiently large $n$. Hence
\begin{equation}\label{eq:gpd_mle_rate}
\widehat\sigma_n-\sigma=O_{\mathbb P_n}\left(\frac{1}{\sqrt n}\right)
\dand
 \widehat\xi_n-\xi_n=O_{\mathbb P_n}\left(\frac{1}{\sqrt n}\right).
\end{equation}
\end{proof}

\subsubsection{Regularized Score and Truncation Error}

For $t_i(s)\coloneqq Y_i/s$, define
\[
 A_n(s,\eta)
 \coloneqq\frac1n\sum_{i=1}^n
   \frac{(1+\eta)t_i(s)}{1+\eta t_i(s)}
 \dand
 B_n(s,\eta)
 \coloneqq\frac1n\sum_{i=1}^n
   \frac{\log{1+\eta t_i(s)}}{\eta}.
\]
Let $S_s(s,\eta)\coloneqq\partial_s\ell_n^p(s,\eta)$ and $S_\eta(s,\eta)\coloneqq\partial_\eta\ell_n^p(s,\eta)$. Differentiating \eqref{eq:gpd_profile_lik} gives
\begin{equation}\label{eq:gpd_profile_scores}
\frac{S_s(s,\eta)}n=\frac{A_n(s,\eta)-1}{s}
\dand
 \frac{S_\eta(s,\eta)}n=\frac{B_n(s,\eta)-A_n(s,\eta)}{\eta}.
\end{equation}
Consequently,
\begin{equation}\label{eq:gpd_root_preserving}
\bq_n(s,\eta)\coloneqq
 \begin{pmatrix}B_n(s,\eta)-1\\[9pt]
 \dfrac{\left(B_n(s,\eta)-A_n(s,\eta)\right)}{\eta}\end{pmatrix}
 =
 \begin{pmatrix}s&\eta\\0&1\end{pmatrix}
 \frac1n\begin{pmatrix}S_s(s,\eta)\\S_\eta(s,\eta)\end{pmatrix}.
\end{equation}
The matrix in \eqref{eq:gpd_root_preserving} has determinant $s>0$, so $\bq_n=\bzero$ and the profile score equation have the same interior roots. We use the continuous extensions of these expressions at $\eta=0$.

The surrogate used in \cref{sub:GPD} is the first-order Taylor truncation in the shape coordinate: in each component of $\bq_n$, it keeps the constant and linear powers $\eta^0$ and $\eta^1$ and discards powers $\eta^2$ and above.

\begin{lemm}\label{lemm:gpd_truncation}
Let $\widetilde\bq_n\coloneqq(\widetilde q_{1,n},\widetilde q_{2,n})^{\mathsf T}$, where
\begin{equation}\label{eq:gpd_surrogate}
\begin{split}
 \widetilde q_{1,n}(s,\eta)
 &\coloneqq
 \frac{M_{1,n}}{s}-1-\frac{\eta M_{2,n}}{(2s^2)},\\
 \widetilde q_{2,n}(s,\eta)
 &\coloneqq
 -\frac{M_{1,n}}{s}+\frac{M_{2,n}}{(2s^2)}
 +\eta\left(\frac{M_{2,n}}{s^2}-\frac{2M_{3,n}}{(3s^3)}\right).
 \end{split}
\end{equation}
For all $s>0$ and $\eta\geq0$,
\begin{equation}\label{eq:gpd_trunc_bound1}
\left|q_{1,n}(s,\eta)-\widetilde q_{1,n}(s,\eta)\right|
 \leq\frac{\eta^2M_{3,n}}{3s^3},
\end{equation}
and
\begin{equation}\label{eq:gpd_trunc_bound2}
\left|q_{2,n}(s,\eta)-\widetilde q_{2,n}(s,\eta)\right|
 \leq
 \eta^2\left(\frac{M_{3,n}}{s^3}
              +\frac{3M_{4,n}}{4s^4}\right).
\end{equation}
Thus, uniformly for $s$ in a fixed compact neighborhood of $\sigma$ and $\eta=O(\xi_n)$,
\begin{equation}\label{eq:gpd_uniform_remainder}
\norm{\widetilde\bq_n(s,\eta)-\bq_n(s,\eta)}
 =O_{\mathbb P_n}(\xi_n^2).
\end{equation}
\end{lemm}

\begin{proof}
It suffices to consider $\eta>0$, since the bounds at $\eta=0$ follow from the continuous extensions. 
For $x\geq0$, put $r_1(x)\coloneqq\log{1+x}-x+x^2/2$. Since
\begin{equation}\label{eq:gpd_log_bound_derivatives}
r_1(0)=0
\dand
r_1'(x)=\frac{x^2}{1+x}
\dand
\left(\frac{x^3}{3}-r_1(x)\right)'=\frac{x^3}{1+x},
\end{equation}
both derivatives in \eqref{eq:gpd_log_bound_derivatives} are nonnegative. Since both $r_1(x)$ and $x^3/3-r_1(x)$ vanish at $x=0$, these identities give
\begin{equation}\label{eq:gpd_log_bound}
0\leq\log{1+x}-x+\frac{x^2}{2}\leq\frac{x^3}{3}.
\end{equation}
Applying \eqref{eq:gpd_log_bound} with $x=\eta t_i(s)$, dividing by $\eta$, and averaging gives
\[
0\leq
q_{1,n}(s,\eta)-\widetilde q_{1,n}(s,\eta)
\leq
\frac{\eta^2}{3n}\sum_{i=1}^nt_i(s)^3
=\frac{\eta^2M_{3,n}}{3s^3},
\]
which is \eqref{eq:gpd_trunc_bound1}.

For the second component, let
\[
h(x)\coloneqq\log{1+x}-\frac{x}{1+x}
\dand
r_2(x)\coloneqq h(x)-\frac{x^2}{2}+\frac{2x^3}{3}.
\]
The identities
\begin{equation}\label{eq:gpd_h_bound_derivatives}
r_2(0)=0
\dand
r_2'(x)=\frac{x^3(3+2x)}{(1+x)^2}
\dand
\left(\frac{3x^4}{4}-r_2(x)\right)'
=\frac{x^4(4+3x)}{(1+x)^2}
\end{equation}
show that
\begin{equation}\label{eq:gpd_h_bound}
h(x)=\frac{x^2}{2}-\frac{2x^3}{3}+r_2(x)
 \dand 0\leq r_2(x)\leq\frac{3x^4}{4}.
\end{equation}
Moreover,
\begin{equation}\label{eq:gpd_g2_identity}
 \frac{B_n-A_n}{\eta}
 =
 \frac1n\sum_{i=1}^n
 \left[
 \frac{h\left(\eta t_i(s)\right)}{\eta^2}
 -\frac{t_i(s)}{1+\eta t_i(s)}
 \right]
\end{equation}
and
\begin{equation}\label{eq:gpd_rational_identity}
 -\frac{t}{1+\eta t}
 =-t+\eta t^2-\frac{\eta^2t^3}{1+\eta t}
\end{equation}
hold exactly. Substituting \eqref{eq:gpd_h_bound} and \eqref{eq:gpd_rational_identity} into \eqref{eq:gpd_g2_identity} gives
\begin{equation}\label{eq:gpd_g2_remainder}
q_{2,n}(s,\eta)-\widetilde q_{2,n}(s,\eta)
=\frac1n\sum_{i=1}^n
\left[
\frac{r_2\left(\eta t_i(s)\right)}{\eta^2}
-\frac{\eta^2t_i(s)^3}{1+\eta t_i(s)}
\right].
\end{equation}
Using \eqref{eq:gpd_h_bound} in \eqref{eq:gpd_g2_remainder}, together with $1+\eta t_i(s)\geq1$, gives
\[
\left|q_{2,n}(s,\eta)-\widetilde q_{2,n}(s,\eta)\right|
\leq
\eta^2\left(
\frac{M_{3,n}}{s^3}+\frac{3M_{4,n}}{4s^4}
\right),
\]
which proves \eqref{eq:gpd_trunc_bound2}. Finally, \cref{lemm:gpd_moments} shows that $M_{3,n}$ and $M_{4,n}$ are bounded in probability. Choose $0<\varepsilon<\sigma<M<\infty$ so that the fixed compact scale neighborhood in the statement is contained in $[\varepsilon,M]$. Then the powers of $s^{-1}$ are uniformly bounded there. Equations \eqref{eq:gpd_trunc_bound1} and \eqref{eq:gpd_trunc_bound2}, together with $\eta=O(\xi_n)$, therefore give \eqref{eq:gpd_uniform_remainder}.
\end{proof}

\subsubsection{Explicit Surrogate Root and Population Discrepancy}

\begin{lemm}\label{lemm:gpd_surrogate}
Let
\[
 C_n\coloneqq6M_{1,n}M_{2,n}+8M_{3,n}
 \dand
 E_n\coloneqq8M_{1,n}M_{3,n}-3M_{2,n}^2.
\]
With probability tending to one, the root of $\widetilde\bq_n(s,\eta)=\bzero$ whose scale coordinate converges to $\sigma$ is
\begin{equation}\label{eq:gpd_explicit_root}
 \widetilde\sigma_n
 \coloneqq
 \frac{C_n-\sqrt{C_n^2-48M_{2,n}E_n}}{24M_{2,n}}
 \dand
 \widetilde\xi_n
 \coloneqq
 \frac{2\widetilde\sigma_n
       (M_{1,n}-\widetilde\sigma_n)}{M_{2,n}}.
\end{equation}
If $\sqrt n\,\xi_n\to\infty$, this root lies in $(0,\infty)\times(0,\bar\xi)$ with probability tending to one and
\begin{equation}\label{eq:gpd_surrogate_rate}
\widetilde\sigma_n-\sigma
 =O_{\mathbb P_n}(n^{-\frac{1}{2}}+\xi_n^2)
 \dand
 \widetilde\xi_n-\xi_n
 =O_{\mathbb P_n}(n^{-\frac{1}{2}}+\xi_n^2).
\end{equation}
Replacing $M_{k,n}$ by $m_k(\xi)$ in \eqref{eq:gpd_explicit_root} defines the population branch $(\sigma_*(\xi),\xi_*(\xi))$, which satisfies
\begin{equation}\label{eq:gpd_population_bias}
\sigma_*(\xi)
 =\sigma\left(1+8\xi^2-68\xi^3+O(\xi^4)\right) \dand
 \xi_*(\xi)=\xi-10\xi^2+100\xi^3+O(\xi^4).
\end{equation}
\end{lemm}

\begin{proof}
Setting $\widetilde q_{1,n}(s,\eta)=0$ in \eqref{eq:gpd_surrogate} gives
\begin{equation}\label{eq:gpd_shape_from_scale}
 \eta=\frac{2s(M_{1,n}-s)}{M_{2,n}}.
\end{equation}
Substituting \eqref{eq:gpd_shape_from_scale} into $\widetilde q_{2,n}(s,\eta)=0$ in \eqref{eq:gpd_surrogate} and multiplying by $6M_{2,n}s^2$ yields
\begin{equation}\label{eq:gpd_quadratic}
12M_{2,n}s^2-(6M_{1,n}M_{2,n}+8M_{3,n})s
 +(8M_{1,n}M_{3,n}-3M_{2,n}^2)=0,
\end{equation}
so the quadratic formula gives two candidate scale coordinates, where the corresponding shape coordinates are determined by \eqref{eq:gpd_shape_from_scale}. At $(M_1,M_2,M_3)=(\sigma,2\sigma^2,6\sigma^3)$, the two scale roots are $\sigma$ and $3\sigma/2$, and the discriminant is $144\sigma^6>0$. Thus the branch with the smaller scale coordinate is the unique branch converging to $(\sigma,0)$. Selecting the minus sign in the quadratic formula and then applying \eqref{eq:gpd_shape_from_scale} gives \eqref{eq:gpd_explicit_root}.

For the population branch, substitute
\[
 m_k(\xi)=\frac{k!\sigma^k}{\prod_{j=1}^k(1-j\xi)}
\]
into \eqref{eq:gpd_explicit_root}. The resulting series are
\begin{equation}\label{eq:gpd_population_intermediate}
 \begin{split}
 \frac{C(\xi)}{\sigma^3}&=60+336\xi+1332\xi^2+4632\xi^3+O(\xi^4),\\
 \frac{E(\xi)}{\sigma^4}&=36+264\xi+1260\xi^2+4992\xi^3+O(\xi^4),\\
 \frac{\left(C(\xi)^2-48m_2(\xi)E(\xi)\right)^{1/2}}{\sigma^3}
 &=12+192\xi+612\xi^2+6024\xi^3+O(\xi^4).
 \end{split}
\end{equation}
Substitution of \eqref{eq:gpd_population_intermediate} into \eqref{eq:gpd_explicit_root} gives \eqref{eq:gpd_population_bias}.

Let $\bm m_n\coloneqq(M_{1,n},M_{2,n},M_{3,n})^{\mathsf T}$ and $\bm m(\xi_n)\coloneqq\left(m_1(\xi_n),m_2(\xi_n),m_3(\xi_n)\right)^{\mathsf T}$. By \cref{lemm:gpd_moments},
\[
    \bm m_n-\bm m(\xi_n)=O_{\mathbb P_n}\left(\frac{1}{\sqrt n}\right).
\]
The discriminant in \eqref{eq:gpd_explicit_root} converges in probability to its value $144\sigma^6>0$ at the exponential moment vector. Thus the root is real with probability tending to one, and both coordinates in \eqref{eq:gpd_explicit_root} are smooth functions of the first three moments in a fixed neighborhood of that vector. The mean value theorem therefore gives
\[
\begin{pmatrix}\widetilde\sigma_n\\ \widetilde\xi_n\end{pmatrix}
-
\begin{pmatrix}\sigma_*(\xi_n)\\ \xi_*(\xi_n)\end{pmatrix}
=O_{\mathbb P_n}\left(\frac{1}{\sqrt n}\right).
\]
Together with \eqref{eq:gpd_population_bias}, this proves \eqref{eq:gpd_surrogate_rate}. Under $\sqrt n\,\xi_n\to\infty$,
\[
\frac{\widetilde\xi_n}{\xi_n}
=
\frac{\xi_*(\xi_n)}{\xi_n}
+O_{\mathbb P_n}\left(\frac{1}{\sqrt n\,\xi_n}\right)
\xrightarrow{\mathbb P_n}1,
\]
where the deterministic term converges to one by \eqref{eq:gpd_population_bias}. Because $\xi_n>0$, $\xi_n\to0$, and $\bar\xi>0$ is fixed, it follows that $\widetilde\xi_n\in(0,\bar\xi)$ with probability tending to one. Moreover,
\[
\widetilde\sigma_n
=\sigma_*(\xi_n)+O_{\mathbb P_n}(n^{-1/2})
=\sigma+o_{\mathbb P_n}(1),
\]
so $\widetilde\sigma_n>0$ with probability tending to one. On the complementary event, we use the completion convention of \cref{sec:defs}.
\end{proof}

\subsection{Uniform Root Stability}

For the root comparison, use the chart
\[
 \bphi(s,\eta)\coloneqq
 \begin{pmatrix}\dfrac{\eta}{s}\\[9pt]
    \dfrac{\eta}{s^2}\end{pmatrix}
 \dand
 U\coloneqq\left\{(u_2,u_3):u_2,u_3>0,\ \frac{u_2^2}{u_3}<\bar\xi\right\}.
\]
This map is one-to-one from $(0,\infty)\times(0,\bar\xi)$ onto $U$, with
\begin{equation}\label{eq:gpd_inverse_chart}
\bF(u_2,u_3)\coloneqq\bphi\inv(u_2,u_3)
 =
 \begin{pmatrix}\dfrac{u_2}{u_3}\\[9pt]
    \dfrac{u_2^2}{u_3}\end{pmatrix}.
\end{equation}
At $(\sigma,\xi_n)$,
\begin{equation}\label{eq:gpd_inverse_jacobian}
\bu_n\coloneqq
 \begin{pmatrix}\dfrac{\xi_n}{\sigma}\\[9pt]\dfrac{\xi_n}{\sigma^2}\end{pmatrix}\dand
 \Jb_{\bF}(\bu_n)=
 \begin{pmatrix}
 \dfrac{\sigma^2}{\xi_n}&-\dfrac{\sigma^3}{\xi_n}\\[9pt]
 2\sigma&-\sigma^2
 \end{pmatrix}.
\end{equation}

\begin{prop}\label{prop:gpd_root_stability}
Assume $\sqrt n\,\xi_n\to\infty$, and let $(\widehat\sigma_n,\widehat\xi_n)$ be the restricted global MLE from \cref{prop:gpd_global_mle}. Then
\begin{equation}\label{eq:gpd_original_approx}
\norm{\begin{pmatrix}\widetilde\sigma_n\\\widetilde\xi_n\end{pmatrix}
       -\begin{pmatrix}\widehat\sigma_n\\\widehat\xi_n\end{pmatrix}}
 =O_{\mathbb P_n}(\xi_n^2)=o_{\mathbb P_n}(1).
\end{equation}
If additionally $\sqrt n\,\xi_n^2\to0$, then
\begin{equation}\label{eq:gpd_rootn_equivalence}
\sqrt n
 \begin{pmatrix}
 \widetilde\sigma_n-\widehat\sigma_n\\
 \widetilde\xi_n-\widehat\xi_n
 \end{pmatrix}
 =o_{\mathbb P_n}(1).
\end{equation}
\end{prop}

\begin{proof}
Set $c_n\coloneqq n^{-1/2}+\xi_n^2$. By \eqref{eq:gpd_mle_rate} and \eqref{eq:gpd_surrogate_rate}, both original parameter vectors differ from $(\sigma,\xi_n)^{\mathsf T}$ by $O_{\mathbb P_n}(c_n)$. The assumption $\sqrt n\,\xi_n\to\infty$, together with $\xi_n\to0$, implies
\[
\frac{c_n}{\xi_n}
=\frac{1}{\sqrt n\,\xi_n}+\xi_n
\longrightarrow0.
\]
Consequently, uniformly along either line segment joining these vectors to $(\sigma,\xi_n)^{\mathsf T}$,
\[
s-\sigma=O_{\mathbb P_n}(c_n)
\dand
\eta-\xi_n=O_{\mathbb P_n}(c_n)
\dand
\frac{s}{\sigma}=1+O_{\mathbb P_n}(c_n)
\dand
\frac{\eta}{\xi_n}=1+O_{\mathbb P_n}\left(\frac{c_n}{\xi_n}\right).
\]
In particular, $s$ is bounded away from zero with probability tending to one. Direct multiplication gives
\[
\Jb_{\bF}(\bu_n)\Jb_{\bphi}(s,\eta)
=
\begin{pmatrix}
\dfrac{\sigma^2\eta(2\sigma-s)}{\xi_ns^3}
&\dfrac{\sigma^2(s-\sigma)}{\xi_ns^2}\\[9pt]
\dfrac{2\sigma\eta(\sigma-s)}{s^3}
&\dfrac{2\sigma}{s}-\dfrac{\sigma^2}{s^2}
\end{pmatrix}
=\Ib_2+o_{\mathbb P_n}(1).
\]
Indeed, the four entries on the right are, respectively,
\[
1+O_{\mathbb P_n}\left(\frac{c_n}{\xi_n}\right),
\qquad
O_{\mathbb P_n}\left(\frac{c_n}{\xi_n}\right),
\qquad
O_{\mathbb P_n}(\xi_nc_n),
\qquad
1+O_{\mathbb P_n}(c_n^2),
\]
uniformly along either segment. For
\[
\btheta_n^\dagger
\in
\left\{
(\widehat\sigma_n,\widehat\xi_n)^{\mathsf T},
(\widetilde\sigma_n,\widetilde\xi_n)^{\mathsf T}
\right\} \dand
\bu_n^\dagger\coloneqq\bphi(\btheta_n^\dagger),
\]
the fundamental theorem of calculus gives
\[
\begin{aligned}
\Jb_{\bF}(\bu_n)(\bu_n^\dagger-\bu_n)
&=
\left(
\int_0^1
\Jb_{\bF}(\bu_n)
\Jb_{\bphi}\left(
(\sigma,\xi_n)^{\mathsf T}
+t\left[\btheta_n^\dagger-(\sigma,\xi_n)^{\mathsf T}\right]
\right)
\dif t
\right)
\left(\btheta_n^\dagger-(\sigma,\xi_n)^{\mathsf T}\right)\\
&=
\left\{\Ib_2+o_{\mathbb P_n}(1)\right\}
O_{\mathbb P_n}(c_n)
=O_{\mathbb P_n}(c_n).
\end{aligned}
\]
Thus
\begin{equation}\label{eq:gpd_transformed_root_localization}
\Jb_{\bF}(\bu_n)(\widehat\bu_n-\bu_n)
=O_{\mathbb P_n}(c_n)
\dand
\Jb_{\bF}(\bu_n)(\widetilde\bu_n-\bu_n)
=O_{\mathbb P_n}(c_n).
\end{equation}
Choose a deterministic sequence $t_n\to\infty$ slowly enough that $t_nc_n=o(\xi_n)$. The deterministic convex sets
\begin{equation}\label{eq:gpd_localization_sets}
D_n\coloneqq
\left\{
\bu:
\norm{\Jb_{\bF}(\bu_n)(\bu-\bu_n)}\leq t_nc_n
\right\}
\end{equation}
then contain both transformed roots with probability tending to one by \eqref{eq:gpd_transformed_root_localization}.

We next verify that $D_n\subseteq U$ eventually and record the uniform bounds needed below. For $\bu\in D_n$, write
\[
\bv\coloneqq\Jb_{\bF}(\bu_n)(\bu-\bu_n)
=(v_1,v_2)^{\mathsf T}
\dand
    r\coloneqq\frac{u_2}{u_{2,n}}
\dand
q\coloneqq\frac{u_3}{u_{3,n}}.
\]
Since $\Jb_{\bF}(\bu_n)^{-1}=\Jb_{\bphi}(\sigma,\xi_n)$, direct substitution gives
\begin{equation}\label{eq:gpd_coordinate_ratios}
r=1-\frac{v_1}{\sigma}+\frac{v_2}{\xi_n}
\dand
q=1-\frac{2v_1}{\sigma}+\frac{v_2}{\xi_n}
\dand
q-r=-\frac{v_1}{\sigma}.
\end{equation}
Thus $r=1+o(1)$, $q=1+o(1)$, and $q-r=o(\xi_n)$ uniformly on $D_n$. Moreover,
\begin{equation}\label{eq:gpd_inverse_chart_on_localization_set}
\bF(\bu)
=
\begin{pmatrix}
\dfrac{\sigma r}{q}\\[9pt]
\xi_n \dfrac{r^2}{q}
\end{pmatrix}
=
\begin{pmatrix}
\sigma+o(1)\\
\xi_n\left(1+o(1)\right)
\end{pmatrix}
\end{equation}
uniformly on $D_n$. Consequently, $D_n\subseteq U$ eventually. An explicit calculation using \eqref{eq:gpd_coordinate_ratios} also gives
\begin{equation}\label{eq:gpd_normalized_inverse_jacobian}
\Jb_{\bF}(\bu)\Jb_{\bF}(\bu_n)^{-1}
=
\begin{pmatrix}
-\dfrac{1}{q}+\dfrac{2r}{q^2} & \dfrac{(q-r)}{(u_{2,n}q^2)}\\[9pt]
2u_{2,n}r\dfrac{(r-q)}{q^2} & \dfrac{2r}{q}-\dfrac{r^2}{q^2}
\end{pmatrix}
=\Ib_2+o(1)
\end{equation}
uniformly on $D_n$. Next, let $\bQ_n\coloneqq\bq_n\circ\bF$ and $\widetilde\bQ_n\coloneqq\widetilde\bq_n\circ\bF$, and put
\[
 \Hb_n(s,\eta)\coloneqq\Jb_{\widetilde\bq_n}(s,\eta)
 =
 \begin{pmatrix}
 -\dfrac{M_{1,n}}{s^2}+\dfrac{\eta M_{2,n}}{s^3}&
 -\dfrac{M_{2,n}}{(2s^2)}\\[9pt]
 \dfrac{M_{1,n}}{s^2}-\dfrac{M_{2,n}}{s^3}
 +\eta\left(-\dfrac{2M_{2,n}}{s^3}+\dfrac{2M_{3,n}}{s^4}\right)&
 \dfrac{M_{2,n}}{s^2}-\dfrac{2M_{3,n}}{(3s^3)}
 \end{pmatrix}.
\]
The sample moment limits in \eqref{eq:gpd_sample_moments} give
\begin{equation}\label{eq:gpd_G0}
\Hb_n(\sigma,\xi_n)\xrightarrow{\mathbb P_n}
 \Gb_0\coloneqq
 \begin{pmatrix}-\dfrac{1}{\sigma}&-1\\[9pt]
    -\dfrac{1}{\sigma}&-2\end{pmatrix},
\end{equation}
whose determinant is $1/\sigma$. Set
\[
 \Jb_n\coloneqq\Hb_n(\sigma,\xi_n)\Jb_{\bF}(\bu_n).
\]
Equations \eqref{eq:gpd_sample_moments} and \eqref{eq:gpd_inverse_chart_on_localization_set}, together with the displayed formula for $\Hb_n$, imply
\begin{equation}\label{eq:gpd_normalized_jacobians}
\sup_{\bu\in D_n}
\norm{\Hb_n\left(\bF(\bu)\right)-\Hb_n(\sigma,\xi_n)}
=o_{\mathbb P_n}(1)
\dand
\sup_{\bu\in D_n}
\norm{\Jb_{\bF}(\bu)\Jb_{\bF}(\bu_n)^{-1}-\Ib_2}
=o(1).
\end{equation}
In particular,
\begin{equation}\label{eq:gpd_composed_jacobians}
\Jb_{\widetilde\bQ_n}(\bu)
=\Hb_n\left(\bF(\bu)\right)\Jb_{\bF}(\bu)
\dand
\Jb_{\bF}(\bu_n)\Jb_n^{-1}
=\Hb_n(\sigma,\xi_n)^{-1}.
\end{equation}
The normalized Jacobian difference required by \cref{prop:matrix_root_stability} is therefore
\[
\Hb_n(\sigma,\xi_n)^{-1}
\Hb_n\left(\bF(\bu)\right)
\Jb_{\bF}(\bu)\Jb_{\bF}(\bu_n)^{-1}-\Ib_2.
\]
Combining the two uniform approximations in \eqref{eq:gpd_normalized_jacobians} with $\Hb_n(\sigma,\xi_n)\inP\Gb_0$ gives
\begin{equation}\label{eq:gpd_matrix_jacobian_condition}
\sup_{\bu\in D_n}
 \norm{
 \Jb_{\bF}(\bu_n)\Jb_n^{-1}
 \left(\Jb_{\widetilde\bQ_n}(\bu)-\Jb_n\right)
 \Jb_{\bF}(\bu_n)^{-1}}
 =o_{\mathbb P_n}(1).
\end{equation}
By \eqref{eq:gpd_inverse_chart_on_localization_set}, every point in $\bF(D_n)$ has $s$ in a fixed compact neighborhood of $\sigma$ and $\eta=O(\xi_n)$. Hence \cref{lemm:gpd_truncation} gives
\begin{equation}\label{eq:gpd_matrix_equation_condition}
\sup_{\bu\in D_n}
 \norm{
 \Jb_{\bF}(\bu_n)\Jb_n^{-1}
 \left(\widetilde\bQ_n(\bu)-\bQ_n(\bu)\right)}
 =
 \sup_{\bu\in D_n}
 \norm{
 \Hb_n(\sigma,\xi_n)^{-1}
 [\widetilde\bq_n\left(\bF(\bu)\right)-\bq_n\left(\bF(\bu)\right)]}
 =O_{\mathbb P_n}(\xi_n^2).
\end{equation}
To retain the order in \eqref{eq:gpd_matrix_equation_condition}, define
\[
\Ab_n
\coloneqq\int_0^1
\Jb_{\widetilde\bQ_n}
\bigl[\widehat\bu_n+t(\widetilde\bu_n-\widehat\bu_n)\bigr]
\dif t
\]
on the event that both roots belong to $D_n$, and put
\[
\Bb_n
\coloneqq
\Jb_{\bF}(\bu_n)\Jb_n^{-1}
(\Ab_n-\Jb_n)
\Jb_{\bF}(\bu_n)^{-1}.
\]
Equation \eqref{eq:gpd_matrix_jacobian_condition} gives $\Bb_n=o_{\mathbb P_n}(1)$. The root equations and the fundamental theorem of calculus then yield
\[
\begin{aligned}
&\Jb_{\bF}(\bu_n)(\widetilde\bu_n-\widehat\bu_n)
=
-(\Ib_2+\Bb_n)^{-1}
\Jb_{\bF}(\bu_n)\Jb_n^{-1}
\left(\widetilde\bQ_n(\widehat\bu_n)-\bQ_n(\widehat\bu_n)\right).
\end{aligned}
\]
Together with \eqref{eq:gpd_matrix_equation_condition}, this proves
\begin{equation}\label{eq:gpd_chart_matrix_rate}
\Jb_{\bF}(\bu_n)(\widetilde\bu_n-\widehat\bu_n)
 =O_{\mathbb P_n}(\xi_n^2).
\end{equation}
To translate this bound through the inverse chart, apply the fundamental theorem of calculus once more. Equations \eqref{eq:gpd_normalized_inverse_jacobian} and \eqref{eq:gpd_chart_matrix_rate} give
\begin{equation}\label{eq:gpd_inverse_chart_linearization}
\begin{aligned}
\bF(\widetilde\bu_n)-\bF(\widehat\bu_n)
&=
\left(
\int_0^1
\Jb_{\bF}\bigl(\widehat\bu_n
+t(\widetilde\bu_n-\widehat\bu_n)\bigr)
\Jb_{\bF}(\bu_n)^{-1}\dif t
\right)
\Jb_{\bF}(\bu_n)
(\widetilde\bu_n-\widehat\bu_n)\\
&=\left(\Ib_2+o_{\mathbb P_n}(1)\right)
\Jb_{\bF}(\bu_n)
(\widetilde\bu_n-\widehat\bu_n)
=O_{\mathbb P_n}(\xi_n^2),
\end{aligned}
\end{equation}
which is \eqref{eq:gpd_original_approx}. If $\sqrt n\,\xi_n^2\to0$, multiplying \eqref{eq:gpd_inverse_chart_linearization} by $\sqrt n$ gives \eqref{eq:gpd_rootn_equivalence}.
\end{proof}

\subsection{Limiting Distributions}

\begin{lemm}\label{lemm:gpd_score_clt}
Let
\[
V_{n,i}\coloneqq\frac{X_{n,i}-\mu}{\sigma}
\]
and define the score contribution for $(\sigma,\xi_n)$, with the location held fixed at $\mu$, by
\begin{equation}\label{eq:gpd_score_contributions}
\bs_{\xi_n}(V_{n,i})
\coloneqq
\begin{pmatrix}
\dfrac{V_{n,i}-1}{\sigma(1+\xi_nV_{n,i})}\\[9pt]
\dfrac{\log{1+\xi_nV_{n,i}}}{\xi_n^2}
-\dfrac{(1+\xi_n)V_{n,i}}{\xi_n(1+\xi_nV_{n,i})}
\end{pmatrix}.
\end{equation}
Then $\E_{\mathbb P_n}\left[\bs_{\xi_n}(V_{n,i})\right]=\bzero$ and
\begin{equation}\label{eq:gpd_information}
\Var_{\mathbb P_n}\left(\bs_{\xi_n}(V_{n,i})\right)
=
\mathcal I(\sigma,\xi_n)
 =
 \frac1{1+2\xi_n}
 \begin{pmatrix}
 \dfrac{1}{\sigma^2}&\dfrac{1}{\sigma(1+\xi_n)}\\[9pt]
 \dfrac{1}{\sigma(1+\xi_n)}&\dfrac{2}{(1+\xi_n)}
 \end{pmatrix}.
\end{equation}
Moreover,
\begin{equation}\label{eq:gpd_score_clt}
 \frac1{\sqrt n}\mathcal I(\sigma,\xi_n)^{-\frac{1}{2}}
 \sum_{i=1}^n\bs_{\xi_n}(V_{n,i})
 \inD\mathcal N(\bzero,\Ib_2).
\end{equation}
\end{lemm}

\begin{proof}
For $V=(X-\mu)/\sigma$, the transformation
\[
W=\frac{1}{\xi}\log{1+\xi V}
\]
has the $\mathrm{Exp}(1)$ distribution. Writing the two components of \eqref{eq:gpd_score_contributions} as $s_\sigma$ and $s_\xi$, respectively, we have
\[
 s_\sigma
 =\frac{1-(1+\xi)e^{-\xi W}}{\sigma\xi}
 \dand
 s_\xi
 =\frac W\xi-\frac{(1+\xi)(1-e^{-\xi W})}{\xi^2}.
\]
For $a>-1$ and $k>-1$, the elementary exponential moments
\[
 \E\left[e^{-aW}\right]=\frac1{1+a}
 \dand
 \E\left[W^ke^{-aW}\right]=\frac{\Gamma(k+1)}{(1+a)^{k+1}},
\]
show that both score components have mean zero and give
\[
 \E\left[s_\sigma^2\right]=\frac1{\sigma^2(1+2\xi)}
 \dand
 \E\left[s_\sigma s_\xi\right]=
 \frac1{\sigma(1+\xi)(1+2\xi)}
 \dand
 \E\left[s_\xi^2\right]=\frac2{(1+\xi)(1+2\xi)}.
\]
This proves \eqref{eq:gpd_information}. Direct inversion gives
\begin{equation}\label{eq:gpd_information_inverse}
\mathcal I(\sigma,\xi)^{-1}
 =
 \begin{pmatrix}
 2\sigma^2(1+\xi)&-\sigma(1+\xi)\\
 -\sigma(1+\xi)&(1+\xi)^2
 \end{pmatrix}
 \longrightarrow_{\xi\to0}
 \begin{pmatrix}2\sigma^2&-\sigma\\-\sigma&1\end{pmatrix}.
\end{equation}

It remains to verify the central limit theorem. For $x\geq0$,
\[
1-e^{-x}
=x-x^2\int_0^1(1-t)e^{-tx}\dif t.
\]
Consequently, the two score components satisfy
\begin{equation}\label{eq:gpd_score_moment_bounds}
\begin{aligned}
|s_\sigma|
&=\left|\frac{1-e^{-\xi W}}{\sigma\xi}
-\frac{e^{-\xi W}}{\sigma}\right|
\leq\frac{W+1}{\sigma},\\
|s_\xi|
&=\left|-W+(1+\xi)W^2
\int_0^1(1-t)e^{-t\xi W}\dif t\right|
\leq W+\frac{1+\bar\xi}{2}W^2.
\end{aligned}
\end{equation}
The same inequalities hold at $\xi=0$ by continuity. Since $W\sim\mathrm{Exp}(1)$, \eqref{eq:gpd_score_moment_bounds} implies
\begin{equation}\label{eq:gpd_uniform_score_third_moment}
\sup_n
\E_{\mathbb P_n}\left[\norm{\bs_{\xi_n}(V_{n,1})}^3\right]
<\infty.
\end{equation}
The matrices in \eqref{eq:gpd_information} are uniformly nonsingular for all sufficiently large $n$. Hence, for every fixed unit vector $\ba\in\R^2$, the variables
\[
T_{n,i}
\coloneqq\ba^{\mathsf T}\mathcal I(\sigma,\xi_n)^{-\frac{1}{2}}
\bs_{\xi_n}(V_{n,i})
\]
have mean zero, variance one, and uniformly bounded third absolute moments. Therefore,
\begin{equation}\label{eq:gpd_score_lyapunov_ratio}
\frac{\sum_{i=1}^n\E_{\mathbb P_n}\left[|T_{n,i}|^3\right]}
{\left(\sum_{i=1}^n\Var_{\mathbb P_n}(T_{n,i})\right)^{\frac{3}{2}}}
=\frac{\E_{\mathbb P_n}\left[|T_{n,1}|^3\right]}{\sqrt n}
=O(n^{-\frac{1}{2}}).
\end{equation}
Lyapunov's theorem and the Cram\'er--Wold device now give \eqref{eq:gpd_score_clt}.
\end{proof}

\begin{thm}\label{thm:gpd_limits}
Under $\sqrt n\,\xi_n\to\infty$,
\begin{equation}\label{eq:gpd_exact_clt}
\sqrt n
 \begin{pmatrix}
 \widehat\sigma_n-\sigma\\
 \widehat\xi_n-\xi_n
 \end{pmatrix}
 \inD
 \mathcal N\left(
 \bzero,\begin{pmatrix}2\sigma^2&-\sigma\\-\sigma&1\end{pmatrix}
 \right).
\end{equation}
If also $\sqrt n\,\xi_n^2\to0$, the same limit holds with $(\widehat\sigma_n,\widehat\xi_n)$ replaced by $(\widetilde\sigma_n,\widetilde\xi_n)$. Separately,
\begin{equation}\label{eq:gpd_location_limit}
\widehat\mu_n=\widetilde\mu_n=X_{(1)}
\dand
 \frac n\sigma(\widehat\mu_n-\mu)
 =
 \frac n\sigma(\widetilde\mu_n-\mu)
 \inD\mathrm{Exp}(1).
\end{equation}
\end{thm}

\begin{proof}
Put
\[
\btheta_n\coloneqq(\sigma,\xi_n)^{\mathsf T}
\dand
\widehat\btheta_n\coloneqq(\widehat\sigma_n,\widehat\xi_n)^{\mathsf T}.
\]
Let $\bS_n^p$ and $\Hb_n^p$ denote the profile score and its Jacobian. The first bound in \eqref{eq:gpd_profile_derivative_comparison} and \cref{lemm:gpd_score_clt} give
\begin{equation}\label{eq:gpd_profile_score_comparison}
\frac1{\sqrt n}\bS_n^p(\btheta_n)
=
\frac1{\sqrt n}\sum_{i=1}^n\bs_{\xi_n}(V_{n,i})
+o_{\mathbb P_n}(1).
\end{equation}

Equations \eqref{eq:gpd_uniform_hessian}, \eqref{eq:gpd_profile_derivative_comparison}, and \eqref{eq:gpd_mle_rate}, together with continuity of \eqref{eq:gpd_information}, imply
\begin{equation}\label{eq:gpd_profile_hessian_segment}
\sup_{0\leq t\leq1}
\norm{
-\frac1n\Hb_n^p\left(\btheta_n+t(\widehat\btheta_n-\btheta_n)\right)
-\mathcal I(\sigma,\xi_n)
}
=o_{\mathbb P_n}(1).
\end{equation}
By \cref{prop:gpd_global_mle}, $\widehat\btheta_n$ is an interior root of the profile score with probability tending to one. Applying the fundamental theorem of calculus to the profile score and using \eqref{eq:gpd_profile_hessian_segment}, we obtain
\begin{equation}\label{eq:gpd_mle_asymptotic_linearization}
\sqrt n
 (\widehat\btheta_n-\btheta_n)
=
\mathcal I(\sigma,\xi_n)^{-1}
\frac1{\sqrt n}\sum_{i=1}^n\bs_{\xi_n}(V_{n,i})
+o_{\mathbb P_n}(1).
\end{equation}
Equations \eqref{eq:gpd_score_clt}, \eqref{eq:gpd_information_inverse}, and \eqref{eq:gpd_mle_asymptotic_linearization} prove \eqref{eq:gpd_exact_clt}. Under $\sqrt n\,\xi_n^2\to0$, by Slutsky's theorem \eqref{eq:gpd_rootn_equivalence} transfers this limit to the approximate MLE. Finally, \eqref{eq:gpd_location_limit} follows from \eqref{eq:gpd_min_limit}.
\end{proof}

\subsection{Simulation Study}\label{app:gpd_simulation}

We used the same sample size grid and number of replications as in the gamma and zeta studies. We simulated from the restricted GPD with $\mu\coloneqq0$, $\sigma\coloneqq1$, $\xi_n\coloneqq0.3n^{-\gamma}$, and $\bar\xi\coloneqq0.16$, using $\gamma\in\{0.20,0.30,0.55\}$ and $n\in\{250,500,1000,2500,5000,10000,25000\}$. The three exponents distinguished the relevant regimes: $\gamma=0.20$ violated $\sqrt n\,\xi_n^2\to0$, $\gamma=0.30$ satisfied both conditions in $n^{-1/2}\ll\xi_n\ll n^{-1/4}$, and $\gamma=0.55$ violated $\sqrt n\,\xi_n\to\infty$.

As a closed-form benchmark, we used the probability weighted moment estimator of \citet[Section~3.3]{hoskingParameterQuantileEstimation1987a}. Their shape parameter $k$ corresponds to $-\xi$ in our convention. To translate their estimator, we used the fact that, for a nonnegative GPD excess $Y$ and $V\coloneqq1-F(Y)$, we have $V\sim\mathrm{Unif}(0,1)$ and $Y=\sigma(V^{-\xi}-1)/\xi$. Hence
\[
    \alpha_r\coloneqq\E\left[Y\left(1-F(Y)\right)^r\right]
 =\frac{\sigma}{\xi}\int_0^1(v^{r-\xi}-v^r)\dif v
 =\frac{\sigma}{(r+1)(r+1-\xi)}.
\]
In particular,
\[
 \alpha_0=\frac{\sigma}{1-\xi}
 \dand
 \alpha_1=\frac{\sigma}{2(2-\xi)}.
\]
Writing $Y_{(1)}\leq\cdots\leq Y_{(n)}$ for the ordered excesses over $X_{(1)}$, we defined
\[
 a_0\coloneqq\frac1n\sum_{j=1}^nY_{(j)}
 \dand
 a_1\coloneqq\frac1n\sum_{j=1}^n\frac{n-j}{n-1}Y_{(j)}.
\]
Replacing $(\alpha_0,\alpha_1)$ by $(a_0,a_1)$ and solving the two moment equations then gave
\[
 \widehat\xi_{\mathrm{PWM}}
 \coloneqq2-\frac{a_0}{a_0-2a_1}
 \dand
 \widehat\sigma_{\mathrm{PWM}}
 \coloneqq\frac{2a_0a_1}{a_0-2a_1}.
\]
\cite{hoskingParameterQuantileEstimation1987a} derived these formulas for a two-parameter GPD with a known threshold. Applying them to excesses over $X_{(1)}$ was our plug-in adaptation to the three-parameter model.

The exact restricted MLE was computed by comparing the interior likelihood root with the two shape boundary solutions; this calculation had no numerical failures. For our estimator, we used the selected raw surrogate root whenever it belonged to the restricted parameter space, and the shifted exponential boundary estimator $(M_{1,n},0)$ otherwise, exactly as specified in the construction above. Consequently, every comparison involving our constrained estimator used all $2000$ replications. The probability-weighted moment estimator was not constrained in this way, so its RMSE ratios used the replications in which it belonged to the parameter space, with the exact MLE evaluated on the same subset. We reported both the raw root replacement frequency and the PWM exclusion frequency.

\begin{figure}[ht]
\centering
\includegraphics[width=0.90\textwidth]{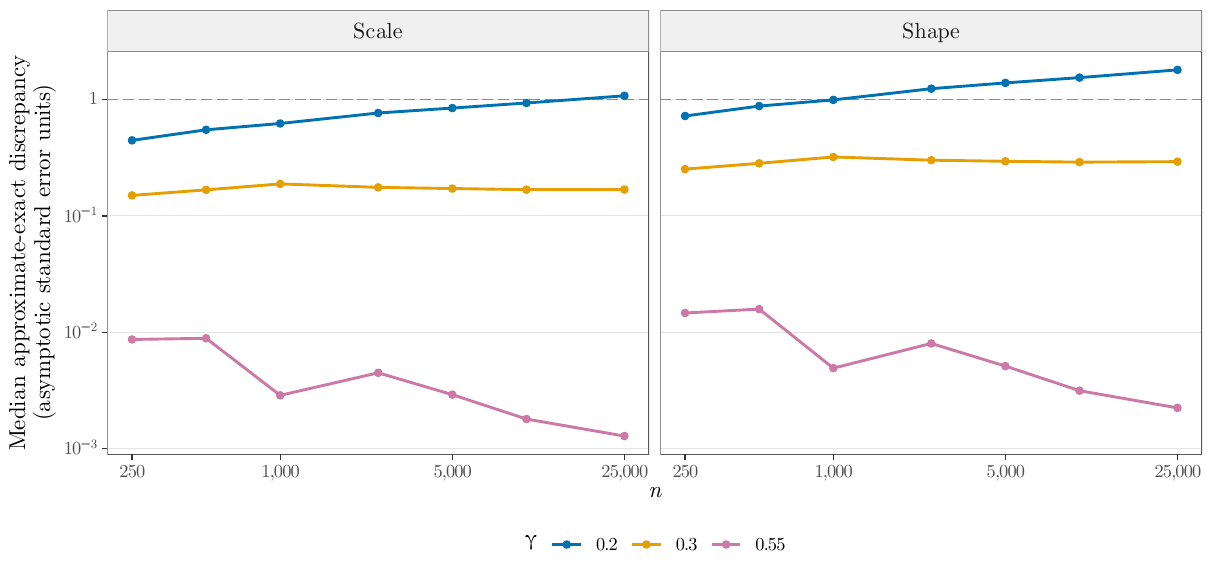}
\caption{Median absolute discrepancy between the constrained approximate and exact GPD estimators, measured in units of the limiting standard errors $\sqrt{2/n}$ for scale and $1/\sqrt n$ for shape. Each point uses all $2000$ replications.}
\label{fig:gpd-discrepancy}
\end{figure}

The discrepancy patterns in \cref{fig:gpd-discrepancy} must be interpreted together with the boundary replacement. When $\gamma=0.20$, the median discrepancy rises from $0.447$ to $1.080$ standard errors for scale and from $0.724$ to $1.805$ for shape. When $\gamma=0.30$, it remains below $0.19$ for scale and $0.33$ for shape throughout the grid; after $n=2500$, the medians lie between $0.168$ and $0.176$ for scale and between $0.290$ and $0.302$ for shape. The smaller values at the first design points partly reflect more frequent use of the boundary replacement, whose frequency falls from $24.0\%$ at $n=250$ to $1.25\%$ at $n=25000$. For $\gamma=0.55$, the very small unconditional discrepancies do not indicate interior root localization: the raw surrogate root is outside the parameter space in about $42$--$45\%$ of replications, and the constrained estimator then coincides with the same shifted-exponential boundary form often selected by the exact MLE.

\begin{figure}[!htbp]
\centering
\includegraphics[width=0.90\textwidth]{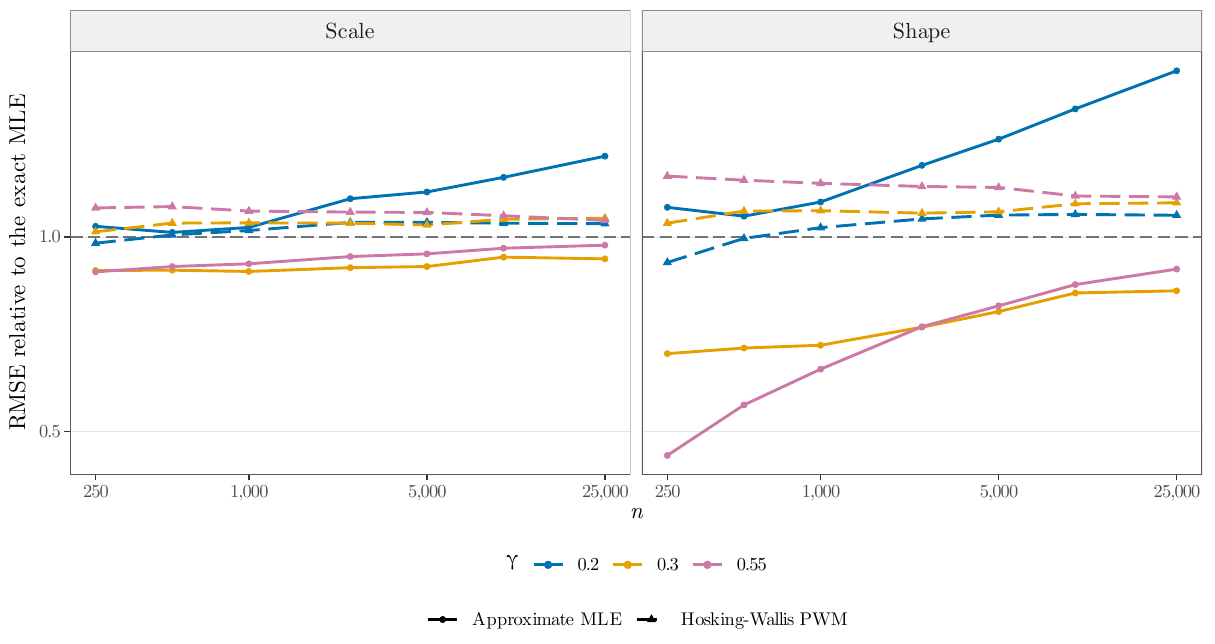}
\caption{RMSE of each explicit GPD estimator relative to the exact restricted MLE. The constrained approximate estimator ratio uses all $2000$ replications. The PWM ratio uses only replications in which the PWM estimator belongs to the restricted parameter space, with the denominator calculated from the same replications.}
\label{fig:gpd-rmse}
\end{figure}

\cref{fig:gpd-rmse} compares the two explicit estimators with the exact MLE. For $\gamma=0.20$, the unconditional RMSE ratio of our constrained estimator grows to $1.334$ for scale and $1.808$ for shape, whereas the probability weighted moment estimator remains within about 9\% of the exact MLE. For $\gamma=0.30$, the RMSE ratios of our estimator range from $0.885$ to $0.931$ for scale and from $0.660$ to $0.826$ for shape. For $\gamma=0.55$, the corresponding ranges are $0.884$--$0.972$ and $0.460$--$0.893$; these values again reflect frequent boundary replacement rather than validity of the interior approximation. No replication was excluded from the constrained estimator ratios. The PWM exclusion rate falls from $28.9\%$ at $n=250$ to $2.30\%$ at $n=25000$ when $\gamma=0.30$.

\FloatBarrier

\begin{figure}[!htbp]
\centering
\includegraphics[width=0.90\textwidth]{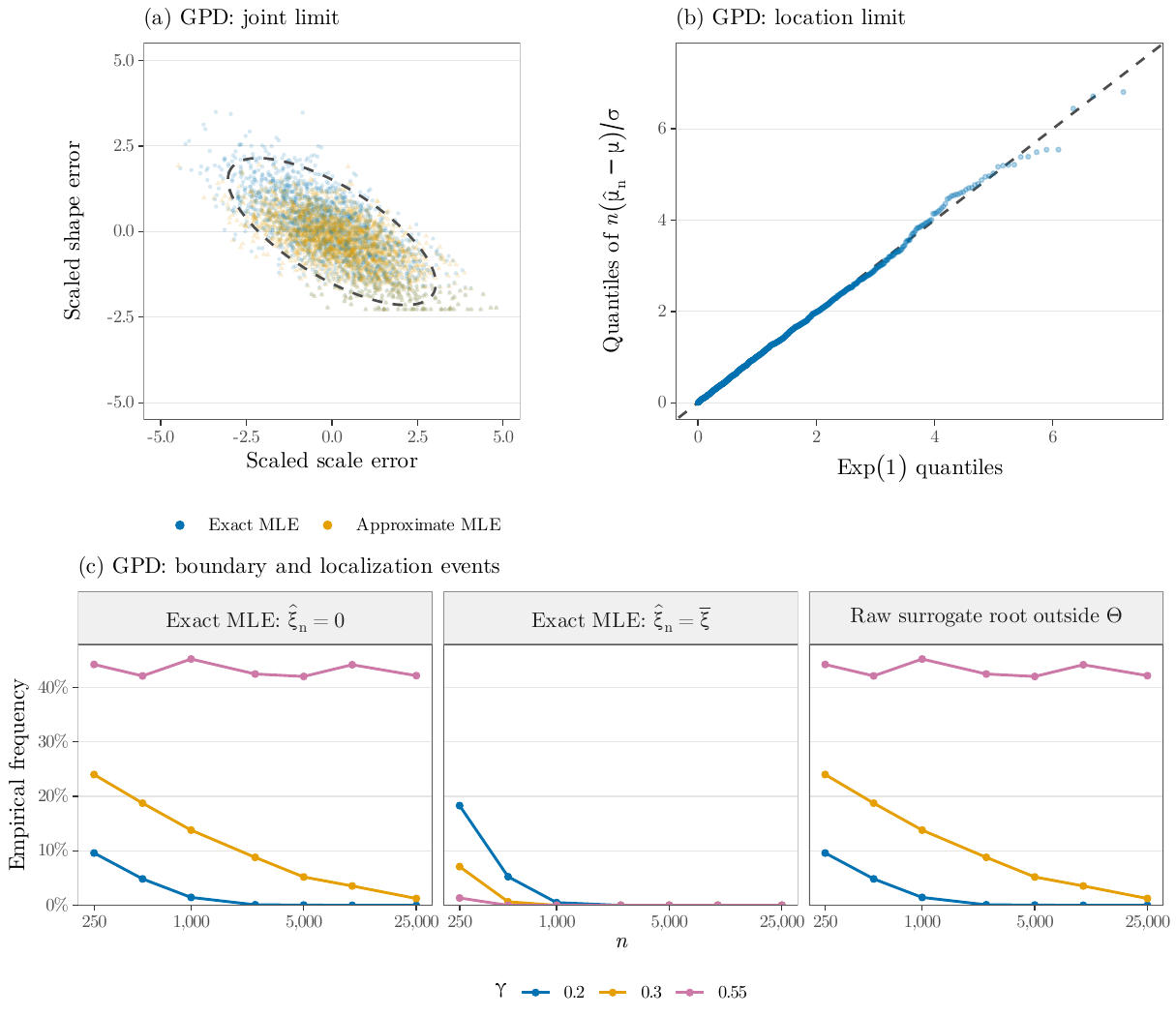}
\caption{GPD diagnostics. Panels (a) and (b) use $\gamma=0.30$ and $n=25000$. The dashed ellipse in panel (a) contains 90\% probability under the Gaussian limit in \eqref{eq:gpd_exact_clt}; the constrained approximate errors use all replications. Panel (b) compares the location error with its $\mathrm{Exp}(1)$ limit. Panel (c) reports the frequencies with which the exact MLE lies at either shape boundary and the raw surrogate root falls outside the restricted parameter space.}
\label{fig:gpd-diagnostics}
\end{figure}

At the representative setting $\gamma=0.30$ and $n=25000$, the empirical scale variance, shape variance, and covariance of the scaled exact MLE are $2.055$, $1.052$, and $-1.049$, close to the limiting values $2$, $1$, and $-1$. The mean of the scaled constrained-estimator error is $(0.251,-0.317)^\intercal$, compared with the leading population displacement $(0.261,-0.327)^\intercal$. The empirical mean and variance of $n(X_{(1)}-\mu)/\sigma$ are $0.997$ and $0.974$, and its Kolmogorov--Smirnov distance from $\mathrm{Exp}(1)$ is $0.0133$. These comparisons appear in panels (a) and (b) of \cref{fig:gpd-diagnostics}. Panel (c) shows that, when $\gamma=0.55$, the exact MLE remains at $\xi=0$ and the raw surrogate root leaves the parameter space in about 42--45\% of replications. At $\gamma=0.30$, both empirical frequencies decrease toward zero. These are boundary and localization events rather than optimization failures.

\FloatBarrier

\end{document}